\documentclass{article}
\usepackage[applemac]{inputenc}
\usepackage[T1]{fontenc}
\usepackage{lmodern}
\usepackage[a4paper]{geometry}
\usepackage[english]{babel}
\usepackage{graphicx}
\usepackage{onlyamsmath}
\usepackage{amsfonts}
\usepackage{amssymb}
\usepackage{ntheorem}

\newcommand{\e}{\epsilon}
\newcommand{\se}{\sqrt{\epsilon}}
\newcommand{\ga}{\gamma}
\newcommand{\gaa}{\lvert\gamma\rvert}
\newcommand{\sga}{\sqrt{\gamma}}
\newcommand{\sgap}{\sqrt{\gamma'}}
\newcommand{\sgaa}{\sqrt{\lvert\gamma\rvert}}
\newcommand{\D}{\mathcal{D}}
\newcommand{\E}{\mathbb{E}}
\newcommand{\esp}{\mathcal{H}^\omega}

\newcommand{\ko}[1]{k^2 (\omega_{#1})}
\newcommand{\N}[1]{N (\omega{#1})}

\newcommand{\Bh}[2]{\beta_{#1} (\omega{#2})}

\newcommand{\M}{\widehat{p} (\omega,x,z)}

\newcommand{\p}[1]{\widehat{p}_{#1} (\omega,z)}

\newcommand{\ha}[1]{\widehat{a}_{#1} (\omega,z)}

\newcommand{\hb}[1]{\widehat{b}_{#1} (\omega,z)}

\newcommand{\hf}{\widehat{f}(\omega)}
\newcommand{\1}{\textbf{1}}
\newcommand{\dz}{\frac{d}{dz}}

\newcommand{\Ns}{N'}
\newcommand{\Nss}{N''}
\newcommand{\Nsss}{N'''}

\newtheorem{thm}{Theorem}[section]
\newtheorem{prop}{Proposition}[section]

\newtheorem{lem}{Lemma}[section]
\theoremstyle{nonumberplain}
\theorembodyfont{\upshape}
\newtheorem{preuve}{Proof}

\title{Wave Propagation in Shallow-Water Acoustic Random Waveguides.}

\author{Christophe Gomez\thanks{Laboratoire de Probabilit\'es et Mod\`eles Al\'eatoires, Universit\'e
Paris 7, Boîte courrier 7012, 75251 PARIS Cedex 05 (FRANCE)  (chgomez@math.jussieu.fr).}}

\begin{document}

\maketitle

\begin{abstract} 

In shallow-water waveguides a propagating field can be decomposed over three kinds of modes: the propagating modes, the radiating modes and the evanescent modes. In this paper we consider the propagation of a wave in a randomly perturbed waveguide and we analyze the coupling between these three kinds of modes using an asymptotic analysis based on a separation of scales technique. Then, we derive the asymptotic form of the distribution of the mode amplitudes and the coupled power equation for propagating modes. From this equation, we show that the total energy carried by the propagating modes decreases exponentially with the size of the random section and we give an expression of the decay rate. Moreover, we show that the mean propagating mode powers converge to the solution of a diffusion equation in the high-frequency regime. 

\end{abstract}

\begin{flushleft}
\end{flushleft}

\begin{flushleft}
\end{flushleft}

\pagestyle{myheadings} \thispagestyle{plain} \markboth{CHRISTOPHE
GOMEZ}{C. GOMEZ \hspace{6cm}SHALLOW-WATER PROPAGATION}

\section*{Introduction.}

Acoustic wave propagation in shallow-water waveguides has been studied for a long time because of its numerous domains of applications. One of the most important applications is submarine detection with active or passive sonars, but it can also be used in underwater communication, mines or archaeological artifacts detection, and to study the ocean's structure or ocean biology. Shallow-waters are complicated media because they have indices of refraction with spatial and time dependences. However, the sound speed in water, which is about 1500 m/s, is sufficiently large with respect to the motions of water masses that we can consider this medium as being time independent. Moreover, the presence of spatial inhomogeneities in the water produces a mode coupling which can induce significant effects over large propagation distances.

\begin{figure}\begin{center}
\begin{tabular}{cc}
\includegraphics*[scale=0.3]{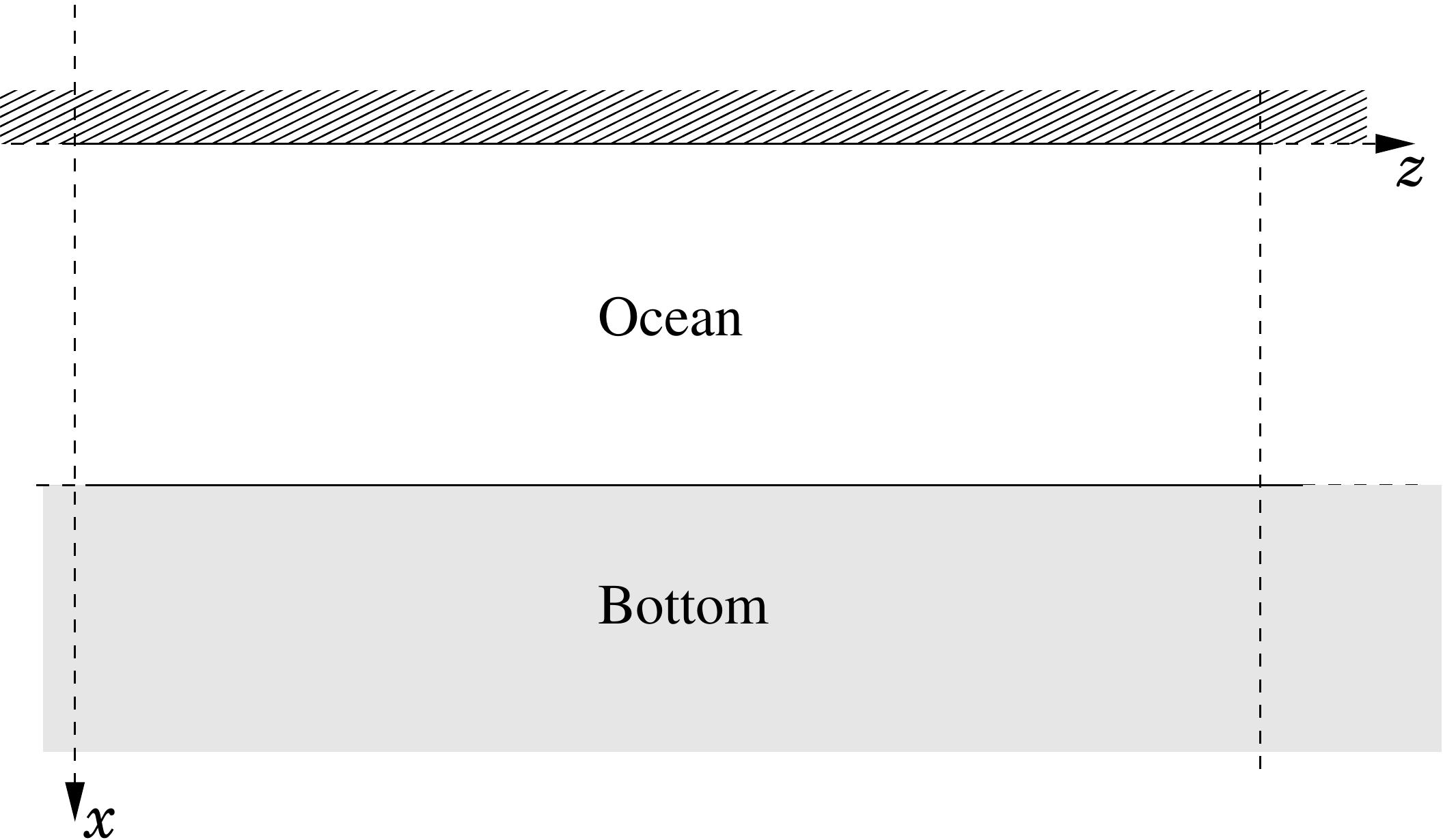} & \includegraphics*[scale=0.35]{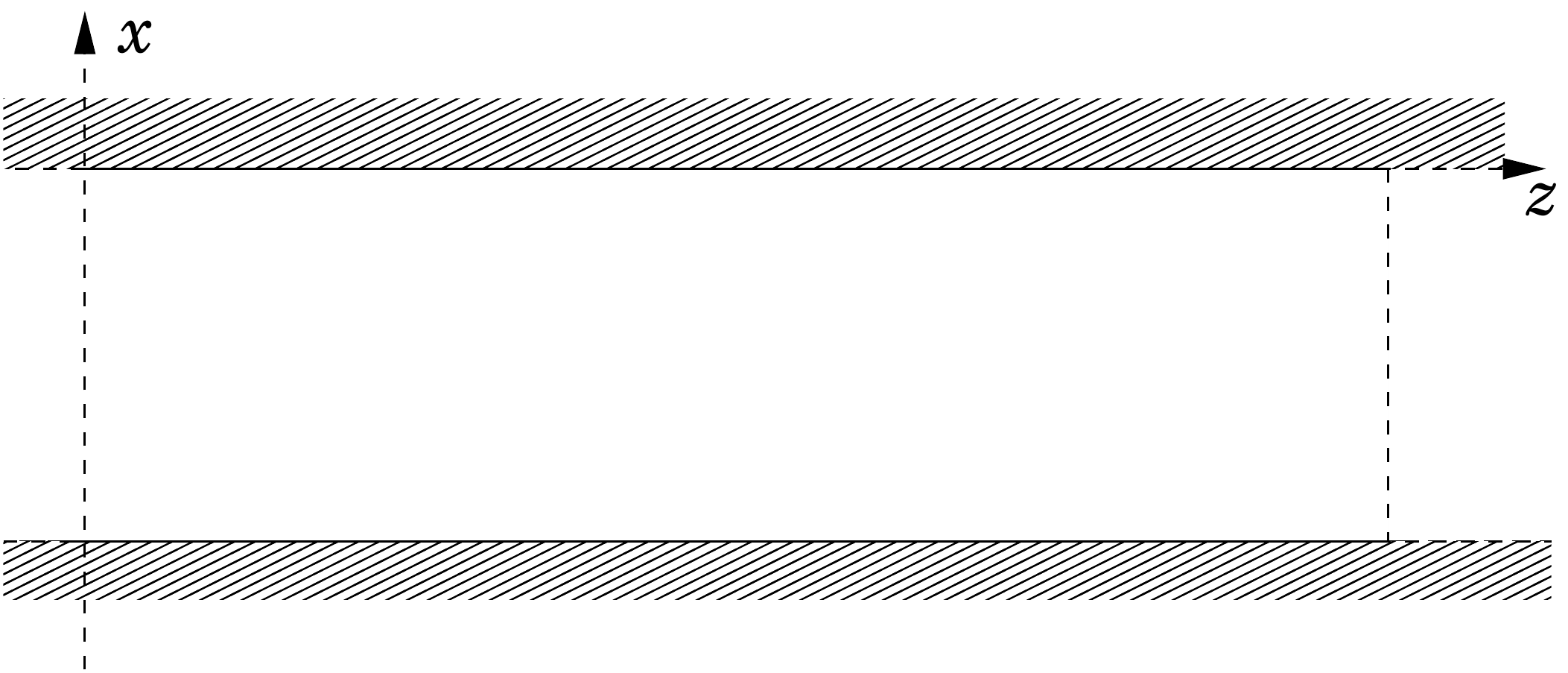}\\
$(a)$ & $(b)$
\end{tabular}
\end{center}
\caption{\label{wavegmod} Illustration of two kinds of waveguides. In $(a)$ we represent a shallow-water  waveguide model with an unbounded cross-section. In $(b)$ we represent a waveguide with a bounded cross-section.}
\end{figure}

In shallow-water waveguides the transverse section can be represented as a semi-infinite interval (see Figure \ref{wavegmod}) and then a wave field can be decomposed over three kinds of modes: the propagating modes which propagate over long distances, the evanescent modes which decrease exponentially with the propagation distance, and the radiating modes representing modes which penetrate under the bottom of the water.  The main purpose of this paper is to analyze how the propagating mode powers are affected by the radiating and evanescent modes. This analysis is carried out using an asymptotic analysis based on a separation of scale technique, where the wavelength and the correlation lengths of the inhomogeneities, which are of the same order, are small compared to the propagation distance, and the fluctuations of the medium are small compared to the wavelength. In the terminology of \cite{book} this is the so-called weakly heterogeneous regime.       

Wave propagation in random waveguides with a bounded cross-section and Dirichlet boundary conditions (see Figure \ref{wavegmod}) has been studied in \cite[Chapter 20]{book} or \cite{papa} for instance. In this case we have only two kinds of modes, the propagating and the evanescent modes. In such a model an asymptotic analysis of the mode powers show total energy conservation and an equipartition of the energy carried by the propagating modes. In \cite{papa} coupled power equations are derived under the assumption that evanescent modes are negligible. In \cite{garniereva} the role of evanescent modes is studied in the absence of radiating modes. In this paper we take into account the influence of the radiating and the evanescent modes on the coupled power equations. In this case we show a mode-dependent and frequency-dependent attenuation on the propagating modes in Theorem \ref{expdeccpeP2}, that is, the total energy carried by the propagating modes decreases exponentially with the size of the random section and we give an expression of the decay rate. Moreover, in the high-frequency regime, we show  in Theorems \ref{hfapproxP2} and \ref{hfapprox0P2} that the propagating mode powers converge to the solution of a diffusion equation. All the results of this paper are also valid for electromagnetic wave propagation in dielectric waveguides and optical fibers \cite{magnanini,marcuse,perrey,rowe,wilcox2}.   

The organization of this paper is as follows. In Section \ref{sect1P2} we present the waveguide model, and in Section \ref{sect2P2} we present the mode decomposition associated to that model and studied in detail in \cite{wilcox}. In Section \ref{sect3P2} we study the mode coupling when the three kinds of modes are taken into account. In the same spirit as in \cite[Chapter 20]{book}, we derive the coupled mode equations, we study the energy flux for the propagating and the radiating modes, and the influence of the evanescent modes on the two other kinds of modes. In Section \ref{coupledprocP2}, under the forward scattering approximation, we study the asymptotic form of the joint distribution of the propagating and radiating mode amplitudes. We apply this result in Section \ref{sect5P2} to derive the coupled power equations for the propagating modes, which was already obtained in \cite{papanicolaou} or \cite{marcuse} for instance. In this section, we study the influence of the radiating and evanescent modes on the mean propagating mode powers. We show that the total energy carried by the propagating modes decreases exponentially with the size of the random section and we give an expression of the decay rate. In other words, the radiating modes induce a mode-dependent attenuation on the propagating modes, that is why these modes are sometimes called dissipative modes. Moreover, under the assumption that nearest-neighbor coupling is the main power transfer mechanism, we show, in the high-frequency regime or in the limit of large number of propagating modes, that the mean propagating mode powers converge to the solution of a diffusion equation. We can refer to \cite{papanicolaou,marcuse} for further references and discussions about diffusion models. In that regime, we can also observe the exponential decay behavior caused by the radiative loss.   

\section{Waveguide Model}\label{sect1P2}

We consider a two-dimensional linear acoustic wave model. The conservation equations of mass and linear momentum are given by
\begin{equation}\label{conservationP2}
\begin{split}
\rho ( x,z)\frac{\partial \textbf{u}}{\partial t} + \nabla p &=\textbf{F}, \\
\frac{1}{K (x,z)} \frac{\partial p}{\partial t} + \nabla . \textbf{u} &=0, 
\end{split}
\end{equation}
where $p$ is the acoustic pressure, $\textbf{u}$ is the acoustic velocity, $\rho$ is the density of the medium, $K$ is the bulk modulus, and the source is modeled by the forcing term $\textbf{F} (t,x,z)$. The third coordinate $z$ represents the propagation axis along the waveguide. The transverse section of the waveguide is the semi-infinite interval $[0,+\infty)$, and $x \in [0,+\infty)$ represents the transverse coordinate. 
Let $d>0$, we assume that the medium parameters are given by
\[\begin{split}
\frac{1}{K(x,z)} & =  \left\{ \begin{array}{ccl} 
                                            \frac{1}{\bar{K}}\left( n^2(x)+\sqrt{\e} V(x,z) \right) & \text{ if }  &  x\in [0,d],\quad z\in [0,L/ \e] \\
                                             \frac{1}{\bar{K}}n^2(x) & \text{ if }  & \left\{\begin{array}{l} x\in[0,+\infty),\,z\in (-\infty,0)\cup(L/\e,+\infty)\\ \text{or}\\ x\in (d,+\infty),\,z\in(-\infty,+\infty). \end{array} \right.
                                          \end{array} \right. \\
\rho(x,z)&=  \bar{\rho}\quad \text{ if }\quad  x\in [0,+\infty),\,z\in\mathbb{R}. \\
\end{split}\] 
In this paper we consider the Pekeris waveguide model.  This kind of model has been studied for half a century \cite{pekeris} and in this model the index of refraction $n(x)$ is given by
\begin{equation*}
n(x)=\left\{\begin{array}{lcl} n_1>1& \text{if}& x\in[0,d)\\
                                                                 1& \text{if}& x\in[d,+\infty).
\end{array}\right.\end{equation*}
\begin{figure}\begin{center}
\includegraphics*[scale=0.4]{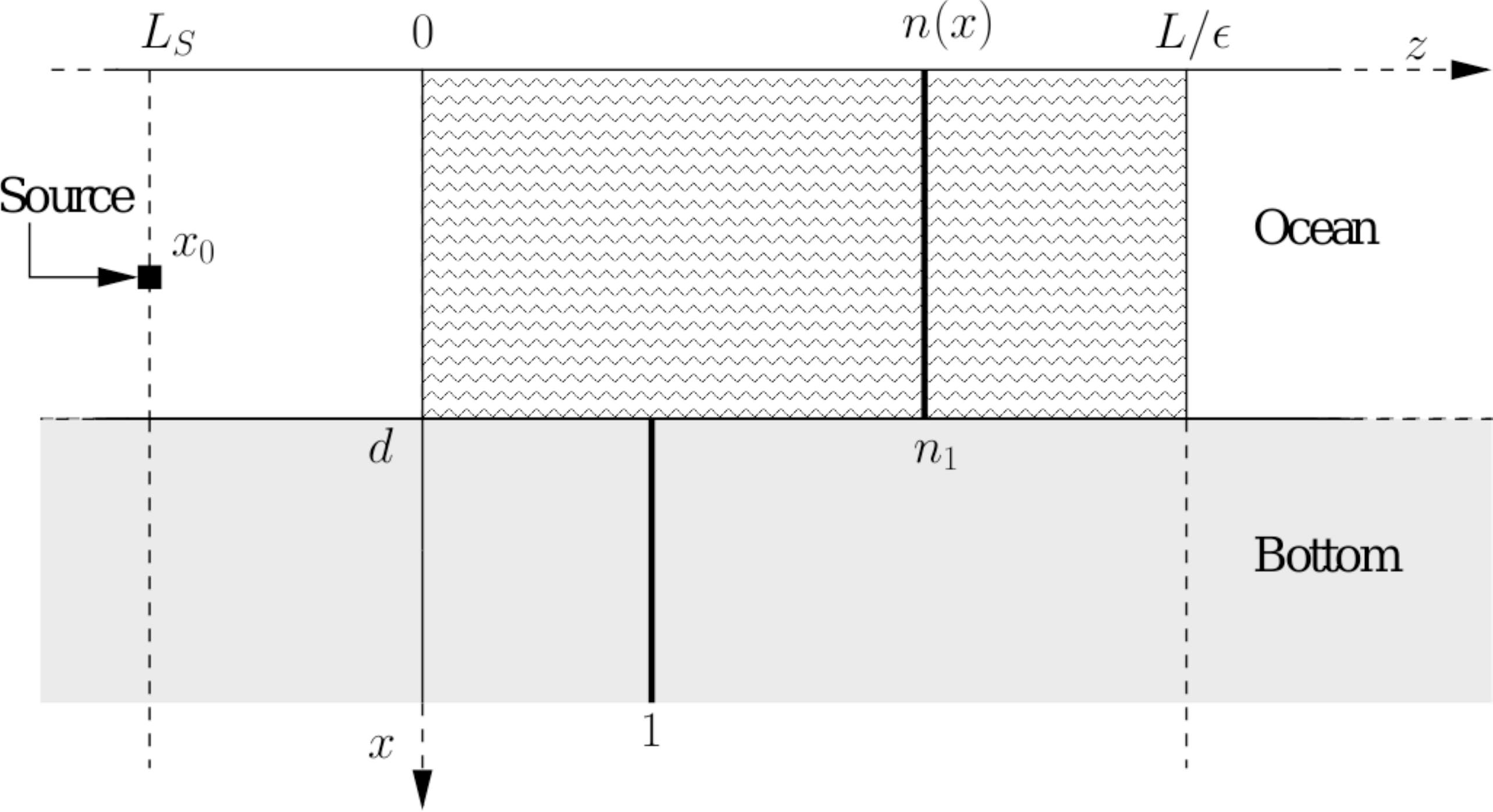}
\end{center}
\caption{\label{figureP2} Illustration of the shallow-water waveguide model.}
\end{figure}

This profile can model an ocean with a constant sound speed. Such conditions can be found during the winter in Earth's mid latitudes and in water shallower than about $30$ meters. The Pekeris profile leads us to simplified algebra but it underestimates the complexity of the medium. However, the analysis that we present in this paper can be extended to more general profiles $n(x)$ with general boundary conditions. In the Pekeris model that we consider $n_1$ represents the index of refraction of the ocean section $[0,d]$, where $d$ is the depth of the ocean, and we consider that the index of refraction of the bottom of the ocean is equal to $1$. This model can also be used to study the propagation of electromagnetic waves in a dielectric slab and an optical fiber  with randomly perturbed index of refraction \cite{magnanini,marcuse,rowe,wilcox2}.

We consider a source that emits a signal in the $z$-direction, which is localized in the plane $z=L_S$. 
\begin{equation}\label{sourcetermP2}  \textbf{F} (t,x,z)=\Psi(t,x)\delta(z-L_S)\textbf{e}_{z}.\end{equation}
$\Psi(t,x)$ represents the profile of the source and $\textbf{e}_{z}$ is the unit vector pointing in the $z$-direction. $L_S<0$ is the location of the source on the propagating axis.

The random process $(V(x,z), x \in [0,d], z\geq 0)$ that we consider, and which represents the spatial inhomogeneities is presented in Section \ref{grfield}.  However, one can remark that the process $V$ is unbounded. This fact implies that the bulk modulus can take negative values. In order to avoid this situation, we can work on the event 
\[ \Big( \forall (x,z) \in [0,d]\times [0,L/\e], \,  n_1+\sqrt{\e}V(x,z) >0 \Big).\]
In fact, the property \eqref{cg1} implies
\[ \begin{split}
\lim_{\e \to 0}\mathbb{P} \Big( \exists (x,z) \in [0,d]&\times [0,L/\e] : \,  n_1+\sqrt{\e}V(x,z) \leq 0 \Big)\\
& \leq \lim_{\e \to 0}\mathbb{P} \Big( \sqrt{\e}\sup_{z\in[0,L]} \sup_{ x \in [0,d]} \left\lvert V\left(x,\frac{z}{\e } \right) \right\rvert \geq n_1 \Big)= 0. \end{split}\]

\section{Wave Propagation in a Homogeneous Waveguide}\label{sect2P2}

In this section, we assume that the medium parameters are given by
\[\rho(x,z)=\bar{\rho} \text{ and } K(x,z)= \frac{\bar{K}}{n^2(x)}, \quad \forall ( x,z)\in [0,+\infty) \times \mathbb{R}.\]
From the conservation equations \eqref{conservationP2}, we can derive the wave equation for the pressure field,
\begin{equation}\label{waveeqP2}
\Delta p - \frac{1}{ c(x) ^2 }\frac{\partial ^2  p }{\partial t^2} = \nabla .\textbf{F},
\end{equation}
where $c(x)=c/n(x)$ with $c= \sqrt{\frac{\bar{K}}{\bar{\rho}}}$, and $\Delta = \partial ^2 _x + \partial ^2 _{z}$.

In underwater acoustics the density of air is very small compared to the density of water, then it is natural to use a pressure-release condition. The pressure is very weak outside the waveguide, and by continuity, the pressure is zero at the free surface $x=0$. This consideration leads us to consider the  Dirichlet boundary conditions
\begin{equation*} 
p (t,0,z)=0 \quad \forall (t,z)\in [0,+\infty)  \times \mathbb{R}.
\end{equation*}

Throughout this manuscript, we consider linear models of propagation. Therefore, the pressure $p(t,x,z)$ can be expressed as the superposition of monochromatic waves by taking its Fourier transform. Here, the Fourier transform and the inverse Fourier transform, with respect to time, are defined by
\[ \hf =\int f(t)e^{i \omega t} dt, \quad f(t)=\frac{1}{2 \pi} \int \hf e^{-i \omega t} d\omega.\] 

In the half-space $z>L_S$ (resp., $z<L_S$), taking the Fourier transform in \eqref{waveeqP2},  we get that $\M$ satisfies the time-harmonic wave equation without source term
\begin{equation}\label{helmo2P2} \partial ^2 _z\M+\partial ^2 _x\M +\ko{}n^2 (x)\M=0, \end{equation}
where $k(\omega)=\frac{\omega}{c}$ is the wavenumber, and with Dirichlet boundary conditions $\widehat{p}(\omega,0,z)=0$ $\forall z$. 
The source term implies the following jump conditions for the pressure field across the plane $z=L_S$
\begin{equation}\label{jumpscondP2}\begin{array}{ccl}
\widehat{p}(\omega,x,L_S^+)-\widehat{p}(\omega,x,L_S^-)&=& \widehat{\Psi}(\omega,x),\\
\partial_z \widehat{p}(\omega,x,L_S^+)-\partial_z \widehat{p}(\omega,x,L_S^-)&=&0.
\end{array}
\end{equation}

\subsection{Spectral Decomposition in Unperturbed Waveguides}\label{spectralP2}

This section is devoted to the presentation of the spectral decomposition of the Pekeris operator $\partial ^2 _x +\ko{}n^2 (x)$. The spectral analysis of this operator is carried out in \cite{wilcox}. Throughout this paper we are interested in solutions of \eqref{helmo2P2} such that
\begin{equation*}\begin{split}
\widehat{p}(\omega,.,.)\textbf{1}_{(L_S,+\infty)}(z) &\in \mathcal{C}^0 \Big((L_S,+\infty), H^1 _0 (0,+\infty)\cap H^2(0,+\infty) \Big)\cap \mathcal{C}^2 \Big((L_S,+\infty),H\Big),\\
\widehat{p}(\omega,.,.)\textbf{1}_{(-\infty,L_S)}(z)&\in \mathcal{C}^0 \Big((-\infty,L_S), H^1 _0 (0,+\infty)\cap H^2(0,+\infty) \Big)\cap \mathcal{C}^2 \Big((-\infty,L_S),H\Big),
\end{split}\end{equation*}  
where $H=L^2(0,+\infty)$. $H$ is equipped with the inner product defined by
\[\forall (h_1,h_2)\in H\times H,\quad \big<h_1,h_2\big>_H=\int_0^{+\infty}h_1(x)\overline{h_2(x)}dx.\]
Consequently, in the half-space $z>L_S$ (resp., $z<L_S$), we can consider \eqref{helmo2P2} as the operational differential equation
\begin{equation}\label{helmo2opP2}
 \frac{d^2}{d z^2} \widehat{p}(\omega,.,z)+R(\omega) \big(\widehat{p}(\omega,.,z)\big)=0 \end{equation}
in $H$, where $R(\omega)$ is an unbounded operator on $H$ with domain
\[ \mathcal{D}(R(\omega))= H^1 _0 (0,+\infty)\cap H^2(0,+\infty),\]
and defined by
\[R(\omega)(y)=\frac{d^2}{ dx^2} y +\ko{}n^2 (x)y\quad \forall y\in \mathcal{D}(R(\omega)).\]
According to \cite{wilcox}, $R(\omega)$ is a self-adjoint operator on the Hilbert space $H$, and its spectrum is given by
\begin{equation}\label{spectrumRP2}
Sp\big(R(\omega)\big)=\left(-\infty,\ko{}\right]\cup\big\{ \beta^2 _{\N{}}(\omega),\dots,\beta_1^2 (\omega)\big\}.\end{equation}
More precisely, $\forall j\in\big\{1,\dots,\N{}\big\}$, the modal wavenumber $\beta_j (\omega)$is positive and
\[\ko{}<\beta^2 _{\N{}}(\omega)<\cdots<\beta_1^2 (\omega)<n_1^2 \ko{}.\]
Moreover, there exists a resolution of the identity $\Pi_\omega$ of $R(\omega)$ such that $\forall y\in H$, $\forall r\in\mathbb{R}$,

\begin{equation*}\begin{split}
 \Pi_\omega(r,+\infty)(y) (x) =&  \sum_{j=1}^{\N{}} \big<y,\phi_j(\omega,.)\big>_H \phi_j(\omega,x)\textbf{1}_{(r,+\infty)}\left(\Bh{j}{}^2\right) \\
                                                                      & + \int_{r}^{\ko{}}  \big<y,\phi_\ga(\omega,.) \big>_H \phi_\ga (\omega,x)d\ga \textbf{1}_{\left(-\infty,\ko{} \right)}(r),
  \end{split}\end{equation*}
and $\forall y\in \mathcal{D}(R(\omega))$, $\forall r\in\mathbb{R}$,  
  \begin{equation*}\begin{split}                                                                    
 \Pi_\omega(r,+\infty)(R(\omega)(y)) (x) =&  \sum_{j=1}^{\N{}}\Bh{j}{}^2 \big<y,\phi_j(\omega,.)\big>_H \phi_j(\omega,x)\textbf{1}_{(r,+\infty)}\left(\Bh{j}{}^2\right) \\
                                                                      & + \int_{r}^{\ko{}}  \ga \big<y,\phi_\ga(\omega,x) \big>_H \phi_\ga (\omega,x)d\ga \textbf{1}_{\left(-\infty,\ko{} \right)}(r).
\end{split}\end{equation*}   
Let us describe these decompositions.

\paragraph{Discrete part of the decomposition}

$\forall j\in \big\{1,\dots,\N{}\big\}$, the $j$th eigenvector is given by \cite{wilcox}
\[\phi_j(\omega, x)=\left\{ \begin{array}{ccl}
A_j(\omega)\sin(\sigma_j(\omega) x/d) & \mbox{ if } & 0\leq x \leq d \\
A_j(\omega)\sin(\sigma_j(\omega)  )e^{-\zeta_j (\omega) \frac{x-d}{d}}& \mbox{ if } & d\leq x,  \end{array} \right.\]
where
\[\sigma_j (\omega)=d \sqrt{n_1 ^2 \ko{}-\beta^2 _j( \omega)}, \quad \zeta_j(\omega) =d \sqrt{\Bh{j}{}^2-\ko{}},\]
and 
\begin{equation}\label{coefajP2}
A_j(\omega)=\sqrt{\frac{2/d}{1+\frac{\sin^2 (\sigma_j (\omega))}{\zeta_j (\omega) }-\frac{\sin(2 \sigma_j (\omega)  )}{2\sigma_j (\omega) } }}.\end{equation}
According to \cite{wilcox},  $\sigma_1 (\omega),\dots,\sigma_{\N{}}(\omega)$ are the solutions on $(0,n_1k(\omega)d\theta)$ of the equation 
\begin{equation}\label{eqvpP2}
\tan(y )=-\frac{y }{\sqrt{( n_1 k d  \theta) ^2-y^2}},
\end{equation}
such that $0<\sigma_1 (\omega)<\cdots<\sigma_{\N{}}(\omega)< n_1 k(\omega)d \theta$,  and with  $\theta =\sqrt{1-1/n_1 ^2}$. This last equation admits exactly one solution over each interval of the form $\big(\pi/2+(j-1)\pi, \pi/2 +j\pi \big)$ for $ j \in \{1,\dots, N(\omega) \}$, where 
\[\N{}=\left[ \frac{n_1  k(\omega)d}{\pi} \theta \right],\]
and $[\cdot]$ stands for the integer part.  
From \eqref{eqvpP2}, we get the following results about the localization of the solutions which is used to show the main result of Section \ref{sechighfreq}.
\begin{lem}\label{coefgP2}
Let $\alpha>1/3$, we have as $\N{} \to +\infty$
\begin{equation*} 
\sup_{j\in \{1,\dots, \N{}-[\N{}^\alpha] -1 \}} \left\lvert \sigma_{j+1}(\omega)-\sigma_j(\omega) -\pi\right\rvert =\mathcal{O}\left( \N{}^{\frac{1}{2}-\frac{3}{2}\alpha } \right).
\end{equation*}

\begin{equation*} 
\sup_{j\in \{1,\dots,\N{}-[\N{}^\alpha ]-2 \}} \left\lvert \sigma_{j+2}(\omega)-2\sigma_{j+1}(\omega)+\sigma_{j}(\omega)\big)\right\rvert =\mathcal{O}\left( \N{}^{1-3\alpha } \right).
\end{equation*}
\end{lem}

\paragraph{Continuous part of the decomposition}

For $\ga \in (-\infty, \ko{})$, we have \cite{wilcox}
\[
\begin{split}
\phi_\ga& (\omega, x)=\\ 
&\left\{ 
\begin{array}{ccl}
A_\gamma(\omega) \sin(\eta(\omega)  x/d ) & \mbox{ if } & 0\leq x \leq d \\
A_\gamma(\omega) \left(\sin(\eta(\omega)  )\cos\big(\xi(\omega) \frac{x-d}{d}\big)+\frac{\eta(\omega) }{\xi(\omega) }\cos(\eta(\omega) )\sin\big(\xi(\omega) \frac{x-d}{d}\big)\right)& \mbox{ if } & d\leq x, \end{array}
 \right.
\end{split}
\]
where
\[\eta (\omega) =d\sqrt{n_1 ^2 \ko{}-\gamma }, \quad \xi(\omega)  =d\sqrt{\ko{}-\gamma },\]
and
\begin{equation}\label{coefagaP2}A_\gamma(\omega) =\sqrt{\frac{d \xi(\omega) }{\pi\big(\xi ^2(\omega) \sin^{2}(\eta(\omega))+\eta ^2(\omega) \cos^2 (\eta(\omega))\big)}}.\end{equation}
It is easy to check that the function $\ga\mapsto A_\ga(\omega)$ is continuous on $\big(-\infty, \ko{}\big)$ and
\begin{equation}\label{coefaga2P2}
A_\ga (\omega) \underset{\ga\to-\infty}{\sim} \frac{1}{\sqrt{\pi} \lvert\ga\rvert ^{1/4}}.
\end{equation}
We can remark that $\phi_\ga(\omega, .)$ does not belong to $H$. Then, $\big<y,\phi_\ga (\omega,.)\big>_H$ is not defined in the classical way. In fact,
\[ \big<y,\phi_\ga (\omega,.) \big>_H=\lim_{M\to+\infty}\int_{0}^M y(x) \phi_\ga (\omega, x)dx\quad  \text{ on }L^2\big(-\infty,\ko{}\big). \]
Moreover, we have $\forall y\in H$
\[ \|y \|^2 _H = \sum_{j=1}^{\N{}}\big \lvert \big<y,\phi_j(\omega,.) \big>_H \big \rvert^2+\int_{-\infty}^{\ko{}}\big\lvert \big<y,\phi_\ga (\omega,.) \big>_H\big\rvert^2 d\ga.\]
Then,
\[
\begin{array}{rcc}
\Theta_\omega:H&\longrightarrow &\mathcal{H}^\omega\\
y& \longrightarrow & \Big(\big(\big<y,\phi_j(\omega,.)\big>_H\big)_{ j=1,...,\N{}},\big(\big<y,\phi_\ga (\omega,.)\big>_H\big)_{\ga\in(-\infty,\ko{})}\Big)
\end{array}
\]
is an isometry, from $H$ onto $\esp=\mathbb{C}^{\N{}}\times L^2\big(-\infty,\ko{}\big)$.

\subsection{Modal Decomposition}

In this section we apply the spectral decomposition introduced in Section \ref{spectralP2} on a solution $\M$ of the equation \eqref{helmo2opP2}. Consequently, we get the modal decomposition for $\M$ in the half-space $z>L_S$, 
\[\M=\sum_{j=1}^{\N{}} \widehat{p}_j (\omega,z )\phi_j(\omega, x)+\int_{-\infty}^{\ko{}}\widehat{p}_\ga (\omega,z)\phi_\ga (\omega,x)d \ga, \]
where $\widehat{p}(\omega,z)=\Theta_\omega (\widehat{p}(\omega,.,z))$. 
For $j\in \big\{1,\dots,\N{}\big\}$, $\Theta_\omega\circ \Pi_\omega(\{j\})$ represents the projection onto the $j$th propagating mode, and $\widehat{p}_j (\omega,z )$ is the amplitude of the $j$th propagating mode. $\Theta_\omega\circ \Pi_\omega(0,\ko{})$ represents the projection onto the radiating modes, and $\widehat{p}_\ga (\omega,z)$ is the amplitude of the $\ga$th radiating mode for almost every $\ga\in(0,\ko{})$. Finally, $\Theta_\omega\circ \Pi_\omega(-\infty,0)$ represents the projection onto the evanescent modes and $\widehat{p}_\ga (\omega,z)$ is the amplitude of the  $\ga$th evanescent mode for almost every $\ga\in(-\infty,0)$.

Consequently, $\widehat{p}(\omega,z)$ satisfies
\begin{equation*}\begin{split}
\frac{d^2}{dz^2} \widehat{p}_j (\omega,z)+\beta^2 _j(\omega)\widehat{p}_j (\omega,z)&=0,\\
\frac{d^2}{dz^2} \widehat{p}_\ga (\omega,z)+\ga\,\,\widehat{p}_\ga (\omega,z)&=0
\end{split}\end{equation*}
in $\esp$ and the pressure field can be written as an expansion over the complete set of modes 
\begin{equation}\label{coefP2}
\begin{split}
\M &=\left[ \sum_{j=1}^{\N{}} \frac{ \widehat{a}_{j,0} (\omega) }{\sqrt{\Bh{j}{}}} e^{ i \Bh{j}{} z} \phi_j(\omega,x)+\int_{0}^{\ko{}} \frac{\widehat{a}_{\ga,0} (\omega)}{\ga^{1/4}} e^{ i \sga z}\phi_\ga(\omega,x)d \ga\right.\\
&\quad\quad\left.+\int_{-\infty}^{0} \frac{\widehat{c}_{\ga,0} (\omega)}{\gaa^{1/4}} e^{ - \sgaa z}\phi_\ga(\omega,x)d\ga
\right]\1 _ {(L_S,+\infty)}(z)\\
&+\left[ \sum_{j=1}^{\N{}} \frac{ \widehat{b}_{j,0} (\omega) }{\sqrt{\Bh{j}{}}} e^{ -i \Bh{j}{} z} \phi_j(\omega,x)+\int_{0}^{\ko{}} \frac{\widehat{b}_{\ga,0} (\omega)}{\ga^{1/4}} e^{ -i \sga z}\phi_\ga(\omega,x)d\ga\right.\\
&\quad \quad\left.+\int_{-\infty}^{0} \frac{\widehat{d}_{\ga,0} (\omega)}{\gaa^{1/4}} e^{ \sgaa z}\phi_\ga(\omega,x)d\ga
\right]\1 _ {(-\infty,L_S)}(z),
\end{split}
\end{equation}
under the assumption that the coefficients $\big(  \widehat{c}_{\ga,0} (\omega) e^{ - \sgaa L_S} /\gaa^{1/4}\big)_\ga$ and $\big(  \widehat{d}_{\ga,0} (\omega) e^{ \sgaa L_S} /\gaa^{1/4}\big)_\ga$ belong to  $L^2(-\infty,0)$. 

In the previous decomposition, $ \widehat{a}_{j,0} (\omega)$ (resp., $ \widehat{b}_{j,0} (\omega)$) is the amplitude of the $j$th right-going (resp., left-going) mode propagating in the right half-space $z>L_S$ (resp., left half-space $z<L_S$),  $ \widehat{a}_{\ga,0} (\omega)$ (resp., $ \widehat{b}_{\ga,0} (\omega)$) is the amplitude of the $\ga$th right-going (resp., left-going) mode radiating in the right half-space $z>L_S$ (resp., left half-space $z<L_S$), and $\widehat{c}_{\ga,0} (\omega)$ (resp., $\widehat{d}_{\ga,0} (\omega)$) is the amplitude of the $\ga$th right-going (resp., left-going) evanescent mode in the right half-space $z>L_S$ (resp., left half-space $z<L_S$). 

We assume that the profile $\Psi(t,x)$ of the source term \eqref{sourcetermP2} is given, in the frequency domain, by
\begin{equation}\label{profsourceP2}
 \widehat{\Psi}(\omega,x)=\widehat{f}(\omega)\left[\sum_{j=1}^{\N{}}\phi_j(\omega,x_0)\phi_j(\omega,x)+\int_{(-S,-\xi)\cup(\xi,\ko{})}\phi_\ga(\omega,x_0)\phi_\ga(\omega,x)d \ga\right],
 \end{equation}
where $x_0\in (0,d)$. The bound $S$ in the spectral decomposition of the source profile was introduced to have $\widehat{\Psi}(\omega,.)\in H$, and $\xi$ was introduced for technical reasons. Note that $S$ can be arbitrarily large and $\xi$ can be arbitrarily small. Therefore, the spatial profile in \eqref{profsourceP2} is an approximation of a Dirac distribution at $x_0$, which models a point source at $x_0$. 
  
Applying $\Theta_\omega$ on \eqref{jumpscondP2} and using \eqref{coefP2}, we get
\begin{equation*}
\begin{split}
\widehat{a}_{j,0} (\omega)  &= -\overline{\widehat{b}_{j,0} (\omega)} = \frac{\sqrt{\Bh{j}{}}}{2} \widehat{f}(\omega)\phi_j(\omega,x_0)e^{-i\Bh{j}{}L_S} \quad\forall j \in \big\{1,\dots,\N{} \big\},\\
\widehat{a}_{\ga,0} (\omega) & =-\overline{\widehat{b}_{\ga,0} (\omega)}  =\left\{ \begin{array}{ccc} \frac{\ga ^{1/4}}{2} \widehat{f}(\omega)\phi_\ga(\omega,x_0) e^{-i \sga L_S}& \text{for almost every} & \ga \in(\xi, \ko{})\\
0 & \text{for almost every} & \ga \in(0,\xi),\end{array}\right.\end{split}\end{equation*}
\begin{equation*}
\widehat{c}_{\ga,0} (\omega) = -\frac{\ga ^{1/4}}{2} \widehat{f}(\omega)\phi_\ga(\omega,x_0) e^{\sgaa L_S},\quad \widehat{d}_{\ga,0} (\omega) = \frac{\ga ^{1/4}}{2} \widehat{f}(\omega)\phi_\ga(\omega,x_0) e^{-\sgaa L_S} \end{equation*}
for almost every $\ga \in(-S,-\xi)$, and
\[\widehat{c}_{\ga,0} (\omega)= \widehat{d}_{\ga,0} (\omega) =0\]
for almost every $\ga \in(-\infty,-S)\cup(-\xi,0)$.

\section{Mode Coupling in Random Waveguides} \label{sect3P2}

In this section we study the expansion of $\M$ when a random section $[0,L/ \e]$ is inserted between two homogeneous waveguides (see Figure \ref{figureP2}). In this section the medium parameters are given by
\[\begin{split}
\frac{1}{K(x,z)} & =  \left\{ \begin{array}{ccl} 
                                            \frac{1}{\bar{K}}\left( n^2(x)+\sqrt{\e} V(x,z) \right) & \text{ if }  &  x\in [0,d],\quad z\in [0,L/ \e] \\
                                             \frac{1}{\bar{K}}n^2(x) & \text{ if }  & \left\{\begin{array}{l} x\in[0,+\infty),\,z\in (-\infty,0)\cup(L/\e,+\infty)\\ \text{or}\\ x\in (d,+\infty),\,z\in(-\infty,+\infty). \end{array} \right.
                                          \end{array} \right. \\
\rho(x,z)&=  \bar{\rho}\quad \text{ if }\quad  x\in [0,+\infty),\,z\in\mathbb{R}. \\
\end{split}\] 

In the perturbed section, the pressure field can be decomposed using the resolution of the identity $\Pi_\omega$ of the unperturbed waveguide:
\[\M=\sum_{j=1}^{\N{}} \widehat{p}_j (\omega,z )\phi_j(\omega, x)+\int_{-\infty}^{\ko{}}\widehat{p}_\ga (\omega,z)\phi_\ga (\omega,x)d \ga,\]
where $\widehat{p}(\omega,z)=\Theta_\omega (\widehat{p}(\omega,.,z))$. In what follows, we shall consider solutions of the form 
\[\M=\sum_{j=1}^{\N{}} \widehat{p}_j (\omega,z )\phi_j(\omega, x)+\int_{(-\infty,-\xi)\cup(\xi,\ko{})}\widehat{p}_\ga (\omega,z)\phi_\ga (\omega,x)d \ga.\]
This assumption leads us to simplified algebra in the proof of Theorem \ref{thasymptP21}. In such a decomposition, the radiating and the evanescent part are separated by the small band $(-\xi,\xi)$ with $\xi\ll 1$. The goal is to isolate the transition mode $0$ between the radiating and the evanescent part of the spectrum $Sp\big(R(\omega)\big)$ given by \eqref{spectrumRP2}. Moreover, we assume that $\e\ll\xi$ and therefore we have two distinct scales. Let us remark that in this paper, we consider in a first time the asymptotic $\e$ goes to $0$ and in a second time the asymptotic $\xi$ goes to $0$.

\subsection{Coupled Mode Equations}

In this section we give the coupled mode equations, which describe the coupling mechanism between the amplitudes of the three kinds of modes. 

In the random section $[0,L/\e]$ the pressure field $\widehat{p}(\omega,z)$ satisfies the following coupled equations in $\esp$:
\begin{equation}\label{eqdiff2P2}\begin{split}
\frac{d^2}{dz^2}\p{j}+\beta^2 _j(\omega) \p{j}&+\se\ko{}\sum_{l=1}^{\N{}}C^\omega_{jl}(z)\p{l}\\
            &+\se\ko{} \int_{(-\infty,-\xi)\cup(\xi,\ko{})}C^\omega_{j\ga'}(z)\p{\ga'}d\ga'=0,\\
\frac{d^2}{dz^2}\p{\ga}+\ga \,\,\p{\ga}&+\se\ko{}\sum_{l=1}^{\N{}}C^\omega_{\ga l}(z)\p{l}\\
            &+\se\ko{} \int_{(-\infty,-\xi)\cup(\xi,\ko{})} C^\omega_{\ga \ga'}(z)\p{\ga'}d\ga'=0,            
\end{split}\end{equation}
where
\begin{equation}\label{coefcouplVP2}\begin{split}
C^\omega_{jl}(z)&=\big<\phi_j(\omega,.),\phi_l(\omega,.) V(.,z)\big>_{H}=\int_0 ^d\phi_{j}(\omega, x) \phi_l(\omega, x) V(x,z) dx,\\
C^\omega_{j\ga}(z)&=C_{\ga j}(z)=\big<\phi_j(\omega,.),\phi_\ga(\omega,.) V(.,z)\big>_{H} =\int_0 ^d\phi_{j}(\omega, x) \phi_\ga(\omega, x) V(x,z) dx,\\
C^\omega_{\ga \ga'}(z)&=\big<\phi_\ga(\omega,.),\phi_{\ga'}(\omega,.) V(.,z)\big>_{H}=\int_0 ^d\phi_{\ga}(\omega, x) \phi_{\ga'}(\omega, x) V(x,z) dx.
\end{split}\end{equation}
We recall that $\widehat{p}(\omega,.,.)\in \mathcal{C}^0 \big((0,+\infty), H^1 _0 (0,+\infty)\cap H^2(0,+\infty) \big)\cap \mathcal{C}^2 \big((0,+\infty),H\big)$, then
\begin{equation} \label{chypP2}\int_{-\infty}^{-\xi} \ga^2 \lvert \p{\ga} \rvert ^2 d\ga<+\infty.\end{equation}
In the previous coupled equation the coefficients $C^\omega(z)$ represent the coupling between the three kinds of modes, which are the propagating, radiating and evanescent modes.
 
Next, we introduce the amplitudes of the  generalized right- and left-going modes 
$\widehat{a}(\omega,z)$ and $\widehat{b}(\omega,z)$, which are given by
\begin{equation*}\begin{split}
\p{j}&=\frac{1}{\sqrt{\Bh{j}{}}}\Big( \ha{j} e^{i\Bh{j}{}z} +\hb{j}e^{-i\Bh{j}{}z} \Big), \\
\frac{ d }{dz}\p{j} &= i \sqrt{\Bh{j}{}} \Big( \ha{j} e^{i\Bh{j}{}z} - \hb{j}e^{-i\Bh{j}{}z} \Big),\\
\p{\ga}&= \frac{1}{\ga^{1/4}}\Big( \ha{\ga} e^{i\sga z} +\hb{\ga}e^{-i\sga z} \Big),\\
\frac{ d }{dz}\p{\ga} &= i \ga^{1/4} \Big( \ha{\ga} e^{i\sga z} - \hb{\ga}e^{-i\sga z} \Big)
\end{split}\end{equation*}
$\forall j \in \big\{1,\dots,\N{}\big\}$ and almost every $\ga \in (\xi,\ko{})$. 
Let 
\[\esp_\xi=\mathbb{C}^{\N{}}\times L^2(\xi,\ko{}).\]
From \eqref{eqdiff2P2}, we obtain the coupled mode equation in $\esp_\xi \times \esp_\xi\times L^2(-\infty,-\xi)$ for the amplitudes $\big(\widehat{a}(\omega,z),\widehat{b}(\omega,z),\widehat{p}(\omega,z)\big)$:
\begin{equation} \label{mcP21}
\begin{split}
\dz &\widehat{a} _{j}(\omega, z)=\se\frac{i\ko{}}{2} \sum_{l=1}^{\N{}}\frac{C^\omega_{jl}(z)}{ \sqrt{{\beta} _{j}{\beta} _{l}}} \left( \widehat{a} _{l}(\omega,z)e^{i \left({\beta} _{l}-{\beta} _{j} \right)z} + \widehat{b} _{l}(\omega, z)e^{-i \left( {\beta} _{l}+{\beta} _{j}\right)z} \right) \\
&+\se\frac{i\ko{}}{2}\int_{\xi}^{\ko{}} \frac{C^\omega_{j\ga'}(z)}{\sqrt{\beta_j\sgap} }\left(  \widehat{a} _{\ga'}(\omega,z)e^{i \left(\sgap-{\beta} _{j} \right)z} + \widehat{b} _{\ga'} (\omega,z)e^{-i \left(\sgap+{\beta} _{j} \right)z} \right) d\ga'\\
&+\se\frac{i\ko{}}{2}\int_{-\infty}^{-\xi} \frac{C^\omega_{j\ga'}(z)}{\sqrt{\beta_j}}\p{\ga'}d\ga' e^{-i \beta_j z},
\end{split}\end{equation}

\begin{equation}\label{mcP22}\begin{split}
\dz& \widehat{a}_{\ga}(\omega,z)=\se\frac{i \ko{}}{2} \sum_{l=1}^{\N{}}\frac{C^\omega_{\ga l}(z)}{\sqrt{\sga {\beta}_{l}}} \left( \widehat{a} _{l}(\omega,z)e^{i \left({\beta} _{l}-\sga\right)z} + \widehat{b} _{l}(\omega,z)e^{-i \left({\beta} _{l}+\sga\right)z} \right)\\
&+\se\frac{i\ko{}}{2}\int_{\xi}^{\ko{}} \frac{C^\omega_{\ga\ga'}(z)}{\ga^{1/4}{\ga'}^{1/4} }\left(  \widehat{a} _{\ga'}(\omega,z)e^{i \left(\sgap-\sga \right)z} + \widehat{b} _{\ga'}(\omega,z) e^{-i \left(\sgap+\sga \right)z} \right) d\ga'\\
&+\se\frac{i\ko{}}{2}\int_{-\infty}^{-\xi} \frac{C^\omega_{\ga\ga'}(z)}{\ga^{1/4}}\p{\ga'}d\ga' e^{-i \ga z},
\end{split}
\end{equation}

\begin{equation}\label{mcP23}\begin{split}
\dz& \widehat{b}_{j}(\omega,z)=-\se\frac{i \ko{}}{2} \sum_{l=1}^{\N{}}\frac{C^\omega_{jl}(z)}{\sqrt{{\beta} _{j} {\beta}_{l}}} \left( \widehat{a} _{l}(\omega,z)e^{i \left({\beta} _{l}+{\beta} _{j}\right)z} + \widehat{b} _{l}(\omega,z)e^{-i \left({\beta} _{l}-{\beta} _{j}\right)z} \right)\\
&-\se\frac{i\ko{}}{2}\int_{\xi}^{\ko{}} \frac{C^\omega_{j\ga'}(z)}{\sqrt{\beta_j\sgap} }\left(  \widehat{a} _{\ga'}(\omega,z)e^{i \left(\sgap+{\beta} _{j} \right)z} + \widehat{b} _{\ga'}(\omega,z) e^{-i \left(\sgap-{\beta} _{j} \right)z} \right) d\ga'\\
&-\se\frac{i\ko{}}{2}\int_{-\infty}^{-\xi} \frac{C^\omega_{j\ga'}(z)}{\sqrt{\beta_j}}\p{\ga'}d\ga' e^{-i \beta_j z},
\end{split}
\end{equation}

\begin{equation}\label{mcP24}\begin{split}
\dz &\widehat{b} _{\ga}(\omega,z)=-\se\frac{i\ko{}}{2} \sum_{l=1}^{\N{}}\frac{C^\omega_{\ga l}(z)}{ \sqrt{\sga{\beta} _{l}}} \left( \widehat{a} _{l}(\omega,z)e^{i \left({\beta} _{l}+\sga \right)z} + \widehat{b} _{l}(\omega,z)e^{-i \left( {\beta} _{l}-\sga\right)z} \right) \\
&-\se\frac{i\ko{}}{2}\int_{\xi}^{\ko{}} \frac{C^\omega_{\ga\ga'}(z)}{\ga^{1/4}{\ga'}^{1/4}}\left(  \widehat{a} _{\ga'}(\omega,z)e^{i \left(\sgap+\sga\right)z} + \widehat{b} _{\ga'}(\omega,z) e^{-i \left(\sgap-\sga \right)z} \right) d\ga'\\
&-\se\frac{i\ko{}}{2}\int_{-\infty}^{-\xi} \frac{C^\omega_{\ga\ga'}(z)}{\ga^{1/4}}\p{\ga'}d\ga' e^{-i \sga z},
\end{split}\end{equation}

\begin{equation}\label{eq5P2}
\frac{d^2}{dz^2}\p{\ga}+\ga \,\,\p{\ga}+\se g_\gamma (\omega, z)=0,           
\end{equation}
where
\begin{equation}\label{ggammaP2}\begin{split} g_\gamma (\omega, z)&= \ko{}\sum_{l=1}^{\N{}}\frac{C^\omega_{\ga l}(z)}{\sqrt{\beta_l}}\left(   \widehat{a} _{l}(\omega,z)e^{i \beta_l z} + \widehat{b} _{l}(\omega,z) e^{-i \beta_l z} \right)\\
&+\ko{} \int_{\xi}^{\ko{}} \frac{C^\omega_{\ga \ga'}(z)}{{\ga'}^{1/4}}\left(   \widehat{a} _{\ga'}(\omega,z)e^{i \sgap z} + \widehat{b} _{\ga'}(\omega,z) e^{-i \sgap z} \right)d\ga'\\
&+\ko{}\int_{-\infty}^{-\xi}C^\omega_{\ga \ga'}(z)\p{\ga'}d\ga'.\end{split}\end{equation}
Let us note that in absence of random perturbations, the amplitudes $\widehat{a}(\omega,z)$ and $\widehat{b}(\omega,z)$ are constant. 

We assume that a pulse is emitted at the source plane $L_S<0$ and propagates toward the randomly perturbed slab $[0,L/\e]$. Using the previous section, the form of this incident field at $z=0$ is given by
\begin{equation}\label{initfieldP2}\widehat{p}(\omega ,x,0)= \sum_{j=1}^{\N{}} \frac{ \widehat{a}_{j,0} (\omega) }{\sqrt{\Bh{j}{}}}\phi_j(\omega,x)+\int_{\xi}^{\ko{}} \frac{\widehat{a}_{\ga,0} (\omega)}{\ga^{1/4}} \phi_\ga(\omega,x)d\ga+\int_{-S}^{-\xi} \frac{\widehat{c}_{\ga,0} (\omega)}{\gaa^{1/4}} \phi_\ga(\omega,x)d\ga.
\end{equation}
Consequently, by the continuity of the pressure field across the interfaces $z=0$ and $z=L/\e$, the coupled mode system is complemented with the boundary conditions
\[ \widehat{a}(\omega,0)=\widehat{a}_{0}(\omega)\quad \text{ and } \quad\widehat{b}\left(\omega,\frac{L}{\e}\right)=0\]
in $\esp_\xi$. 
For $j\in \big\{1,\dots,\N{}\big\}$, $\widehat{a}_{j,0}(\omega)$ represents the initial amplitude of the $j$th propagating mode, and for $\ga \in(\xi, \ko{})$, $\widehat{a}_{\ga,0}(\omega)$ represents the initial amplitude of the $\ga$th radiating mode at $z=0$. Moreover, for $\ga \in(-S, -\xi)$, $\widehat{c}_{\ga,0}(\omega)$ represents the initial amplitude of the $\ga$th evanescent mode at $z=0$.
The second condition implies that no wave comes from the right homogeneous waveguide.

\subsection{Energy Flux for the Propagating and Radiating Modes}\label{globalP2}

In this section we study the energy flux for the propagating and radiating modes, and the influence of the evanescent modes on this flux.

We begin this section by introducing the radiation condition for the evanescent modes
\begin{equation*}
\lim_{z\to+\infty}\big\| \Pi_\omega(-\infty,-\xi)\big(\widehat{p}(\omega,.,z)\big)\big\|^2 _{H}=0.
\end{equation*}
This condition means, in the homogeneous right half-space, that the energy carried by the evanescent modes decay as the propagation distance becomes large. From the radiation condition and \eqref{eq5P2}, we get for almost every $\ga \in(-\infty,-\xi)$
\begin{equation}\label{pgammaP2}\begin{split}
\p{\ga}&=\frac{\se}{2\sgaa}\int_0^{z\wedge L/\e} g_\ga (\omega,u)e^{\sgaa (u-z)}du+ \frac{\se}{2\sgaa}\int_{z\wedge L/\e}^{L/\e} g_\ga (\omega,u)e^{\sgaa (z-u)}du \\
& +  \phi_{\ga}(\omega,x_0)e^{-\sgaa(z-L_S)}\textbf{1}_{(-S,-\xi)}(\ga)
\end{split}\end{equation} 
$\forall z\in[0,+\infty)$.
According to \eqref{coefP2}, the relation \eqref{pgammaP2} can be viewed as a perturbation of the form of the evanescent mode without a random perturbation. Using the same arguments as in \cite[Chapter 20]{book}, we get $\forall z\in [0,L/\e]$,
\begin{equation*}
\frac{d}{dz}\big( \|\widehat{a}(\omega,z)\|^2_{\esp_\xi} -\|\widehat{b}(\omega,z)\|^2_{\esp_\xi}\big)=-\se Im\left( \int_{-\infty}^{-\xi} \overline{g_\ga (\omega, z)}\widehat{p}_\ga(\omega,z)d\ga \right),
\end{equation*}
and
\begin{equation}\label{localflxP2} \begin{split}
\|\widehat{a}(\omega,z)\|^2_{\esp_\xi} -\|\widehat{b}(\omega,z)\|^2_{\esp_\xi}&= \|\widehat{a}_0(\omega)\|^2_{\esp_\xi}- \|\widehat{b}_0(\omega)\|^2_{\esp_\xi}-\frac{\e}{2}\int_{-\infty}^{-\xi}\frac{G_\ga(\omega,z)}{\sgaa}d\ga\\
&-\sqrt{\e}\int_{-S}^{-\xi}\phi_\ga (\omega,x_0)e^{\sgaa L_S}\int_{0}^z Im\big(\overline{g_\ga(\omega,u)}\big)e^{-\sgaa u}du\,d\ga,
\end{split} 
\end{equation}
where
\[G_\ga (\omega,z)=\int_{0}^{z}\int_{z}^{L/\e}Im\big(\overline{g_\ga (\omega, u)}g_\ga (\omega, v)\big)e^{\sgaa (u-v)}dvdu.\]
Consequently, for $z=L/\e$, we get 
\[\begin{split} \|\widehat{a}(\omega,L/\e)\|^2_{\esp_\xi} + \|\widehat{b}(\omega,0)\|^2_{\esp_\xi}&= \|\widehat{a}_0(\omega)\|^2_{\esp_\xi}\\
&-\sqrt{\e}\int_{-S}^{-\xi}\phi_\ga (\omega,x_0)e^{\sgaa L_S}\int_{0}^{L/\e} Im\big(\overline{g_\ga(\omega,u)}\big)e^{-\sgaa u}du\,d\ga. \end{split}\]
The second term on the right side of the previous relation has the factor $\phi_\ga (\omega,x_0)e^{\sgaa L_S}$ which is  the form of the evanescent mode at $z=0$ without a random perturbation. Therefore, if $L_S$ is far away from $0$ and whatever the source (evanescent modes decay exponentially from $L_S$ to $0$) or if there is no excitation of modes $\ga \in (-\infty,-\xi)$ by the source (that is when $S=\xi$), we can get the conservation of the global energy flux for the propagating and radiating modes:
\[\|\widehat{a}(\omega,L/\e)\|^2_{\esp_\xi} + \|\widehat{b}(\omega,0)\|^2_{\esp_\xi}= \|\widehat{a}_0(\omega)\|^2_{\esp_\xi}.\]
However, from \eqref{localflxP2} and even if there is no evanescent modes in \eqref{initfieldP2}, the local energy flux is not conserved. The energy related to the evanescent modes is given by the last two terms on the right side in \eqref{localflxP2}. Let us estimate these two quantities. First,
\[ \begin{split}\sup_{z\in[0,L/\e]}\left\vert\int_{-\infty}^{-\xi}\frac{G_\ga(\omega,z)}{\sgaa}d\ga\right\vert& \leq K(\xi,d)  \sup_{z\in[0,L]}\sup_{x\in[0,d]}\left\lvert V\left(x,\frac{z}{\e}\right) \right\rvert^2\\
&\times \sup_{z\in[0,L/\e]}\|\widehat{a}(\omega,z) \|^2_{\esp_\xi}+\|\widehat{b}(\omega,z) \|^2_{\esp_\xi}+\| \widehat{p}(\omega,z) \|^2 _{L^1(-\infty,-\xi)}.
\end{split}\]
Second,
\[\begin{split}
\sup_{z\in[0,L/\e]}\Big\lvert\int_{-S}^{-\xi}\phi_\ga (\omega,x_0)&e^{\sgaa L_S}\int_{0}^{z} Im\big(\overline{g_\ga(\omega,u)}\big)e^{-\sgaa u}du\,d\ga\Big\rvert\\
&\leq K(\xi,d)  \sup_{z\in[0,L]}\sup_{x\in[0,d]}\left\lvert V\left(x,\frac{z}{\e}\right) \right\rvert\\
&\quad \times \sup_{z\in[0,L/\e]}\|\widehat{a}(\omega,z) \|_{\esp_\xi}+\|\widehat{b}(\omega,z) \|_{\esp_\xi}+\| \widehat{p}(\omega,z) \| _{L^1(-\infty,-\xi)}.
\end{split}\]
In the two previous inequalities $K(\xi,d)$ represents a constant which can change between the different relations. However, it is difficult to get good a priori  estimates about 
\begin{equation}\label{normP2}\sup_{z\in[0,L/\e]}\|\widehat{a}(\omega,z) \|^2_{\esp_\xi}+\|\widehat{b}(\omega,z) \|^2_{\esp_\xi}+\| \widehat{p}(\omega,z) \|^2 _{L^1(-\infty,-\xi)}.\end{equation}
For this reason, let us introduce the stopping "time"
\[L^\e =\inf\left(L>0, \quad \sup_{z\in[0,L/\e]}\|\widehat{a}(\omega,z) \|^2_{\esp_\xi}+\|\widehat{b}(\omega,z) \|^2_{\esp_\xi}+\| \widehat{p}(\omega,z) \|^2 _{L^1(-\infty,-\xi)} \geq \frac{1}{\sqrt{\e}}\right).\]
The role of this stopping "time" is to limit the size of the random section to ensure that the quantity \eqref{normP2} is not too large.
Consequently, the energy carried by the evanescent modes over the section $[0,L/\e]$ for $L\leq L^\e$, is at most of order $\mathcal{O}\big(\e^{1/4}\sup_{z\in[0,L/\e]}\sup_{x\in[0,d]}\!\left\lvert V\left(x,z\right) \right\rvert^2\big)$, and according to \eqref{cg1} the local energy flux for the propagating and the radiating modes is conserved in the asymptotic $\e \to 0$. More precisely, we can show that $\forall \eta >0$,
\begin{equation}\label{conspropP2} \lim_{\e\to0}\mathbb{P}\left(\sup_{z\in[0,L/\e]}\left\lvert   \|\widehat{a}(\omega,z)\|^2_{\esp_\xi} -\|\widehat{b}(\omega,z)\|^2_{\esp_\xi}- \|\widehat{a}_0(\omega)\|^2_{\esp_\xi}+ \|\widehat{b}_0(\omega)\|^2_{\esp_\xi}  \right\rvert>\eta,\, L\leq L^\e  \right)=0.\end{equation}
In Section \ref{coupledprocP2}, we shall see, under the forward scattering approximation, that the condition $L\leq L^\e$ is readily fulfilled in the limit $\e\to 0$, that is we have $\lim_{\e \to 0}\mathbb{P}(L^\e \leq L)=0$.

\subsection{Influence of the Evanescent Modes on the Propagating and Radiating Modes}\label{iemprP2}

We analyze, in this section, the influence of the evanescent modes on the coupling mechanism between the propagating and the radiating modes.

First of all, we recall that $\Theta_\omega\circ\Pi_\omega(-\infty,-\xi)\big(\widehat{p}(\omega,.,z)\big)$ represents the evanescent part of the pressure field $\widehat{p}(\omega,.,z)$, where $\Theta_\omega$ and $\Pi_\omega$ are defined in Section \ref{spectralP2}. In this section we consider $F=L^1 (-\infty,-\xi)$
equipped with the norm
\[\|y\|_F=\int_{-\infty}^{-\xi}\lvert y_\ga\rvert d\ga,\]
which is a Banach space. 
Substituting \eqref{ggammaP2} into \eqref{pgammaP2}, we get 
\begin{equation}\label{eqevmP2}
(Id-\se \Phi^\omega)\Big(\Theta_\omega\circ\Pi_\omega(-\infty,-\xi)\big(\widehat{p}(\omega,.,.)\big)\Big)=\se\,\tilde{p}(\omega,.)+\tilde{p}_0 (\omega,.).
\end{equation}
This equation holds in the Banach space $\big(\mathcal{C}\big([0,+\infty),F\big),\|.\|_{\infty,F}\big)$,
where
\[ \|y\|_{\infty,F}=\sup_{z\geq 0} \|y(z)\|_F \quad  \forall y\in \mathcal{C}\big([0,+\infty),F\big).\]
In \eqref{eqevmP2}, $\Phi^\omega$ is a linear bounded operator, from $\big(\mathcal{C}\big([0,+\infty),F\big),\|.\|_{\infty,F}\big)$ to itself,
defined by
\[\begin{split}
\Phi^\omega_\ga (y)(z)&=\frac{\ko{}}{2\sgaa}\int_{0}^{z\wedge L/\e}\int_{-\infty}^{-\xi} C^\omega_{\ga \ga'}(u)y_{\ga'} (u)d\ga' e^{\sgaa (u-z)}du\\
  &+\frac{\ko{}}{2\sgaa}\int_{z\wedge L/\e}^{L/\e}\int_{-\infty}^{-\xi} C^\omega_{\ga \ga'}(u)y_{\ga'} (u)d\ga' e^{\sgaa (z-u)}du
\end{split}\]
$\forall z\in [0,+\infty)$,
and for almost every $\ga \in (-\infty,-\xi)$
\[\begin{split}
\tilde{p}_\ga(\omega,z)&=\frac{\ko{}}{2\sgaa}\int_0^{z\wedge L/\e} \Big[ \sum_{l=1}^{\N{}}\frac{C^\omega_{\ga l}(u)}{\sqrt{\beta_l}}\big(\widehat{a}_l (\omega,u)e^{i\beta_l u}+\widehat{b}_l (\omega,u)e^{-i\beta_l u}\big)\\
&\quad \quad\quad+\int_{\xi}^{\ko{}}\frac{C^\omega_{\ga \ga'}(u)}{{\ga'}^{1/4}}\big(\widehat{a}_{\ga'} (\omega,u)e^{i\sgap u}+\widehat{b}_{\ga'} (\omega, u)e^{-i\sgap u}\big)\Big]d\ga'e^{\sgaa(u-z)}du\\
&+\frac{\ko{}}{2\sgaa}\int_{z\wedge L/\e}^{L/\e} \Big[ \sum_{l=1}^{\N{}}\frac{C^\omega_{\ga l}(u)}{\sqrt{\beta_l}}\big(\widehat{a}_l (\omega,u)e^{i\beta_l u}+\widehat{b}_l (\omega,u)e^{-i\beta_l u}\big)\\
&\quad \quad\quad+\int_{\xi}^{\ko{}}\frac{C^\omega_{\ga \ga'}(u)}{{\ga'}^{1/4}}\big(\widehat{a}_{\ga'} (\omega, u)e^{i\sgap u}+\widehat{b}_{\ga'} (\omega, u)e^{-i\sgap u}\big)\Big]d\ga'e^{\sgaa(z-u)}du
\end{split}
\]
$\forall z \in [0,+\infty)$.
Finally, for almost every $\ga \in (-\infty,-\xi)$ and $\forall z\in [0,+\infty)$,
\[ \tilde{p}_{\ga,0} (\omega,z)= \phi_{\ga}(\omega,x_0)e^{-\sgaa(z-L_S)}\textbf{1}_{(-S,-\xi)}(\ga). \]

We remark that $\Theta_\omega\circ\Pi_\omega(-\infty,-\xi)\big(\widehat{p}(\omega,.,.)\big)\in \mathcal{C}\big([0,+\infty),F\big)$ thanks to $\eqref{chypP2}$. Moreover, $\tilde{p}(\omega,.)\in \mathcal{C}\big([0,+\infty),F\big)$ since $\int_{-\infty}^{-\xi} \frac{A_\ga (\omega)}{\lvert\ga\rvert}d\ga<+\infty$, where $A_\ga (\omega)$ is defined by \eqref{coefagaP2} and satisfies \eqref{coefaga2P2}.
We can check that the norm of the operator $\Phi^\omega$ is bounded by
\[ \| \Phi^\omega\|\leq K(\xi,d)\sup_{z\in [0,L/\e]}\sup_{x\in[0,d]}\left\lvert V\left(x,z\right)\right\lvert. \]
Consequently, using \eqref{cg1}, $\lim_{\e\to 0}\mathbb{P}(Id- \se \,\Phi^\omega \text{ is invertible})=1$. Then, the condition ($Id- \se \,\Phi^\omega$ is invertible) is satisfied in the asymptotic $\e \to 0$. On the event ($Id- \se \,\Phi^\omega$ is invertible), we have 
\begin{equation}\label{presP2}\begin{split}\Theta_\omega\circ\Pi_\omega(-\infty,-\xi)\big(\widehat{p}(\omega,.,.)\big)&=\big(Id -\se\, \Phi^\omega\big)^{-1}(\se\, \tilde{p}(\omega, .)+ \tilde{p}_{0} (\omega,.))\\
&=\se\, \tilde{p}(\omega, .)+ \tilde{p}_{0} (\omega,.)+\sqrt{\e}\Phi^\omega(\tilde{p}_{0} (\omega,.))\\
& \quad+\sum_{j=1}^{+\infty}(\se\Phi^\omega)^j \big(\sqrt{\e}\tilde{p}(\omega,.)+\sqrt{\e}\Phi^\omega(\tilde{p}_{0} (\omega,.)) \big).
\end{split}\end{equation}
Moreover, 
\[\begin{split}
 \| \Theta_\omega\circ\Pi_\omega(-\infty,-\xi)&\big(\widehat{p}(\omega,.,.)\big)-\se\, \tilde{p}(\omega, .)- \tilde{p}_{0} (\omega,.)-\sqrt{\e}\Phi^\omega(\tilde{p}_{0} (\omega,.)) \|_{\infty,F}\\
 &\leq 2 \e \| \Phi^\omega \|\, \|\tilde{p}(\omega,.)\|_{\infty,F}+2 \e \| \Phi^\omega \|^2\,\|\tilde{p}_0(\omega,.)\|_{\infty,F}\\
 &\leq K(\xi,d) \,\e \sup_{z\in [0,L/\e]}\sup_{x\in[0,d]}\left\lvert V\left(x,z\right)\right\lvert^2\sup_{z\in[0,L /\e]}\|\widehat{a}(\omega,z) \|_{\esp_\xi}+\|\widehat{b}(\omega,z) \|_{\esp_\xi},
 \end{split} \]
and therefore
\[\begin{split}\Theta_\omega\circ\Pi_\omega(-\infty,-\xi)&\big(\widehat{p}(\omega,.,.)\big)=  \se \tilde{p}(\omega, .) + \tilde{p}_{0} (\omega,.)+\sqrt{\e}\Phi^\omega(\tilde{p}_{0} (\omega,.))\\
&+\mathcal{O}\Big(\e\sup_{z\in [0,L/\e]}\sup_{x\in[0,d]}\left\lvert V\left(x,z\right)\right\lvert^2 \sup_{z\in[0,L /\e]}\|\widehat{a}(\omega,z) \|_{\esp_\xi}+\|\widehat{b}(\omega,z) \|_{\esp_\xi}\Big)
\end{split}\]
in  $\mathcal{C}\big([0,+\infty),F\big)$.
Now, we consider 
\[\begin{split}
\tilde{p}_{\ga,2}&(\omega,z)=\frac{\ko{}}{2\sgaa}\int_0^{z\wedge L/\e} \Big[ \sum_{l=1}^{\N{}}\frac{C_{\ga l}(u)}{\sqrt{\beta_l}}\big(\widehat{a}_l (\omega,z\wedge L/\e)e^{i\beta_l u}+\widehat{b}_l (\omega,z\wedge L/\e)e^{-i\beta_l u}\big)\\
&\quad+\int_{\xi}^{\ko{}}\frac{C_{\ga \ga'}(u)}{{\ga'}^{1/4}}\big(\widehat{a}_{\ga'} (\omega,z\wedge L/\e)e^{i\sgap u}+\widehat{b}_{\ga'} (\omega, z\wedge L/\e)e^{-i\sgap u}\big)\Big]d\ga' e^{\sgaa(u-z)}du\\
&+\frac{\ko{}}{2\sgaa}\int_{z\wedge L/\e}^{L/\e} \Big[ \sum_{l=1}^{\N{}}\frac{C_{\ga l}(u)}{\sqrt{\beta_l}}\big(\widehat{a}_l (\omega,z\wedge L/\e)e^{i\beta_l u}+\widehat{b}_l (\omega,z\wedge L/\e)e^{-i\beta_l u}\big)\\
&\quad +\int_{\xi}^{\ko{}}\frac{C_{\ga \ga'}(u)}{{\ga'}^{1/4}}\big(\widehat{a}_{\ga'} (\omega, z\wedge L/\e)e^{i\sgap u}+\widehat{b}_{\ga'} (\omega, z\wedge L/\e)e^{-i\sgap u}\big)\Big]d\ga' e^{\sgaa(z-u)}du
\end{split}
\]
$\forall z \in [0,+\infty)$.
Using \eqref{mcP21}, \eqref{mcP22}, \eqref{mcP23}, \eqref{mcP24}, and \eqref{presP2}, we get 
\[\begin{split}
 \| \tilde{p}(\omega,.)- \tilde{p}_2(\omega, .) \|_{\infty,F}\leq K(\xi,d) \,\se \,&\sup_{z\in [0,L/\e]}\sup_{x\in[0,d]}\left\lvert V\left(x,z\right)\right\lvert^2\\
 &\times \Big( \sup_{z\in[0,L /\e]}\|\widehat{a}(\omega,z) \|_{\esp_\xi}+\|\widehat{b}(\omega,z) \|_{\esp_\xi} +\|\widehat{p}(\omega,z)\|_{F} \Big)
 \end{split}\]
and then
\[\begin{split} \Theta_\omega \circ\Pi_\omega &(-\infty,-\xi)\big(\widehat{p}(\omega,.,.)\big)= \se\, \tilde{p}_2(\omega, .) + \tilde{p}_{0} (\omega,.)+\sqrt{\e}\Phi^\omega(\tilde{p}_{0} (\omega,.))\\
&+\mathcal{O}\Big( \e \sup_{z\in [0,L/\e]}\sup_{x\in[0,d]}\left\lvert V\left(x,z\right)\right\lvert^2\sup_{z\in[0,L /\e]}\|\widehat{a}(\omega,z) \|_{\esp_\xi}+\|\widehat{b}(\omega,z) \|_{\esp_\xi} +\|\widehat{p}(\omega,z)\|_{F} \Big)\end{split}\]
in  $\mathcal{C}\big([0,+\infty),F\big)$.
Consequently, we can rewrite \eqref{mcP21}, \eqref{mcP22}, \eqref{mcP23}, and \eqref{mcP24} in a closed  form in $\esp_\xi \times \esp_\xi$. $\forall z\in[0,L/\e]$, we get
\begin{equation*}\begin{split}
\dz \widehat{a}&(\omega, z)=\se\, \textbf{H}^{aa}(\omega,z)\big(\widehat{a}(\omega, z)\big)+\se\, \textbf{H}^{ab}(\omega,z)\big(\widehat{b}(\omega, z)\big)+\sqrt{\e}\,\textbf{R}^{a,L_S}(\omega,z)\\
&\hspace{1.5cm}+\e\, \textbf{G}^{aa}(\omega,z)\big(\widehat{a}(\omega, z)\big)+\e\, \textbf{G}^{ab}(\omega,z)\big(\widehat{b}(\omega, z)\big)+\e\, \tilde{\textbf{R}}^{a,L_S}(\omega,z)\\
&+\mathcal{O}\Big( \e^{3/2}  \sup_{z\in [0,L/\e]}\sup_{x\in[0,d]}\left\lvert V\left(x,z\right)\right\lvert^2\sup_{z\in[0,L /\e]}\|\widehat{a}(\omega,z) \|_{\esp_\xi}+\|\widehat{b}(\omega,z) \|_{\esp_\xi} +\|\widehat{p}(\omega,z)\|_{F}\Big),
\end{split}\end{equation*}
\begin{equation*}\begin{split}
\dz \widehat{b}(\omega, z)&=\se\, \textbf{H}^{ba}(\omega,z)\big(\widehat{a}(\omega, z)\big)+\se\, \textbf{H}^{bb}(\omega,z)\big(\widehat{b}(\omega, z)\big)+\sqrt{\e}\,\textbf{R}^{b,L_S}(\omega,z)\\
&+\e\, \textbf{G}^{ba}(\omega,z)\big(\widehat{a}(\omega, z)\big)+\e\, \textbf{G}^{bb}(\omega,z)\big(\widehat{b}(\omega, z)\big)+\e \tilde{\textbf{R}}^{b,L_S}(\omega,z)\\
&+\mathcal{O}\Big(\e^{3/2} \sup_{z\in [0,L/\e]}\sup_{x\in[0,d]}\left\lvert V\left(x,z\right)\right\lvert^2\sup_{z\in[0,L /\e]}\|\widehat{a}(\omega,z) \|_{\esp_\xi}+\|\widehat{b}(\omega,z) \|_{\esp_\xi} +\|\widehat{p}(\omega,z)\|_{F}\Big).
\end{split}\end{equation*}
Let us recall that these equations hold on the event $\big(Id- \se \,\Phi^\omega \text{ is invertible}\big)$ which satisfies the condition $\lim_{\e\to 0}\mathbb{P}(Id- \se \,\Phi^\omega \text{ is invertible})=1$. In these equations, $\textbf{H}^{aa}(\omega,z)$, $ \textbf{H}^{ab}(\omega,z)$, $ \textbf{H}^{ba}(\omega,z)$, $\textbf{H}^{bb}(\omega,z)$, $\textbf{G}^{aa}(\omega,z)$, $ \textbf{G}^{ab}(\omega,z)$, $ \textbf{G}^{ba}(\omega,z)$  and $\textbf{G}^{bb}(\omega,z)$ are operators from $\esp_\xi$ to itself defined by:
\begin{equation}\label{haajP2}\begin{split}
\textbf{H}^{aa}_{j}(\omega,z)(y)=\overline{\textbf{H}^{bb}_{j}(\omega, z)}(y)&=\frac{i\ko{}}{2}\Big[\sum_{l=1}^{\N{}} \frac{C^\omega_{jl}(z)}{\sqrt{\beta_j(\omega) \beta_l(\omega)}}y_l e^{i(\beta_l(\omega) -\beta_j(\omega))z}\\
&+ \int_{\xi}^{\ko{}}\frac{C^\omega_{j\ga'}(z)}{\sqrt{\beta_j(\omega) \sgap}}y_{\ga'}e^{i(\sgap -\beta_j(\omega))z}d\ga' \Big],
\end{split}\end{equation}
\begin{equation}\label{haagaP2}\begin{split}
\textbf{H}^{aa}_{\ga}(\omega, z)(y)=\overline{\textbf{H}^{bb}_{\ga}(\omega, z)}(y)&=\frac{i\ko{}}{2}\Big[\sum_{l=1}^{\N{}} \frac{C^\omega_{\ga l}(z)}{\sqrt{\sga \beta_l(\omega)}}y_l e^{i(\beta_l(\omega) -\sga)z}\\
&+ \int_{\xi}^{\ko{}}\frac{C^\omega_{\ga \ga'}(z)}{\ga^{1/4}{\ga'}^{1/4}}y_{\ga'}e^{i(\sgap -\sga)z}d\ga' \Big],
\end{split}\end{equation}
\begin{equation}\label{habjP2}\begin{split}
\textbf{H}^{ab}_{j}(\omega,z)(y)=\overline{\textbf{H}^{ba}_{j}(\omega, z)}(y)&=\frac{i\ko{}}{2}\Big[\sum_{l=1}^{\N{}} \frac{C^\omega_{jl}(z)}{\sqrt{\beta_j(\omega) \beta_l(\omega)}}y_l e^{-i(\beta_l(\omega) +\beta_j(\omega))z}\\
&+ \int_{\xi}^{\ko{}}\frac{C^\omega_{j\ga'}(z)}{\sqrt{\beta_j(\omega) \sgap}}y_{\ga'}e^{-i(\sgap +\beta_j(\omega))z}d\ga' \Big],
\end{split}\end{equation}
\begin{equation}\label{habgaP2}\begin{split}
\textbf{H}^{ab}_{\ga}(\omega,z)(y)=\overline{\textbf{H}^{ba}_{\ga}(\omega,z)}(y)&=\frac{i\ko{}}{2}\Big[\sum_{l=1}^{\N{}} \frac{C^\omega_{\ga l}(z)}{\sqrt{\sga \beta_l(\omega)}}y_l e^{-i(\beta_l(\omega) + \sga)z}\\
&+ \int_{\xi}^{\ko{}}\frac{C_{\ga \ga'}(z)}{\ga^{1/4}{\ga'}^{1/4}}y_{\ga'}e^{-i(\sgap +\sga)z}d\ga' \Big],
\end{split}\end{equation}
\begin{equation}\label{gaajP2}\begin{split}
\textbf{G}^{aa}_{j}&(\omega,z)(y)=\overline{\textbf{G}^{bb}_{j}(\omega,z)}(y)=\\
&\frac{ik^4(\omega)}{4}\Big[\sum_{l=1}^{\N{}}\int_{-\infty}^{-\xi}\Big[ \int_{0}^{z}\frac{C^\omega_{j\ga'}(z)C^\omega_{\ga' l}(u)}{\sqrt{\beta_j(\omega)\lvert \ga'\rvert \beta_l(\omega)}}e^{i\beta_l(\omega) u -\sqrt{\lvert \ga'\rvert}(z-u)}du\\
&\quad \quad\quad+ \int_{z}^{L/\e}\frac{C^\omega_{j\ga'}(z)C^\omega_{\ga' l}(u)}{\sqrt{\beta_j(\omega)\lvert \ga'\rvert \beta_l(\omega)}}e^{i\beta_l(\omega) u -\sqrt{\lvert \ga'\rvert}(u-z)}du\Big]d\ga'e^{-i\beta_j(\omega) z}y_l\Big]\\
&+\frac{ik^4(\omega)}{4}\Big[\int_{\xi}^{\ko{}}\int_{-\infty}^{-\xi}\Big[ \int_{0}^{z}\frac{C^\omega_{j\ga'}(z)C^\omega_{\ga' \ga''}(u)}{\sqrt{\beta_j(\omega)\lvert \ga'\rvert \sqrt{\ga''}}}e^{i\sqrt{\ga''} u -\sqrt{\lvert \ga'\rvert}(z-u)}du\\
&\quad \quad\quad+ \int_{z}^{L/\e}\frac{C^\omega_{j\ga'}(z)C^\omega_{\ga' \ga''}(u)}{\sqrt{\beta_j(\omega)\lvert \ga'\rvert \sqrt{\ga''}}}e^{i\sqrt{\ga''} u -\sqrt{\lvert \ga'\rvert}(u-z)}du\Big]d\ga'e^{-i\beta_j(\omega) z}y_{\ga''}d\ga''\Big],
\end{split}\end{equation}
\begin{equation}\label{gaagaP2}\begin{split}
\textbf{G}^{aa}_{\ga}&(\omega,z)(y)=\overline{\textbf{G}^{bb}_{\ga}(\omega,z)}(y)=\\
&\frac{ik^4(\omega)}{4}\Big[\sum_{l=1}^{\N{}}\int_{-\infty}^{-\xi}\Big[ \int_{0}^{z}\frac{C^\omega_{\ga\ga'}(z)C^\omega_{\ga' l}(u)}{\sqrt{\sga\lvert \ga'\rvert \beta_l(\omega)}}e^{i\beta_l(\omega) u -\sqrt{\lvert \ga'\rvert}(z-u)}du\\
&\quad \quad\quad+ \int_{z}^{L/\e}\frac{C^\omega_{\ga\ga'}(z)C^\omega_{\ga' l}(u)}{\sqrt{\sga\lvert \ga'\rvert \beta_l(\omega)}}e^{i\beta_l(\omega) u -\sqrt{\lvert \ga'\rvert}(u-z)}du\Big]d\ga'e^{-i\sga z}y_l\Big]\\
&+\frac{ik^4(\omega)}{4}\Big[\int_{\xi}^{\ko{}}\int_{-\infty}^{-\xi}\Big[ \int_{0}^{z}\frac{C^\omega_{\ga\ga'}(z)C^\omega_{\ga' \ga''}(u)}{\sqrt{\sga\lvert \ga'\rvert \sqrt{\ga''}}}e^{i\sqrt{\ga''} u -\sqrt{\lvert \ga'\rvert}(z-u)}du\\
&\quad \quad\quad+ \int_{z}^{L/\e}\frac{C^\omega_{\ga\ga'}(z)C^\omega_{\ga' \ga''}(u)}{\sqrt{\sga\lvert \ga'\rvert \sqrt{\ga''}}}e^{i\sqrt{\ga''} u -\sqrt{\lvert \ga'\rvert}(u-z)}du\Big]d\ga'e^{-i\sga z}y_{\ga''}d\ga''\Big],
\end{split}\end{equation}
\begin{equation}\label{gabjP2}\begin{split}
\textbf{G}^{ab}_{j}&(\omega,z)(y)=\overline{\textbf{G}^{ba}_{j}(\omega,z)}(y)=\\
&\frac{ik^4(\omega)}{4}\Big[\sum_{l=1}^{\N{}}\int_{-\infty}^{-\xi}\Big[ \int_{0}^{z}\frac{C^\omega_{j\ga'}(z)C^\omega_{\ga' l}(u)}{\sqrt{\beta_j(\omega)\lvert \ga'\rvert \beta_l(\omega)}}e^{-i\beta_l(\omega) u -\sqrt{\lvert \ga'\rvert}(z-u)}du\\
&\quad \quad\quad+ \int_{z}^{L/\e}\frac{C^\omega_{j\ga'}(z)C^\omega_{\ga' l}(u)}{\sqrt{\beta_j(\omega)\lvert \ga'\rvert \beta_l(\omega)}}e^{-i\beta_l(\omega) u -\sqrt{\lvert \ga'\rvert}(u-z)}du\Big]d\ga'e^{-i\beta_j(\omega) z}y_l\Big]\\
&+\frac{ik^4(\omega)}{4}\Big[\int_{\xi}^{\ko{}}\int_{-\infty}^{-\xi}\Big[ \int_{0}^{z}\frac{C^\omega_{j\ga'}(z)C^\omega_{\ga' \ga''}(u)}{\sqrt{\beta_j(\omega)\lvert \ga'\rvert \sqrt{\ga''}}}e^{-i\sqrt{\ga''} u -\sqrt{\lvert \ga'\rvert}(z-u)}du\\
&\quad \quad\quad+ \int_{z}^{L/\e}\frac{C^\omega_{j\ga'}(z)C^\omega_{\ga' \ga''}(u)}{\sqrt{\beta_j(\omega)\lvert \ga'\rvert \sqrt{\ga''}}}e^{-i\sqrt{\ga''} u -\sqrt{\lvert \ga'\rvert}(u-z)}du\Big]d\ga'e^{-i\beta_j(\omega) z}y_{\ga''}d\ga''\Big],
\end{split}\end{equation}
\begin{equation}\label{gabgaP2}\begin{split}
\textbf{G}^{ab}_{\ga}&(\omega,z)(y)=\overline{\textbf{G}^{ba}_{\ga}(\omega,z)}(y)=\\
&\frac{ik^4(\omega)}{4}\Big[\sum_{l=1}^{\N{}}\int_{-\infty}^{-\xi}\Big[ \int_{0}^{z}\frac{C^\omega_{\ga\ga'}(z)C^\omega_{\ga' l}(u)}{\sqrt{\sga\lvert \ga'\rvert \beta_l(\omega)}}e^{-i\beta_l(\omega) u -\sqrt{\lvert \ga'\rvert}(z-u)}du\\
&\quad \quad\quad+ \int_{z}^{L/\e}\frac{C^\omega_{\ga\ga'}(z)C^\omega_{\ga' l}(u)}{\sqrt{\sga\lvert \ga'\rvert \beta_l(\omega)}}e^{-i\beta_l(\omega) u -\sqrt{\lvert \ga'\rvert}(u-z)}du\Big]d\ga'e^{-i\sga z}y_l\Big]\\
&+\frac{ik^4(\omega)}{4}\Big[\int_{\xi}^{\ko{}}\int_{-\infty}^{-\xi}\Big[ \int_{0}^{z}\frac{C^\omega_{\ga\ga'}(z)C^\omega_{\ga' \ga''}(u)}{\sqrt{\sga\lvert \ga'\rvert \sqrt{\ga''}}}e^{-i\sqrt{\ga''} u -\sqrt{\lvert \ga'\rvert}(z-u)}du\\
&\quad \quad\quad+ \int_{z}^{L/\e}\frac{C^\omega_{\ga\ga'}(z)C^\omega_{\ga' \ga''}(u)}{\sqrt{\sga\lvert \ga'\rvert \sqrt{\ga''}}}e^{-i\sqrt{\ga''} u -\sqrt{\lvert \ga'\rvert}(u-z)}du\Big]d\ga'e^{-i\sga z}y_{\ga''}d\ga''\Big].
\end{split}\end{equation}
The operators $\textbf{H}^{aa}(\omega,z)$ and $\textbf{H}^{ab}(\omega,z)$ represent the coupling between the propagating and the radiating modes with themselves, while the operators $\textbf{G}^{aa}(\omega,z)$ and $\textbf{G}^{ab}(\omega,z)$ represent the coupling between the evanescent modes with the propagating and the radiating modes. Moreover, $\textbf{R}^{a,L_S}(\omega,z)$, $\tilde{\textbf{R}}^{a,L_S}(\omega,z)$, $\textbf{R}^{b,L_S}(\omega,z)$, and $\tilde{\textbf{R}}^{b,L_S}(\omega,z)$ represent the influence of the evanescent modes produced by the source term on the propagating and the radiating modes. These terms are defined by 
\begin{equation}\label{RajP2}
\textbf{R}^{a,L_S}_j (\omega,z)=\overline{\textbf{R}^{b,L_S}_j (\omega,z)}=\frac{i \ko{}}{2}\int_{-S}^{-\xi}\frac{C^\omega_{j\ga'}(z)}{\sqrt{\beta_j(\omega)}}\phi_{\ga'}(\omega, x_0)e^{-\sqrt{\lvert \ga' \rvert}(z-L_S)}d\ga' e^{-i \beta_j(\omega) z},
\end{equation}
\begin{equation}\label{RagaP2}
\textbf{R}^{a,L_S}_\ga (\omega,z)=\overline{\textbf{R}^{b,L_S}_\ga (\omega,z)}=\frac{i \ko{}}{2}\int_{-S}^{-\xi}\frac{C^\omega_{\ga\ga' }(z)}{\gaa ^{1/4}}\phi_{\ga'}(\omega, x_0)e^{-\sqrt{\lvert \ga' \rvert}(z-L_S)}d\ga' e^{-i \sga z},
\end{equation}
\begin{equation}\label{RatjP2}\begin{split}
\tilde{\textbf{R}}&^{a,L_S}_j (\omega,z)=\overline{\tilde{\textbf{R}}^{b,L_S}_j (\omega,z)}=\\ 
&\frac{ik^{4}(\omega)}{4}\int_{-\infty}^{-\xi}\int_{-S}^{-\xi}\left[ \int_0^z \frac{C^\omega_{j\ga'}(z)C^\omega_{\ga' \ga''}(u)}{\sqrt{\beta_j(\omega) \lvert\ga' \rvert}}\phi_{\ga''}(\omega,x_0)e^{-\sqrt{\lvert \ga'' \rvert}(u-L_S)}e^{-\sqrt{\lvert \ga' \rvert}(z-u)}du\right.\\
&\quad +\left.\int_z^{L/\e} \frac{C^\omega_{j\ga'}(z)C^\omega_{\ga' \ga''}(u)}{\sqrt{\beta_j(\omega) \lvert\ga' \rvert}}\phi_{\ga''}(\omega,x_0)e^{-\sqrt{\lvert \ga'' \rvert}(u-L_S)}e^{-\sqrt{\lvert \ga' \rvert}(u-z)}du\right]d\ga''d\ga'\,e^{-i\beta_j(\omega) z},
\end{split}\end{equation}
\begin{equation}\label{RatgaP2}\begin{split}
\tilde{\textbf{R}}&^{a,L_S}_\ga (\omega,z)=\overline{\tilde{\textbf{R}}^{b,L_S}_\ga (\omega,z)}=\\ 
&\frac{ik^{4}(\omega)}{4}\int_{-\infty}^{-\xi}\int_{-S}^{-\xi}\left[ \int_0^z \frac{C^\omega_{\ga\ga'}(z)C^\omega_{\ga' \ga''}(u)}{\sqrt{\sqrt{\ga} \lvert\ga' \rvert}}\phi_{\ga''}(\omega,x_0)e^{-\sqrt{\lvert \ga'' \rvert}(u-L_S)}e^{-\sqrt{\lvert \ga' \rvert}(z-u)}du\right.\\
&\quad +\left.\int_z^{L/\e} \frac{C^\omega_{\ga\ga'}(z)C^\omega_{\ga' \ga''}(u)}{\sqrt{\sqrt{\ga} \lvert\ga' \rvert}}\phi_{\ga''}(\omega,x_0)e^{-\sqrt{\lvert \ga'' \rvert}(u-L_S)}e^{-\sqrt{\lvert \ga' \rvert}(u-z)}du\right]d\ga''d\ga'\,e^{-i\sqrt{\ga} z}.
\end{split}\end{equation}

\subsection{Forward Scattering Approximation}\label{fsaP2}

In this section we introduce the forward scattering approximation, which is widely used in the literature. In this approximation the coupling between forward- and backward-propagating modes is assumed to be negligible compared to the coupling between the forward-propagating modes. We refer to \cite{papa, gomez} for justifications on the validity of this approximation. 

The justification of this approximation is as follows. The coupling between a right-going propagating mode and a left-going propagating mode involves a coefficient of the form
\[ \int_{0}^{+\infty}\mathbb{E}[C^\omega_{jl}(0)C^\omega_{jl}(z)]\cos\big((\Bh{l}{}+\Bh{j}{})z\big)dz,\] 
and the coupling between two right-going propagating modes or two left-going propagating modes involves a coefficient of the form
\[ \int_{0}^{+\infty}\mathbb{E}[C^\omega_{jl}(0)C^\omega_{jl}(z)]\cos\big((\Bh{l}{}-\Bh{j}{})z\big)dz\]
$\forall (j,l) \in \big\{1,\dots,\N{}\big\}^2$. Therefore, if we assume that 
\[ \int_{0}^{+\infty}\mathbb{E}[C^\omega_{jl}(0)C^\omega_{jl}(z)]\cos\big((\Bh{l}{}+\Bh{j}{})z\big)dz=0 \quad \forall (j,l) \in \big\{1,\dots,\N{}\big\}^2,\]
then there is no coupling between right-going and left-going propagating modes, which justifies the forward scattering approximation, but there is still coupling between right-going propagating modes which will be described in Section \ref{coupledprocP2}.

In our context the operator $R(\omega)$, introduced in Section \ref{spectralP2}, has a continuous spectrum and it becomes technically complex to apply a limit theorem for the rescaled process $(\widehat{a}(\omega,z/\e),\widehat{b}(\omega,z/\e))$. The reason is the following. This process is not bounded and the stopping times which are the first exit times of closed balls are not lower semicontinuous for the topology of $\mathcal{C}([0,L],\esp_{\xi,w})$, where $\esp_{\xi,w}$ stands for $\esp_\xi$ equipped with the weak topology. In our context the continuous part $(\xi,\ko{})$ of the spectrum imposes us to use the norm $\|.\|_{\esp_{\xi}}$ to control some quantities.  Moreover, according to Theorem \ref{thasymptP21}, in which the energy of the limit process is not conserved, it seems not possible to show a limit theorem on $\mathcal{C}\big([0,L],(\esp_{\xi},\|.\|_{\esp_{\xi}})\big)$ in view of \eqref{conspropP2}. In \cite{book} and \cite{papa} there is a finite number of propagating modes, so that the weak topology and the strong topology are the same. In \cite{gomez}  the number of propagating modes increases as $\e$ goes to $0$. However, in this last case, the problem can be corrected by considering the first exit times of a closed ball related to the weak topology and by considering the process in an appropriate finite-dimensional dual space.  

In our context if we forget these technical problems, according to \cite{book, papa} the forward scattering approximation should be valid in the asymptotic $\e\to 0$ under the assumption that the power spectral density of the process $V$, i.e. the Fourier transform of its $z$-autocorrelation function, possesses a cut-off wavenumber. In other words, we can consider the case where
\[ \int_{0}^{+\infty}\mathbb{E}[C^\omega_{jl}(0)C^\omega_{jl}(z)]\cos\big((\Bh{l}{}+\Bh{j}{})z\big)dz=0 \quad \forall (j,l) \in \big\{1,\dots,\N{}\big\}^2.\]
Let us remark that the continuous part $(0,\ko{})$ of the spectrum, which corresponds to the radiating modes, does not play any role in the previous assumption. The reason is that the radiating part of the process plays no role in the coupling mechanism as we can see in Theorems \ref{thasymptP21} and \ref{thasympP2} below and therefore remains constant.

Finally, we shall consider the simplified equation on $[0,L/\e]$,
\begin{equation*}\begin{split}
\dz \widehat{a}(\omega, z)&=\se\, \textbf{H}^{aa}(\omega,z)\left(\widehat{a}(\omega, z)\right)+\sqrt{\e}\,\textbf{R}^{a,L_S}(\omega,z)\\&
+\e\, \textbf{G}^{aa}(\omega,z)\left(\widehat{a}(\omega, z)\right)+\e\,\tilde{\textbf{R}}^{a,L_S}(\omega,z)\\
&+\mathcal{O}\Big( \e^{3/2} \sup_{z\in [0,L/\e]}\sup_{x\in[0,d]}\left\lvert V\left(x,z\right)\right\lvert^2 \sup_{z\in [0,L/\e]}\| \widehat{a}(\omega,z) \|_{\esp_\xi}+\|\widehat{p}(\omega,z)\|_{F}\Big)
\end{split}\end{equation*}
in $\esp_\xi$. 
We shall see in Section \ref{coupledprocP2}, under the forward scattering approximation, that 
\[ \lim_{\e\to 0}\mathbb{P}\left(L^\e \leq L\right) =0\quad \forall\, L>0,\]
where
\[L^\e=\inf\Big(L>0, \quad \sup_{z\in[0,L/\e]}\|\widehat{a}(\omega,z) \|^2_{\esp_\xi}+\| \widehat{p}(\omega,z) \|^2 _{F} \geq \frac{1}{\sqrt{\e}} \Big).\]
Consequently, we can show that $\forall \,\eta>0$
\[ \lim_{\e\to0}\mathbb{P}\Big(\sup_{z\in[0,L/\e]}\left\lvert   \|\widehat{a}(\omega,z)\|^2_{\esp_\xi}- \|\widehat{a}_0(\omega)\|^2_{\esp_\xi} \right\rvert>\eta  \Big)=0.\]
This result means that the local energy flux for the propagating and the radiating modes is conserved in the asymptotic $\e\to 0$.

\section{Coupled Mode Processes}\label{coupledprocP2} 

In this section, we study the asymptotic behavior, as $\e\to 0$ in first and $\xi \to 0$ in second, of the statistical properties of the coupling mechanism in terms of a diffusion process. 

Let us define the rescaled process
\[\widehat{a} ^\e (\omega ,z)=\widehat{a}\left(\omega,\frac{z}{\e}\right)\quad \forall z\in[0,L].\]
This scaling corresponds to the size of the random section $[0,L/ \e]$. This process satisfies the rescaled coupled mode equations on $[0,L]$
\begin{equation}\label{mccP22}\begin{split}
\dz \widehat{a}^{\e}(\omega, z)&=\frac{1}{\se} \textbf{H}^{aa}\left(\omega,\frac{z}{\e}\right)\left(\widehat{a}^\e(\omega, z)\right)+\frac{1}{\sqrt{\e}}\textbf{R}^{a,L_S}\left(\omega,\frac{z}{\e}\right)\\
&+ \textbf{G}^{aa}\left(\omega,\frac{z}{\e}\right)\left(\widehat{a}^{\e}(\omega, z)\right)+\tilde{\textbf{R}}^{a,L_S}\left(\omega,\frac{z}{\e}\right)\\
&+\mathcal{O}\Big( \se\sup_{z\in [0,L/\e]}\sup_{x\in[0,d]}\left\lvert V\left(x,z\right)\right\lvert^2 \sup_{z\in[0,L]}\|\widehat{a}^\e (\omega,z) \|_{\esp_\xi}+\| \widehat{p}(\omega,z/\e) \| _{F}\Big)
\end{split}\end{equation}
in $\esp_\xi$, with the initial condition $\widehat{a}^\e(\omega, 0)=\widehat{a}_0(\omega)$. 
We shall see that under the forward scattering approximation the condition $L^\e>L$ is readily fulfilled in the asymptotic $\e$ goes to $0$. 

\begin{prop} \label{prop1P2}
$\forall L>0$,
\[\lim_{\e\to 0}\mathbb{P}\left(L^\e\leq L\right) =0,\]
where
\[L^\e=\inf\Big(L>0, \quad \sup_{z\in[0,L/\e]}\|\widehat{a}(\omega,z) \|^2_{\esp_\xi}+\| \widehat{p}(\omega,z) \|^2 _{F} \geq \frac{1}{\sqrt{\e}} \Big),\]
and 
\[\lim_{M\to+\infty}\overline{\lim_{\e \to 0}}\,\,\mathbb{P}\left( \sup_{z\in[0,L]}\|\widehat{a}^\e(\omega, z)\|^2 _{\esp_{\xi}}\geq M \right)=0.\]
\end{prop}
This result means that the amplitude $\widehat{a}^\e(\omega, z)$ is asymptotically uniformly bounded in the limit $\e\to0$ on $[0,L]$. More precisely, according to Section \ref{globalP2}, we have $\forall \,\eta>0$
\[ \lim_{\e\to0}\mathbb{P}\Big(\sup_{z\in[0,L]}\left\lvert   \|\widehat{a}^\e(\omega,z)\|^2_{\esp_\xi}- \|\widehat{a}_0(\omega)\|^2_{\esp_\xi} \right\rvert>\eta  \Big)=0,\]
that is the local energy flux for the propagating and the radiating modes is conserved in the asymptotic $\e\to 0$.

\begin{preuve}
Using Gronwall's inequality, $\forall L>0$ we get 
\[\lim_{M\to+\infty}\overline{\lim_{\e \to 0}}\,\,\mathbb{P}\left( \sup_{z\in[0,L]}\|\widehat{a}^\e(\omega, z)\|^2 _{\esp_{\xi}}\geq M,\, L\leq L^\e  \right)=0.\]
This result means that the process $\widehat{a}^\e(\omega, .)$ is asymptotically uniformly bounded on $[0,L]$ and then $L^\e $ is large compared to $L$ in the asymptotic $\e\to0$. In fact, $\forall L>0$ and $\forall M>0$
\[\begin{split}
\mathbb{P}\left(L^\e \leq L\right)&\leq \mathbb{P}\left(L^\e \leq L,\,\, \sup_{z\in[0,L\wedge L^\e ]}\|\widehat{a}^\e(\omega, z)\|^2 _{\esp_{\xi}}\leq M \right)\\
&\quad+\mathbb{P}\left( \sup_{z\in[0, L\wedge L^\e ]}\|\widehat{a}^\e(\omega, z)\|^2 _{\esp_{\xi}}\geq M \right).
\end{split}\]
Moreover, 
\[\mathbb{P}\left(L^\e \leq L,\,\, \sup_{z\in[0,L\wedge L^\e ]}\|\widehat{a}^\e(\omega, z)\|^2_{\esp_{\xi}}\leq M \right)=0\]
for $\e$ small enough, since for $L^\e \leq L$
\[\begin{split}
\e^{-1/2} \leq\sup_{z\in[0,L^\e ]}& \|\widehat{a}^\e(\omega, z)\|^2 _{\esp_{\xi}}+\|\widehat{p}(\omega,.)\|^2 _{F}\\
&\leq M+K(\xi,d)\e \sup_{z\in [0,L/\e]}\sup_{x\in[0,d]}\left\lvert V\left(x,z\right)\right\lvert^2 M+2\|\tilde{p}_0(\omega,.)\|^2_{\infty,F}\end{split}\]
according to \eqref{presP2}. 
$\blacksquare$
\end{preuve}

Let us introduce $\widehat{a}^\e _{1}(\omega,.)$ the unique solution of the differential equation on $[0,L]$
\begin{equation}\label{mccP23}
\dz \widehat{a}^\e _1(\omega, z)=\frac{1}{\se}\textbf{H}^{aa}\left(\omega,\frac{z}{\e}\right)\left(\widehat{a}^\e _1(\omega, z)\right)+\textbf{G}^{aa}\left(\omega,\frac{z}{\e}\right)\left(\widehat{a}^\e_1(\omega, z)\right)
\end{equation}
in $\esp_\xi$, with initial condition $\widehat{a}^\e _1(\omega, 0)=\widehat{a}_0(\omega)$.
Using Gronwall's inequality and \eqref{cg2} we can state that
\[ \lim_{M\to+\infty}\overline{\lim_{\e \to0}}\,\,\mathbb{P}\left( \sup_{z\in[0,L]}\|\widehat{a}^\e_1(\omega, z)\|_{\esp_{\xi}}\geq M \right)=0. \] 
The relation between the solution of the full system \eqref{mccP22} and the one of the simplified system \eqref{mccP23} is given by the following proposition.
\begin{prop}\label{propevlossP2}
\[\forall \eta>0 \text{ and } \forall \mu > 0, \quad\lim_{\e \to 0}  \mathbb{P}\left(\sup_{z\in[\mu,L]} \|\widehat{a}^\e(\omega,z)-\widehat{a}^\e _1(\omega,z)\|_{\esp_{\xi}}>\eta\right)=0.\]
\end{prop}
Proposition \ref{propevlossP2} means that the information about the evanescent part of the source profile is lost in the asymptotic $\e$ goes to $0$. In fact, the coupling mechanism described by the system  \eqref{mccP22} implies that the information about the evanescent part of the source profile is transmitted to the propagating modes through the coefficients $\textbf{R}^{a,L_S}(\omega,z)$ and $\tilde{\textbf{R}}^{a,L_S}(\omega,z)$ defined by \eqref{RajP2}, \eqref{RagaP2}, \eqref{RatjP2} and \eqref{RatgaP2}. In these expressions we have the term $\phi_{\ga'}(\omega,x)e^{-\sqrt{\lvert \ga'\rvert} (z-L_S)}$ which comes from the right-hand side of \eqref{pgammaP2} and which is the form of evanescent modes without a random perturbation. This term is responsible for the loss of information about the evanescent part of the source profile because of its exponentially decreasing behavior.
\begin{preuve}
We begin by proving that $\forall L>0$, $\forall \eta >0$ and $\forall \mu >0$
\[ \lim_{\e \to 0}\mathbb{P}\left( \sup_{z\in[\mu,L]}\|\widehat{a}^\e(\omega, z)-\widehat{a}^\e _1(\omega, z)\|^2_{\esp_{\xi}}>\eta,\,L\leq L^\e \right)=0.\]
In fact, $\textbf{R}^{a,L_S}(\omega,z)$ decreases exponentially fast with the propagation distance.
Moreover, $\tilde{\textbf{R}}^{a,L_S}(\omega,z)$ can be treated as $\textbf{G}^{aa}$ in the proof of Theorem \ref{thasymptP21} because $e^{-\sqrt{\lvert \ga' \rvert} (u-L_S)}$ cannot be compensated by $e^{-i\Bh{j}{}z}$ nor by $e^{-i\sga z}$. 
Moreover, using Proposition \ref{prop1P2} we get the result.
$\blacksquare$
\end{preuve}

Finally, we introduce the transfer operator $\textbf{T}^{\xi,\e} (\omega,z)$ from $\esp_\xi$ to itself, which is the unique operator solution of the differential equation
\begin{equation}\label{dercTP2}
\dz \textbf{T}^{\xi,\e} (\omega,z)=\frac{1}{\sqrt{\e}} \textbf{H}^{aa} \left(\omega,\frac{z}{\e}\right)
\textbf{T}^{\xi,\e} (\omega,z)+ \textbf{G}^{aa} \left(\omega,\frac{z}{\e}\right)
\textbf{T}^{\xi,\e} (\omega,z)
\end{equation}
with $\textbf{T}^{\xi,\e} (\omega,0)=  Id$.
Then,
\[\forall z\in [0,L],\quad \widehat{a}_1(\omega,z)=\textbf{T}^{\xi,\e}(\omega,z)(\widehat{a}_0 (\omega)),\]
and we get the following result.
\begin{prop}
\[\forall \eta>0 \text{ and } \forall \mu > 0,\quad \lim_{\e \to 0}\mathbb{P}\left( \sup_{z\in[\mu,L]}\|\widehat{a}^\e(\omega, z)-\emph{\textbf{T}}^{\xi,\e}(\omega,z)(\widehat{a}_0 (\omega))\|^2_{\esp_{\xi}} > \eta \right)=0.\]
\end{prop}

\subsection{Limit Theorem}\label{thasyptoticP2}

This section presents the basic theoretical results of this paper. In \cite{papa} and \cite{papanicolaou}, the authors used the limit theorem stated in \cite{papakohler} since the number of propagating modes was fixed. However, in our configuration, in addition to the $\N{}$-discrete propagating modes the wave field consists of a continuum of radiating modes. The two following results are based on a diffusion-approximation result for the solution of an ordinary differential equation with random coefficients. This result is an extension of that stated in \cite{papakohler} to the case of processes with values in a Hilbert space.

\begin{thm}\label{thasymptP21}
$\forall L>0$ and $\forall y \in \esp_{\xi}=\mathbb{C}^{\N{}}\times L^2(\xi,\ko{})$, the family $\big(\emph{\textbf{T}}^{\xi,\e} (\omega,.)(y)\big)_{\e\in(0,1)}$, solution of the differential equation \eqref{dercTP2}, converges in distribution on 
$\mathcal{C}([0,L], \esp_{\xi,w})$ as $\e \to 0$ to a limit denoted by $ \emph{\textbf{T}}^{\xi} ( \omega,.)(y)$. Here $\esp_{\xi,w}$ stands for the Hilbert space $\esp_{\xi}$ equipped with the weak topology. This limit is the unique diffusion process on $\esp_{\xi}$, starting from $y$, associated to the infinitesimal generator
\[\mathcal{L}^\omega_\xi=\mathcal{L}^\omega_1+\mathcal{L}^\omega_{2,\xi}+\mathcal{L}^\omega_{3,\xi},\]
where
\begin{equation*}\begin{split}
\mathcal{L}^\omega_{1}&=\frac{1}{2}\sum_{\substack{j,l =1\\ j\not=l}}^{\N{}}\Gamma^{c}_{jl}(\omega)\left(T_j\overline{T_j}\partial_{T_l}\partial_{\overline{T_l}}+T_l\overline{T_l}\partial_{T_j}\partial_{\overline{T_j}}-T_jT_l\partial_{T_j}\partial_{T_l}-\overline{T_j}\overline{T_l}\partial_{\overline{T_j}}\partial_{\overline{T_l}}\right)\\
&+\frac{1}{2}\sum_{j,l =1}^{\N{}}\Gamma^{1}_{jl}(\omega)\left(T_j\overline{T_l}\partial_{T_j}\partial_{\overline{T_l}}+\overline{T_j}T_l\partial_{\overline{T_j}}\partial_{T_l}-T_jT_l\partial_{T_j}\partial_{T_l}-\overline{T_j}\overline{T_l}\partial_{\overline{T_j}}\partial_{\overline{T_l}}\right)\\
&+\frac{1}{2}\sum_{j=1}^{\N{}} \left( \Gamma^{c}_{jj}(\omega) -\Gamma^{1}_{jj}(\omega)\right)\left(T_j \partial_{T_j}+\overline{T_j}\partial_{\overline{T_j}}\right)+\frac{i}{2}\sum_{j=1}^{\N{}} \Gamma^{s}_{jj}(\omega) \left(T_j \partial_{T_j}-\overline{T_j}\partial_{\overline{T_j}}\right),
\end{split}\end{equation*}
and
\begin{equation*}\begin{split}
\mathcal{L}^\omega_{2,\xi}&=-\frac{1}{2}\sum_{j=1}^{\N{}} \left(\Lambda^{c,\xi}_j(\omega) +i\Lambda^{s,\xi}_j(\omega)\right)    T_j \partial_{T_j}+ \left(\Lambda^{c,\xi}_j(\omega)  - i\Lambda^{s,\xi}_j(\omega)\right)\overline{T_j}\partial_{\overline{T_j}},\\
\mathcal{L}^\omega_{3,\xi}&=i\sum_{j=1}^{\N{}} \kappa^\xi_j(\omega) \left(T_j \partial_{T_j}-\overline{T_j}\partial_{\overline{T_j}}\right).
\end{split}\end{equation*}
\end{thm}
Here, we have considered the classical complex derivative with the following notation: If $v=v_1 +iv_2$, then $\partial_v =\frac{1}{2}\left(\partial_{v_1} -i \partial_{v_2} \right)$ and $\partial_{\overline{v}}=\frac{1}{2}\left(\partial_{v_1} +i \partial_{v_2} \right)$. We have used the following notations. $\forall (j,l)\in\big\{1,\dots,\N{} \big\}^{2}$ and $j\not=l$
\begin{equation*}\begin{split}
\Gamma^{c}_{jl}(\omega)&= \frac{k^4(\omega)}{2\Bh{j}{}\Bh{l}{}}\int_{0}^{+\infty}\mathbb{E}\big[C^\omega_{jl}(0)C^\omega_{jl}(z)\big]\cos\big((\Bh{l}{}-\Bh{j}{})z\big)dz,\\
\Gamma^{c}_{jj}(\omega)&=-\sum_{\substack{l=1\\l\not=j}}^{\N{}}\Gamma^{c}_{jl}(\omega),\\
\Gamma^{s}_{jl}(\omega)&= \frac{k^4(\omega)}{2\Bh{j}{}\Bh{l}{}}\int_{0}^{+\infty}\mathbb{E}\big[C^\omega_{jl}(0)C^\omega_{jl}(z)\big]\sin\big((\Bh{l}{}-\Bh{j}{})z\big)dz, \\
\Gamma^{s}_{jj}(\omega)&=-\sum_{\substack{l=1\\ l\not= j}}^{\N{}}\Gamma^{s}_{jl}(\omega),
\end{split}\end{equation*}
and 
$\forall (j,l)\in\big\{1,\dots,\N{} \big\}^{2}$,
\[\begin{split}
\Gamma^{1}_{jl}(\omega)&= \frac{k^4(\omega)}{2\Bh{j}{}\Bh{l}{}}\int_{0}^{+\infty}\mathbb{E}\big[C^\omega_{jj}(0)C^\omega_{ll}(z)\big]dz,\\
\Lambda^{c,\xi}_{j}(\omega)&= \int_{\xi}^{\ko{}}\frac{k^4 (\omega)}{2\sgap \Bh{j}{}}\int_{0}^{+\infty}\mathbb{E}\big[C^\omega_{j\ga'}(0)C^\omega_{j\ga'}(z)\big]\cos\big((\sgap-\Bh{j}{})z\big)dzd\ga', \\
\Lambda^{s,\xi}_{j}(\omega)&= \int_{\xi}^{\ko{}} \frac{k^4 (\omega)}{2\sgap \Bh{j}{}}\int_{0}^{+\infty}\mathbb{E}\big[C^\omega_{j\ga'}(0)C^\omega_{j\ga'}(z)\big]\sin\big((\sgap-\Bh{j}{})z\big)dzd\ga',\\
\kappa^{\xi}_j (\omega)&=\int_{-\infty}^{-\xi}\frac{k^{4}(\omega)}{2\Bh{j}{} \sqrt{\lvert\ga'\rvert}} \int_{0}^{+\infty}\mathbb{E}\big[C^\omega_{j\ga'}(0)C^\omega_{j\ga'}(z)\big]\cos\big(\Bh{j}{}z\big)e^{-\sqrt{\lvert\ga'\rvert}z}dzd\ga' .
\end{split}\]
The coupling coefficients $C^\omega(z)$ are defined by \eqref{coefcouplVP2}.
We get the following result in the asymptotic $\xi\to 0$. 
\begin{thm}\label{thasympP2}
$\forall L>0$ and $\forall y \in\esp_0=\mathbb{C}^{\N{}}\times L^2(0,\ko{})$, the family $\big(\emph{\textbf{T}}^{\xi} (\omega,.)(y)\big)_{\xi\in(0,1)}$ converges in distribution on 
$\mathcal{C}([0,L], (\esp_{0},\|.\|_{\esp_{0}}))$ as $\xi \to 0$ to a limit denoted by $ \emph{\textbf{T}}^0( \omega,.)(y)$. This limit is the unique diffusion process on $\esp_{0}$, starting from $y$, associated to the infinitesimal generator
\begin{equation*}\mathcal{L}^\omega=\mathcal{L}^\omega_1+\mathcal{L}^\omega_2+\mathcal{L}^\omega_3,\end{equation*}
where
\begin{equation*}\begin{split}
\mathcal{L}^\omega_{2}&=-\frac{1}{2}\sum_{j=1}^{\N{}} \left(\Lambda^{c}_j(\omega) +i\Lambda^{s}_j(\omega)\right)    T_j \partial_{T_j}+ \left(\Lambda^{c}_j(\omega)  - i\Lambda^{s}_j(\omega)\right)\overline{T_j}\partial_{\overline{T_j}},\\
\mathcal{L}^\omega_{3}&=i\sum_{j=1}^{\N{}} \kappa_j(\omega) \left(T_j \partial_{T_j}-\overline{T_j}\partial_{\overline{T_j}}\right).
\end{split}\end{equation*}
\end{thm}
Here, we have $\forall j\in\big\{1,\dots,\N{} \big\}$
\[
\Lambda^{c}_{j}(\omega)=\lim_{\xi\to 0}\Lambda^{c,\xi}_{j}(\omega),\quad
\Lambda^{s}_{j}(\omega)=\lim_{\xi\to 0}\Lambda^{s,\xi}_{j}(\omega),\quad
\kappa _j (\omega)=\lim_{\xi \to 0}\kappa ^{\xi}_j (\omega).
\]

Theorems \ref{thasymptP21} and \ref{thasympP2} describe the asymptotic behavior, as $\e\to 0$ first and $\xi \to 0$ second, of the statistical properties of the transfer operator $\textbf{T}^{\xi,\e}(\omega,L)$, in terms of a diffusion process. In the appendix we give the proofs of Theorems \ref{thasymptP21} and \ref{thasympP2}, which are based on a martingale approach using the perturbed-test-function method. In a first step we show the tightness of the processes, and in a second step we characterize all subsequence limits by mean of a well-posed martingale problem in a Hilbert space. 

The infinitesimal generator $\mathcal{L}^\omega$ is composed of three parts which represent different behaviors on the diffusion process. We can remark that the infinitesimal generator depends only on the $\N{}$-discrete coordinates. The first operator $\mathcal{L}^\omega_1$ describes the coupling between the $\N{}$-propagating modes. This part is of the form of the infinitesimal generator obtained in \cite{book,papa}, and the total energy is conserved.  The second operator $\mathcal{L}^\omega_2$ describes the coupling between the propagating modes with the radiating modes. This part implies a mode-dependent and frequency-dependent attenuation on the $\N{}$-propagating modes that we study in Section \ref{expdecP2}, and a mode-dependent and frequency-dependent phase modulation. The third operator $\mathcal{L}^\omega_3$ describes the coupling between the propagating and the evanescent modes, and implies a mode-dependent and frequency-dependent phase modulation. The purely imaginary part of the operator $\mathcal{L}^\omega$ does not remove energy from the propagating modes but gives an effective dispersion.    

Moreover, let us remark that the convergence in Theorem \ref{thasymptP21} holds also on $\mathcal{C}([0,L], (\esp_{\xi},\|.\|_{\esp_{\xi}}))$ for the $\N{}$-discrete propagating mode amplitudes.

\subsection{Mean Mode Amplitudes}

In this section we study the asymptotic mean mode amplitudes. From Theorem \ref{thasympP2}, we get the following result about the mean mode amplitudes.
\begin{prop}
 $\forall y\in \esp_0$, $\forall z\in[0,L]$, $\forall j\in\big\{1,\dots,\N{} \big\}$ 
\begin{equation}\label{asymptmmatrP2}\begin{split}
\lim_{\xi\to0}\lim_{\e\to0}\mathbb{E}&\Big[\emph{\textbf{T}}^{\xi,\e}_j( \omega,z)(y)\Big]=
\mathbb{E}\Big[\emph{\textbf{T}}^0_j( \omega,z)(y)\Big]\\
&=\exp\left[\left(\frac{\Gamma^{c}_{jj}(\omega)-\Gamma^1 _{jj}(\omega) - \Lambda^{c}_j (\omega)}{2}\right)z+ i\left(\frac{\Gamma^{s}_{jj}(\omega)-\Lambda^{s}_{jj}(\omega)}{2}+k_j(\omega)\right)z \right]y_{j} (\omega).\end{split}\end{equation}
\end{prop}
First, let us remark that the mean amplitude of the radiating part remains constant on $L^2(0,\ko{})$. Second, $\forall j\in \big\{1,\dots,\N{}\big\}$, the coefficient $(\Gamma^1 _{jj}(\omega)+ \Lambda^{c}_j (\omega)-\Gamma^{c}_{jj}(\omega))/2$ is nonnegative. In fact, for $(j,l)\in\{1,\dots,\N{}\}^2$ such that $j\not=l$, $\Gamma^{c}_{jl}(\omega)$ and $\Gamma^{1}_{jj}(\omega)$ are nonnegative because they are proportional to the power spectral density of $C^\omega_{jl}$ and $C^\omega_{jj}$ at $\Bh{l}{}-\Bh{j}{}$ and $0$ frequencies. Therefore, $-\Gamma^c_{jj}(\omega)$ is also nonnegative. Moreover, $\Lambda^{c}_j (\omega)$ is also nonnegative because it is proportional to the integral over $(0,\ko{})$ of the power spectral density of $C^\omega_{j\ga}$ at $\sga-\Bh{j}{}$ frequency.

The exponential decay rate for the mean $j$th-propagating mode is given by
\[\Big\lvert\mathbb{E}\Big[ \textbf{T}^0_j( \omega,L)(y)\Big]\Big\rvert= \big\lvert y_j \big\rvert \exp\left[-\left(\frac{\Gamma^{1}_{jj}(\omega)-\Gamma^c _{jj}(\omega)+ \Lambda^{c}_j (\omega)}{2}\right) L\right], \]
which depends on the effective coupling between the propagating modes, and the coupling between the propagating and the radiating modes. This exponential decay corresponds to a loss of coherence of the transmitted field.

\section{Coupled Power Equations}\label{sect5P2}

\begin{figure}\begin{center}
\includegraphics*[scale=0.65]{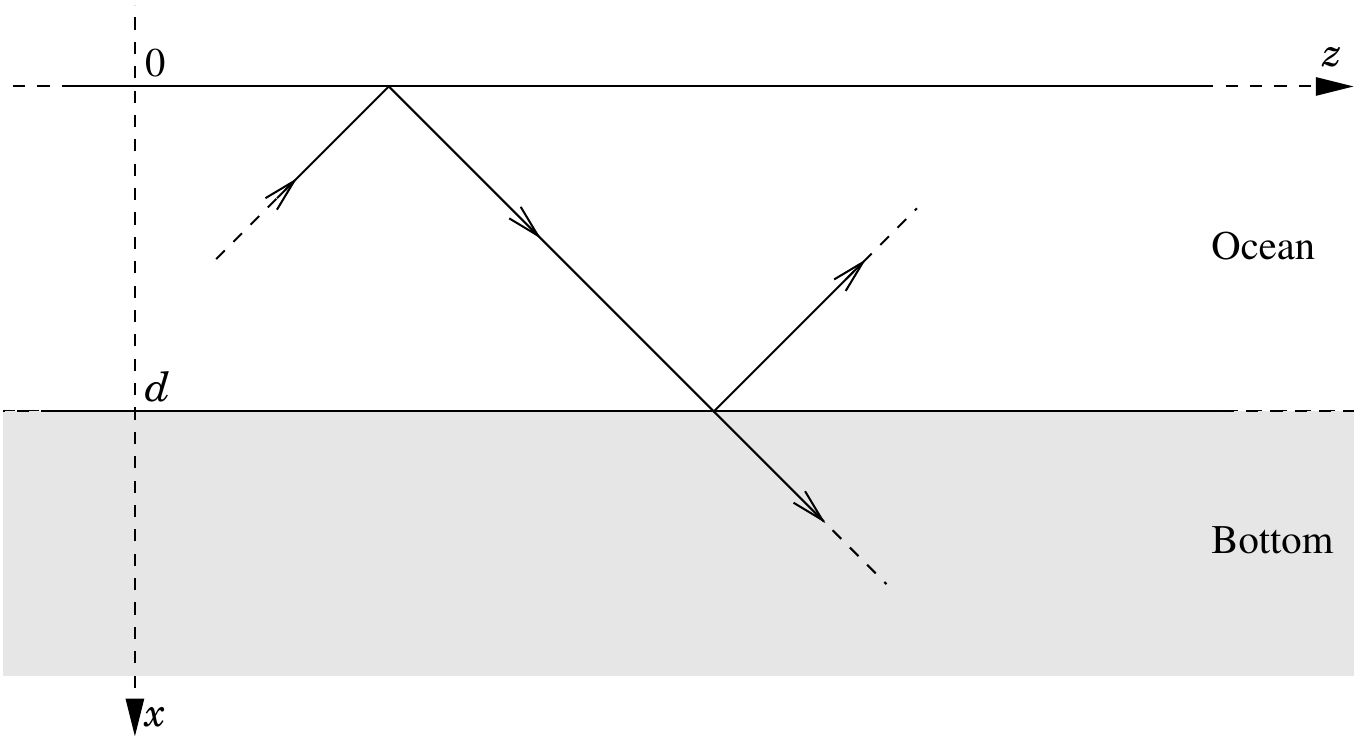}
\end{center}
\caption{\label{perterad1P2} Illustration of the radiative loss in the shallow-water random waveguide model.}
\end{figure}

This section is devoted to the analysis of the asymptotic mean mode powers of the propagating modes. More precisely, we study the asymptotic effects of the coupling between the propagating modes with the radiating modes. Let 
\begin{equation}\label{mmpP2}
\mathcal{T}_{j}^l(\omega,z)=\lim_{\xi\to 0}\lim_{\e\to 0}\mathbb{E}\Big[\big\lvert \textbf{T}^{\xi,\e}_j(\omega,L)(y^l) \big\rvert^2\Big]=\mathbb{E}\Big[ \big\lvert \textbf{T}^0_j( \omega,z)(y^l) \big\rvert^2 \Big],
\end{equation}
be the asymptotic mean mode power of the $j$th propagating modes. $\mathcal{T}^l_j(\omega,L)$ is the expected power of the $j$th propagating mode at the propagation distance $z=L$. Here $y^l \in \esp_0$ is defined by $y^l_j=\delta _{jl}$ and $y^l_\ga=0$ for $\ga \in(0,\ko{})$, and where $\delta_{jl}$ is the Kronecker symbol. The initial condition $y^l$ means that an impulse equal to one charges only the $l$th propagating mode. From Theorem \ref{thasympP2}, we have the coupled power equations:
\begin{equation}\label{cpeP2}
\dz \mathcal{T}_{j}^l(\omega,z)=-\Lambda_j ^c (\omega)\mathcal{T}_{j}^l(\omega,z)+\sum_{\substack{n=1\\n\not=j}}^{\N{}}\Gamma_{nj}^{c}(\omega)\left(\mathcal{T}_{n}^l(\omega,z)-\mathcal{T}_{j}^l(\omega,z)\right),
\end{equation}
with initial conditions $\mathcal{T}_{j}^l(\omega,0)=\delta_{jl}$. These equations describe the transfer of energy between the propagating modes and $\Gamma^{c}(\omega)$ is the energy transport matrix.  In our context, we also have the coefficients $\Lambda^{c}_j(\omega)$ given by the coupling between the propagating modes with the radiating modes. These coefficients, defined in Theorem \ref{thasympP2}, are responsible for the radiative loss of energy in the ocean bottom (see Figure \ref{perterad1P2}). This loss of energy is described more precisely in the following section.

\subsection{Exponential Decay of the Propagating Modes Energy}\label{expdecP2}

In this section, we assume that at least one of the coefficients $\Lambda^{c}(\omega)$ is positive. 
With this assumption, we show that the total energy carried by the propagating modes decays exponentially with the size $L$ of the random section. 
In the opposite situation, that is when there is no radiative loss $\Lambda^{c}(\omega)=0$, it has been shown in \cite{book} and \cite[Chapter 20] {papa} that the energy of the propagating modes is conserved and for large $L$ the asymptotic distribution of the energy becomes uniform over the propagating modes.

Let us defined 
\[  \mathcal{S}^{\N{}}_{+}=\left\{ X\in\mathbb{R}^{\N{}}, \quad X_j\geq 0\quad\forall j\in\{1,\dots,\N{}\}  \text{ and }\|X\|^2_{2,\mathbb{R}^{\N{}}} =\big<X,X\big>_{\mathbb{R}^{\N{}}}=1 \right\}\]
with $\big< X,Y\big>_{\mathbb{R}^{\N{}}}=\sum_{j=1}^{\N{}}X_j Y_j$ for $(X,Y)\in (\mathbb{R}^{\N{}})^2$,
and 
\[\Lambda^c _d(\omega)=diag\big(\Lambda^c _1(\omega),\dots,\Lambda^c _{\N{}}(\omega)\big).\]
\begin{thm}\label{expdeccpeP2} Let us assume that the energy transport matrix $\Gamma^{c}(\omega)$ is irreducible. Then, we have
 \[ \lim_{L\to+\infty}\frac{1}{L}\ln\left[ \sum_{j=1}^{\N{}}\mathcal{T}_{j}^l(\omega,L) \right]=-\Lambda_{\infty}(\omega)\]
with 
\begin{equation}\label{expcoefP2}\Lambda_{\infty}(\omega)  =\inf_{X\in \mathcal{S}^{\N{}}_{+}} \big< \big(-\Gamma^{c}(\omega)+\Lambda^c _d(\omega)\big)X , X \big>_{\mathbb{R}^{\N{}}},
 \end{equation}
which is positive as soon as one of the coefficients $\Lambda^c_j(\omega)$ is positive.
\end{thm}    
This result means that the total energy carried by the expected powers of the propagating modes decays exponentially with the propagation distance, and the decay rate can be expressed in terms of a variational formula over a finite-dimensional space. 
\begin{preuve}
The coupled power equations admit a probabilistic representation in terms of a jump Markov process. 
If we denote by $\big(Y^{\N{}}_{t}\big)_{t\geq0}$ a jump Markov process  with state space $\{1,\dots,\N{}\}$ and intensity matrix $\Gamma^{c}(\omega)$, then we have using the Feynman-Kac formula:  
\begin{equation}\label{prcpeP2} \mathcal{T}_{j}^l(\omega,z)=\mathbb{E}\left[ \exp\left(-\int_{0}^{z}\Lambda^{c}_{Y^{\N{}}_{s}} (\omega)ds\right)\textbf{1}_{\left(  Y^{\N{}}_{z} = j\right)} \Big\vert Y^{\N{}}_{0}=l \right].\end{equation}
Moreover, we have supposed that $\Gamma^{c}(\omega)$ is irreducible. Then, $\big(Y^{\N{}}_{t}\big)_{t\geq0}$ in an ergodic process with invariant measure $\mu_{\N{}}$, which is the uniform distribution over $\{1,\dots,\N{}\}$. That is, $\mu_{\N{}}(j)=1/\N{}$ $\forall j\in \{1,\dots,\N{}\}$. The self-adjoint generator of the jump Markov process $(Y^{\N{}}_{t})_{t\geq0}$ is given by
\[ \mathcal{L}^{\N{}}\phi (j)= \sum_{n=1}^{\N{}}\Gamma^{c}_{nj}(\omega)\left(\phi(n)-\phi(j)\right),\]
for every function $\phi$ from $\{1,\dots,\N{}\}$ to $\mathbb{R}$, and it is easy to check that $ \mathcal{L}^{\N{}}\mu_{\N{}}=0$.
Let us consider the local times
\[l_T(j)=\int_{0}^{T}\textbf{1}_{\left(Y^{\N{}}_{s}=j\right)}ds\]
 for $j\in \{1,\dots,\N{}\}$ and $T>0$, which corresponds to the time spent by the process $\big(Y^{\N{}}_{t}\big)_{t\geq0}$ in the state $j$ during the time interval $[0,T]$. According to \cite{donvar}, we have a large deviation principle for $\frac{1}{T}l_T$ viewed as a random process with values in $\mathcal{M}^{\N{}}_1$ which is the set of probability measures on $\{1,\dots,\N{}\}$. More precisely, we have 
\[\begin{split}
\lim_{L\to+\infty}\frac{1}{L}\ln\mathbb{E}\Big[ \exp\Big(-L\,\,\big< \Lambda^{c} ,\frac{1}{L} l_L&\big>_{\mathbb{R}^{\N{}}}\Big) \Big\vert Y^{\N{}}_{0}=l\Big]\\
&= \lim_{L\to+\infty}\frac{1}{L}\ln\mathbb{E}\Big[ \exp\Big( -\int_0^L \Lambda^c _{Y^{\N{}}_s }ds \Big) \Big\vert Y^{\N{}}_{0}=l\Big]  \\
&=-\inf_{\mu\in \mathcal{M}^{\N{}}_1}\left(I(\mu)+\big<\Lambda^{c}(\omega),\mu\big>\right)
\end{split}\]
 with
\[I(\mu)=\big\|\big(-\Gamma^{c}(\omega)\big)^{1/2}\sqrt{\mu}\big\|_{2,\mathbb{R}^{\N{}}}^2=\big<\big(-\Gamma^c(\omega) \big)\sqrt{\mu},\sqrt{\mu}\big>_{\mathbb{R}^{\N{}}}.\]
Consequently,
 \[ \lim_{L\to+\infty}\frac{1}{L}\ln\left[ \sum_{j=1}^{\N{}}\mathcal{T}_{j}^l(\omega,L) \right]=-\Lambda_{\infty}(\omega).\]
Let us assume that $\Lambda_{\infty}(\omega)=0$. As $\mathcal{S}^{\N{}}_{+}$ is a compact space, there exists $X_0 \in\mathcal{S}^{\N{}}_{+}$ such that 
\[\Lambda_{\infty}(\omega)=\big< \big(-\Gamma^{c}(\omega)+\Lambda^c_d (\omega)\big)X_0, X_0 \big>_{\mathbb{R}^{\N{}}}=0.\]
Moreover, $-\Gamma^{c}(\omega)$ and $\Lambda^c _d(\omega)$ are two nonnegative matrices and $0$ is a simple eigenvalue of $-\Gamma^{c}(\omega)$ by the Perron-Frobenius theorem. 
Then,
\[\big<(- \Gamma^{c}(\omega))X_0, X_0 \big>_{\mathbb{R}^{\N{}}}=0\Leftrightarrow X_0=\sqrt{\mu_{\N{}}},\]
and
\[\big< \Lambda^c _d(\omega)X_0, X_0 \big>_{\mathbb{R}^{\N{}}}=0\Rightarrow \exists j\in \{1,\dots,\N{}\},\quad  X_0(j)=0.\] 
Therefore, 
\[\Lambda_{\infty}(\omega)>0.\]$\blacksquare$
\end{preuve}
The expression \eqref{expcoefP2} of $\Lambda_{\infty}(\omega)$  is not simple. However, we have the following inequalities.
\begin{equation}\label{inedecP2}
\min_{j\in\{1,\dots,\N{}\}}\Lambda^c_{j}(\omega)\quad\leq \quad\Lambda_{\infty}(\omega) \quad\leq\quad \overline{\Lambda}(\omega)=\frac{1}{\N{}}\sum_{j=1}^{\N{}}\Lambda^c_{j}(\omega).
\end{equation}

First, we assume that $\forall j\in\{1,\dots,\N{}\}$, $\Lambda^c _j(\omega)=\Lambda (\omega)>0$. In this case, using \eqref{inedecP2}
\[ \Lambda_{\infty}(\omega) =\Lambda(\omega).\]
This means that if all the coefficients which represent the radiation losses are equal, the decay rate of the total energy of the propagating modes is given by this coefficient. 

Second, we assume that the coupling matrix is small, that is, we replace $\Gamma^{c}(\omega)$ by $\tau \Gamma^{c}(\omega)$ with $\tau \ll 1$. If $\forall j\in\{1,\dots,\N{}\}$, $\Lambda^c _j(\omega)>0$ we have
\[\lim_{\tau\to0}\Lambda^{\tau}_{\infty}(\omega)=\min_{j\in\{1,\dots,\N{}\}}\Lambda^c_{j}(\omega). \]
From \eqref{inedecP2}, it is the smallest value that $\Lambda_{\infty}(\omega)$ can take. This result is consistent with the fact that the coupling process on the transfer of energy between propagating modes is negligible and the decay rate of the energy of a particular propagating mode $j$ is given by its own decay coefficient $\Lambda_j(\omega)$. Then, for the total energy of propagating modes the decay rate is given by the minimum of those decay coefficients. Consequently, if there exists $\Lambda^c_{j_0}(\omega) =0$, we have   
\[\lim_{\tau\to0}\Lambda^{\tau}_{\infty}(\omega)=0.\]
The reason is the energy of the $j_0$th propagating mode stays approximately constant with a weak transfer of energy, and
 \[\lim_{\tau\to0}\frac{1}{\tau}\Lambda^{\tau}_{\infty}(\omega)=\inf_{X\in\tilde{V}} \big< \big(-\Gamma^{c}(\omega)\big)X , X \big>_{\mathbb{R}^{\N{}}} >0,\]
where
 \[\tilde{V}=\left\{ X\in \mathcal{S}^{\N{}}_+,\quad supp{X}\subset\{1,\dots,\N{}\}\setminus supp(\Lambda^c(\omega))  \right\},\]
because $\sqrt{\mu_{\N{}}}\not\in\tilde{V}$.

Now, we assume that the coupling matrix is large, that is we replace $\Gamma^{c}(\omega)$ by $\frac{1}{\tau} \Gamma^{c}(\omega)$ with $\tau \ll 1$. In this case, we have
\[\lim_{\tau\to0}\Lambda^{\tau}_{\infty}(\omega)=\overline{\Lambda}(\omega).\] 
From \eqref{inedecP2}, it is the largest value that $\Lambda_{\infty}(\omega)$ can take. The strong coupling produces a uniform distribution of energy over the propagating modes and the decay rate becomes $\big<\Lambda^{c}(\omega),\mu_{\N{}}\big>_{\mathbb{R}^{\N{}}}=\overline{\Lambda}(\omega)$ for each mode. A more convenient way to get this result is to use a probabilistic representation. In fact, we have  
\[ \begin{split}
\mathcal{T}_{j}^l(\omega,z)&=\mathbb{E}\left[ \exp\left(-\int_{0}^{z}\Lambda^{c}_{Y^{\N{}}_{s/\tau}} (\omega)\right)\textbf{1}_{\left(Y^{\N{}}_{z} = j\right)} \Big\vert Y^{\N{}}_{0}=l \right]\\
&=\mathbb{E}\left[ \exp\left(-z\frac{\tau}{z}\int_{0}^{z/\tau}\Lambda^{c}_{Y^{\N{}}_{s}} (\omega)\right)\textbf{1}_{\left(Y^{\N{}}_{z/\tau} = j\right)} \Big\vert Y^{\N{}}_{0}=l \right],
\end{split}\]
where $\big(Y^{\N{}}_{t}\big)_{t\geq0}$ is a jump Markov process with state space $\{1,\dots,\N{}\}$ and intensity matrix $\Gamma^{c}(\omega)$. Using the ergodic properties of  $\big(Y^{\N{}}_{t}\big)_{t\geq0}$, we get that 
 \[\lim_{\tau\to 0}\mathcal{T}_{j}^{\tau,l}(\omega,L)=\frac{1}{\N{}}\exp\left(-\overline{\Lambda}(\omega)L\right).\]
\begin{figure}\begin{center}
\includegraphics*[scale=0.65]{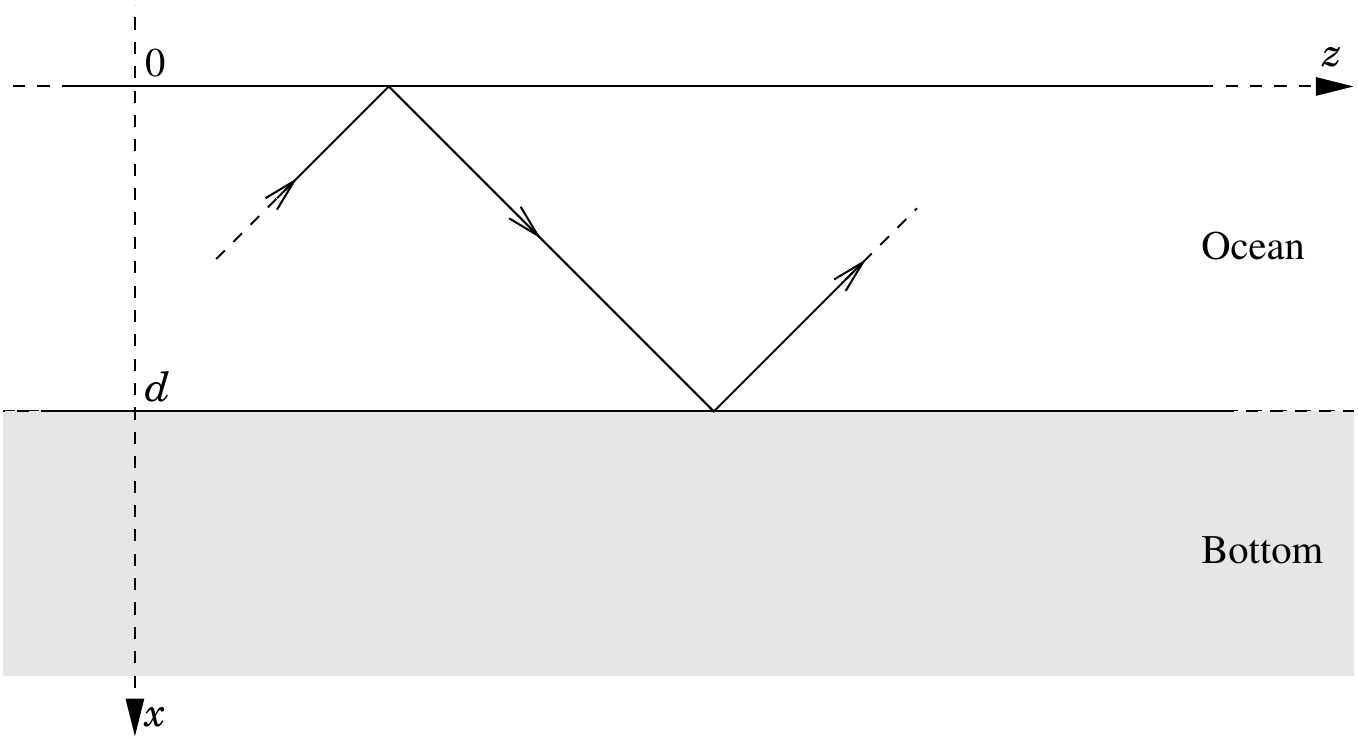}
\end{center}
\caption{\label{perterad2P2} Illustration of negligible radiation losses in the shallow-water random waveguide model.}
\end{figure}
Finally, if we assume that the radiation losses are negligible, that is we replace $\Lambda^{c}(\omega)$ by $\tau \Lambda^{c}(\omega)$ with $\tau \ll 1$, we have
\[\lim_{\tau\to0}\Lambda^{\tau}_{\infty}(\omega)=0.\] 
 In fact, if the radiative loss is negligible, the coupling process becomes dominant, and we can show that 
 \[\forall L>0, \quad \sup_{z\in[0,L]}\| \mathcal{T}_{j}^{\tau,l}(\omega,z)-\mathcal{T}_{j}^{0,l}(\omega,z)\|_{2,\mathbb{R}^{\N{}}}=\mathcal{O}(\tau),\]
 where $\mathcal{T}^{0,l}(\omega,.)$ satisfies \eqref{cpeP2} without the coefficient $\Lambda^{c}(\omega)$. In this situation 
 \[\mathcal{T}_{j}^{0,l}(\omega,L)=\mathbb{P}\left(Y^{\N{}}_L =j\Big\vert Y^{\N{}}_{0}=l\right),\]
 and the total energy is conserved (see Figure \ref{perterad2P2}), and
 \[\lim_{\tau\to0}\frac{1}{\tau}\Lambda^{\tau}_{\infty}(\omega)=\overline{\Lambda}(\omega)>0.\]
 As it was already observed in \cite{papa} the modal energy distribution converges as $L\to +\infty$ to a uniform distribution:
 \[\lim_{L\to+\infty}\mathcal{T}_{j}^{0,l}(\omega,L)=\frac{1}{\N{}}. \]

 \subsection{High-Frequency Approximation to Coupled Power Equations}\label{sechighfreq}

In this section, under the assumption that nearest-neighbor coupling, introduced in Section \ref{bliP2}, is the main power transfer mechanism, we give an approximate solution of the coupled power equations \eqref{cpeP2} in the high-frequency regime or in the limit of large number of propagating modes $\N{}\gg1$. Let us note that the limit of a large number of propagating modes $\N{} \gg 1$ corresponds to the high-frequency regime $\omega \to +\infty$. Next, we analyze the energy carried by the propagating modes in this regime.  

The coupled power equations can be approximated in the high-frequency regime by a diffusion equation. This approximation has been already obtained in \cite{papanicolaou} for instance, in which we can find further references about this topic. We can also refer to \cite{marcuse} for more discussions on this approximation. For an application of such a diffusion model to acoustic propagation in randomly sound channels we refer to \cite{mellen}, and for applications to time reversal of waves we refer to \cite{gomez}.

Using the form of the covariance function \eqref{covfunc}, we find
\[\Gamma^c _{jl}(\omega)=\frac{a k^4(\omega) I_{j,l}(\omega)}{2 \Bh{j}{}\Bh{l}{}(a^2+(\Bh{j}{}-\Bh{l}{})^2 )}\]
 and
 \[\Lambda^c _j (\omega)=\int_{0}^{\ko{}}\frac{a k^4(\omega) I_{j,\gamma }(\omega)}{2\Bh{j}{}\sqrt{\gamma }(a^2+(\Bh{j}{}-\sqrt{\gamma })^2 )} d\gamma ,\]
where
\[\begin{split}
I_{jl}  =& \frac{1}{4}A^{2}_j A^{2}_l \Big[ S\big(\sigma_j -\sigma_l,\sigma_j-\sigma_l\big)+S\big(\sigma_j +\sigma_l,\sigma_j +\sigma_l\big)  \\
 & \hspace{1.5cm}-  S\big(\sigma_j -\sigma_l,\sigma_j +\sigma_l\big)-S\big(\sigma_j +\sigma_l,\sigma_j -\sigma_l\big) \Big], \\
I_{j\gamma }  =&  \frac{1}{4}A^{2}_j A^{2}_\ga \Big[ S\big(\sigma_j -\eta ,\sigma_j -\eta\big)+S\big(\sigma_j +\eta ,\sigma_j+\eta \big) \\
 & \hspace{1.5cm} -  S\big(\sigma_j-\eta ,\sigma_j +\eta\big )-S\big(\sigma_j +\eta ,\sigma_j-\eta\big ) \Big],
\end{split}\]
with
\[S(v_1,v_2)=\int_{0}^{d} \int_{0}^{d} \ga_0(x_1,x_2) \cos\big(\frac{v_1}{d}x_1\big)\cos\big(\frac{v_2}{d}x_2\big)dx_1dx_2,\]
and where $A_j(\omega)$, $A_\ga(\omega)$, $\sigma_j(\omega)$, $\eta(\omega)$, $\phi_j(\omega,x)$, and $\phi_{\ga}(\omega, x)$ are defined in Section \ref{spectralP2}.

\subsubsection{Band-Limiting Idealization}\label{bliP2}

In this section, we introduce a band-limiting idealization hypothesis  
in which the power spectral density of the random fluctuations is assumed to be limited
in both the transverse and the longitudinal directions.

We assume that the support of $S$ lies in the square $\big[- \frac{3\pi}{2}, \frac{3\pi}{2} \big] \times \big[- \frac{3\pi}{2},\frac{3\pi}{2} \big]$.
Then,
\[I_{jl}(\omega)=\left\{ \begin{array}{ccl}
\frac{1}{4}A^{2}_j(\omega) A^2 _l(\omega) S\big(\sigma_j(\omega) -\sigma_l(\omega),\sigma_j (\omega)-\sigma_l(\omega) \big) & \text{ if }& \lvert j-l\rvert=1 \\
0 & \mbox{ otherwise, } & \end{array} \right. \]
and
\[I_{j\gamma }(\omega)=\left\{ \begin{array}{ccl}
\frac{1}{4}A^{2}_j(\omega) A^2 _l (\omega)S\big(\sigma_j(\omega) -\eta(\omega) ,\sigma_j(\omega) -\eta(\omega) \big) & \text{if}& \lvert \sigma_j(\omega)-\eta(\omega) \rvert \leq \frac{3\pi}{2} \\
0 & \text{otherwise.} & \end{array} \right. \]
From this assumption we get  $\forall \, 0<\gamma< \ko{}$ and $j\in\{1,\dots,\N{}-2 \}$,
\[\begin{split}
\eta(\omega) -\sigma_j(\omega) \geq  n_1 k(\omega)d \sqrt{1-\frac{1}{n_1 ^2}} -\sigma_j(\omega) & \geq  n_1 k(\omega)d \theta -(\N{}-2)\pi \\
 & \geq  \pi\underbrace{\left(\frac{n_1 k(\omega) d}{\pi}\theta -\N{}\right)}_{\in [0,1)} +2 \pi.
 \end{split}\]
Then, for $j\in\{1,\dots,\N{}-2\}$, 
\[\inf_{0<\ga < k^2} \eta(\omega) -\sigma _j(\omega)>\frac{3\pi}{2},\]
and
\[\Lambda^c _j(\omega) = 0, \quad \forall j\in \{1,\dots,\N{}-2 \}.\]
Consequently, the coupled power equations \eqref{cpeP2} become
\begin{equation}\label{cpe2P2}\begin{split}
\dz \mathcal{T}_{N} ^l (z)&=-\Lambda^c _{N}\mathcal{T}_{N} ^l (z)+\Gamma^c _{N-1\,N}\left(\mathcal{T}_{N-1} ^l (z)-\mathcal{T}_{N} ^l (z)\right),\\
\dz \mathcal{T}_{N-1} ^l (z)&=-\Lambda^c _{N-1}\mathcal{T}_{N-1} ^l (z)+\Gamma^c _{N-1\,N-2}\left(\mathcal{T}_{N-2} ^l (z)-\mathcal{T}_{N-1} ^l (z)\right)\\
&\hspace{3.1cm}+\Gamma^c _{N-1\,N}\left(\mathcal{T}_{N} ^l (z)-\mathcal{T}_{N-1} ^l (z)\right),\\
\dz \mathcal{T}_{j} ^l (z)&=\Gamma^c _{j-1\,j}\left(\mathcal{T}_{j-1} ^l (z)-\mathcal{T}_{j} ^l (z)\right)+\Gamma^c _{j+1\,j}\left(\mathcal{T}_{j+1} ^l (z)-\mathcal{T}_{j} ^l (z)\right) \text{ for }j\in\{2,\dots,N-2\},\\
\dz \mathcal{T}_{1} ^l (z)&=\Gamma^c _{2\,1}\left(\mathcal{T}_{2} ^l (z)-\mathcal{T}_{1} ^l (z)\right),
\end{split}\end{equation} 
with $\mathcal{T}_{j} ^l (0)=\delta_{jl}$.

The band-limiting idealization hypothesis is tantamount to a nearest-neighbor coupling. More precisely, this assumption implies that  $\forall (j,l)\in \{1,\dots,\N{}\}^2$ the $j$th mode amplitude can exchange informations with the $l$th amplitude mode if they are direct neighbors, that is, if they satisfy $\lvert j-l \rvert \leq 1$.

\subsubsection{High-Frequency Approximation}

The evolution of the mean mode powers of the propagating modes can be described, in the high-frequency regime or in the limit of a large number of propagating modes $\N{}\gg1$, by a diffusion model. This diffusive continuous model is equipped with boundary conditions which take into account the effect of the radiating modes at the bottom and the free surface of the waveguide (see Figure \ref{perterad1P2} page \pageref{perterad1P2}).

Let, $\forall \varphi \in \mathcal{C}^0([0,1])$, $\forall u \in[0,1]$, and $z\geq 0$, 
\[\mathcal{T}^{\N{}}_\varphi(z,u)=\mathcal{T}^{[\N{}u]}_\varphi (\omega,z)=\sum_{j=1}^{\N{}} \varphi\Big(\frac{j}{\N{}}\Big)\mathcal{T}^{[\N{}u]}_j (\omega,z),\]
where $\varphi \mapsto \mathcal{T}^{\N{}}_\varphi (z,.)$ can be extended to an operator from $L^2(0,1)$ to itself. Here, $L^2(0,1)$ is equipped with the inner product defined as follows:
$\forall (\varphi,\psi)\in L^2(0,1)^2$
\[\big<\varphi,\psi\big>_{L^2(0,1)}=\int_0^1 \varphi(v)\psi(v)dv.\]

\begin{thm}\label{hfapproxP2} We have
\begin{enumerate}
\item $\forall \varphi \in L^2(0,1)$ and $\forall z\geq 0$,
\[\lim_{\omega \to +\infty}\mathcal{T}^{\N{}}_\varphi(z,u)=\mathcal{T}_\varphi(z,u)\quad\text{in }L^2(0,1),\]
where
$\mathcal{T}_\varphi (z,u)$ satisfies the partial differential equation : $\forall (z,u)\in (0,+\infty)\times(0,1)$,
\[\frac{\partial}{\partial z} \mathcal{T}_\varphi (z,u)= \frac{\partial}{\partial u}\left(a_{\infty}(\cdot)\frac{\partial}{\partial u}  \mathcal{T}_\varphi\right)(z,u), \]
with the boundary conditions 
\[ \frac{\partial}{\partial u} \mathcal{T}_\varphi (z,0)=0,\quad  \mathcal{T}_\varphi (z,1)=0, \quad \text{and}\quad \mathcal{T}_\varphi (0,u)=\varphi(u), \]
$\forall z>0$.
\item $\forall u\in [0,1]$, $\forall z\geq0$, and $\forall \varphi \in \mathcal{C}^0 ([0,1])$ such that $\varphi(1)=0$, we have
\[\lim_{\omega\to +\infty} \mathcal{T}^{\N{}} _\varphi(z,u)=\mathcal{T}_\varphi (z,u).\]
 \end{enumerate}
 Here,
\[a_{\infty}(u)=\frac{a_0 }{1-\left(1-\frac{\pi^2 }{a^2 d^2} \right)(\theta u)^2},\]
with $a_0  = \frac{ \pi^2  S_{0} }{ 2 a n_1 ^4 d^4 \theta ^2}$, $\theta=\sqrt{1-1/n^2_1}$, $S_0 =\int_0^d\int_0^d \ga_0(x_1,x_2)\cos\big(\frac{\pi}{d}x_1\big)\cos\big(\frac{\pi}{d}x_2\big)dx_1 dx_2$. $n_1$ is the index of refraction in the ocean section $[0,d]$, $1/a$ is the correlation length of the random inhomogeneities in the longitudinal direction, and $\ga_0$ is the covariance function of the random inhomogeneities in the transverse direction.  
\end{thm} 
This theorem is a continuum approximation in the limit of a large number of propagating modes $\N{}\gg1$. This approximation gives us, in the high-frequency regime, a diffusion model for the transfer of energy between the $\N{}$-discrete propagating modes, with a reflecting boundary condition at $x=0$ (the top of the waveguide in Figure \ref{figureP2} page \pageref{figureP2}) and an absorbing boundary condition at $u=1$ (the bottom of the waveguide in Figure \ref{figureP2}) which represents the radiative loss (see Figure \ref{perterad1P2}).

\subsubsection{Exponential Decay in the High-Frequency Regime}\label{edhfrintP2}

In the high-frequency regime, we also observe that the energy carried by the continuum of propagating modes decays exponentially with the propagation distance.
The exponential decay of the energy in the high-frequency regime is given by the following result.
\begin{thm}\label{hfexpdecP2} $\forall \varphi \in L^2(0,1)\setminus\{0\}$ such that $\varphi\geq 0$, and $\forall u \in [0,1)$,
\[\lim_{L\to+\infty} \frac{1}{L}\ln \left[  \mathcal{T}_\varphi (L,u)  \right]=-\Lambda_{\infty},\]
where 
\[\Lambda_{\infty}=\inf_{\varphi \in \mathcal{D}} \int_0^1 a_\infty (v)\varphi'(v)^2dv >0\]
and
\[\mathcal{D}=\left\{ \varphi\in \mathcal{C}^\infty([0,1]),\quad \|\varphi\|_{L^2(0,1)}=1,\quad \frac{\partial}{\partial v}\varphi(0)=0,\quad  \varphi(1)=0\right\}.\]
\end{thm}
This result means that the energy carried by each propagating modes decays exponentially with the propagation propagation, and the decay rate can be expressed in terms of a variational formula. Consequently, the spatial inhomogeneities of the medium and the geometry of the shallow-water waveguide lead us to an exponential decay phenomenon caused by the radiative loss into the ocean bottom.

\begin{preuve}
We can see that the operator $P_\infty=\frac{\partial}{\partial v}\big(a_{\infty}(\cdot)\frac{\partial}{\partial v}  \big)$ on $L^2([0,1])$, with domain 
\[\mathcal{D}(P_\infty)=\left\{\varphi\in H^2(0,1),\quad  \frac{\partial}{\partial v}\varphi(0)=0,\quad  \varphi(1)=0  \right\}\] 
is self-adjoint. $P_\infty$ has a compact resolvent $R_\lambda=(\lambda Id-P_\infty)^{-1}$ because $[0,1]$ is a compact set and then it has a point spectrum $(\lambda_{j})_{j\geq 1}$ with eigenvectors denoted  by $(\phi_{\infty,j})_{j\geq 1}$.  
Moreover, all the eigenspaces are finite-dimensional subspaces of $\mathcal{D}(P_\infty)$
and $\forall \varphi \in \mathcal{D}(P_\infty)\setminus\{0\}$ 
\[\big< P_\infty (\varphi),\varphi \big>_{L^2(0,1)}<0.\]
Let us organize the point spectrum in the nonincreasing way, $\cdots<\lambda_2<\lambda_1<0$. We have
\[ \mathcal{T}_\varphi (L,v)=\sum_{j\geq1}\big<\varphi,\phi_{\infty,j}\big>_{L^2(0,1)}e^{\lambda_j L}\phi_{\infty,j}(v). \]

\begin{lem}\label{Krein}
$\lambda_1$ is a simple eigenvalue and one can choose $\phi_{\infty,1}$ such that $\phi_{\infty,1}(v)>0$ $\forall v\in[0,1)$.

\end{lem}

\begin{preuve}[of Lemma \ref{Krein}]
This lemma is a consequence of the Krein-Rutman theorem, but not its strongest form \cite{schaefer}. Indeed, the set of nonnegative functions in $L^2([0,1])$ has an empty interior. However, using the smoothness of the eigenvector the proof also works in our case as we shall see it. 

Using the maximum principle we know that if $\varphi \in L^2([0,1])$ such that $\varphi \geq 0$, we have $\mathcal{T}_{\varphi}(L,.)\geq0$, and then $R_\lambda(\varphi)\geq 0$. Consequently, applying the Krein-Rutman theorem \cite{schaefer} to the resolvent operator $R_\lambda$ with $\lambda>0$ and which is a compact operator,  the spectral radius $\rho(R_\lambda)$ is an eigenvalue, for which one can associate an eigenvector $\varphi_{\rho(R_{\lambda})}$ such that $\forall v\in [0,1]$, $\varphi_{\rho(R_{\lambda})}(v)\geq 0$. However, we have $\forall v\in [0,1)$, $\varphi_{\rho(R_{\lambda})}(v)> 0$. In fact, let us assume that there exists $v_0\in[0,1)$ such that $\varphi_{\rho(R_{\lambda})}(v_0)=0$, then $R_\lambda(\varphi_{R_\lambda})(v_0)=0$. Moreover, $R_{\lambda}(\varphi_{R_\lambda})=\rho(R_{\lambda})\varphi_{\rho(R_{\lambda})}$ is an eigenvector for $P_\infty$, and  then $\varphi_{\rho(R_{\lambda})}$ is a smooth function on $[0,1]$. Therefore, according to the proof of Theorem \ref{hfapproxP2} we have
\[\begin{split}
R_{\lambda}(\varphi_{R_\lambda})(v_0)&=\int_0^{+\infty}e^{-\lambda t}\mathcal{T}_{\varphi_{\rho(R_\lambda)}}(t,v_0)dt\\
&=\int_{0}^{+\infty}e^{-\lambda t}\mathbb{E}^{\overline{\mathbb{P}}_{v_0}}\big[\varphi_{\rho(R_\lambda)}(x(t))\textbf{1}_{(t<\tau_1)}\big]dt\\
&=\mathbb{E}^{\overline{\mathbb{P}}_{v_0}}\Big[\int_0^{\tau_1}e^{-\lambda t}\varphi_{\rho(R_\lambda)}(\lvert x(t)\rvert) dt\Big]=0,
\end{split}\]
where $\overline{\mathbb{P}}_{v_0}$ is the unique solution of the martingale problem associated to $\mathcal{L}_{\overline{a}_\infty}=\frac{\partial}{\partial v}\left(\overline{a}_\infty(\cdot)\frac{\partial}{\partial v}\right)$ and starting from $v_0$. Here, we have chosen $\overline{a}_\infty$ such that $\forall v\in[0,1]$, $\overline{a}_{\infty}(v)=\overline{a}_{\infty}(-v)=a_{\infty}(v)$, and the martingale problem associated to $\mathcal{L}_{\overline{a}_\infty}$ is well-posed. Moreover, $\tau_1=\inf(t\geq 0,\,\lvert x(t)\rvert\geq 1)$.
Consequently, $\overline{\mathbb{P}}_{v_0}\big(\int_0^{\tau_1} e^{-\lambda t}\varphi_{\rho(R_\lambda)}(\lvert x(t)\rvert)dt=0\big)=1$. However, we know that there exists $v_1\in(0,1)$ such that $\varphi_{\rho(R_\lambda)}(v_1)>0$, and then $v_1<v_0<1$. Therefore, $\overline{\mathbb{P}}_{v_0}(\tau_1<\tau_{v_1})=1$, and by the Markov property
\[\begin{split}
0<\mathbb{E}^{\overline{\mathbb{P}}_{v_0}}\big[e^{-\tau_{v_1}}\textbf{1}_{(\tau_{v_1}<+\infty)}\big]&=\mathbb{E}^{\overline{\mathbb{P}}_{v_0}}\big[e^{-\tau_{v_1}}\textbf{1}_{(\tau_{v_1}<+\infty,\tau_1<\tau_{v_1})}\big]\\
&<\mathbb{E}^{\overline{\mathbb{P}}_{1}}\big[e^{-\tau_{v_1}}\textbf{1}_{(\tau_{v_1}<+\infty)}\big]\\
&<\mathbb{E}^{\overline{\mathbb{P}}_{v_0}}\big[e^{-\tau_{v_1}}\textbf{1}_{(\tau_{v_1}<+\infty)}\big],
\end{split}\]   
which is impossible. Therefore, $\forall v\in[0,1)$, $\varphi_{\rho(R_\lambda)}>0$.
Now, to see that the eigenvalue $\rho(R_\lambda)$ is simple, let $\varphi\in L^2(0,1)\setminus\{0\}$ such that $R_\lambda (\varphi)=\rho(R_\lambda)\varphi$, and let 
\[\begin{array}{rcc}
P_{R_\lambda}:\mathbb{R}&\longrightarrow &\mathcal{C}^0([0,1])\\
t&\longmapsto&\varphi_{R_\lambda}-t\varphi,
\end{array}\]
which is a continuous function. We recall that $\varphi$ is a smooth function on $[0,1]$. Let us show that $\exists t\in\mathbb{R}$ such that $\varphi=t\,\varphi_{\rho(R_\lambda)}$, that is $0\in P_{R_\lambda}(\mathbb{R})$. To do this let us assume that $0\not\in P_{R_\lambda}(\mathbb{R})$. By linearity one can assume that $\exists v_0\in[0,1)$ such that $\varphi(v_0)>0$. Let $\eta>0$ be small enough to have $v_0\in[0,1-\eta]$. Let $K^+_\eta=\big\{\varphi \in\mathcal{C}^0([0,1-\eta]),\forall v\in[0,1-\eta], \varphi(v)\geq 0\big\}$, then the interior of $K^+_\eta$ for the sup norm on $[0,1]$ is $K^{++}_\eta=\big\{\varphi \in\mathcal{C}^0([0,1-\eta]),\forall v\in[0,1-\eta], \varphi(v)> 0\big\}$. Moreover, for $t$ small enough $\varphi_{R_\lambda}-t\varphi \in K^{++}_\eta$, and $\varphi_{R_\lambda}-t\varphi \not\in K^{+}_\eta$  for $t$ large enough. Then $\exists t_0\in \mathbb{R}$ such that $\varphi_{R_\lambda}-t_0\varphi \in K^{+}_\eta\setminus K^{++}_\eta$. However,  $\varphi_{R_\lambda}-t_0\varphi \geq 0$, but  $\varphi_{R_\lambda}-t_0\varphi \not= 0$ because $0\not\in P_{R_\lambda}(\mathbb{R})$. Following the previous work we have 
\[ \rho(R_\lambda)(\varphi_{R_\lambda}-t_0\varphi)=R_\lambda(\varphi_{R_\lambda}-t_0\varphi)\in K^{++}_\eta.\]    
Consequently, $\rho(R_\lambda)=1/(\lambda-\lambda_1)$ implies that $\lambda_1$ is also a simple eigenvalue and one can choose 
\[\phi_{\infty,1}=R_\lambda(\varphi_{R_\lambda})=\rho(R_\lambda)\varphi_{R_\lambda}\in K^{++}_\eta.\]
That concludes the proof of Lemma \ref{Krein}.
$\square$\end{preuve}
As a result, $\forall \varphi \in L^2(0,1)\setminus\{0\}$ such that $\varphi\geq 0$, $\forall v\in[0,1)$ we get
\[\lim_{L\to+\infty} \frac{1}{L}\ln \left[ \mathcal{T}_\varphi (L,v)  \right]=\lambda_1,\] 
and
\[\lambda_1=\sup_{\substack{\varphi \in \mathcal{D}(P_\infty)\\ \|\varphi\|_{L^2([0,1])}= 1}}\big<P_\infty(\varphi),\varphi\big>_{L^2([0,1])}=-\Lambda_{\infty}<0.\]
$\blacksquare$
\end{preuve}
 In Theorem \ref{hfexpdecP2}, we take $\varphi \in L^2(0,1)\setminus\{0\}$ such that $\varphi\geq 0$, which can be consider as being the initial repartition of energy over the continuum of modes. However, the result of Theorem \ref{hfexpdecP2} is also valid for any $\varphi \in L^2(0,1)\setminus\{0\}$ such that $\big<\varphi,\phi_{\infty,1}\big>_{L^2(0,d)} > 0$.

\subsection{High-Frequency Approximation to Coupled Power Equation with Negligible Radiation Losses}\label{sechighfreqnrad}

In the case of negligible radiation losses, we also get a continuous diffusive model for the coupled power equations in the high-frequency regime or in the limit of a large number of propagating modes $\N{}\gg1$. This diffusive continuous model is equipped with boundary conditions which take into account the negligible effect of the radiation losses at the bottom and the free surface of the waveguide (see Figure \ref{perterad2P2} page \pageref{perterad2P2}).

Now, let us assume that the radiation losses are negligible, that is, $\Lambda^{c}(\omega)=\tau \tilde{\Lambda}^{c}(\omega)$ with $\tau \ll 1$. We have already remarked that, if the radiation losses are negligible, then the coupling process is predominant and we have 
 \[\forall L>0, \quad \sup_{z\in[0,L]}\| \mathcal{T}_{j}^{\tau,l}(\omega,z)-\mathcal{T}_{j}^{0,l}(\omega,z)\|_{2,\mathbb{R}^{\N{}}}=\mathcal{O}(\tau),\]
 where $\mathcal{T}^{0,l}(\omega,.)$ satisfies 
\begin{equation*}\begin{split}
\dz \mathcal{T}_{N} ^{0,l} (z)&=\Gamma^c _{N-1\,N}\left(\mathcal{T}_{N-1} ^{0,l} (z)-\mathcal{T}_{N} ^{0,l} (z)\right),\\
\dz \mathcal{T}_{j} ^{0,l} (z)&=\Gamma^c _{j-1\,j}\left(\mathcal{T}_{j-1} ^{0,l} (z)-\mathcal{T}_{j} ^{0,l} (z)\right)+\Gamma^c _{j+1\,j}\left(\mathcal{T}_{j+1} ^{0,l} (z)-\mathcal{T}_{j} ^{0,l} (z)\right) \text{ for }j\in\{2,\dots,N-1\},\\
\dz \mathcal{T}_{1} ^{0,l} (z)&=\Gamma^c _{2\,1}\left(\mathcal{T}_{2} ^{0,l} (z)-\mathcal{T}_{1} ^{0,l} (z)\right),
\end{split}\end{equation*} 
with $\mathcal{T}_{j} ^{0,l} (0)=\delta_{jl}$.

\subsubsection{High Frequency Approximation}

Let, $\forall \varphi \in \mathcal{C}^0([0,1])$, $\forall u \in[0,1]$, and $z\geq 0$, 
\[\mathcal{T}^{\N{}}_\varphi(z,u)=\mathcal{T}^{[\N{}u]}_\varphi (z)=\sum_{j=1}^{\N{}} \varphi\Big(\frac{j}{\N{}}\Big)\mathcal{T}^{[\N{}u]}_j (z),\]
where $\varphi \mapsto \mathcal{T}^{\N{}}_\varphi(z,.)$ can be extended into an operator from $L^2(0,1)$ to itself. 
\begin{thm}\label{hfapprox0P2} 
We have
\begin{enumerate}
\item $\forall \varphi \in L^2(0,1)$ and $\forall z\geq 0$,
\[\lim_{\omega \to +\infty}\mathcal{T}^{\N{}}_\varphi(z,u)=\mathcal{T}_\varphi(z,u)\quad\text{in }L^2(0,1),\]
where
$\mathcal{T}_\varphi (z,u)$ satisfies the partial differential equation : $\forall (z,u)\in (0,+\infty)\times(0,1)$,
\[\frac{\partial}{\partial z} \mathcal{T}_\varphi (z,u)= \frac{\partial}{\partial u}\left(a_{\infty}(\cdot)\frac{\partial}{\partial u}  \mathcal{T}_\varphi\right)(z,u), \]
with the boundary conditions 
\[ \frac{\partial}{\partial u} \mathcal{T}_\varphi (z,0)=0,\quad  \frac{\partial}{\partial v} \mathcal{T}_\varphi (z,1)=0, \quad \text{and}\quad \mathcal{T}_\varphi (0,u)=\varphi(u), \]
$\forall z>0$.
\item $\forall u\in [0,1)$, $\forall z\geq0$, and $\forall \varphi \in \mathcal{C}^0 ([0,1])$ such that $\varphi(1)=0$, we have
\[\lim_{\omega\to +\infty} \mathcal{T}^{\N{}} _\varphi(z,u)=\mathcal{T}_\varphi (z,u).\]
 \end{enumerate}
Here,
\[a_{\infty}(u)=\frac{a_0 }{1-\left(1-\frac{\pi^2 }{a^2 d^2} \right)(\theta u)^2},\]
with $a_0  = \frac{ \pi^2  S_{0} }{ 2 a n_1 ^4 d^4 \theta ^2}$, $\theta=\sqrt{1-1/n^2_1}$, $S_0 =\int_0^d\int_0^d \ga_0(x_1,x_2)\cos\big(\frac{\pi}{d}x_1\big)\cos\big(\frac{\pi}{d}x_2\big)dx_1 dx_2$. $n_1$ is the index of refraction in the ocean section $[0,d]$, $1/a$ is the correlation length of the random inhomogeneities in the longitudinal direction, and $\ga_0$ is the covariance function of the random inhomogeneities in the transverse direction.  \end{thm} 

This theorem is a continuum approximation in the limit of a large number of propagating modes in the case where the radiation losses are negligible.This approximation gives us, in the high-frequency regime, a diffusion model for the transfer of energy between the $\N{}$-discrete propagating modes, with two reflecting boundary conditions at $u=0$ (the top of the waveguide in Figure \ref{figureP2} page \pageref{figureP2}) and $u=1$ (the bottom of the waveguide in Figure \ref{figureP2}). Here, the two reflecting boundary conditions mean that there is no radiative loss anymore (see Figure \ref{perterad2P2}).

\subsubsection{Asymptotic behavior of $\mathcal{T}(L,v)$ as $L\to +\infty$}

In the case where the radiation losses are negligible, we have seen in Section \ref{expdecP2} that the decay rate satisfies $\lim_{\tau\to 0}\Lambda^\tau _\infty (\omega)=0$ and $\mathcal{T}^{0,l}(\omega,L)$ converge to the uniform distribution over $\{1,\dots,\N{}\}$ as $L\to +\infty$  \cite{papa}. In the high-frequency regime we have the following continuous version.
\begin{thm}
$\forall \varphi \in L^2(0,1)$ and $\forall u\in[0,1]$,
\[\lim_{L\to+\infty}  \mathcal{T}_\varphi (L,u)=\int_0^1\varphi(v) dv,\]
that is, the energy carried by the continuum of propagating modes converges exponentially fast to the uniform distribution over $[0,1]$ as $L\to +\infty$. 
\end{thm}
As a result, the energy is conserved and the modal energy distribution converges to a uniform distribution as $L\to+\infty$. 
\begin{preuve}
We can see that the operator $P_\infty=\frac{\partial}{\partial v}\big(a_{\infty}(\cdot)\frac{\partial}{\partial v}  \big)$ on $L^2([0,1])$, with domain 
\[\mathcal{D}(P_\infty)=\left\{\varphi\in H^2(0,1),\quad  \frac{\partial}{\partial v}\varphi(0)=0,\quad \frac{\partial}{\partial v} \varphi(1)=0  \right\}\] 
is self-adjoint. Moreover, $P_\infty$ has a compact resolvant because $[0,1]$ is a compact set and then it has a point spectrum $(\lambda_{j})_{j\geq 0}$ with eigenvectors denoted by $(\phi_{\infty,j})_{j\geq 0}$.  
Moreover, all the eigenspaces are finite-dimensional subspaces of $\mathcal{D}(P_\infty)$
and $\forall \varphi \in \mathcal{D}(P_\infty)\setminus\{0\}$ 
\[\big< P_\infty (\varphi),\varphi \big>_{L^2(0,1)}\leq 0.\]
Let us remark that $\lambda_0=0$ is a simple eigenvalue with eigenvector $\phi_{\infty,0}=1$. Then, the spectrum is include in $(-\infty,0]$ and we have the following decomposition
\[ \mathcal{T}_\varphi (z,v)= \int_0^1 \varphi(v)dv+ \sum_{j\geq1}\big<\varphi,\phi_{\infty,j}\big>_{L^2(0,1)}e^{\lambda_j z}\phi_{\infty,j}(v).\]
Therefore, $\forall u\in[0,1]$,
\[\lim_{L\to+\infty} \mathcal{T}_\varphi (L,u)= \int_0^1 \varphi(v)dv,\]
with exponential rate $\lambda_1<0$. 
$\blacksquare$
\end{preuve}

\section*{Conclusion}

In this paper we have analyzed the propagation of waves in a shallow-water acoustic waveguide with random perturbations. In such a waveguide, the wave field can be decomposed into three kinds of modes, which are the propagating, the radiating, and the evanescent modes, and the random perturbations produce a coupling between these modes. 

We have shown that the evolution of the propagating mode amplitudes can be described as a diffusion process (Theorems \ref{thasymptP21} and \ref{thasympP2}). This diffusion takes into account the main coupling mechanisms: The coupling with the evanescent modes induces a mode-dependent and frequency-dependent phase modulation on the propagating modes, the coupling with the radiating modes, in addition to a mode-dependent and frequency-dependent phase modulation, induces a mode-dependent and frequency-dependent attenuation on the propagating modes. In other words, the propagating modes lose energy in the form of radiation into the bottom of the waveguide and their total energy decays exponentially with the propagation distance. We can express the decay rate in terms of a variational formula over a finite-dimensional space (Theorem \ref{expdeccpeP2}).

Under the assumption that nearest-neighbor coupling is the main power transfer mechanism, the evolution of the mean mode powers of the propagating modes can be described, in the high frequency regime or in the limit of a large number of propagating modes, by a continuous diffusive model with boundary conditions which take into account the effect of the radiation losses at the bottom and the free surface of the waveguide. In this regime, we observe that the energy carried by the continuum of propagating modes also decay exponentially with the propagation distance. The exponential decay rate can be express in terms of a variational formula (Theorem \ref{hfexpdecP2}).

The diffusive systems obtained in this paper can be used to analyze pulse propagation and refocusing during time-reversal experiments in underwater acoustics \cite{gomez}.

\section{Appendix}

\subsection{Gaussian Random Field}\label{grfield}

This section is a short remainder about some properties of Gaussian random fields that we use in the proofs of Theorems \ref{thasymptP21}, and in Sections \ref{globalP2} and \ref{iemprP2}. All the results exposed in this section can be shown using the standard properties of Gaussian random fields presented in \cite{adler} and \cite{adlertaylor} for instance.

In this paper, the random perturbations of the medium parameters are modeled using a random process denoted by $(V(x,t), x \in [0,d], t\geq 0)$. Throughout this paper the process $V$ is a continuous real-valued zero-mean Gaussian field with a covariance function given by
\begin{equation}\label{covfunc} \mathbb{E}\left[V(x,z_1)V(y,z_2)\right]=\gamma_0(x,y) e^{-a \lvert z_1-z_2 \rvert} \quad \forall(x,y)\in[0,d]^2\text{ and }\forall (z_1,z_2)\in[0,+\infty)^2.\end{equation} 
Here, $a>0$; $\gamma_0:[0,d]\times[0,d] \rightarrow\mathbb{R}$ is a Lipschitz function, which is the kernel of a nonnegative operator, that is, there exists a nonnegative operator $Q_{\gamma_0}$ from $L^2(0,d)$ to itself such that $\forall(\varphi,\psi)\in L^2(0,d)^2$ 
\[\big<Q_{\gamma_0}(\varphi),\psi\big>_{L^2(0,d)}=\int_0^d\int_0^d \gamma_0(x,y)\varphi(x)\psi(y)dxdy.\]
Consequently, one can consider the process $(V(.,t))_{t\geq0}$ as being a continuous zero-mean Gaussian field with values in $L^2(0,d)$ and covariance operator $Q_{\gamma_0}$. In other words, $\forall n\in\mathbb{N}^\ast$, $\forall (\varphi_1,\dots,\varphi_n) \in L^2(0,d)^n$, and $\forall (t_1\dots,t_n)\in[0,+\infty)^n$ 
\[\big(V_{\varphi_1}(t_1),\dots,V_{\varphi_n}(t_n)\big)=\big(\big<V(.,t_1),\varphi_1\big>_{L^2(0,d)},\dots,\big<V(.,t_n),\varphi_n\big>_{L^2(0,d)}\big)\] 
is a real-valued zero-mean Gaussian vector such that $\forall(j,l)\in\{1,\dots,n\}^2$
\begin{equation}\label{autofunct}
\mathbb{E}\left[V_{\varphi_j}(t_j)V_{\varphi_l}(t_l)\right]=\big<Q_{\gamma_0}(\varphi_j),\varphi_l\big>_{L^2(0,d)}e^{-a\lvert t_j-t_l\rvert}.
\end{equation}
With this point of view we have the following proposition.
\begin{prop} We have
\begin{enumerate}
\item $(V(.,t))_{t\geq0}$ is a continuous zero-mean stationary Gaussian field with values in $L^2(0,d)$ and autocorrelation function given by \eqref{autofunct}. Then, we have $\forall n\in \mathbb{N}^{\ast}$ and $\forall t\geq 0$,
\begin{equation}\label{cg2}
\mathbb{E}\left[\left(\int_0^d \big\lvert V\big(x,t\big) \big\rvert^2 dx \right)^n\right]=\mathbb{E}\left[\left(\int_0^d \big\lvert V\big(x,0\big) \big\rvert^2 dx \right)^n\right]<+\infty.
\end{equation}
\item We have the following Markov property. Let 
\[\mathcal{F}_t=\sigma(V(.,s),s\leq t)\] be the $\sigma$-algebra generated by $(V(.,s), s\leq t)$. We have 
\[\Big( V(.,t+h) \Big\vert \,\mathcal{F}_t\Big) = \Big( V(.,t+h)  \Big\vert\,\sigma(V(.,t) )\Big),\]
where the equality holds in law, and this law is the one of a Gaussian field with mean 
\[\mathbb{E}\big[ V(.,t+h ) \vert \mathcal{F}_t\big]=e^{-ah}V(.,t)\] 
and covariance, $\forall(\varphi,\psi)\in L^2(0,d)^2$,
\[\begin{split}
\mathbb{E}\Big[ V_\varphi(t+h )V_\psi(t+h)-  \mathbb{E}\big[ V_\varphi(t+h ) \vert \mathcal{F}_t\big]&\mathbb{E}\big[ V_\psi(t+h ) \vert \mathcal{F}_t\big]  \Big\vert \mathcal{F}_t\Big]\\
&=\big<Q_{\gamma_0}(\varphi),\psi\big>_{L^2(0,d)}\left(1-e^{-2ah} \right).
\end{split}\]
\end{enumerate}
\end{prop}
The Markov property of the random process $(V(.,t))_{t\geq0}$ is a direct consequence of the exponential form of the autocorrelation function \eqref{autofunct} with respect to the variable $t$ \cite{adler}. This property will be used in the proof of Theorems \ref{thasymptP21}, which are based on the perturbed-test-function method.

Now, we are interested in some estimation on the supremum of $V(x,t)$ with respect to the two variables $x$ and $t$. To this end, let us introduce some notations \cite{adlertaylor}. Let $\e>0$ be a small parameter and $L>0$. We consider the following pseudo-metric on the square $[0,d]\times[0,L/\e]$ defined by
\[m\big((x,t),(y,s)\big)=\mathbb{E}\left[(V(x,t) -V(y,s))^2\right]^{1/2}\leq K_{\gamma_0}\big[\lvert t-s\rvert+\lvert x-y\rvert \big].\] 
Let us remark that $[0,d]\times[0,L/\e]$ associated to the pseudo-metric $m$ is a compact set. From Theorem 1.3.3 in \cite{adlertaylor}, we have
\[\begin{split}
\mathbb{E}\left[\sup_{\substack{x\in[0,d]\\t\in[0,L/\e]}}\big\lvert V(x,t)\big\rvert \right]&\leq K \int_0^{diam([0,d]\times [0,L/\e])/2}H^{1/2}(r)dr\\
&\leq K_1\int_0^{\sup_{x\in[0,d]}\gamma_0(x,x)} \sqrt{\ln\left( K_2\frac{dL}{r^2\e}  \right)}dr,
\end{split} \]
where $H(r)=\ln(N(r))$, and $N(r)$ denotes the smallest number of balls, for the pseudo-metric $m$, with radius $r$ to cover the square $[0,d]\times[0,L/\e]$. Here, $diam$ stands for the diameter with respect to the pseudo-metric $m$. Consequently, we have the following proposition. 
\begin{prop}
$\forall \mu>0$ and $\forall K>0$,
\begin{equation}\label{cg1}
\lim_{\e\to 0} \mathbb{P}\Big(\e^\mu\sup_{x\in[0,d]}\sup_{t\in[0,L/\e]}\big\lvert V(x,t)\big\rvert \geq K\Big)=0.\end{equation}
\end{prop}
Moreover, according to Theorem 2.1.1 in \cite{adlertaylor}, one can show that 
 the limit \eqref{cg1} is obtained exponentially fast as $\e \to 0$.

\subsection{Proof of Theorem \ref{thasymptP21}}\label{proofthasympt21}

The proof of this theorem is in two parts. The process $\big(\textbf{T}^{\xi,\e}(z)\big)_{z \geq 0}$ is not adapted with respect to the filtration $\mathcal{F}^\e _z=\mathcal{F}_{z/\e}$. Then, the first part of the proof consists in simplifying the problem and introducing a new process for which the martingale approach can be used. The first part of the proof follows the ideas of \cite{khasminskii}. The second part of proof of this theorem is based on a martingale approach using the perturbed-test-function method and  follows the ideas developed in \cite{carmona}.

Then, let us introduce $ \tilde{\textbf{T}}^{\xi,\e} (.)$ the unique solution of the differential equation

\begin{equation}\label{transferP22}
\dz \tilde{\textbf{T}}^{\xi,\e} (z)=\frac{1}{\sqrt{\e}} \textbf{H}^{aa} \left(\frac{z}{\e}\right)
\tilde{\textbf{T}}^{\xi,\e} (z)+ \big<\textbf{G}^{aa}\big> \tilde{\textbf{T}}^{\xi,\e} (z),
\end{equation}
with $\textbf{T}^{\xi,\e} (0)=  Id$ and where $\big<\textbf{G}^{aa}\big>$ is defined, $\forall y \in \mathcal{H}_\xi$, by
\[
 \big<\textbf{G}^{aa}\big> _{j}(y)=\int_{-\infty}^{-\xi}\frac{i k^{4}}{2\beta_j \sqrt{\lvert\ga'\rvert}} \int_{0}^{+\infty}\mathbb{E}\big[C_{j\ga'}(0)C_{j\ga'}(z)\big]\cos\big(\beta_j z\big)e^{-\sqrt{\lvert\ga'\rvert}z}dzd\ga' y_j  
\]
 $\forall j\in \big\{1,\dots,N\big\}$ and $\big<\textbf{G}^{aa}\big> _{\ga}(y)=0 $ for $\ga\in(\xi,k^2)$.
We have the following proposition that describes the relation between the two processes $\textbf{T}^{\xi,\e}(z)$ and $\tilde{\textbf{T}}^{\xi,\e}(z)$.
\begin{prop}\label{propoP2}
\[\forall y \in\mathcal{H}_{\xi} \text{ and } \forall \eta>0,\quad \lim_{\e \to 0}\mathbb{P}\left( \sup_{z\in[0,L]}\|\emph{\textbf{T}}^{\xi,\e}(z)(y)-\tilde{\emph{\textbf{T}}}^{\xi,\e}(z)(y)\|^2_{\mathcal{H}_{\xi}} > \eta \right)=0.\]
\end{prop}  
Let us remark that the new process $\big(\tilde{\textbf{T}}^{\xi,\e}(z)\big)_{z \geq 0}$ is adapted to the filtration $\mathcal{F}^\e _z$ and 
\[\|\tilde{\textbf{T}}^{\xi,\e} (z)(y) \|^2_{\mathcal{H}_{\xi}}=\|y\|^2_{\mathcal{H}_{\xi}} \quad \forall z\geq 0. \]
Let $r_y=\|y\|_{\mathcal{H}_{\xi}}$,
\[\mathcal{B}_{r_y, \mathcal{H}_{\xi}}=\left\{\lambda \in\mathcal{H}_\xi, \|\lambda\|_{\mathcal{H}_\xi}=\sqrt{\left<\lambda,\lambda\right>_{\mathcal{H}_\xi}} \leq r_y\right\}\] 
the closed ball with radius $r_y$, and $\{g_n, n\geq 1\}$ a dense subset of $\mathcal{B}_{r_y,\mathcal{H}_\xi}$. We equip  $\mathcal{B}_{r_y, \mathcal{H}_\xi}$ with the distance $d_{\mathcal{B}_{r_y,\mathcal{H}_\xi}}$ defined by
\[d_{\mathcal{B}_{r_y,\mathcal{H}_\xi}}(\lambda, \mu)=\sum_{j=1}^{+\infty}\frac{1}{2^j}\left\lvert\big<\lambda-\mu,g_n\big>_{\mathcal{H}_\xi}\right\rvert\]
$\forall (\lambda,\mu)\in({\mathcal{B}_{r_y,\mathcal{H}_\xi}})^2$, and then $(\mathcal{B}_{\mathcal{H}_\xi} ,d_{\mathcal{B}_{r_y,\mathcal{H}_\xi}})$ is a compact metric space.

Using a particular tightness criterion, we prove the tightness of the family $(\tilde{\textbf{T}}^{\xi,\e} (.))_{\e \in (0,1)}$ on $\mathcal{C}([0,+\infty),(\mathcal{B}_{r_y, \mathcal{H}_\xi},d_{\mathcal{B}_{r_y,\mathcal{H}_\xi}} ))$, which is a polish space. We have chosen such a space to be able to apply the Portmanteau theorem. In a second part, we shall characterize all subsequence limits as solutions of a well-posed martingale problem in the Hilbert space $\mathcal{H}_{\xi}$. 

We have the following version of the Arzelà-Ascoli theorem \cite{billingsley,inf2} for processes with values in a complete separable metric space. 
\begin{thm}
A set $B\subset \mathcal{C}([0,+\infty),(\mathcal{B}_{r_y,\mathcal{H}_\xi},d_{\mathcal{B}_{r_y,\mathcal{H}_\xi}} ))$ has a compact closure if and only if 
\[\forall T>0, \quad \lim_{\eta \to 0}\sup_{g \in A} m_T (g, \eta)=0,  \]
with 
\[m_T (g,\eta)=\sup_{\substack{(s,t)\in [0,T]^2\\ \lvert t-s \rvert \leq\eta}} d_{\mathcal{B}_{r_y,\mathcal{H}_\xi}}( g(s) ,g(t)). \] 
\end{thm}
From this result, we obtain the classical tightness criterion. 
\begin{thm}
A family of probability measure $\big(\mathbb{P}^\e\big)_{\e\in(0,1)}$ on $\mathcal{C}([0,+\infty),(\mathcal{B}_{r_y,\mathcal{H}_\xi},d_{\mathcal{B}_{r_y, \mathcal{H}_\xi}} ))$ is tight if and only if
\[\forall T>0, \eta'>0 \quad \lim_{\eta \to 0}\sup_{\e \in (0,1) }\mathbb{P}^\e \big(g\,; \,\,m_T (g, \eta)>\eta'\big)=0.\]
\end{thm}
From the definition of the metric $d_{\mathcal{B}_{r_y,\mathcal{H}_\xi}}$, the tightness criterion becomes the following.
\begin{thm}\label{crit}
A family of processes $(X^\e)_{\e\in(0,1)}$ is tight in $\mathcal{C}([0,+\infty),(\mathcal{B}_{r_y,\mathcal{H}_\xi},d_{\mathcal{B}_{r_y, \mathcal{H}_\xi}} ))$ if and only if $\big(\big<X^\e,\lambda\big>_{\mathcal{H}_\xi}\big)_{\e\in(0,1)}$ is tight on $\mathcal{C}([0,+\infty),\mathbb{C})$ $\forall \lambda \in \mathcal{H}_\xi$.
\end{thm}
This last theorem looks like the tightness criterion of Mitoma and Fouque \cite{mitoma,fouque}.

For any $\lambda \in \mathcal{H}_\xi$, we set $\tilde{\textbf{T}}^{\xi,\e}_\lambda (z)(y)=\big<\tilde{\textbf{T}}^{\xi,\e}(z)(y),\lambda\big>_{\mathcal{H}_\xi}$. According to Theorem \ref{crit}, the family $(\tilde{\textbf{T}}^{\xi,\e}(.)(y))_{\e}$ is tight in $\mathcal{C}([0,+\infty), (\mathcal{B}_{r_y,\mathcal{H}_\xi},d_{\mathcal{B}_{r_y,\mathcal{H}_\xi}} ))$ if and only if the family $(\tilde{\textbf{T}}^{\xi,\e}_\lambda(.)(y))_{\e }$ is tight on $\mathcal{C}([0,+\infty),\mathbb{C} )$  $\forall \lambda \in \mathcal{H}_\xi$. Furthermore, $(\tilde{\textbf{T}}^{\xi,\e}(.)(y))_{\e}$ is a family of continuous processes. Then, it is sufficient to prove that $\forall\lambda \in \mathcal{H}_\xi$, $(\tilde{\textbf{T}}^{\xi,\e}_\lambda(.)(y))_{\e}$ is tight in $\mathcal{D}([0,+\infty),\mathbb{C})$, which is the set of cad-lag functions with values in $\mathbb{C}$ equipped with the Skorokhod topology. 

\begin{preuve}[of Proposition \ref{propoP2}]
Differentiating the square norm and using the fact that $\textbf{H}^{aa}(z)$ is skew Hermitian, we get
\[\begin{split}
 \|\textbf{T}^{\xi,\e}(z)(y) &-\tilde{\textbf{T}}^{\xi,\e}(z)(y) \|^2_ {\mathcal{H}_\xi} \\
 &\leq 2 \Big\lvert \int_0 ^z \Big<\big(\textbf{G}^{aa}\left(\frac{z}{\e}\right)-\big<\textbf{G}^{aa}\big>\big) \textbf{T}^{\xi,\e}(z)(y),\textbf{T}^{\xi,\e}(z)(y)- \tilde{\textbf{T}}^{\xi,\e}(z)(y)\Big>_{\mathcal{H}_\xi}   du\Big\rvert\\
 &\quad + 2 \big\| \big<\textbf{G}^{aa}\big>\big\|\int_0^z\| \textbf{T}^{\xi,\e}(u)(y)- \tilde{\textbf{T}}^{\xi,\e}(u)(y)\|^2_{\mathcal{H}_\xi}.
\end{split}\]
Let $\eta' >0$, we will split the interval $[0,z/\e]$ into intervals of length $\eta'/\sqrt{\e}$. The idea is that over these intervals the fast dynamic of $\textbf{G}^{aa}$ averages out while $\textbf{T}^{\xi,\e}$ does not move significantly. We have
\[\begin{split}
\Big\lvert\e \int_0 ^{z/\e} & \Big<\big(\textbf{G}^{aa}(u))-\big<\textbf{G}^{aa}\big>\big)\textbf{T}^{\xi,\e}(\e u)(y),\textbf{T}^{\xi,\e}(\e u)(y)- \tilde{\textbf{T}}^{\xi,\e}(\e u)(y)\Big>_{\mathcal{H}_\xi}   du\Big\rvert \\
& \leq\Big\lvert \e \int_0 ^{\big[\frac{z}{\sqrt{\e}\eta'} \big]\frac{\eta'}{\sqrt{\e}}} \Big<\big(\textbf{G}^{aa}(u)-\big<\textbf{G}^{aa}\big>\big) \textbf{T}^{\xi,\e}( \e u)(y),\textbf{T}^{\xi,\e}(\e u)(y)- \tilde{\textbf{T}}^{\xi,\e}(\e u)(y)\Big>_{\mathcal{H}_\xi}   du\Big\rvert \\
&\quad+\Big\lvert \e \int_{\big[\frac{z}{\sqrt{\e}\eta'} \big]\frac{\eta'}{\sqrt{\e}}} ^{z/\e}  \Big<\big(\textbf{G}^{aa}(u)-\big<\textbf{G}^{aa}\big>\big) \textbf{T}^{\xi,\e}(\e u)(y),\textbf{T}^{\xi,\e}(\e u)(y)- \tilde{\textbf{T}}^{\xi,\e}(\e u)(y)\Big>_{\mathcal{H}_\xi}   du\Big\rvert,
\end{split}\]
with
\[\begin{split}
\Big\lvert \e \int_{\big[\frac{z}{\sqrt{\e}\eta'} \big]\frac{\eta'}{\sqrt{\e}}}^{z/\e} & \Big<\big(\textbf{G}^{aa}(u)-\big<\textbf{G}^{aa}\big>\big) \textbf{T}^{\xi,\e}(\e u)(y),\textbf{T}^{\xi,\e}(\e u)(y)- \tilde{\textbf{T}}^{\xi,\e}(\e u)(y)\Big>_{\mathcal{H}_\xi}   du\Big\rvert\\
& \leq \left[ \e^{1/4}\sqrt{ \eta'}\left(\int_0^L \big\| \textbf{G}^{aa}\left(\frac{u}{\e}\right) \big\|^2 du\right)^{1/2} +\sqrt{\e}\eta' \| \big<\textbf{G}^{aa}\big>\|\right]\\
&\quad\times \sup_{z\in[0,L]}\|\textbf{T}^{\xi,\e}(z)(y)\|_{\mathcal{H}_\xi} \|\textbf{T}^{\xi,\e}(z)(y)-\tilde{\textbf{T}}^{\xi,\e}(z)(y)\|_{\mathcal{H}_\xi}
\end{split}\]
since $0\leq z-\big[\frac{z}{\sqrt{\e}\eta'}\big]\sqrt{\e}\eta' \leq \sqrt{\e}\eta'$,
and
\[\begin{split}
\Big\lvert &\e \int_0 ^{\big[\frac{z}{\sqrt{\e}\eta'} \big]\frac{\eta'}{\sqrt{\e}}} \Big<\big(\textbf{G}^{aa}(u)-\big<\textbf{G}^{aa}\big>\big) \textbf{T}^{\xi,\e}( \e u)(y),\textbf{T}^{\xi,\e}(\e u)(y)- \tilde{\textbf{T}}^{\xi,\e}(\e u)(y)\Big>_{\mathcal{H}_\xi}   du\Big\rvert\\
&\leq \sqrt{\e}\!\!\!\sum_{m=0}^{\big[\frac{L}{\sqrt{\e}\eta'} \big]-1}\hspace{-0.2cm} \Big\lvert \sqrt{\e} \int_{m\frac{\eta'}{\sqrt{\e}}} ^{(m+1)\frac{\eta'}{\sqrt{\e}}}  \Big<\big(\textbf{G}^{aa}(u)-\big<\textbf{G}^{aa}\big>\big) \textbf{T}^{\xi,\e}( \e u)(y),\textbf{T}^{\xi,\e}(\e u)(y)- \tilde{\textbf{T}}^{\xi,\e}(\e u)(y)\Big>_{\mathcal{H}_\xi}   du\Big\rvert. 
\end{split}\]
Moreover,
\[\textbf{T}^{\xi,\e}( \e u)(y)=\textbf{T}^{\xi,\e}(m\eta' \sqrt{\e})(y)+\int_{m\frac{\eta'}{\sqrt{\e}}}^{u} \sqrt{\e}\textbf{H}^{aa}(v)\textbf{T}^{\xi,\e}( \e v)(y) +\e \textbf{G}^{aa}(v)\textbf{T}^{\xi,\e}( \e v)(y)dv\]
and 
\[\tilde{\textbf{T}}^{\xi,\e}( \e u)(y)=\tilde{\textbf{T}}^{\xi,\e}(m\eta' \sqrt{\e})(y)+\int_{m\frac{\eta'}{\sqrt{\e}}}^{u} \sqrt{\e}\textbf{H}^{aa}(v)\tilde{\textbf{T}}^{\xi,\e}( \e v)(y) +\e \big<\textbf{G}^{aa}\big>\tilde{\textbf{T}}^{\xi,\e}( \e v)(y)dv.\]
Therefore, we have 
\[\begin{split}
\sqrt{\e} &\int_{m\frac{\eta'}{\sqrt{\e}}} ^{(m+1)\frac{\eta'}{\sqrt{\e}}}  \Big<\big(\textbf{G}^{aa}(u)-\big<\textbf{G}^{aa}\big>\big)\textbf{T}^{\xi,\e}( \e u)(y),\textbf{T}^{\xi,\e}(\e u)(y)- \tilde{\textbf{T}}^{\xi,\e}(\e u)(y)\Big>_{\mathcal{H}_\xi}   du\\
&\hspace{-0.3cm}= \sqrt{\e} \hspace{-0.1cm}\int_{m\frac{\eta'}{\sqrt{\e}}} ^{(m+1)\frac{\eta'}{\sqrt{\e}}} \hspace{-0.2cm} \Big<\big(\textbf{G}^{aa}(u)-\big<\textbf{G}^{aa}\big>\big) \textbf{T}^{\xi,\e}(m\eta' \sqrt{\e})(y),\textbf{T}^{\xi,\e}(m\eta' \sqrt{\e})(y)- \tilde{\textbf{T}}^{\xi,\e}(m\eta' \sqrt{\e})(y)\Big>_{\mathcal{H}_\xi} \hspace{-0.3cm} du\\
&+\int_{m\frac{\eta'}{\sqrt{\e}}} ^{(m+1)\frac{\eta'}{\sqrt{\e}}} \int_{m\frac{\eta'}{\sqrt{\e}}}^{u} \e \Big<\big(\textbf{G}^{aa}(u)-\big<\textbf{G}^{aa}\big>\big) \textbf{H}^{aa}(v)\textbf{T}^{\xi,\e}(\e v)(y),\textbf{T}^{\xi,\e}(\e u)(y)- \tilde{\textbf{T}}^{\xi,\e}(\e u)(y)\Big>_{\mathcal{H}_\xi} \\
&\hspace{1.5cm}+\e^{3/2} \Big<\big(\textbf{G}^{aa}(u)-\big<\textbf{G}^{aa}\big>\big) \textbf{G}^{aa}(v)\textbf{T}^{\xi,\e}(\e v)(y),\textbf{T}^{\xi,\e}(\e u)(y)- \tilde{\textbf{T}}^{\xi,\e}(\e u)(y)\Big>_{\mathcal{H}_\xi} \\
&\hspace{1.5cm}+\e  \Big<\big(\textbf{G}^{aa}(u)-\big<\textbf{G}^{aa}\big>\big)\textbf{T}^{\xi,\e}(m\eta' \sqrt{\e})(y), \textbf{H}^{aa}(v)\big(\textbf{T}^{\xi,\e}(\e v)(y)- \tilde{\textbf{T}}^{\xi,\e}(\e v)(y)\big)\Big>_{\mathcal{H}_\xi} \\
&+\e^{3/2}  \Big<\big(\textbf{G}^{aa}(u)-\big<\textbf{G}^{aa}\big>\big)\textbf{T}^{\xi,\e}(m\eta' \sqrt{\e})(y), \textbf{G}^{aa}(v)\textbf{T}^{\xi,\e}(\e v)(y)- \big<\textbf{G}^{aa}\big>\tilde{\textbf{T}}^{\xi,\e}(\e v)(y)\Big>_{\mathcal{H}_\xi} du.
\end{split}\]
Consequently, by the Gronwall's inequality 
\[ \sup_{z\in[0,L]}\|\textbf{T}^{\xi,\e}(z)(y)-\tilde{\textbf{T}}^{\xi,\e}(z)(y) \|^2_ {\mathcal{H}_\xi}\leq B(\e,\eta')e^{2\big\|\big<G^{aa}\big> \big\|L},\]
where
\[\begin{split}
B&(\e,\eta')= 2\left[ \e^{1/4}\sqrt{ \eta'}\left(\int_0^L \big\| \textbf{G}^{aa}\left(\frac{u}{\e}\right) \big\|^2 du\right)^{1/2} +\sqrt{\e}\eta' \| \big<\textbf{G}^{aa}\big>\|\right]\\
&\hspace{5cm}\times \sup_{z\in[0,L]}\|\textbf{T}^{\xi,\e}(z)(y)\|_{\mathcal{H}_\xi} \|\textbf{T}^{\xi,\e}(z)(y)-\tilde{\textbf{T}}^{\xi,\e}(z)(y)\|_{\mathcal{H}_\xi}\\
& +2\sqrt{\e}\sum_{m=0}^{\big[\frac{L}{\sqrt{\e}\eta'} \big]-1} \int_{m\frac{\eta'}{\sqrt{\e}}} ^{(m+1)\frac{\eta'}{\sqrt{\e}}} \int_{m\frac{\eta'}{\sqrt{\e}}}^{u}\Big( 2 \e \big[\|\textbf{G}^{aa}(u)\| +\|\big<G^{aa}\big>\|\big]\|\textbf{H}^{aa}(v) \|\\
&+ \e^{3/2}\big[\|\textbf{G}^{aa}(u)\| +\|\big<G^{aa}\big>\|\big]^2\Big) \sup_{z\in[0,L]}\|\textbf{T}^{\xi,\e}(z)(y)\|_{\mathcal{H}_\xi}\|\textbf{T}^{\xi,\e}(z)(y)-\tilde{\textbf{T}}^{\xi,\e}(z)(y)\|_{\mathcal{H}_\xi}dv\,du\\
&\hspace{-0.2cm}+\Big\lvert \sqrt{\e} \hspace{-0.1cm}\int_{m\frac{\eta'}{\sqrt{\e}}} ^{(m+1)\frac{\eta'}{\sqrt{\e}}}\hspace{-0.2cm}  \Big<\big(\textbf{G}^{aa}(u)-\big<\textbf{G}^{aa}\big>\big) \textbf{T}^{\xi,\e}(m\eta' \sqrt{\e})(y),\textbf{T}^{\xi,\e}(m\eta' \sqrt{\e})(y)- \tilde{\textbf{T}}^{\xi,\e}(m\eta' \sqrt{\e})(y)\Big>_{\mathcal{H}_\xi}  \hspace{-0.3cm} du\Big\rvert,
\end{split}\]
and
\[\mathbb{P}\Big( \sup_{z\in[0,L]}\|\textbf{T}^{\xi,\e}(z)(y)-\tilde{\textbf{T}}^{\xi,\e}(z)(y) \|^2_ {\mathcal{H}_\xi} > \eta \Big)\leq \mathbb{P}\Big(B(\e,\eta')\geq\eta e^{-2\big\|\big<G^{aa}\big> \big\|L}\Big).\]
Setting $\eta''=\eta e^{-2\big\|\big<G^{aa}\big> \big\|L}$, we have
\[\begin{split}
\mathbb{P}\Big(B(\e,\eta')\geq \eta''\Big)&\leq \mathbb{P}\Big(B(\e,\eta')\geq \eta'', \sup_{z\in[0,L]}\|\textbf{T}^{\xi,\e}(z)(y) \|^2_ {\mathcal{H}_\xi}\leq M \Big)\\
&\quad+\mathbb{P}\Big(\sup_{z\in[0,L]}\|\textbf{T}^{\xi,\e}(z)(y)\|^2_ {\mathcal{H}_\xi}\geq M \Big).
\end{split}\]
We already know that the process $ \tilde{\textbf{T}}^{\xi,\e}(.)(y)$ is bounded.
Moreover,
\[\mathbb{P}\Big(B(\e,\eta')\geq \eta'', \sup_{z\in[0,L]}\|\textbf{T}^{\xi,\e}(z)(y) \|^2_ {\mathcal{H}_\xi}\leq M \Big)\leq \frac{1}{\eta''}\mathbb{E}\left[B(\e,\eta')\textbf{1}_{\big(\sup_{z\in[0,L]}\|\textbf{T}^{\xi,\e}(z)(y) \|^2_ {\mathcal{H}_\xi}\leq M\big)}\right]\]
with
\[\begin{split}
\mathbb{E}&\left[B(\e,\eta')\textbf{1}_{\big(\sup_{z\in[0,L]}\|\textbf{T}^{\xi,\e}(z)(y) \|^2_ {\mathcal{H}_\xi}\leq M\big)}\right]\leq K\big[\eta'^2+\e^{1/4}\sqrt{\eta'}+ \sqrt{\e}(\eta'+\eta'^2)\big]\\
& \quad+2\sqrt{\e}\sum_{m=0}^{\big[\frac{L}{\sqrt{\e}\eta'} \big]-1}
\mathbb{E}\Big[\textbf{1}_{\big(\sup_{z\in[0,L]}\|\textbf{T}^{\xi,\e}(z)(y) \|^2_ {\mathcal{H}_\xi}\leq M\big)}\\
&\hspace{-0.1cm}\times\Big\lvert \sqrt{\e}\hspace{-0.1cm} \int_{m\frac{\eta'}{\sqrt{\e}}} ^{(m+1)\frac{\eta'}{\sqrt{\e}}} \hspace{-0.2cm} \Big<\big(\textbf{G}^{aa}(u)-\big<\textbf{G}^{aa}\big>\big)\textbf{T}^{\xi,\e}(m\eta' \sqrt{\e})(y),\textbf{T}^{\xi,\e}(m\eta' \sqrt{\e})(y)- \tilde{\textbf{T}}^{\xi,\e}(m\eta' \sqrt{\e})(y)\Big>_{\mathcal{H}_\xi}  \hspace{-0.3cm} du\Big\rvert\Big]\\
&\leq K\big[\eta'^2+\e^{1/4}\sqrt{\eta'}+ \sqrt{\e}(\eta'+\eta'^2)\big]\\
&\quad+ 2\sqrt{\e} K \sum_{m=0}^{\big[\frac{L}{\sqrt{\e}\eta'} \big]-1} \E\left[ \Big\| \sqrt{\e} \int_{m\frac{\eta'}{\sqrt{\e}}} ^{(m+1)\frac{\eta'}{\sqrt{\e}}} \textbf{G}^{aa}(u)-\big<\textbf{G}^{aa}\big>\,du  \Big\|^2 \right]^{1/2},
\end{split}\]
since
$\int_{m\frac{\eta'}{\sqrt{\e}}} ^{(m+1)\frac{\eta'}{\sqrt{\e}}} \int_{m\frac{\eta'}{\sqrt{\e}}}^{u}dv\, du=\frac{\eta'}{\e}$ and 
\begin{equation}\label{ineqGP2}
\mathbb{E}\left[\big\| \textbf{G}^{aa}\left(\frac{u}{\e}\right) \big\|^2 \right]\leq 
K \,\E\left[\Big( \int_0^d\lvert V(x,0) \rvert^2dx \Big)^2\right]\end{equation}
for $u\in[0,L]$. As a result, it remains us to estimate only one term.
\begin{lem}\label{lemtechnP2}
\[\lim_{\e \to 0}\sqrt{\e}  \sum_{m=0}^{\big[\frac{L}{\sqrt{\e}\eta'} \big]-1}\E\left[ \Big\| \sqrt{\e} \int_{m\frac{\eta'}{\sqrt{\e}}} ^{(m+1)\frac{\eta'}{\sqrt{\e}}} \emph{\textbf{G}}^{aa}(u)-\big<\emph{\textbf{G}}^{aa}\big>\,du  \Big\|^2 \right]^{1/2}=0.\]
\end{lem}
\begin{preuve}[of Lemma \ref{lemtechnP2}]
Let us remark that we have the following decomposition. For each $j\in\big\{1,\dots,N\big\}$, almost every $\ga\in(\xi,k^2)$, and $\forall y\in \mathcal{H}_{\xi}$,
\[\begin{split}
\textbf{G}^{aa}_j(z)(y)&=\sum_{l=1}^N \textbf{G}^{aa}_{jl}(z)y_l +\int_\xi ^{k^2} \textbf{G}^{aa}_{j\ga'}(z)y_{\ga'} d\ga',\\
\textbf{G}^{aa}_\ga(z)(y)&=\sum_{l=1}^N \textbf{G}^{aa}_{\ga l}(z)y_l +\int_\xi ^{k^2} \textbf{G}^{aa}_{\ga\ga'}(z)y_{\ga'} d\ga'.
\end{split}\]
Letting 
\[ \textbf{P}=\sqrt{\e} \int_{m\frac{\eta'}{\sqrt{\e}}} ^{(m+1)\frac{\eta'}{\sqrt{\e}}} \textbf{G}^{aa}(u)-\big<\textbf{G}^{aa}\big>\,du, \]
we have $(j,l)^2\in\big\{1,\dots,N\big\}^2$ such that $j\not=l$, and almost every $\ga\in(\xi,k^2)$
\[\begin{split}
\textbf{P}_{jj}&=\sqrt{\e} \int_{m\frac{\eta'}{\sqrt{\e}}} ^{(m+1)\frac{\eta'}{\sqrt{\e}}} \textbf{G}^{aa}_{jj}(u)-\big<\textbf{G}^{aa}\big>_{jj}\,du,\\
\textbf{P}_{jl}&=\sqrt{\e} \int_{m\frac{\eta'}{\sqrt{\e}}} ^{(m+1)\frac{\eta'}{\sqrt{\e}}} \textbf{G}^{aa}_{jl}(u)du,\\
\textbf{P}_{j\ga'}&=\sqrt{\e} \int_{m\frac{\eta'}{\sqrt{\e}}} ^{(m+1)\frac{\eta'}{\sqrt{\e}}} \textbf{G}^{aa}_{j \ga'}(u)du,\\
\textbf{P}_{\ga l}&=\sqrt{\e} \int_{m\frac{\eta'}{\sqrt{\e}}} ^{(m+1)\frac{\eta'}{\sqrt{\e}}} \textbf{G}^{aa}_{\ga l}(u)du,\\
\textbf{P}_{\ga\ga'}&=\sqrt{\e} \int_{m\frac{\eta'}{\sqrt{\e}}} ^{(m+1)\frac{\eta'}{\sqrt{\e}}} \textbf{G}^{aa}_{\ga \ga'}(u)du,
\end{split}\]
and
\[\frac{1}{2}\|\textbf{P}\|^2\leq \sum_{j,l=1}^N\lvert \textbf{P}_{jl} \rvert^2 +\sum_{j=1}^N\int_\xi^{k^2}\lvert\textbf{P}_{j\ga'}\rvert^2 d\ga'+\int_\xi^{k^2}\sum_{l=1}^N\lvert\textbf{P}_{\ga l}\rvert^2 d\ga+\int_\xi^{k^2}\int_\xi^{k^2}\lvert\textbf{P}_{\ga \ga'}\rvert^2 d\ga'\,d\ga. \]
Moreover,
\[\begin{split}
\E\Big[V(x_1,z_1)V(x_2,z_2)V(x_3,z_3)V(x_4,z_4)\Big]=&E\Big[V(x_1,z_1)V(x_2,z_2)\Big]E\Big[V(x_3,z_3)V(x_4,z_4)\Big]\\
&+E\Big[V(x_1,z_1)V(x_3,z_3)\Big]E\Big[V(x_2,z_2)V(x_4,z_4)\Big]\\
&+E\Big[V(x_1,z_1)V(x_4,z_4)\Big]E\Big[V(x_2,z_2)V(x_3,z_3)\Big]\\
&=\ga_0(x_1,x_2)\ga_0(x_3,x_4)e^{-a\lvert z_1-z_2 \rvert}e^{-a\lvert z_3-z_4 \rvert}\\
&+\ga_0(x_1,x_3)\ga_0(x_2,x_4)e^{-a\lvert z_1-z_3 \rvert}e^{-a\lvert z_2-z_4 \rvert}\\
&+\ga_0(x_1,x_4)\ga_0(x_2,x_3)e^{-a\lvert z_1-z_4 \rvert}e^{-a\lvert z_2-z_3 \rvert},
\end{split} \]
which is the fourth order moment of a Gaussian field. To compute the expectation of the square norm of $\textbf{P}$ we must know these moments. Following that decomposition, the square norm of $\textbf{P}$ can be decomposed in three parts. First,  after a long computation, the two parts corresponding to the two last terms of the previous decomposition are dominated by $\sqrt{\e}$ uniformly in $m$. Then, we focus our attention on the part corresponding to the first part of the previous decomposition. For $\E\big[\lvert \textbf{P}_{j\ga'}\rvert^2\big]$, $\E\big[\lvert\textbf{P}_{\ga l}\rvert^2\big]$, and $\E\big[\lvert\textbf{P}_{jl}\rvert^2\big]$ with $j\not=l$, we get after a long computation terms of the form
\[\begin{split}
\e\int_{m\frac{\eta'}{\sqrt{\e}}} ^{(m+1)\frac{\eta'}{\sqrt{\e}}}  \int_{m\frac{\eta'}{\sqrt{\e}}} ^{(m+1)\frac{\eta'}{\sqrt{\e}}} e^{i(\sqrt{\ga'}-\beta_j)u_1}e^{-i(\sqrt{\ga'}-\beta_j)u_2}du_1\,du_2&=\mathcal{O}(\e),\\
\e\int_{m\frac{\eta'}{\sqrt{\e}}} ^{(m+1)\frac{\eta'}{\sqrt{\e}}}  \int_{m\frac{\eta'}{\sqrt{\e}}} ^{(m+1)\frac{\eta'}{\sqrt{\e}}} e^{i(\beta_l-\sqrt{\ga})u_1}e^{-i(\beta_l-\sqrt{\ga})u_2}du_1\,du_2&=\mathcal{O}(\e),\\ \e\int_{m\frac{\eta'}{\sqrt{\e}}} ^{(m+1)\frac{\eta'}{\sqrt{\e}}}  \int_{m\frac{\eta'}{\sqrt{\e}}} ^{(m+1)\frac{\eta'}{\sqrt{\e}}} e^{i(\beta_l-\beta_j)u_1}e^{-i(\beta_l-\beta_j)u_2}du_1\,du_2&=\mathcal{O}(\e).
\end{split}\]
For $\E\big[\lvert\textbf{P}_{\ga \ga'}\rvert^2\big]$ we separate the integral into two parts.
\[ \int_\xi^{k^2}\int_\xi^{k^2}\E\big[\lvert\textbf{P}_{\ga \ga'}\rvert^2 \big]d\ga'\,d\ga=\int_{I_{\geq \mu}} \E\big[\lvert\textbf{P}_{\ga \ga'}\rvert^2\big] d\ga'\,d\ga+\int_{I_{<\mu}} \E\big[\lvert\textbf{P}_{\ga \ga'}\rvert^2\big] d\ga'\,d\ga, \]
where $\mu>0$ and
\[\begin{split}
I_{\geq \mu}&=  \Big\{ (\ga,\ga')\in(\xi,k^2)^2,\quad\big\lvert \sqrt{\ga}-\sqrt{\ga'}\big\rvert\geq \mu \Big\},\\
I_{<\mu}&=\Big\{ (\ga,\ga')\in(\xi,k^2)^2,\quad\big\lvert \sqrt{\ga}-\sqrt{\ga'}\big\rvert < \mu  \Big\}.
\end{split}\]
Consequently,
\[ \int_{I_{<\mu}} \E\big[\lvert\textbf{P}_{\ga \ga'}\rvert^2\big] d\ga'\,d\ga\leq K  \int_{I_{<\mu}}d\ga'\,d\ga, \]
and on $I_{\geq \mu}$ we get terms of the form
\[   \e\int_{I_{\geq \mu}}\int_{m\frac{\eta'}{\sqrt{\e}}} ^{(m+1)\frac{\eta'}{\sqrt{\e}}}  \int_{m\frac{\eta'}{\sqrt{\e}}} ^{(m+1)\frac{\eta'}{\sqrt{\e}}} e^{i(\sqrt{\ga'}-\sqrt{\ga})u_1}e^{-i(\sqrt{\ga'}-\sqrt{\ga})u_2}du_1\,du_2d\ga'd\ga=\mathcal{O}(\e).\]
Now, it remains us to study $\E\big[\lvert\textbf{P}_{jj}\rvert^2\big]$. After a long computation, the terms of order one produced by $\textbf{G}^{aa}_{jj}$ are compensated by the terms of order one given by $\big<\textbf{G}^{aa}\big>_{j}$. Moreover, the other terms are dominated by $\sqrt{\e}$.

As a result, we get
\[\overline{\lim_{\e \to 0}}\,\, \sqrt{\e}\sum_{m=0}^{\big[\frac{L}{\sqrt{\e}\eta'} \big]-1}\E\left[ \Big\| \sqrt{\e} \int_{m\frac{\eta'}{\sqrt{\e}}} ^{(m+1)\frac{\eta'}{\sqrt{\e}}} \textbf{G}^{aa}(u)-\big<\textbf{G}^{aa}\big>\,du  \Big\|^2 \right]^{1/2}\leq K \sqrt{ \int_{I_{<\mu}}d\ga'\,d\ga} \]
and one can conclude the proof of Lemma \ref{lemtechnP2} by letting $\mu\to 0$.$\square$
\end{preuve}

From the previous lemma, we finally get, $\forall \eta' >0$
\[\overline{\lim_{\e\to 0}}\,\,\mathbb{P}\Big( \sup_{z\in[0,L]}\|\textbf{T}^{\xi,\e}(z)(y)-\tilde{\textbf{T}}^{\xi,\e}(z)(y) \|^2_ {\mathcal{H}_\xi} > \eta \Big)\leq \frac{e^{2\big\|\big<G^{aa}\big> \big\|L}}{\eta}K\eta'^2,\]
since using the Gronwall's inequality and \eqref{ineqGP2} we have
\[\lim_{M\to +\infty}\overline{\lim_{\e\to 0}}\,\,\mathbb{P}\Big(\sup_{z\in[0,L]}\|\textbf{T}^{\xi,\e}(z)(y)\|^2_ {\mathcal{H}_\xi}\geq M \Big)=0.\]
Consequently, we conclude the proof of Proposition \ref{propoP2} by letting $\eta'\to0$.
$\blacksquare$
\end{preuve}

According to Proposition \ref{propoP2}, to study the convergence in distribution of the process $\big(\textbf{T}^{\xi,\e}(.)(y)\big)_{\e}$ it suffices to study the convergence for $\big(\tilde{\textbf{T}}^{\xi,\e}(.)(y)\big)_\e$. Moreover, we shall consider the complex case for more convenient manipulations. Letting $\lambda \in\mathcal{H}_\xi$, we consider the equation
\begin{equation*}
\frac{d}{dt} \tilde{\textbf{T}}^{\xi,\e} _\lambda (t)(y)=\frac{1}{\sqrt{\e}}H _{\lambda}\left(\tilde{\textbf{T}}^{\xi,\e} (t)(y),C\left(\frac{t}{\e}\right),\frac{t}{\e}\right)+G_\lambda\left(\tilde{\textbf{T}}^{\xi,\e} (t)(y)\right),
\end{equation*}
with $H_\lambda=\big<H,\lambda\big>_{\mathcal{H}_\xi}$, $G_\lambda=\big<\big<\textbf{G}^{aa}\big>(.),\lambda\big>_{\mathcal{H}_\xi}$,
where, for $j\in\big\{1,\dots,N\big\}$ and almost every $\ga \in (\xi,k^2)$
\begin{equation*}\begin{split}
H _{j}\left(\textbf{T},C,s \right) &=\frac{i k^2}{2}\Big[\sum_{l=1}^N \frac{C_{jl}}{\sqrt{\beta_j \beta_l}}e^{i(\beta_l -\beta_j)s}\textbf{T}_l + \int_{\xi}^{k^2}\frac{C_{j\ga'}}{\sqrt{\beta_j \sgap}}e^{i(\sgap -\beta_j)s}\textbf{T}_{\ga'}d\ga' \Big],\\
H _{\ga}\left(\textbf{T},C,s \right)&=\frac{i k^2}{2}\Big[\sum_{l=1}^N \frac{C_{\ga l}}{\sqrt{\sga \beta_l}} e^{i(\beta_l -\sga)s}\textbf{T}_l+ \int_{\xi}^{k^2}\frac{C_{\ga \ga'}}{\ga^{1/4}{\ga'}^{1/4}}e^{i(\sgap -\sga)s}\textbf{T}_{\ga'}d\ga' \Big].
\end{split}\end{equation*}

The proof of Theorem \ref{thasymptP21} is based on the perturbed-test-function approach. Using the notion of a pseudogenerator, we prove tightness and characterize all subsequence limits.

\subsubsection{Pseudogenerator}\label{pseudogene}

We recall the techniques developed by Kurtz and Kushner \cite{kushner}. Let $\mathcal{M}^\e $ be the set of all $\mathcal{F}^\e$-measurable functions $f(t)$ for which $\sup_{t\leq T} \mathbb{E}\left[\lvert f(t) \rvert \right] <+\infty $ and where $T>0$ is fixed. The $p-\lim$ and the pseudogenerator are defined as follows. Let $f$ and $f^\delta$ in $\mathcal{M}^\e $ $\forall\delta>0$. We say that $f=p-\lim_\delta f^\delta$ if
\[ \sup_{t, \delta }\mathbb{E}[\lvert f^\delta(t)\rvert]<+\infty\quad \text{and}\quad \lim_{\delta\rightarrow 0}\mathbb{E}[\lvert f^\delta (t) -f(t)\rvert]=0 \quad \forall t.\]
The domain of $\mathcal{A}^\e$ is denoted by $\mathcal{D}\left(\mathcal{A}^\e\right)$. We say that $f\in \mathcal{D}\left(\mathcal{A}^\e\right)$ and $\mathcal{A}^\e f=g$ if $f$ and $g$ are in $\mathcal{D}\left(\mathcal{A}^\e\right)$ and 
\begin{equation*}
p-\lim_{\delta \to 0} \left[ \frac{\mathbb{E}^\e _t [f(t+\delta)]-f(t)}{\delta}-g(t) \right]=0,
\end{equation*}
where $\mathbb{E}^\e _t$ is the conditional expectation given $\mathcal{F}^\e _t$ and $\mathcal{F}^\e _t =\mathcal{F}_{t/\e}$. A useful result about $\mathcal{A}^\e$ is given by the following theorem.
\begin{thm}\label{martingaleP2}
Let $f\in \mathcal{D}\left(\mathcal{A}^\e\right)$. Then
\begin{equation*}
M_f ^\e (t)=f(t)-\int _0 ^t  \mathcal{A}^\e f(u)du
\end{equation*}
is an $\left( \mathcal{F}^\e _t \right)$-martingale.
\end{thm}

\subsubsection{Tightness}

We consider the classical complex derivative with the following notation: If $v=\alpha +i\beta$, then $\partial_v =\frac{1}{2}\left(\partial_\alpha -i \partial_ \beta \right)$ and $\partial_{\overline{v}}=\frac{1}{2}\left(\partial_\alpha +i \partial_ \beta \right)$.
\begin{prop}\label{tightth1P2}
$\forall \lambda \in \mathcal{H}_\xi$, the family $\big(\tilde{\emph{\textbf{T}}}^{\xi,\e} _{\lambda} (.)(y) \big)_{\e\in(0,1)}$ is tight on $\mathcal{D}\left([0,+\infty),\mathbb{C} \right)$.
\end{prop}
\begin{preuve}
According to Theorem 4 in \cite{kushner}, we need to show the three following lemmas. Let $\lambda \in \mathcal{H}_\xi$, $f$ be a smooth function, and $f_0 ^\e (t)=f\left(\tilde{\textbf{T}}^{\xi,\e} _\lambda (t)(y) \right)$. We have,
\begin{equation*}
\begin{split}
\mathcal{A}^\e f_0 ^\e (t)&= \partial_v f\left(\tilde{\textbf{T}}^{\xi,\e} _{\lambda} (t)(y) \right)\left[\frac{1}{\sqrt{\e}}H _{\lambda}\left(\tilde{\textbf{T}}^{\xi,\e}  (t)(y),C\left(\frac{t}{\e}\right),\frac{t}{\e} \right)+G_\lambda\left(\tilde{\textbf{T}}^{\xi,\e} (t)(y)\right)\right] \\
&+ \partial_{\overline{v}} f\left(\tilde{\textbf{T}}^{\xi,\e} _{\lambda} (t)(y)\right)\overline{ \left[\frac{1}{\sqrt{\e}}H _{\lambda}\left(\tilde{\textbf{T}}^{\xi,\e}  (t)(y),C\left(\frac{t}{\e}\right),\frac{t}{\e} \right)+G_\lambda\left(\tilde{\textbf{T}}^{\xi,\e} (t)(y)\right)\right]}.
\end{split}
\end{equation*}
Let
\begin{equation*}
\begin{split}
f^\e _1 (t)&=\frac{1}{\sqrt{\e}}\partial_v f\left(\tilde{\textbf{T}}^{\xi,\e} _{\lambda} (t)(y) \right)\int_{t}^{+\infty} \mathbb{E}^\e _t\left[ H _{\lambda}\left(\tilde{\textbf{T}}^{\xi,\e}  (t)(y),C\left(\frac{u}{\e}\right),\frac{u}{\e} \right)\right] du \\
&+\frac{1}{\sqrt{\e}}\partial_{\overline{v}} f\left(\tilde{\textbf{T}}^{\xi,\e} _{\lambda} (t)(y) \right)\int_{t}^{+\infty} \mathbb{E}^\e _t\left[ \overline{H_{\lambda}\left(\tilde{\textbf{T}}^{\xi,\e} (t)(y),C\left(\frac{u}{\e}\right),\frac{u}{\e} \right)}\right]du. \\
\end{split}
\end{equation*}
\begin{lem}\label{f1P2}  $\forall T>0$,
$\lim_{\e} \sup_{0\leq t \leq T} \lvert f^\e _1 (t)\rvert =0$ almost surely, and $\sup_{t\geq 0}\mathbb{E}\left[\lvert f^\e _1 (t) \rvert \right]=\mathcal{O}\left(\sqrt{\e} \right)$.
\end{lem}
\begin{preuve}[of Lemma \ref{f1P2}]
Using the Markov property of the Gaussian field $V$, we have
\[f^{\e}_1(t)= \frac{i k^2 \sqrt{\e}}{2}\partial_{v}f\left( \tilde{\textbf{T}}^{\xi,\e} _\lambda (t)(y)\right)F^\e_{1,\lambda}(t)- \frac{i k^2 \sqrt{\e}}{2}\partial_{\overline{v}}f\left(\tilde{ \textbf{T}}^{\xi,\e} _\lambda (t)(y)\right)\overline{F^\e_{1,\lambda}(t)}\]
with
\begin{equation*}
\begin{split}
F^\e _{1,\lambda} (t)&=\sum_{j=1}^N \left[  \sum_{l=1}^N \frac{C_{jl}\left(\frac{t}{\e}\right)}{ \sqrt{\beta _j\beta _l}}e^{i(\beta _l-\beta _j)\frac{t}{\e}}\right. \tilde{\textbf{T}}^{\xi,\e}_{l} (t)(y)\frac{a+i(\beta _l-\beta _j)}{a^2+(\beta _l-\beta _j)^2}\\
&\left.\hspace{1cm}+\int_\xi ^{k^2} \frac{C_{j \ga'}\left(\frac{t}{\e}\right)}{ \sqrt{\beta^\e _j \sgap}}e^{i(\sgap-\beta _j)\frac{t}{\e}} \tilde{\textbf{T}}^{\xi,\e}_{\ga' } (t)(y)\frac{a+i(\sgap-\beta _j)}{a^2+(\sgap-\beta _j)^2}d\ga'\right]\overline{\lambda_{j}}\\
&\quad+\int_\xi ^{k^2}\left[  \sum_{l=1}^N \frac{C_{\ga l}\left(\frac{t}{\e}\right)}{ \sqrt{\sga\beta _l}}e^{i(\beta _l-\sga)\frac{t}{\e}}\right. \tilde{\textbf{T}}^{\xi,\e}_{l} (t)(y)\frac{a+i(\beta _l-\sga)}{a^2+(\beta _l-\sga)^2}\\
&\left.\hspace{1cm}+\int_\xi ^{k^2} \frac{C_{\ga \ga'}\left(\frac{t}{\e}\right)}{ \ga^{1/4} {\ga'}^{1/4}}e^{i(\sgap-\sga)\frac{t}{\e}} \tilde{\textbf{T}}^{\xi,\e}_{\ga' } (t)(y)\frac{a+i(\sgap-\sga)}{a^2+(\sgap-\sga)^2}d\ga'\right]\overline{\lambda_{\ga}}d\ga.
\end{split}
\end{equation*}
Using \eqref{cg2}, we easily get
\[\mathbb{E}\left[\lvert f^\e _1 (t) \rvert \right] \leq \sqrt{\e}K(f,\lambda).\]
and
\[ \lvert f^\e _1 (t) \rvert \leq K(\lambda,f) \sqrt{\e} \sup_{0\leq t \leq T/\e} \sup_{x\in [0,d]} \left\lvert V\left(x,t \right)\right\rvert.\]
Then, we can conclude with \eqref{cg1}.$\square$
\end{preuve}
\begin{lem}\label{A1P2}
$\left\{\mathcal{A}^\e \left(f^\e _0 +f^\e _1\right)(t), \e \in(0,1), 0\leq t\leq T\right\}$ is uniformly integrable $\forall T>0$.
\end{lem}
\begin{preuve}[of Lemma \ref{A1P2}]
A computation gives us
\[\mathcal{A}^\e \left(f^\e _0 +f^\e _1\right)(t)=\tilde{F}_{\lambda}\left(\tilde{\textbf{T}}^{\xi,\e} (t)(y),C\left(\frac{t}{\e}\right)\otimes C\left(\frac{t}{\e}\right),\frac{t}{\e} \right),\] 
where 
\[C(T)\otimes C(T)_{q_1\,q_2\,q_3\,q_4}=C_{q_1\,q_2}(T)C_{q_3\,q_4}(T)\]
for $(q_1,q_2,q_3,q_4)\in\big(\{1,\dots,N\}\cup(\xi,k^2)\big)^4$,
with
\begin{equation*}\begin{split}
\tilde{F}_{\lambda}\left(\textbf{T},C,s \right)&= \partial_{v}f(\textbf{T})\left[\tilde{F}^{1,\e} _\lambda (\textbf{T},C,s) +G_\lambda\left(\textbf{T}\right)\right]+\partial_{\overline{v}}f(\textbf{T})\overline{\left[\tilde{F}^{1} _\lambda (\textbf{T},C,s)+G_\lambda\left(\textbf{T}\right)\right]} \\
&+ \partial ^2 _{v}f(\textbf{T})\tilde{F}^{2} _\lambda (\textbf{T},C,s) +\partial ^2 _{\overline{v}}f(\textbf{T})\overline{\tilde{F}^{2} _\lambda (\textbf{T},C,s)} \\
&+ \partial _{\overline{v}}\partial _{v}f(\textbf{T})\tilde{F}^{3} _\lambda (\textbf{T},C,s) +\partial _{v}\partial _{\overline{v}}f(\textbf{T})\overline{\tilde{F}^{3} _\lambda (\textbf{T},C,s)},
\end{split}\end{equation*}
and
\[\begin{split}
\tilde{F}^{1} _\lambda (\textbf{T},C,s)=&\\
& -\frac{k^4}{4}\sum_{j=1}^N \left[\sum_{l ,l'=1}^{N}\frac{C_{jlll'}}{\sqrt{\beta _j \beta _l ^2 \beta _{l'}}}e^{i (\beta _{l'}-\beta _{j})s} \textbf{T}_{l'} \frac{a+i(\beta _l-\beta_j)}{a^2+ (\beta _l-\beta_j) ^2}\right.\\
&+\sum_{l=1}^N\int_{\xi}^{k^2} \frac{C_{jll\ga''}}{ \sqrt{\beta _j \beta _l^2  \sqrt{\ga''}}}e^{i (\sqrt{\ga''} -\beta _j)s} \textbf{T}_{\ga'' } \frac{a+i(\beta_l-\beta_j)}{a^2+ (\beta _l-\beta_j) ^2}d\ga'' \\
&+\int_\xi^{k^2}\sum_{l'=1}^N\frac{C_{j\ga'\ga' l'}}{ \sqrt{\beta _j \ga' \beta _{l' }}}e^{i (\beta _{l'} -\beta_{j})s} \textbf{T}_{l'} \frac{a+i(\sgap-\beta_j)}{a^2+ (\sgap-\beta_j) ^2}d\ga' \\
&\left.+\int_\xi^{k^2}\int_\xi^{k^2}\frac{C_{j\ga'\ga'\ga''}}{ \sqrt{\beta _j \ga' \sqrt{\ga''}}}e^{i (\sqrt{\ga''} -\beta _{j})s} \textbf{T}_{\ga'' }  \frac{a+i(\sgap -\beta_j)}{a^2+ (\sgap-\beta_j) ^2}d\ga' d\ga'' \right]\overline{\lambda_{j}}\end{split}\]
\[\begin{split}
\hspace{1.5cm}&-\frac{k^4}{4}\int_\xi^{k^2} \left[\sum_{l ,l'=1}^{N}\frac{C_{\ga lll'}}{\sqrt{\sga \beta _l ^2 \beta _{l'}}}e^{i (\beta _{l'}-\sga)s} \textbf{T}_{l'} \frac{a+i(\beta _l-\sga)}{a^2+ (\beta _l-\sga) ^2}\right.\\
&+\sum_{l=1}^N\int_{\xi}^{k^2} \frac{C_{\ga ll\ga''}}{ \sqrt{\sga \beta _l^2  \sqrt{\ga''}}}e^{i (\sqrt{\ga''} -\sga)s} \textbf{T}_{\ga'' } \frac{a+i(\beta_l-\sga)}{a^2+ (\beta _l-\sga) ^2}d\ga'' \\
&+\int_\xi^{k^2}\sum_{l'=1}^N\frac{C_{\ga \ga'\ga' l'}}{ \sqrt{\sga \ga' \beta _{l' }}}e^{i (\beta _{l'} -\sga)s} \textbf{T}_{l'} \frac{a+i(\sgap-\sga)}{a^2+ (\sgap-\sga) ^2}d\ga' \\
&\left.+\int_\xi^{k^2}\int_\xi^{k^2}\frac{C_{\ga\ga'\ga'\ga''}}{ \sqrt{\sga \ga' \sqrt{\ga''}}}e^{i (\sqrt{\ga''} -\sga)s} \textbf{T}_{\ga'' }  \frac{a+i(\sgap -\beta_j)}{a^2+ (\sgap-\sga) ^2}d\ga' d\ga'' \right]\overline{\lambda_{\ga}}d\ga,
\end{split}\]

\[\begin{split}
\tilde{F}^{2} _\lambda &(\textbf{T},C,s)\\
&= -\frac{k^4}{4}\sum_{j,j'=1}^N \left[\sum_{l ,l'=1}^{N}\frac{C_{jlj'l'}}{\sqrt{\beta _j \beta _l \beta_{j'} \beta _{l'}}}e^{i (\beta _{l}-\beta _{j}+\beta _{l'}-\beta _{j'})s} \textbf{T}_{l}\textbf{T}_{l'} \frac{a+i(\beta _l-\beta_j)}{a^2+ (\beta _l-\beta_j) ^2}\right.\\
&+\sum_{l=1}^N\int_{\xi}^{k^2} \frac{C_{jlj'\ga'_2}}{ \sqrt{\beta _j \beta _l \beta_{j'}  \sqrt{\ga'_2}}}e^{i(\beta _{l}-\beta _{j}+\sqrt{\ga'_2}-\beta _{j'})s} \textbf{T}_{l}\textbf{T}_{\ga'_2} \frac{a+i(\beta_l-\beta_j)}{a^2+ (\beta _l-\beta_j) ^2}d\ga'_2 \\
&+\int_\xi^{k^2}\sum_{l'=1}^N\frac{C_{j\ga'_1j' l'}}{ \sqrt{\beta _j \sqrt{\ga'_1}\beta_{j'} \beta _{l' }}}e^{i (\sqrt{\ga'_1}-\beta_j+\beta _{l'} -\beta_{j'})s} \textbf{T}_{\ga'_1}\textbf{T}_{l'} \frac{a+i(\sqrt{\ga'_1}-\beta_j)}{a^2+ (\sqrt{\ga'_1}-\beta_j) ^2}d\ga'_1 \\
&\left.+\int_\xi^{k^2}\int_\xi^{k^2}\frac{C_{j\ga'_1j'\ga'_2}}{ \sqrt{\beta _j \sqrt{\ga'_1} \beta_{j'}\sqrt{\ga'_2}}}e^{i (\sqrt{\ga'_1} -\beta _{j}+\sqrt{\ga'_2}-\beta_{j'})s} \textbf{T}_{\ga'_1}\textbf{T}_{\ga'_2}  \frac{a+i(\sqrt{\ga'_1} -\beta_j)}{a^2+ (\sqrt{\ga'_1}-\beta_j) ^2}d\ga'_1 d\ga'_2\right]\overline{\lambda_{j}\lambda_{j'}}\\
& -\frac{k^4}{4}\sum_{j=1}^N \int_\xi ^{k^2}\left[\sum_{l ,l'=1}^{N}\frac{C_{jl \ga_2 l'}}{\sqrt{\beta _j \beta _l \sqrt{\ga_2} \beta _{l'}}}e^{i (\beta _{l}-\beta _{j}+\beta _{l'}-\sqrt{\ga_2})s} \textbf{T}_{l}\textbf{T}_{l'} \frac{a+i(\beta _l-\beta_j)}{a^2+ (\beta _l-\beta_j) ^2}\right.\\
&+\sum_{l=1}^N\int_{\xi}^{k^2} \frac{C_{jl\ga_2\ga'_2}}{ \sqrt{\beta _j \beta _l   \sqrt\ga_2 {\ga'_2}}}e^{i(\beta _{l}-\beta _{j}+\sqrt{\ga'_2}-\sqrt{\ga_2})s} \textbf{T}_{l}\textbf{T}_{\ga'_2} \frac{a+i(\beta_l-\beta_j)}{a^2+ (\beta _l-\beta_j) ^2}d\ga'_2 \\
&+\int_\xi^{k^2}\sum_{l'=1}^N\frac{C_{j\ga'_1\ga_2 l'}}{ \sqrt{\beta _j \sqrt{\ga'_1\ga_2} \beta _{l' }}}e^{i (\sqrt{\ga'_1}-\beta_j+\beta _{l'} -\sqrt{\ga_2})s} \textbf{T}_{\ga'_1}\textbf{T}_{l'} \frac{a+i(\sqrt{\ga'_1}-\beta_j)}{a^2+ (\sqrt{\ga'_1}-\beta_j) ^2}d\ga'_1 \\
&\left. +\int_\xi^{k^2}\int_\xi^{k^2}\frac{C_{j\ga'_1\ga_2 \ga'_2}}{ \sqrt{\beta _j \sqrt{\ga'_1\ga_2\ga'_2}}}e^{i (\sqrt{\ga'_1} -\beta _{j}+\sqrt{\ga'_2}-\sqrt{\ga_2})s} \textbf{T}_{\ga'_1}\textbf{T}_{\ga'_2}  \frac{a+i(\sqrt{\ga'_1} -\beta_j)}{a^2+ (\sqrt{\ga'_1}-\beta_j) ^2}d\ga'_1 d\ga'_2\right]\overline{\lambda_{j}\lambda_{\ga_2}}d\ga_2\\
&-\frac{k^4}{4}\int_\xi ^{k^2}\sum_{j'=1}^N \left[\sum_{l ,l'=1}^{N}\frac{C_{\ga_1 lj'l'}}{\sqrt{\sqrt{\ga_1} \beta _l \beta_{j'} \beta _{l'}}}e^{i (\beta _{l}-\sqrt{\ga_1}+\beta _{l'}-\beta _{j'})s} \textbf{T}_{l}\textbf{T}_{l'} \frac{a+i(\beta _l-\sqrt{\ga_1})}{a^2+ (\beta _l-\sqrt{\ga_1}) ^2}\right.\\
&+\sum_{l=1}^N\int_{\xi}^{k^2} \frac{C_{\ga_1 lj'\ga'_2}}{ \sqrt{\sqrt{\ga_1} \beta _l \beta_{j'}  \sqrt{\ga'_2}}}e^{i(\beta _{l}-\sqrt{\ga_1}+\sqrt{\ga'_2}-\beta _{j'})s} \textbf{T}_{l}\textbf{T}_{\ga'_2} \frac{a+i(\beta_l-\sqrt{\ga_1})}{a^2+ (\beta _l-\sqrt{\ga_1}) ^2}d\ga'_2 \\
&+\int_\xi^{k^2}\sum_{l'=1}^N\frac{C_{\ga_1\ga'_1j' l'}}{ \sqrt{ \sqrt{\ga_1\ga'_1}\beta_{j'} \beta _{l' }}}e^{i (\sqrt{\ga'_1}-\sqrt{\ga_1}+\beta _{l'} -\beta_{j'})s} \textbf{T}_{\ga'_1}\textbf{T}_{l'} \frac{a+i(\sqrt{\ga'_1}-\sqrt{\ga_1})}{a^2+ (\sqrt{\ga'_1}-\sqrt{\ga_1}) ^2}d\ga'_1 \\
&\left.+\int_\xi^{k^2}\int_\xi^{k^2}\frac{C_{\ga_1 \ga'_1j'\ga'_2}}{ \sqrt{ \sqrt{\ga_1\ga'_1} \beta_{j'}\sqrt{\ga'_2}}}e^{i (\sqrt{\ga'_1} -\sqrt{\ga_1}+\sqrt{\ga'_2}-\beta_{j'})s} \textbf{T}_{\ga'_1}\textbf{T}_{\ga'_2}  \frac{a+i(\sqrt{\ga'_1} -\sqrt{\ga_1})}{a^2+ (\sqrt{\ga'_1}-\sqrt{\ga_1}) ^2}d\ga'_1 d\ga'_2\right]\overline{\lambda_{\ga_1}\lambda_{j'}}
\end{split}\]
\[\begin{split}
\,& -\frac{k^4}{4}\int_\xi^{k^2} \int_\xi ^{k^2}\left[\sum_{l ,l'=1}^{N}\frac{C_{\ga_1l \ga_2 l'}}{\sqrt{\sqrt{\ga_1} \beta _l \sqrt{\ga_2} \beta _{l'}}}e^{i (\beta _{l}-\sqrt{\ga_1}+\beta _{l'}-\sqrt{\ga_2})s} \textbf{T}_{l}\textbf{T}_{l'} \frac{a+i(\beta _l-\sqrt{\ga_1})}{a^2+ (\beta _l-\sqrt{\ga_1}) ^2}\right.\\
&+\sum_{l=1}^N\int_{\xi}^{k^2} \frac{C_{\ga_1 l\ga_2\ga'_2}}{ \sqrt{\sqrt{\ga_1} \beta _l   \sqrt{\ga_2 \ga'_2}}}e^{i(\beta _{l}-\sqrt{\ga_1}+\sqrt{\ga'_2}-\sqrt{\ga_2})s} \textbf{T}_{l}\textbf{T}_{\ga'_2} \frac{a+i(\beta_l-\sqrt{\ga_1})}{a^2+ (\beta _l-\sqrt{\ga_1}) ^2}d\ga'_2 \\
&+\int_\xi^{k^2}\sum_{l'=1}^N\frac{C_{\ga_1\ga'_1\ga_2 l'}}{ \sqrt{ \sqrt{\ga_1 \ga'_1\ga_2} \beta _{l' }}}e^{i (\sqrt{\ga'_1}-\sqrt{\ga_1}+\beta _{l'} -\sqrt{\ga_2})s} \textbf{T}_{\ga'_1}\textbf{T}_{l'} \frac{a+i(\sqrt{\ga'_1}-\sqrt{\ga_1})}{a^2+ (\sqrt{\ga'_1}-\sqrt{\ga_1}) ^2}d\ga'_1 \\
&\left. +\int_\xi^{k^2}\int_\xi^{k^2}\frac{C_{\ga_1\ga'_1\ga_2 \ga'_2}}{ (\ga'_1\ga_2\ga'_2)^{1/4}}e^{i (\sqrt{\ga'_1} -\sqrt{\ga_1}+\sqrt{\ga'_2}-\sqrt{\ga_2})s} \textbf{T}_{\ga'_1}\textbf{T}_{\ga'_2}  \frac{a+i(\sqrt{\ga'_1} -\sqrt{\ga_1})}{a^2+ (\sqrt{\ga'_1}-\sqrt{\ga_1}) ^2}d\ga'_1 d\ga'_2\right]\overline{\lambda_{\ga_1}\lambda_{\ga_2}}d\ga_1d\ga_2,
\end{split}\]

\begin{equation*}\begin{split}
\tilde{F}^{3}_\lambda& (\textbf{T},C,s) \\
&= \frac{k^4}{4}\sum_{j,j'=1}^N \left[\sum_{l ,l'=1}^{N}\frac{C_{jlj'l'}}{\sqrt{\beta _j \beta _l \beta_{j'} \beta _{l'}}}e^{i (\beta _{l}-\beta _{j}-\beta _{l'}+\beta _{j'})s} \textbf{T}_{l}\overline{\textbf{T}_{l'}} \frac{a+i(\beta _l-\beta_j)}{a^2+ (\beta _l-\beta_j) ^2}\right.\\
&+\sum_{l=1}^N\int_{\xi}^{k^2} \frac{C_{jlj'\ga'_2}}{ \sqrt{\beta _j \beta _l \beta_{j'}  \sqrt{\ga'_2}}}e^{i(\beta _{l}-\beta _{j}-\sqrt{\ga'_2}+\beta _{j'})s} \textbf{T}_{l}\overline{\textbf{T}_{\ga'_2}} \frac{a+i(\beta_l-\beta_j)}{a^2+ (\beta _l-\beta_j) ^2}d\ga'_2 \\
&+\int_\xi^{k^2}\sum_{l'=1}^N\frac{C_{j\ga'_1j' l'}}{ \sqrt{\beta _j \sqrt{\ga'_1}\beta_{j'} \beta _{l' }}}e^{i (\sqrt{\ga'_1}-\beta_j-\beta _{l'} +\beta_{j'})s} \textbf{T}_{\ga'_1}\overline{\textbf{T}_{l'}} \frac{a+i(\sqrt{\ga'_1}-\beta_j)}{a^2+ (\sqrt{\ga'_1}-\beta_j) ^2}d\ga'_1 \\
&\left.+\int_\xi^{k^2}\int_\xi^{k^2}\frac{C_{j\ga'_1j'\ga'_2}}{ \sqrt{\beta _j \sqrt{\ga'_1} \beta_{j'}\sqrt{\ga'_2}}}e^{i (\sqrt{\ga'_1} -\beta _{j}-\sqrt{\ga'_2}+\beta_{j'})s} \textbf{T}_{\ga'_1}\overline{\textbf{T}_{\ga'_2}}  \frac{a+i(\sqrt{\ga'_1} -\beta_j)}{a^2+ (\sqrt{\ga'_1}-\beta_j) ^2}d\ga'_1 d\ga'_2\right]\overline{\lambda_{j}}\lambda_{j'}\\
&+\frac{k^4}{4}\sum_{j=1}^N \int_\xi ^{k^2}\left[\sum_{l ,l'=1}^{N}\frac{C_{jl \ga_2 l'}}{\sqrt{\beta _j \beta _l \sqrt{\ga_2} \beta _{l'}}}e^{i (\beta _{l}-\beta _{j}-\beta _{l'}+\sqrt{\ga_2})s} \textbf{T}_{l}\overline{\textbf{T}_{l'}} \frac{a+i(\beta _l-\beta_j)}{a^2+ (\beta _l-\beta_j) ^2}\right.\\
&+\sum_{l=1}^N\int_{\xi}^{k^2} \frac{C_{jl\ga_2\ga'_2}}{ \sqrt{\beta _j \beta _l   \sqrt{\ga_2 \ga'_2}}}e^{i(\beta _{l}-\beta _{j}-\sqrt{\ga'_2}+\sqrt{\ga_2})s} \textbf{T}_{l}\overline{\textbf{T}_{\ga'_2}} \frac{a+i(\beta_l-\beta_j)}{a^2+ (\beta _l-\beta_j) ^2}d\ga'_2 \\
&+\int_\xi^{k^2}\sum_{l'=1}^N\frac{C_{j\ga'_1\ga_2 l'}}{ \sqrt{\beta _j \sqrt{\ga'_1\ga_2} \beta _{l' }}}e^{i (\sqrt{\ga'_1}-\beta_j-\beta _{l'}+\sqrt{\ga_2})s} \textbf{T}_{\ga'_1}\overline{\textbf{T}_{l'}} \frac{a+i(\sqrt{\ga'_1}-\beta_j)}{a^2+ (\sqrt{\ga'_1}-\beta_j) ^2}d\ga'_1 \\
&\left. +\int_\xi^{k^2}\int_\xi^{k^2}\frac{C_{j\ga'_1\ga_2 \ga'_2}}{ \sqrt{\beta _j \sqrt{\ga'_1\ga_2\ga'_2}}}e^{i (\sqrt{\ga'_1} -\beta _{j}-\sqrt{\ga'_2}+\sqrt{\ga_2})s} \textbf{T}_{\ga'_1}\overline{\textbf{T}_{\ga'_2}}  \frac{a+i(\sqrt{\ga'_1} -\beta_j)}{a^2+ (\sqrt{\ga'_1}-\beta_j) ^2}d\ga'_1 d\ga'_2\right]\overline{\lambda_{j}}\lambda_{\ga_2}
d\ga_2\\
&+\frac{k^4}{4}\int_\xi ^{k^2}\sum_{j'=1}^N \left[\sum_{l ,l'=1}^{N}\frac{C_{\ga_1 lj'l'}}{\sqrt{\sqrt{\ga_1} \beta _l \beta_{j'} \beta _{l'}}}e^{i (\beta _{l}-\sqrt{\ga_1}-\beta _{l'}+\beta _{j'})s} \textbf{T}_{l}\overline{\textbf{T}_{l'}} \frac{a+i(\beta _l-\sqrt{\ga_1})}{a^2+ (\beta _l-\sqrt{\ga_1}) ^2}\right.\\
&+\sum_{l=1}^N\int_{\xi}^{k^2} \frac{C_{\ga_1 lj'\ga'_2}}{ \sqrt{\sqrt{\ga_1} \beta _l \beta_{j'}  \sqrt{\ga'_2}}}e^{i(\beta _{l}-\sqrt{\ga_1}-\sqrt{\ga'_2}+\beta _{j'})s} \textbf{T}_{l}\overline{\textbf{T}_{\ga'_2}} \frac{a+i(\beta_l-\sqrt{\ga_1})}{a^2+ (\beta _l-\sqrt{\ga_1}) ^2}d\ga'_2 \\
&+\int_\xi^{k^2}\sum_{l'=1}^N\frac{C_{\ga_1\ga'_1j' l'}}{ \sqrt{ \sqrt{\ga_1\ga'_1}\beta_{j'} \beta _{l' }}}e^{i (\sqrt{\ga'_1}-\sqrt{\ga_1}-\beta _{l'}+\beta_{j'})s} \textbf{T}_{\ga'_1}\overline{\textbf{T}_{l'}} \frac{a+i(\sqrt{\ga'_1}-\sqrt{\ga_1})}{a^2+ (\sqrt{\ga'_1}-\sqrt{\ga_1}) ^2}d\ga'_1 \\
&\left.+\int_\xi^{k^2}\int_\xi^{k^2}\frac{C_{\ga_1 \ga'_1j'\ga'_2}}{ \sqrt{ \sqrt{\ga_1\ga'_1} \beta_{j'}\sqrt{\ga'_2}}}e^{i (\sqrt{\ga'_1} -\sqrt{\ga_1}-\sqrt{\ga'_2}+\beta_{j'})s} \textbf{T}_{\ga'_1}\overline{\textbf{T}_{\ga'_2}}  \frac{a+i(\sqrt{\ga'_1} -\sqrt{\ga_1})}{a^2+ (\sqrt{\ga'_1}-\sqrt{\ga_1}) ^2}d\ga'_1 d\ga'_2\right]\overline{\lambda_{\ga_1}}\lambda_{j'}
\end{split}\]
\[\begin{split}
\,&+\frac{k^4}{4}\int_\xi^{k^2} \int_\xi ^{k^2}\left[\sum_{l ,l'=1}^{N}\frac{C_{\ga_1l \ga_2 l'}}{\sqrt{\sqrt{\ga_1} \beta _l \sqrt{\ga_2} \beta _{l'}}}e^{i (\beta _{l}-\sqrt{\ga_1}-\beta _{l'}+\sqrt{\ga_2})s} \textbf{T}_{l}\overline{\textbf{T}_{l'}} \frac{a+i(\beta _l-\sqrt{\ga_1})}{a^2+ (\beta _l-\sqrt{\ga_1}) ^2}\right.\\
&+\sum_{l=1}^N\int_{\xi}^{k^2} \frac{C_{\ga_1 l\ga_2\ga'_2}}{ \sqrt{\sqrt{\ga_1} \beta _l   \sqrt{\ga_2 \ga'_2}}}e^{i(\beta _{l}-\sqrt{\ga_1}-\sqrt{\ga'_2}+\sqrt{\ga_2})s} \textbf{T}_{l}\overline{\textbf{T}_{\ga'_2}} \frac{a+i(\beta_l-\sqrt{\ga_1})}{a^2+ (\beta _l-\sqrt{\ga_1}) ^2}d\ga'_2 \\
&+\int_\xi^{k^2}\sum_{l'=1}^N\frac{C_{\ga_1\ga'_1\ga_2 l'}}{ \sqrt{ \sqrt{\ga_1 \ga'_1\ga_2} \beta _{l' }}}e^{i (\sqrt{\ga'_1}-\sqrt{\ga_1}-\beta _{l'}+\sqrt{\ga_2})s} \textbf{T}_{\ga'_1}\overline{\textbf{T}_{l'}} \frac{a+i(\sqrt{\ga'_1}-\sqrt{\ga_1})}{a^2+ (\sqrt{\ga'_1}-\sqrt{\ga_1}) ^2}d\ga'_1 \\
&\left. +\int_\xi^{k^2}\int_\xi^{k^2}\frac{C_{\ga_1\ga'_1\ga_2 \ga'_2}}{ (\ga_1\ga'_1\ga_2\ga'_2)^{1/4}}e^{i (\sqrt{\ga'_1} -\sqrt{\ga_1}-\sqrt{\ga'_2}+\sqrt{\ga_2})s} \textbf{T}_{\ga'_1}\overline{\textbf{T}_{\ga'_2}}  \frac{a+i(\sqrt{\ga'_1} -\sqrt{\ga_1})}{a^2+ (\sqrt{\ga'_1}-\sqrt{\ga_1}) ^2}d\ga'_1 d\ga'_2\right]\overline{\lambda_{\ga_1}}\lambda_{\ga_2}d\ga_1d\ga_2.
\end{split}\]
This expression combined with \eqref{cg2} gives us, $\sup_{\e,t} \mathbb{E}\big[\left \lvert \mathcal{A}^\e \left(f^\e _0 +f^\e _1\right)(t)\right \rvert ^2 \big] < + \infty$. $\square$ 
\end{preuve}
\begin{lem}\label{borneP2}
\[\lim_{M\to+\infty }\varlimsup_{\e \to 0} \mathbb{P}\left( \sup_{0 \leq t \leq T} \lvert \tilde{\emph{\textbf{T}}}^{\xi,\e}_\lambda (t)(y) \rvert \geq M \right)=0.\]
\end{lem}
\begin{preuve}[of Lemma \ref{borneP2}]
We recall that   $\|  \tilde{\textbf{T}}^{\xi,\e}  (t)(y) \| _{\mathcal{H}_\xi}= \| y \| _{\mathcal{H}_\xi}$ and then
\[\lvert \tilde{\textbf{T}}^{\xi,\e} _\lambda (t)(y) \rvert \leq \|  \tilde{\textbf{T}}^{\xi,\e}  (t)(y) \| _{\mathcal{H}_\xi}\| \lambda \| _{\mathcal{H}_\xi}=\| y \| _{\mathcal{H}_\xi}\| \lambda \| _{\mathcal{H}_\xi}.  \] 
$\square$ \\
This last lemma completes the proof of Proposition \ref{tightth1P2}. $\blacksquare$
\end{preuve}
\end{preuve}

\subsubsection{Martingale problem}

In this section, we shall characterize all subsequence limits by showing they are solution of a well-posed martingale problem. To do that, we consider a converging subsequence of $(\tilde{\textbf{T}}^{\xi,\e}  (.)(y))_{\epsilon \in (0,1)}$ which converges to a limit $\textbf{T}^{\xi} (.)(y)$. For the sake of simplicity we denote by $(\tilde{\textbf{T}}^{\xi,\e}  (.)(y))_{\epsilon \in (0,1)}$ the subsequence.

\paragraph{Convergence Result}

\begin{prop}\label{propmgP2}
$\forall \lambda \in \mathcal{H}_\xi$ and $\forall f$ smooth test function,
\begin{equation*}
\begin{split}
f\big(\emph{\textbf{T}}^\xi_{\lambda}&(t)(y)\big)\\
&-\int_0 ^t  \partial_{v}f\big(\emph{\textbf{T}}^\xi_{\lambda}(s)(y)\big)\left<J^\xi(\emph{\textbf{T}}^\xi(s)(y)),\lambda\right>_{\mathcal{H}_\xi} +\partial_{\overline{v}}f\big(\emph{\textbf{T}}^\xi_{\lambda}(s)(y)\big) \overline{\left<J^\xi(\emph{\textbf{T}}^\xi(s)(y)),\lambda\right>_{\mathcal{H}_\xi}}\\
+&\partial ^2 _{v} f\big(\emph{\textbf{T}}^\xi_{\lambda}(s)(y)\big)\left<K \big(\emph{\textbf{T}}^\xi(s)(y)\big) (\lambda),\lambda \right>_{\mathcal{H}_\xi}+\partial ^2_{\overline{v}}f\big(\emph{\textbf{T}}^\xi_{\lambda}(s)(y)\big) \overline{\left<K \big(\emph{\textbf{T}}^\xi(s)(y)\big) (\lambda),\lambda \right>_{\mathcal{H}_\xi}}\\
+&\partial  _{\overline{v}} \partial  _{v}f\big(\emph{\textbf{T}}^\xi_{\lambda}(s)(y)\big)\left<L\big(\emph{\textbf{T}}^\xi(s)(y)\big) (\lambda),\lambda \right>_{\mathcal{H}_\xi}+ \partial  _{v} \partial  _{\overline{v}}f\big(\emph{\textbf{T}}^\xi_{\lambda}(s)(y)\big)\overline{\left<L \big(\emph{\textbf{T}}^\xi(s)(y)\big) (\lambda),\lambda \right>_{\mathcal{H}_\xi}}ds
\end{split}
\end{equation*}
is a martingale, where
\begin{equation*}
\begin{split}
J^\xi(\emph{\textbf{T}})_{j}&= \left[\frac{\Gamma^c _{jj}-\Gamma^1_{jj}-\Lambda^{c,\xi}_j}{2}+i\left(\frac{\Gamma^{s}_{jj}-\Lambda^{s,\xi}_j}{2}+\kappa^\xi_j\right)\right]\emph{\textbf{T}}_j,\\
K (\emph{\textbf{T}}) (\lambda)_{j}&= - \frac{1}{2}\sum_{l=1}^N \Gamma^1_{jl} \emph{\textbf{T}}_j \emph{\textbf{T}}_l\overline{\lambda_l}- \frac{1}{2}\sum_{\substack{l=1\\l\not=j}}^N\left( \Gamma^{c}_{jl}+i\Gamma^{s}_{jl} \right)\emph{\textbf{T}}_j \emph{\textbf{T}}_l\overline{\lambda_l},\\
L (\emph{\textbf{T}}) (\lambda)_{j}&=\frac{1}{2}\sum_{l=1}^N \Gamma^1_{jl} \emph{\textbf{T}}_j\overline{\emph{\textbf{T}}_l}\lambda_l+ \frac{1}{2}\sum_{\substack{l=1\\l\not=j}}^N \Gamma^c_{jl} \emph{\textbf{T}}_l \overline{\emph{\textbf{T}}_l}\lambda_j,
\end{split}
\end{equation*}
and
\[J^\xi(\emph{\textbf{T}})_{\ga}=K (\emph{\textbf{T}}) (\lambda)_{\ga}=L (\emph{\textbf{T}}) (\lambda)_{\ga}=0\]
for almost every $\ga \in (\xi,k^2)$, and for $(\emph{\textbf{T}},\lambda)\in \mathcal{H}_\xi^2$.
\end{prop} 
\begin{preuve}[of Proposition \ref{propmgP2}]
Let 
\begin{equation*}
\begin{split}
f^\e _2 (t)&=\int_t ^{+\infty} \mathbb{E}_t ^\e \left[ \tilde{F}_{\lambda}\left(\tilde{\textbf{T}}^{\xi,\e} (t)(y),C\left(\frac{u}{\e}\right)\otimes C\left(\frac{u}{\e}\right),\frac{u}{\e} \right)\right]\\
&\quad-\tilde{F}_{\lambda}\left(\tilde{\textbf{T}}^{\xi,\e} (t)(y),\mathbb{E}[C(0)\otimes C(0)],\frac{u}{\e} \right)du.
\end{split}
\end{equation*} 
\begin{lem}\label{f2P2}
\[\sup_{t\geq 0}\mathbb{E}\left[\lvert f^\e _2 (t) \rvert \right]=\mathcal{O}\left(\e \right)\]
and
\[\mathcal{A}^\e \left( f_0 ^\e + f_1^\e +f_2 ^\e \right)(t)=\tilde{F}_{\lambda}\left(\tilde{\emph{\textbf{T}}}^{\xi,\e} (t)(y),\mathbb{E}[C(0)\otimes C(0)],\frac{t}{\e} \right) + A(\e,t),\]
where $\sup_{t\geq 0}\mathbb{E}\left[ \lvert A(\e,t) \rvert \right]=\mathcal{O}(\sqrt{\e})$.
\end{lem}
\begin{preuve}[of Lemma \ref{f2P2}]
Using a change of variable we get $f^\e _2 (t)=\e B(\e ,t)$ with
\begin{equation*}\begin{split}
  B(\e ,t)=&\int_0 ^{+\infty} \mathbb{E}_t ^\e \left[ \tilde{F}_{\lambda}\left(\tilde{\textbf{T}}^{\xi,\e} (t)(y),C\left(u+\frac{t}{\e}\right)\otimes C\left(u+\frac{t}{\e}\right),u+\frac{t}{\e} \right)\right]\\
  & -\tilde{F}_{\lambda}\left(\tilde{\textbf{T}}^{\xi,\e} (t)(y),\mathbb{E}[C(0)\otimes C(0)],u+\frac{t}{\e} \right)du.
 \end{split}\end{equation*}
By a computation, we get that $\sup_{\e, t\geq 0} \mathbb{E}\left[ \lvert B(\e,t) \rvert \right]<+\infty$, and after a long but straightforward computation we get the second part of the lemma. $\square$
\end{preuve}
Next, let $\tilde{G}_{\lambda}\big( \tilde{\textbf{T}}^{\xi,\e} (t)(y),\frac{t}{\e} \big)=\tilde{F}_{\lambda}\big(\tilde{\textbf{T}}^{\xi,\e} (t)(y),\mathbb{E}[C(0)\otimes C(0)],\frac{t}{\e} \big)$ and 
\[f^\e _3 (t)=-\int_0 ^t \Big[\tilde{G}_{\lambda}\big( \tilde{\textbf{T}}^{\xi,\e} (t)(y),\frac{u}{\e} \big)-\lim_{T\rightarrow+\infty} \frac{1}{T} \int_{0}^{T} \tilde{G}_{\lambda}\big(\tilde{ \textbf{T}}^{\xi,\e} (t)(y),s \big)ds \Big]du.\]
\begin{lem}\label{f3P2}$\forall T'>0$, we have
\[\lim_{\e\to 0}\sup_{0\leq t\leq T'} \mathbb{E}\left[\lvert f^\e _3 (t) \rvert \right]=0.\]
\end{lem}
\begin{preuve}[of Lemma \ref{f3P2}]
Using a change of variable, we get
\[f^\e _3 (t)=-\e \int_0 ^{\frac{t}{\e}}  \left[\tilde{G}_{\lambda}\left(\tilde{\textbf{T}}^{\xi,\e} (t)(y),u \right)-\lim_{T\rightarrow+\infty} \frac{1}{T} \int_{0}^{T} \tilde{G}_{\lambda}\left( \tilde{\textbf{T}}^{\xi,\e} (t)(y),s \right)ds \right]du.\] 
Let $\mu>0$, we have
\[\begin{split}
\Big\lvert \int_0 ^{\frac{t}{\e}}&  \left[\tilde{G}_{\lambda}\left( \textbf{T}^{\xi,\e} (t)(y),u \right)-\lim_{T\rightarrow+\infty} \frac{1}{T} \int_{0}^{T} \tilde{G}_{\lambda}\left( \textbf{T}^{\xi,\e} (t),s \right)ds \right]du \Big \rvert\\
& \hspace{3cm}\leq K(\mu,T',\xi,y)+\frac{K(T',\xi,y)}{\e} \sum_{j=1}^4 \int_{I^j_{< \mu}}d\ga_1\dots d\ga_j,
\end{split}\]
where for $j\in\{1,2,3,4\}$,
\[\begin{split}
I^j_{< \mu}=\Big\{& (\ga_l)_{l\in\{1,\dots,j\}}\in(\xi,k^2)^{j}, \exists (q_l)_{l\in\{1,\dots,4-j\}} \in\{\beta_1,\dots,\beta_N\}^{4-j}\\
&\hspace{2cm}\text{ and }(\mu_l)_{l\in\{1,\dots,4\}}\in\{-1,1\}^4,\text{ with } \Big\lvert \sum_{l=1}^j\mu_l\sqrt{\ga_l}+\sum_{l=1}^{4-j}\mu_{l+j}q_l \Big\rvert<\mu\Big\}.
\end{split}\]
Finally,
\[\begin{split}
\varlimsup_{\e \to 0}\sup_{0\leq t\leq T'}\E&\left[\e\left\lvert \int_0 ^{\frac{t}{\e}}  \left[\tilde{G}_{\lambda}\left( \textbf{T}^{\xi,\e} (t)(y),u \right)-\lim_{T\rightarrow+\infty} \frac{1}{T} \int_{0}^{T} \tilde{G}_{\lambda}\left( \textbf{T}^{\xi,\e} (t)(y),s \right)ds \right]du \right \rvert \right]\\
&\hspace{4cm}\leq K(T',\xi,y) \sum_{j=1}^4 \int_{I^j_{< \mu}}d\ga_1\dots d\ga_j,\end{split}\]
and then by letting $\mu\to 0$ we get the announced result. 
 $\square$ 
\end{preuve}
Let $f^\e (t)=f^\e _0 (t) + f^\e _1 (t)+ f^\e _2 (t) + f^\e _3 (t) $. A computation gives
\[\begin{split}
\mathcal{A}^\e f^\e (t)&=\lim_{T\rightarrow+\infty} \frac{1}{T} \int_{0}^{T} \tilde{G}_{\lambda}\left(\tilde{\textbf{T}}^{\xi,\e} (t)(y),s \right)ds +C(\e,t)\\
&=\tilde{G}^\infty_{\lambda}\left(\tilde{\textbf{T}}^{\xi,\e} (t)(y) \right) +C(\e,t),
\end{split}\] 
where, $\forall \mu>0$, 
\[\varlimsup_{\e\to 0}\sup_{0\leq t \leq T'} \mathbb{E}\left[ \lvert C(\e,t) \rvert \right] \leq K(T',\xi,y) \sum_{j=1}^4 \int_{I^j_{< \mu}}d\ga_1\dots d\ga_j,\]
using the boundness condition \eqref{cg2}. Moreover, for $(\textbf{T},\lambda)\in \mathcal{H}_\xi^2$, $\tilde{G}^\infty$ is defined as follow
\[\begin{split}
\tilde{G}^\infty _\lambda(\textbf{T})=&\partial_{v}f\left(\textbf{T}\right)\left<J^\xi(\textbf{T}),\lambda\right>_{\mathcal{H}_\xi} +\partial_{\overline{v}}f\left(\textbf{T}\right) \overline{\left<J^\xi(\textbf{T}),\lambda\right>_{\mathcal{H}_\xi} }\\
&+\partial ^2 _{v} f\left(\textbf{T}\right)\left<\tilde{K} \left(\textbf{T}\right) (\lambda),\lambda \right>_{\mathcal{H}_\xi}+\partial ^2_{\overline{v}}f\left(\textbf{T}\right) \overline{\left<\tilde{K} \left(\textbf{T}\right) (\lambda),\lambda \right>_{\mathcal{H}_\xi}}\\
&+\partial  _{\overline{v}} \partial  _{v}f\left(\textbf{T}\right)\left<\tilde{L} \left(\textbf{T}\right) (\lambda),\lambda \right>_{\mathcal{H}_\xi}+ \partial  _{v} \partial  _{\overline{v}}f\left(\textbf{T}\right)\overline{\left<\tilde{L} \left(\textbf{T}\right) (\lambda),\lambda \right>_{\mathcal{H}_\xi}}, 
\end{split}\]
where
\[
\begin{split}
\tilde{K} (\textbf{T}) (\lambda)_{j}&= -\frac{k^4}{4}\sum_{\beta _{j}+\beta _{j'}=\beta _{l}+\beta _{l'}}\frac{\textbf{C}_{jlj'l'}}{\sqrt{\beta _j \beta _l \beta_{j'} \beta _{l'}}} \textbf{T}_{l}\textbf{T}_{l'} \frac{a+i(\beta _l-\beta_j)}{a^2+ (\beta _l-\beta_j) ^2}\overline{\lambda_{j'}}\\
\tilde{L} (\textbf{T}) (\lambda)_{j}&=\frac{k^4}{4}\sum_{\beta _{j}-\beta _{j'}=\beta _{l}-\beta _{l'}}\frac{\textbf{C}_{jlj'l'}}{\sqrt{\beta _j \beta _l \beta_{j'} \beta _{l'}}} \textbf{T}_{l}\overline{\textbf{T}_{l'}} \frac{a+i(\beta _l-\beta_j)}{a^2+ (\beta _l-\beta_j) ^2}\lambda_{j'}
\end{split}
\]
for $j\in\{1,\dots,N\}$, with
\[\textbf{C}=\E\big[C(0)\otimes C(0)\big],\]
and
\[\tilde{K} (\textbf{T}) (\lambda)_{\ga}=\tilde{L} (\textbf{T}) (\lambda)_{\ga}=0\]
for almost every $\ga \in (\xi,k^2)$.

We assume that the following nondegeneracy condition holds. The wavenumbers $\beta_j$ are distinct along with their sums and differences. This assumption is also considered in \cite{book}, \cite{papa} and \cite{papanicolaou}. As a result we get 
\begin{equation}\label{geneP2}
\begin{split}
\tilde{G}^\infty_{\lambda}\left(\tilde{\textbf{T}}^{\xi,\e} (t)(y) \right) &=\partial_{v}f\left(\tilde{\textbf{T}}^{\xi,\e} _{\lambda}(t)(y)\right)\left<J^\xi(\tilde{\textbf{T}}^{\xi,\e} (t)(y)),\lambda\right>_{\mathcal{H}_\xi}\\ & +\partial_{\overline{v}}f\left(\tilde{\textbf{T}}^{\xi,\e} _{\lambda}(t)(y)\right) \overline{\left<J^\xi(\tilde{\textbf{T}}^{\xi,\e} (t)(y)),\lambda\right>_{\mathcal{H}_\xi} }\\
&+\partial ^2 _{v} f\left(\tilde{\textbf{T}}^{\xi,\e} _{\lambda}(t)(y)\right)\left<K \left(\tilde{\textbf{T}}^{\xi,\e} (t)(y)\right) (\lambda),\lambda \right>_{\mathcal{H}_\xi}\\ &+\partial ^2_{\overline{v}}f\left(\tilde{\textbf{T}}^{\xi,\e} _{\lambda}(t)(y)\right) \overline{\left<K \left(\tilde{\textbf{T}}^{\xi,\e} (t)(y)\right) (\lambda),\lambda \right>_{\mathcal{H}_\xi}}\\
&+\partial  _{\overline{v}} \partial  _{v}f\left(\tilde{\textbf{T}}^{\xi,\e} _{\lambda}(t)(y)\right)\left<L \left(\tilde{\textbf{T}}^{\xi,\e} (t)(y)\right) (\lambda),\lambda \right>_{\mathcal{H}_\xi}\\ &+ \partial  _{v} \partial  _{\overline{v}}f\left(\tilde{\textbf{T}}^{\xi,\e} _{\lambda}(t)(y)\right)\overline{\left<L \left(\tilde{\textbf{T}}^{\xi,\e} (t)(y)\right) (\lambda),\lambda \right>_{\mathcal{H}_\xi}}.
\end{split}
\end{equation}
By Theorem \ref{martingaleP2}, $\big(M^\e _{f^\e} (t) \big)_{t\geq 0}$ is an $\left( \mathcal{F}^\e _t \right)$-martingale. Then, for every bounded continuous function $h$, every sequence $0< s_1<\cdots <s_n \leq s <t$, and every family $(\lambda_j)_{j\in \{1,\dots,n\}}$ with values in $\mathcal{H}_\xi ^n$ we have
\[\mathbb{E}\left[ h\left( \tilde{\textbf{T}}^{\xi,\e} _{\lambda_j}(s_j)(y),1\leq j \leq n \right)\left( f^\e (t) - f^\e (s)-\int _s ^t  \mathcal{A}^\e f^\e (u)du \right) \right]=0 .\]
Finally, using  \eqref{geneP2} and Lemmas \ref{f1P2}, \ref{f2P2}, and \ref{f3P2}, we can conclude the proof of Proposition \ref{propmgP2}. $\blacksquare$
 \end{preuve}

\paragraph{Uniqueness}

In order to prove uniqueness, we decompose  $\textbf{T}^\xi (.)(y)$ into real and imaginary parts. Then, let us consider the new process
\[\textbf{Y} ^\xi(t)=\begin{bmatrix} \textbf{Y}^{1,\xi} (t) \\ \textbf{Y}^{2,\xi}(t) \end{bmatrix}, \text{ where }   \textbf{Y}^{1,\xi} (t)=Re\big( \textbf{T}^\xi(t) (y)\big) \text{ and } \textbf{Y}^{2,\xi} (t)=Im\big( \textbf{T}^\xi(t)(y) \big).\]
This new process takes its values in $\mathcal{G}_\xi\times \mathcal{G}_\xi$, where $\mathcal{G}_\xi=\mathbb{R}^N\times L^2((\xi,k^2),\mathbb{R})$. $\mathcal{G}_\xi\times \mathcal{G}_\xi$ is equipped with the inner product defined by
\[\left< \textbf{T},\textbf{S} \right>_{\mathcal{G}_\xi\times \mathcal{G}_\xi}=\sum_{j= 1}^N \textbf{T}^1_{j}\textbf{S}^1_{j} + \textbf{T}^2_{j}\textbf{S}^2_{j}+\int_\xi^{k^2}\textbf{T}^1_\ga \textbf{S}^1_\ga+\textbf{T}^2_\ga\textbf{S}^2_\ga d\ga\] 
$\forall (\textbf{T},\textbf{S})\in \mathcal{G}_\xi\times \mathcal{G}_\xi$. We also use the notation $\textbf{Y}^\xi_{\lambda}(t)=\left<\textbf{Y}^\xi(t),\lambda \right>_{\mathcal{G}_\xi\times \mathcal{G}_\xi}$ with $\lambda \in\mathcal{G}_\xi\times \mathcal{G}_\xi $. 
We introduce the operator $\Upsilon$ on $\mathcal{G}_\xi\times \mathcal{G}_\xi$ given by
\begin{equation*}
\begin{split}
\Upsilon:\,&\mathcal{G}_\xi\times \mathcal{G}_\xi\longrightarrow \mathcal{G}_\xi\times \mathcal{G}_\xi,\\
&\begin{bmatrix} \textbf{T}^1 \\ \textbf{T}^2 \end{bmatrix} \longmapsto \begin{bmatrix} \textbf{T}^2 \\ -\textbf{T}^1 \end{bmatrix}.
\end{split}
\end{equation*} 
By Proposition \ref{propmgP2}, we get the following result.
\begin{prop}\label{propmg2P2}
$\forall \lambda \in \mathcal{G}_\xi\times \mathcal{G}_\xi$, $\forall f \in \mathcal{C}_b^{\infty}(\mathbb{R})$
\begin{equation*}
\begin{split}
f\big(\emph{\textbf{Y}}^\xi_{\lambda}(t)\big)-\int_0 ^t & \big< B^\xi(\emph{\textbf{Y}}^\xi(s)),\lambda\big>_{\mathcal{G}_\xi\times \mathcal{G}_\xi} f'\big(\emph{\textbf{Y}}^\xi_{\lambda}(s)\big) \\
&+\frac{1}{2}\big<A\big(\emph{\textbf{Y}}^\xi(s)\big)(\lambda),\lambda\big>_{\mathcal{G}_\xi \times \mathcal{G}_\xi} f^{''}\big(\emph{\textbf{Y}}^\xi_{\lambda}(s)\big)ds
\end{split}
\end{equation*} 
is a martingale, where
\begin{equation*}
A(\emph{\textbf{Y}})(\lambda)=A_1(\emph{\textbf{Y}})(\lambda)+A_2(\emph{\textbf{Y}})(\lambda)+A_3(\emph{\textbf{Y}})(\lambda),
\end{equation*}
with, for $j\in\{1,\dots,N\}$,
\[\begin{split}
B^\xi(\emph{\textbf{Y}})_j&=\left[\frac{\Gamma^c _{jj}-\Lambda^{c,\xi}_j}{2}\right]\emph{\textbf{Y}}_j-\left[\frac{\Gamma^s_{jj}-\Lambda^{s,\xi}_j}{2}+\kappa^{\xi}_j\right]\Upsilon_j(\emph{\textbf{Y}})\\
A_1(\emph{\textbf{Y}})(\lambda)_j&=\Upsilon_j(\emph{\textbf{Y}})\sum_{l=1}^N\Gamma^1_{jl}\big[ \Upsilon^1_l(\emph{\textbf{Y}})\lambda^1_l+ \Upsilon^2_l(\emph{\textbf{Y}})\lambda^2_l\big]\\
A_2(\emph{\textbf{Y}})(\lambda)_j&=-\emph{\textbf{Y}}_j\sum_{\substack{l=1\\ l\not=j}}^N \Gamma^c_{jl}\big[\emph{\textbf{Y}}^1_l\lambda^1_l+\emph{\textbf{Y}}^2_l\lambda^2_l\big]+\Upsilon_j(\emph{\textbf{Y}})\sum_{\substack{l=1\\ l\not=j}}^N \Gamma^c_{jl}\big[\Upsilon^1_l(\emph{\textbf{Y}})\lambda^1_l+\Upsilon^2_l(\emph{\textbf{Y}})\lambda^2_l\big]\\
A_3(\emph{\textbf{Y}})(\lambda)_j&=\lambda_j\sum_{\substack{l=1\\l\not=j}}^N\Gamma^c_{jl}\big[(\emph{\textbf{Y}}^1_l)^2+(\emph{\textbf{Y}}^2_l)^2\big],
\end{split}\]
and
\[B^\xi_\ga(\emph{\textbf{Y}})=A_\ga(\emph{\textbf{Y}})(\lambda)=A_\ga(\emph{\textbf{Y}})(\lambda)=A_\ga(\emph{\textbf{Y}})(\lambda)=0\]
for almost every $\ga\in(\xi,k^2)$, and for $(\emph{\textbf{Y}},\lambda)\in \big(\mathcal{G}_\xi\times \mathcal{G}_\xi\big)^2$.
\end{prop}
\begin{preuve}[of Proposition \ref{propmg2P2}]
Using Proposition \ref{propmgP2},
\begin{equation*}
\begin{split}
f\big(\textbf{Y}^\xi_{\lambda}(t)\big)-\int_0 ^t &  Re\big(\big< J^\xi(\textbf{T}^\xi(s)(y)),\lambda\big>_{\mathcal{H}_\xi} \big)f'\big(\textbf{Y}^\xi _{\lambda}(s)\big) \\
& +\frac{1}{2}Re\big(\big<(L+K) \big(\textbf{T}^\xi(s)(y)\big) (\lambda),\lambda \big>_{\mathcal{H}_\xi}\big)f^{''}\big(\textbf{Y}^\xi_{\lambda}(s)\big)ds
\end{split}
\end{equation*} 
is a martingale. Let us remark that we also denote by $\lambda$ the function $\lambda^1+i \lambda^2$, and 
\begin{equation*}
Re\Big( \big<\textbf{T}^\xi (t)(y),\lambda\big>_{\mathcal{H}_\xi} \Big)= \big<\textbf{Y}^\xi(t),\lambda\big>_{\mathcal{G}_\xi \times \mathcal{G}_\xi} \text{ and }
Im\Big( \big<\textbf{T}^\xi (t)(y) ,\lambda\big>_{\mathcal{H}_\xi} \Big)= \big<\Upsilon(\textbf{Y}^\xi(t)),\lambda\big>_{\mathcal{G}_\xi \times \mathcal{G}_\xi}.
\end{equation*}
Then, we have  
\[\begin{split}
Re\big(\big< J^\xi(\textbf{T}^\xi(s)(y)),\lambda\big>_{\mathcal{H}_\xi} \big)&=\big< B^\xi(\textbf{Y}^\xi(s)),\lambda\big>_{\mathcal{G}_\xi\times \mathcal{G}_\xi} \\
Re\big(\big<(L+K) (\textbf{T}^\xi(s)(y)) (\lambda),\lambda \big>_{\mathcal{H}_\xi}\big)&=\big<A(\textbf{Y}^\xi (s))(\lambda),\lambda\big>_{\mathcal{G}_\xi\times \mathcal{G}_\xi}.
\end{split}\] 
$\blacksquare$
\end{preuve}
As a consequence of Proposition \ref{propmg2P2}, $\forall \lambda \in \mathcal{G}\times \mathcal{G}$, letting successively $f\in\mathcal{C}^\infty _b(\mathbb{R})$ such that $f(s)=s$ and  $f(s)=s^2$ if $\lvert s\rvert\leq r_y \|\lambda \|_{ \mathcal{G}\times \mathcal{G}}$, we get that
\begin{equation*}
\big< M^\xi (t), \lambda\big>_{\mathcal{G}_\xi\times\mathcal{G}_\xi}=M^\xi_\lambda (t)= \big< \textbf{Y}^\xi ( t)-\int_0 ^t B^\xi(\textbf{Y}^\xi(s))ds ,\lambda\big>_{\mathcal{G}_\xi\times \mathcal{G}_\xi}
\end{equation*}
is a continuous martingale with quadratic variation given by
\[<M^\xi_\lambda >(t)=\int _0 ^t \big< A(\textbf{Y}^\xi(s))(\lambda),\lambda\big>_{\mathcal{G}_\xi\times \mathcal{G}_\xi}ds.\]
\begin{prop}\label{uniquemgP2}
$\forall f \in \mathcal{C}^2 _b( \mathcal{G}_\xi\times \mathcal{G}_\xi)$,
\begin{equation}\label{mguniciteP2}
M^\xi_f(t)=f(\emph{\textbf{Y}}^\xi(t))-\int_0 ^t L^\xi f(\emph{\textbf{Y}}^\xi(s))ds
\end{equation}
is a continuous martingale, where $\forall \emph{\textbf{Y}} \in\mathcal{G}_\xi\times \mathcal{G}_\xi$
\[ L^\xi f(\emph{\textbf{Y}})=\frac{1}{2}trace\left( A(\emph{\textbf{Y}})D^2 f(\emph{\textbf{Y}}) \right) +\left< B^\xi(\emph{\textbf{Y}}), Df(\emph{\textbf{Y}})\right>_{\mathcal{G}_\xi\times \mathcal{G}_\xi}.\]
Moreover, the martingale problem associated to the generator $L^\xi$ is well-posed.
\end{prop}
\begin{preuve}[of Proposition \ref{uniquemgP2}]
We begin with the following lemma.
\begin{lem}
\begin{equation*}\begin{split}
A:\,&\mathcal{G}_\xi\times \mathcal{G}_\xi\longrightarrow L^+ _1\left(\mathcal{G}_\xi\times \mathcal{G}_\xi\right),\\
B^\xi:\,&\mathcal{G}_\xi\times \mathcal{G}_\xi\longrightarrow \mathcal{G}_\xi\times \mathcal{G}_\xi,
\end{split}\end{equation*}
where $L^+  _1\left(\mathcal{G}_\xi\times \mathcal{G}_\xi\right)$ is a set of nonnegative operators with finite trace. 
\end{lem}
\begin{preuve}
$\forall (\textbf{Y},\lambda) \in (\mathcal{G}_\xi\times \mathcal{G}_\xi)^2$, we have
\[ \begin{split}
\left< A(\textbf{Y})(\lambda),\lambda\right>_{\mathcal{G}_\xi\times \mathcal{G}_\xi}&=Re\big(\big<(L+K) (\textbf{T}) (\lambda),\lambda \big>_{\mathcal{H}_\xi}\big)\\
&=Re\Big(\sum_{j,l=1}^N \Gamma^1_{jl}\big[\textbf{T}_j\overline{\lambda_j}-\overline{\textbf{T}_j}\lambda_j  \big]\overline{\big[\textbf{T}_l\overline{\lambda_l}-\overline{\textbf{T}_l}\lambda_l  \big]}\Big)\\
&\quad+\sum_{\substack{j,l=1\\j\not=l}}^N\Gamma^c_{jl}\Big\lvert \textbf{T}_j\overline{\lambda_l}-\overline{\textbf{T}_l}\lambda_j  \Big\rvert^2.
\end{split}\]
with $\textbf{T}=\textbf{Y}^1+i\textbf{Y}^2$ and $\lambda=\lambda^1+i\lambda^2$.
First, $\forall (j,l)\in\{1,\dots,N\}^2$ such that $j\not=l$, $\Gamma^c_{jl}$ is nonnegative because it is proportional to the power spectral density of $C_{jl}$ at $\beta_l-\beta_j$ frequency. Second, the matrix $\Gamma^1$ is nonnegative since $\forall X\in\mathbb{C}^N$, we have
\[^{t}X \Gamma ^{1}X =\frac{k^4}{2}\sum_{j,l=1}^N \int_{0}^{+\infty }\mathbb{E}[C_{jj}(0)C_{ll}(z)]dz\tilde{X }_{j}\tilde{X}_l =\frac{k^4}{2}\int_0^{+\infty}\E[C_{\tilde{X}}(0)C_{\tilde{X}}(z)]dz\geq0\]
because it is proportional to the power spectral density of $C_{\tilde{X}}(z)=\sum_{j}C_{jj}(z)\tilde{X} _{j}$ at 0 frequency, and with $\tilde{X}_j=X_j/\beta_j$, $\forall j\in\{1,\dots,N\}$. Moreover,
\begin{equation*}
trace(A(\textbf{Y}))=\sum_{j=1}^N \Gamma^1_{jj}\big[(\textbf{Y}^1_j)^2+(\textbf{Y}^2_j)^2\big] \leq \sup_{j\in\{1,\dots,N\}}\Gamma^1_{jj}\,\, \| \textbf{Y}\| ^2 _{\mathcal{G}_\xi\times \mathcal{G}_\xi}.
\end{equation*}
$\square$
\end{preuve}
Consequently, following the proof of Theorem 4.1.4 in \cite{inf}, \eqref{mguniciteP2} is a martingale. However, $B^\xi$ and $A$ are not bounded functions but this problem can be compensated by the fact that the process $\textbf{Y}^\xi(.)$ takes its values in $\mathcal{B}_{r_y,\mathcal{G}_\xi\times\mathcal{G}_\xi}$.

Moreover, from this lemma there exists a linear operator $\sigma$ from $\mathcal{G}_\xi \times \mathcal{G}_\xi$ to $L_2(\mathcal{G}_\xi \times \mathcal{G}_\xi)$, which is the set of Hilbert-Schmidt operators from $\mathcal{G}_\xi \times \mathcal{G}_\xi$ to itself, such that $A(\textbf{Y})=\sigma(\textbf{Y})\circ\sigma^{\ast}(\textbf{Y})$. According to Theorem 3.2.2 and 4.4.1 in \cite{inf}, the martingale problem associated to $L^\xi$ is well-posed because $\forall \textbf{Y} \in \mathcal{G}_\xi\times \mathcal{G}_\xi$  
\[\| \sigma(\textbf{Y})\|  \leq K(N) \| \textbf{Y} \| _{\mathcal{G}_\xi\times \mathcal{G}_\xi}. \]
$\blacksquare$
\end{preuve}

Let us recall that the process $\textbf{Y}^\xi(.)$ is an element of $\mathcal{C}([0,+\infty),(\mathcal{B}_{r_y,\mathcal{G}_\xi\times \mathcal{G}_\xi},d_{\mathcal{B}_{r_y,\mathcal{G}_\xi\times \mathcal{G}_\xi}}))$, and we cannot assert that $\textbf{Y}^\xi(.)$ is uniquely determined. In fact, we need to know if its law is supported by $\mathcal{C}([0,+\infty),(\mathcal{G}_\xi\times\mathcal{G}_\xi,\|.\|_{\mathcal{G}_\xi\times\mathcal{G}_\xi}))$. Letting 
\[f(\textbf{Y})=\|\Pi(\xi,k^2)\otimes\Pi(\xi,k^2)(\textbf{Y}-y)\|^2_{\mathcal{G}_\xi\times\mathcal{G}_\xi},\]
where
\begin{equation*}
\begin{split}
\Pi(\xi,k^2)\otimes\Pi(\xi,k^2):\,&\mathcal{G}_\xi\times \mathcal{G}_\xi\longrightarrow \mathcal{G}_\xi\times \mathcal{G}_\xi,\\
&\begin{bmatrix} \textbf{Y}^1 \\ \textbf{Y}^2 \end{bmatrix} \longmapsto \begin{bmatrix} \Pi(\xi,k^2)(\textbf{Y}^1) \\ \Pi(\xi,k^2)(\textbf{Y}^2) \end{bmatrix}.
\end{split}
\end{equation*}
As $\textbf{Y}^\xi(.)$ is a solution on $\mathcal{C}([0,+\infty),(\mathcal{B}_{r_y,\mathcal{G}_\xi\times \mathcal{G}_\xi},d_{\mathcal{B}_{r_y,\mathcal{G}_\xi\times \mathcal{G}_\xi}}))$ of the martingale associated to $L^\xi$, we get
\[\E[f(\textbf{Y}^\xi(t))]=0 \quad\forall t\geq0,\] 
and therefore $\Pi(\xi,k^2)\otimes\Pi(\xi,k^2)(\textbf{Y}^\xi(.))=\Pi(\xi,k^2)\otimes\Pi(\xi,k^2)(Re(y),Im(y))$.
Consequently, the process $\textbf{Y}^\xi(.)$ is strongly continuous since the weak and the strong topologies are the same on $\mathbb{R}^N$. Finally, $\textbf{Y}^\xi(.)$ is uniquely characterized as being the unique solution of the martingale problem associated to $L^\xi$ and starting from $(Re(y),Im(y))$, and that concludes the proof of Theorem \ref{thasymptP21}.

\subsection{Proof of Theorem \ref{thasympP2}}

Let $\mathcal{H}_0=\mathbb{C}^N\times L^2(0,k^2)$ and $y\in\mathcal{H}_0$. We begin by showing the tightness of the process $(\textbf{T}^\xi(.)(y^\xi))_\xi$, which is the unique solution of the martingale problem associated to $\mathcal{L}_\xi$ and starting from $y^\xi=\Pi(\xi,+\infty)(y)$. As the radiating part $\Pi(0,k^2)(\textbf{T}^\xi(.)(y^\xi))$ of the process $\textbf{T}^\xi(.)(y^\xi)$ is constant equal to $\Pi(\xi,k^2)(y^\xi)$, to prove the tightness of $(\textbf{T}^\xi(.)(y^\xi))_\xi$ is suffices to show the tightness of the finite-dimensional process $(\Pi(k^2,+\infty)(\textbf{T}^\xi(.)(y^\xi)))_\xi$. Let $\E^\xi_t$ be the conditional expectation given the $\sigma$-algebra $\sigma(\textbf{T}^\xi(u)(y^\xi),0\leq u\leq t)$. Then, $\forall t\geq 0$, $\forall h\in(0,1)$ and $\forall s\in[0,h]$, we have
\[\begin{split}
\E^\xi_t\Big[\| \textbf{T}^\xi(t+s)(y^\xi)-&\textbf{T}^\xi(t)(y^\xi) \|^2_{\mathbb{C}^N} \Big]\\
&\leq \E^\xi_t\Big[\| \textbf{Y}^{1,\xi}(t+s)-\textbf{Y}^{1,\xi}(t) \|^2_{\mathbb{R}^N} \Big]+ \E^\xi_t\Big[\| \textbf{Y}^{2,\xi}(t+s)-\textbf{Y}^{2,\xi}(t) \|^2_{\mathbb{R}^N}\Big]\\
&\leq \sum_{\substack{j=1\\l=1,2}}^N\E^\xi_t\Big[ (\textbf{Y}^{l,\xi}_j(t+s)-\textbf{Y}^{l,\xi}_j(t))^2 \Big]\\
&\leq  \sum_{\substack{j=1\\l=1,2}}^N\E^\xi_t\Big[\Big(\int_t^{t+s} L^\xi f_j^l(\textbf{Y}^\xi(u)) du\Big)^2  \Big]  +\E^\xi_t\Big[\Big(M^\xi_{f^l_j}(t+s)-M^\xi_{f^l_j}(t)\Big)^2 \Big],
\end{split}\]
with $\forall \textbf{Y}\in \mathcal{G}_0\times\mathcal{G}_0$, $f^l_j(\textbf{Y})=\textbf{Y}^l_j$. Therefore, using that the process $\textbf{T}^\xi(.)(y^\xi)$ takes its values in $\mathcal{B}_{r_y,\mathcal{H}_\xi}$, we first get 
\[\E^\xi_t\Big[\Big(\int_t^{t+s} L^\xi f_j^l(\textbf{Y}^\xi(u)) du\Big)^2  \Big]  \leq K\,h^2,\] 
and second,
\[\E^\xi_t\Big[\Big(M^\xi_{f^l_j}(t+s)-M^\xi_{f^l_j}(t)\Big)^2 \Big]=\E^\xi_t\Big[<M^\xi_{f^l_j}>_{t+s}-<M^\xi_{f^l_j}>_{t} \Big]\leq K\,h\]
with
\[<M^\xi_{f^l_j}>_{t}=\int_0^t L^\xi (f^l_j)^2(\textbf{Y}^\xi(u)) -2f^l_j(\textbf{Y}^\xi(u))L^\xi f^l_j(\textbf{Y}^\xi(u))\,du.\]
Consequently, the process $(\textbf{T}^\xi(.)(y^\xi))_\xi$ is tight on $\mathcal{C}([0,+\infty),(\mathcal{H}_0,\|.\|_{\mathcal{H}_0}))$. Now, to characterize all limits of converging subsequences, let us denote by $\textbf{T}^0(.)(y)$ such a limit point. First, for every smooth function $f$ on $\mathcal{H}_0$, for every bounded continuous function $h$, and every sequence $0< s_1<\cdots <s_n \leq s <t$, we have
\[\mathbb{E}\left[ h\left( \textbf{T}^{\xi}(s_j)(y^\xi),1\leq j \leq n \right)\left( f (\textbf{T}^{\xi}(t)(y^\xi)) -f (\textbf{T}^{\xi}(s)(y^\xi))-\int _s ^t  \mathcal{L}_\xi f (\textbf{T}^{\xi}(u)(y^\xi))du \right) \right]=0 .\]
Second,
\[\sup_{\textbf{T}\in \mathcal{B}_{r_y,\mathcal{H}_0}}\Big\rvert \mathcal{L}f(\textbf{T})-\mathcal{L}_\xi f(\textbf{T})\Big\lvert\leq K\sup_{j\in\{1,\dots,N\}} \big\lvert\Lambda^{c,\xi}_j-\Lambda^{c}_j \big\rvert+\big\lvert\Lambda^{s,\xi}_j-\Lambda^{s}_j \big\rvert+\big\lvert\kappa^\xi_j-\kappa_j \big\rvert. \]
Consequently, $\textbf{T}^0(.)(y)$ is a solution of the martingale problem associated to $\mathcal{L}$ and starting from $y$. However, following the proof of the uniqueness in Theorem \ref{thasymptP21}, this martingale problem is well-posed and therefore $\textbf{T}^\xi(.)(y^\xi)$ converges in distribution to the unique solution of the martingale problem associated to $\mathcal{L}$ and starting from $y$.

\subsection{Proof of Theorem \ref{hfapproxP2}}

The proof of this theorem follows ideas developed in \cite[Chapter 11]{stroock}. In order to prove this theorem we use a probabilistic representation of $\mathcal{T}_{j} ^l (\omega, z)$ by using the Feynman-Kac formula. To this end, we introduce the jump Markov process $\big( X^{N}_t \big)_{t\geq 0}$ with state space $\big\{-(N-1)/N, \dots,0,\dots, (N-1)/N\big\}$ and generator given by
\[{\mathcal{L}}^{N}\phi \left(\frac{l}{N}\right)= \Gamma^c _{l\,l+1}\left( \phi \left(\frac{l-1}{N}\right) -\phi \left(\frac{ l}{N}\right) \right)+\Gamma ^c_{ l +2\,l +1}\left( \phi \left(\frac{ l +1}{N}\right) -\phi \left(\frac{ l }{N}\right) \right)\]
for $l\in \{1,\dots,N-2 \} $,
\[{\mathcal{L}}^{N}\phi \left(\frac{l}{N}\right)=\Gamma^c _{\lvert l\rvert +2\,\lvert l\rvert +1} \left( \phi \left(\frac{l-1}{N}\right) -\phi \left(\frac{l}{N}\right) \right)+\Gamma^c _{\lvert l\rvert \,\lvert l\rvert +1}\left( \phi \left(\frac{l+1}{N}\right) -\phi \left(\frac{l}{N}\right) \right)\]
for $l\in \{-(N-2),\dots,-1 \} $,
\[{\mathcal{L}}^{N}\phi (0)=\frac{\Gamma^c _{2\,1}}{2}\left( \phi \left(\frac{1}{N}\right) -\phi (0) \right)+\frac{\Gamma^{c} _{2\,1}}{2}\left( \phi \left(\frac{-1}{N}\right) -\phi (0) \right),\]
and
\[{\mathcal{L}}^{N}\phi \left(\frac{\pm (N-1)}{N}\right)=\Gamma^{c} _{N-1\,N}\left( \phi \left(\frac{\pm (N-2)}{N}\right) -\phi \left(\frac{\pm (N-1)}{N}\right) \right).\]
Using the Feynman-Kac formula, we get for $(j,l)\in\{1,\dots,\N{}\} ^2$
\[\mathcal{T}^l _j ( \omega , L ) =  \E_{\frac{l-1}{N}}\left[  e^{-\Lambda^c_N \int_0 ^L \textbf{1}_{\left(\lvert X^{N}_v \rvert = \frac{N-1}{N}\right)}dv  -\Lambda^c_{N-1} \int_0 ^L \textbf{1}_{\left(\lvert X^{N}_v \rvert= \frac{N-2}{N}\right)}dv } \textbf{1}_{\left( \lvert X^{N} _L \rvert + \frac{1}{N}= \frac{j}{N} \right)}  \right]. \]
Let $f$ be a bounded continuous function on $[0,1]$, we consider $\mathcal{T}^l ( \omega , L )$ as a family of bounded measures on $[0,1]$ by setting 
\[ \mathcal{T}^l _f (\omega, L)=\E_{\frac{l-1}{N}}\left[  e^{-\Lambda^c_N \int_0 ^L \textbf{1}_{\left(\lvert X^N_v  \rvert= \frac{N-1}{N}\right)}dv  -\Lambda^c_{N-1} \int_0 ^L \textbf{1}_{\left(\lvert X^N_v \rvert= \frac{N-2}{N}\right)}dv } f\left( \lvert X^N _L\rvert +\frac{1}{N} \right)  \right]. \] 

In the first part of the proof, we consider the case $v\in[0,1)$ and in a second part we shall treat the case $v=1$.
  
Let $u\in[0,1)$ such that $l(N)/N \to u$. We begin by introducing some notations. Throughout the proof we  denote by $\tau^{(l)}_{j/N}$ the $l$th passage in $j/N$, for $j\in\{-(N-1),\dots,N-1 \}$. To avoid the unboundness in ${\mathcal{L}}^{N}$ of the reflecting barriers ${\mathcal{L}}^{N}\phi(\pm (N-1)/N)$, we introduce the stopping time
\[ \tau^\alpha _{N}=\tau^{(1)}_{(N-[N^\alpha])/N}\wedge \tau^{(1)}_{- (N-[N^\alpha])/N}\]
with $\alpha \in (0,1)$.
Let $X^{N,\tau}_t = X^{N}_{t\wedge  \tau^\alpha _{N}}$, $\forall t\geq 0$, be
the stopped process and $d(N)=(l(N)-1)/N$. We denote by $\mathbb{P}^{N} _{d(N)}$ the law of $( X^N_t )_{t\geq 0}$ starting from $d(N)$ and by $\mathbb{P}^{N,\tau} _{d(N)}$ the law of $( X^{N,\tau}_t )_{t\geq 0}$ starting from $d(N)$. 
Let
\[{\mathcal{L}_{\overline{a}_\infty}}=\frac{\partial  }{\partial  v}\left( \overline{a}_\infty(\cdot) \frac{\partial  }{\partial  v} \right),\]
where $\overline{a}_\infty(\cdot) \in \mathcal{C}^1 (\mathbb{R})$ is an extension over $\mathbb{R}$ of $a_\infty(\cdot)$, which is defined on $[-1,1]$, and such that the martingale problem associated to $\mathcal{L}_{\overline{a}_\infty}$ and starting from $u$ is well posed. We denote by $\overline{\mathbb{P}}_u$ this unique solution. Let $\varphi \in \mathcal{C}^{\infty}_0 (\mathbb{R})$,
\[M_\varphi(t) = \varphi(x(t)) - \varphi(x(0))-\int_{0}^{t} {\mathcal{L}_{\overline{a}_\infty}}\varphi(x(s))ds,\]
and $\tau_{r} =\inf (u\geq 0, \, \lvert x(t)\rvert \geq r  )$ for $r\in(0,1)$.

\begin{lem}\label{generatorP2} 
$\forall\varphi\in {\mathcal{C}}^{\infty  }_{0}(\mathbb{R})$, $\forall \alpha \in (2/3,1)$,
\[\lim_{N \to +\infty } \sup_{v\in  \left[-\frac{N-[N^\alpha ]}{N},-\frac{1}{N}\right]\cup \left[\frac{1}{N},\frac{N-[N^\alpha ]}{N}\right]}\lvert {\mathcal{L}}^{N} \varphi(v)-{\mathcal{L}_{\overline{a}_\infty}}\varphi(v) \rvert=0,\]
where $ {\mathcal{L}}^{N} \varphi(v)$ is defined as follows. $\forall j\in \{1,\dots,N-2\}$,
\[\begin{split}
{\mathcal{L}}^{N} \varphi (v) &= \Gamma^c _{j\,j+1}\left( \varphi \left(\frac{j-1}{N}\right) -\varphi\left(\frac{j}{N}\right) \right) \\
                                          & +  \Gamma^c _{j+1\,j+2}\left( \varphi \left( \frac{j+1}{N}\right) -\varphi\left(\frac{j}{N}\right) \right),
\end{split}\] 
for $v\in[j/N, (j+1)/N)$, and
\[\begin{split}
{\mathcal{L}}^{N} \varphi (v) & =  \Gamma^c _{j\,j+1}\left( \varphi \left( \frac{-j+1}{N}\right) -\varphi\left( \frac{-j}{N}\right) \right) \\
                                          & +  \Gamma^c _{j+1\,j+2}\left( \varphi \left( \frac{-j-1}{N}\right) -\varphi\left( \frac{-j}{N}\right) \right),
\end{split}\] 
if $v\in(- (j+1)/N,-j/N]$.                     
\end{lem}
\begin{preuve}[of Lemma \ref{generatorP2}]
We shall restrict the proof of this lemma to the proof of 
\[\lim_{N \to +\infty } \sup_{v\in \left[\frac{1}{N},\frac{N-[N^\alpha ]}{N}\right]}\lvert {\mathcal{L}}^{N} \varphi(v)-{\mathcal{L}_{\overline{a}_\infty}}\varphi(v) \rvert=0,\]
since the other case is completely similar by symmetry. We start the proof of this technical lemma by proving Lemma \ref{coefgP2} page \pageref{coefgP2}.
\begin{preuve}[of Lemma \ref{coefgP2}]
Let $h_e(v)=\frac{v}{\sqrt{\left( n_1 k d \theta  \right) ^2-v ^2}}$ and $g(v)=\arctan(v)$.
We recall that $\forall j\in \big\{1,\dots,N\big\}$, $\tan(\sigma_j )=-h_e(\sigma_j)$ 
First,
\[\begin{split}
\big\lvert \sigma_{j+1}-\sigma_j -\pi\big\rvert & \leq\big\lvert g\big( \tan \big(\sigma_{j+1} -(j+1)\pi\big) \big)- g\big( \tan \big(\sigma_j -j\pi\big) \big) \big\rvert \\
& \leq  K \lvert \tan\left(  \sigma_{j+1}  \right)-\tan\left(  \sigma_{j} \right) \rvert \\
& \leq  K \lvert h_e( \sigma_{j+1} )- h_e( \sigma_{j} ) \rvert \\
& \leq  K \sup_{v\in[\sigma_{j}, \sigma_{j+1}]} h'_e(v) 
\end{split}\]
where $h'_e(v)=\frac{(n_1 k d \theta )^2}{\left( (n_1 k d \theta )^2-v^2 \right)^{\frac{3}{2}}}$
which is a positive and increasing function.
Moreover, 
\[\sigma_{N-[N^\alpha] }\leq  (N-[N^\alpha]) \pi \]
and then
\[\sup_{j\in \{1,\dots, N-[N^\alpha]  \}} \big\lvert \sigma_{j+1}-\sigma_j -\pi\big\rvert ={\cal{O}}\left( N^{\frac{1}{2}-\frac{3}{2}\alpha } \right).\]
Second, in the same way we have
\[\begin{split}
\sigma_{j+2}-&2\sigma_{j+1}+\sigma_j   \\
&=   g\big( \tan \big(\sigma_{j+2} -(j+2)\pi\big)\big)- 2 g\big( \tan \big(\sigma_{j+1} -(j+1)\pi\big) \big) +g\big( \tan \big(\sigma_j -j\pi\big) \big)  \\
& = - \big( g( h_e(\sigma_{j+2}))- 2 g( h_e(\sigma_{j+1})) +g(h_e(\sigma_j) )\big)
\end{split}\]
and 
\[\begin{split}
g( h_e(\sigma_{j+2}))- 2 g( h_e(\sigma_{j+1})) +g(h_e(\sigma_j) )&=\big[g'(h_e(\sigma_{j+1}))-g'(h_e(\sigma_j))\big].\big[h_e(\sigma_{j+2})-h_e(\sigma_{j+1})\big]\\
&+g'(h_e(\sigma_{j}))\big[h_e(\sigma_{j+2})-2h_e(\sigma_{j+1})+h_e(\sigma_{j})\big]\\
&+\int_{h_e(\sigma_{j+1})}^{h_e(\sigma_{j+2})} \left(h_e(\sigma_{j+2})-t\right)g''(t)\,\,dt\\
&- \int_{h_e(\sigma_{j})}^{h_e(\sigma_{j+1})} \left(h_e(\sigma_{j+1})-t\right)g''(t)\,\,dt.
\end{split}\]
Moreover,
\[\big\lvert g'(h_e(\sigma_{j+1}))-g'(h_e(\sigma_j))\big\rvert.\big\lvert h_e(\sigma_{j+2})-h_e(\sigma_{j+1}) \big\rvert\leq  K\, N^{1-3\alpha}.\]
\[\begin{split}
h_e(\sigma_{j+2})-2h_e(\sigma_{j+1})+h_e(\sigma_{j})=&\int_{\sigma_{j+1}}^{\sigma_{j+2}}  h'_e(t)-h'_e(t-\pi)\,dt\\
&+\int_{\sigma_{j+1}}^{\sigma_{j+2}} h'_e(t-\pi)\,dt-\int_{\sigma_{j}}^{\sigma_{j+1}} h'_e(t)\,dt,
\end{split}\]
with
\[\left\lvert\int_{\sigma_{j+1}}^{\sigma_{j+2}}  h'_e(t)-h'_e(t-\pi)\,dt\right\rvert \leq K \, h''_e(\sigma_{N-[N^\alpha]})=\mathcal{O}(N^{\frac{1}{2}-\frac{5}{2} \alpha}),\]
because $h''_{e}(v)=\frac{3(n_1 k d \theta)^2 v}{((n_1 k d \theta)^2-v^2)^{5/2}}$,
and
\[\begin{split}
\Big\lvert \int_{\sigma_{j+1}}^{\sigma_{j+2}} h'_e(t-\pi)\,dt&-\int_{\sigma_{j}+\pi}^{\sigma_{j+1}+\pi} h'_e(t-\pi)\,dt\Big\rvert \\
&\leq h'_e(\sigma_{N-[N^\alpha]})\left(\left\lvert \sigma_{j+2}-\sigma_{j+1}-\pi\right\rvert + \left\lvert \sigma_{j+2}-\sigma_{j+1}-\pi\right\rvert\right)\\
&\leq K N^{1-3\alpha}.
\end{split}\]
Finally,
\[\begin{split}
\Big\lvert \int_{h_e(\sigma_{j+1})}^{h_e(\sigma_{j+2})} \left(h_e(\sigma_{j+2})-t\right)g''(t)\,\,dt
&- \int_{h_e(\sigma_{j})}^{h_e(\sigma_{j+1})} \left(h_e(\sigma_{j+1})-t\right)g''(t)\,\,dt\Big\rvert\\
&\leq K \left((h_e(\sigma_{j+2})-h_e(\sigma_{j+1}))^2+(h_e(\sigma_{j+1})-h_e(\sigma_{j}))^2\right)\\
&\leq K\, N^{1-3\alpha},
\end{split}\]
and
\[\sup_{j\in\{1,\dots,N-[N^\alpha ]-2\}}\lvert \sigma_{j+2}-2\sigma_{j+1}+\sigma_{j} \rvert=\mathcal{O}(N^{1-3\alpha}).\]
This completes the proof of Lemma \ref{coefgP2} since we can take $\alpha>1/3$ and we have $N^{\frac{1}{2}-\frac{5}{2}\alpha}\leq N^{1-3\alpha}$.
$\square$
\end{preuve}
From this lemma, we immediately get
\[ \begin{split}
\sup_{j\in\{1,\dots,N-[N^\alpha ]-1\}}\big\lvert S(\sigma_{j+1}-\sigma _j &,\sigma_{j+1}-\sigma _j)-S( \pi,\pi)\big\rvert \\
&\leq K\sup_{j\in\{1,\dots,N-[N^\alpha ]-1\}}\big\lvert \sigma_{j+1}-\sigma _j-\pi \big\rvert=\mathcal{O}\left(N^{\frac{1}{2}-\frac{3}{2}\alpha }\right).
\end{split} \]
Before showing that  
\[ \sup_{j\in \{1,\dots,N-[N^{\alpha }]\}}\Big\lvert A^2 _j -\frac{2}{d} \Big\rvert=\mathcal{O}(N^{\alpha-1}),\]
where $A _j$ is defined by \eqref{coefajP2}, we prove that 
\[\sup_{j\in \{ 1,\dots, [N^{ \alpha }] \}}\big\lvert \sigma_{j} - j\pi \big\rvert=\mathcal{O}\left( \frac{1}{N^{1- \alpha }}\right).\]
In fact, $\forall j\in \big\{1,\dots,N^{\alpha}\big\}$
\[\begin{split}
\big\lvert \sigma_{j} -j\pi\big\rvert & = \big\lvert  g \big( \tan \big(  \sigma_{j}-j\pi \big) \big)- \arctan ( \tan (0))\Big\rvert\\
& \leq  K \big\lvert \tan ( \sigma_{j})\big\rvert \\
& \leq K \,h_e(\sigma_{[N^{ \alpha }]}).
\end{split}\]
Moreover, $\sigma_{[N^{\alpha }]}\leq N^{ \alpha }\pi$
and then
\[h(\sigma_{[N^{ \alpha }]})=\mathcal{O}\left(\frac{1}{N^{1-\alpha }}\right).\]
Consequently,
\[
 \sup_{j\in \{1,\dots,N-[N^{\alpha }]\}}\Big\lvert A^2 _j -\frac{2}{d} \Big\rvert  \leq K  \sup_{j\in \{1,\dots,N-[N^{\alpha }]\}} \left\lvert \frac{\sin^2(\sigma_j )}{\zeta_j}-\frac{\sin(2\sigma_j )}{2\sigma_j}  \right\rvert\leq \frac{K}{N^{1-\alpha}}\]
 because
 \[\begin{split}
 \sup_{j\in \{1,\dots,N-[N^{\alpha }]\}} &\Big\lvert \frac{\sin^2(\sigma_j )}{\zeta_j}-\frac{\sin(2\sigma_j )}{2\sigma_j}  \Big\rvert\leq \frac{1}{\sqrt{(n_1 k d \theta)^2-\sigma^2_{N-[N^\alpha]}}} \\
&+\sup_{j\in \{[N^{\alpha}]+1,\dots,N-[N^{\alpha }]\}}\left\lvert \frac{\sin(2\sigma_j )}{2\sigma_j} \right\rvert+\sup_{j\in \{ 1,\dots, [N^{ \alpha }] \}}\left\lvert \frac{\sin(2\sigma_j )}{2\sigma_j} \right\rvert\\
&\leq \frac{K}{N^{1/2+\alpha/2}}+\frac{1}{2\sigma_{[N^\alpha]+1}}+ \sup_{j\in \{ 1,\dots, [N^{ \alpha }] \}}\left\lvert \frac{\sin(2\sigma_j )-\sin(2j \pi)}{2\sigma_j} \right\rvert\\
&\leq K\left(  \frac{1}{N^{1/2+\alpha/2}} +\frac{1}{N^{\alpha}}+ \frac{1}{\sigma_1}\sup_{j\in \{ 1,\dots, [N^{ \alpha }] \}}\big\lvert\sigma_j -j\pi\big\rvert \right).
\end{split}\]

Now, let us introduce 
\[
B^N _j= \frac{a}{2n_1^2 \sqrt{1-\frac{\sigma_{j} ^2}{n_1 ^2 k^2 d^2} }\sqrt{1-\frac{\sigma_{j+1} ^2}{n_1 ^2 k^2 d^2}}} \frac{\frac{1}{4} A_{j} ^2 A^2 _{j+1}S(\sigma _{j+1}-\sigma _{j},\sigma _{j+1}-\sigma _{j})}{ a^2+\left( \beta _{j}-\beta _{j+1} \right)^2 }. 
\]
Then, for $j\in\big\{1,\dots,N-2\big\}$ and $v\in\left[\frac{j}{N}, \frac{j+1}{N}\right]$. 
\[\begin{split}
\mathcal{L}^N \varphi(v) =\left(\frac{n_1 k d\theta}{N\pi} \right)^2 \left(\frac{N\pi}{n_1 d \theta }\right)^2
& \Big[ \Big(\varphi \Big(\frac{[Nv]+1}{N}\Big) -\varphi\Big(\frac{[Nv]}{N}\Big)\Big)B^N_{j+1}\\
&+\Big( \varphi \Big(\frac{[Nv]-1}{N}\Big) -\varphi\Big(\frac{[Nv]}{N}\Big)\Big) B^N_{j} \Big].
\end{split}\]
Consequently, from the following decomposition
\[\begin{split}
N^2 \Big[ \Big(\varphi& \Big(\frac{[Nv]+1}{N}\Big) -\varphi\Big(\frac{[Nv]}{N}\Big)\Big)B^N_{j+1}+\Big( \varphi \Big(\frac{[Nv]-1}{N}\Big) -\varphi\Big(\frac{[Nv]}{N}\Big)\Big) B^N_{j}  \Big]\\
&\hspace{10cm}-\frac{n_1 ^2 d^2 \theta ^2}{\pi^2}\mathcal{L}_{a_\infty}\varphi(v)\\
&=N^2\Big[ \varphi \Big(\frac{[Nv]+1}{N}\Big)-2\varphi\Big(\frac{[Nv]}{N}\Big)+\varphi \Big(\frac{[Nv]-1}{N}\Big)\Big]\Big[B^N_j -\frac{n_1 ^2 d^2 \theta ^2}{\pi^2}a_\infty(v)\Big]\\
&+N\Big[ \varphi \Big(\frac{[Nv]+1}{N}\Big) -\varphi\Big(\frac{[Nv]}{N}\Big)\Big]\Big[N(B^N_{j+1}  -B^N_j)-\frac{n_1 ^2 d^2 \theta ^2}{\pi^2}\frac{d}{dv}a_\infty(v) \Big]\\
&+\frac{n_1 ^2 d^2 \theta ^2}{\pi^2}a_\infty(v)\Big[N^2\Big( \varphi \Big(\frac{[Nv]+1}{N}\Big)-2\varphi\Big(\frac{[Nv]}{N}\Big)+\varphi \Big(\frac{[Nv]-1}{N}\Big)\Big)-\varphi''(v)\Big]\\
&+\frac{n_1 ^2 d^2 \theta ^2}{\pi^2}\frac{d}{dv}a_\infty(v) \Big[N\Big( \varphi \Big(\frac{[Nv]+1}{N}\Big) -\varphi\Big(\frac{[Nv]}{N}\Big)\Big)-\varphi'(v)\Big],
\end{split}\]
and because it is easy to show that 
\[\begin{split}
\sup_{v\in\left[\frac{1}{N},\frac{N-[N^\alpha ]}{N}\right]}&\Big\lvert N\Big( \varphi \Big(\frac{[Nv]+1}{N}\Big) -\varphi\Big(\frac{[Nv]}{N}\Big)\Big)-\varphi'(v)\Big\rvert =\mathcal{O}\Big(\frac{1}{N}\Big)\\
\sup_{v\in\left[\frac{1}{N},\frac{N-[N^\alpha ]}{N}\right]}&\Big\lvert N^2\Big( \varphi \Big(\frac{[Nv]+1}{N}\Big)-2\varphi\Big(\frac{[Nv]}{N}\Big)+\varphi \Big(\frac{[Nv]-1}{N}\Big)\Big)-\varphi''(v)\Big\rvert=\mathcal{O}\Big(\frac{1}{N}\Big),
\end{split}\]
it suffices to show the two following points
\begin{itemize}
\item \[\lim_{N} \sup_{j\in\{1 ,\dots,N-[N^\alpha]-1\}} \sup_{v \in \left[\frac{j}{N},\frac{j+1}{N}\right]} \Big\lvert B^N_j -\frac{n_1 ^2 d^2 \theta ^2}{\pi^2}a_\infty(v)\Big\rvert=0.\]
\item\[\lim_{N} \sup_{j\in\{1 ,\dots,N-[N^\alpha]-1\}} \sup_{v \in \left[\frac{j}{N},\frac{j+1}{N}\right]} \Big\lvert N\Big( B^N_j -B^N_{j+1}\Big) - \frac{n_1 ^2 d^2 \theta ^2}{\pi^2}\frac{d}{dv}a_\infty(v)\Big\rvert =0.\]
\end{itemize}
We decompose the proof of these two points into two sublemmas.

\begin{lem}\label{point1P2}
\[\lim_{N}\sup_{j\in\{1 ,\dots,N-[N^\alpha]-1\}} \sup_{v \in \left[\frac{j}{N},\frac{j+1}{N}\right]}\Big\lvert B^N_j  -\frac{n_1 ^2 d^2 \theta ^2}{\pi^2}a_\infty(v)\Big\rvert=0\]
\end{lem}
\begin{preuve}[of Lemma \ref{point1P2}]
$\forall j\in \big\{1,\dots,N-[N^\alpha]-1\big\}$ and $\forall v\in \Big[\frac{j}{N}, \frac{j+1}{N}\Big]$, we have the following inequalities, 
\[\frac{1}{ 1-\frac{\sigma _{j} ^2}{n_1 ^2 k^2d^2}}\leq \frac{1}{\sqrt{ \left( 1-\frac{\sigma_{j}^2}{n_1 ^2 k^2d^2}\right) \left( 1-\frac{\sigma_{j+1}^2}{n_1 ^2 k^2 d^2}\right)}} 
\leq \frac{1}{ 1-\frac{\sigma _{j+1} ^2}{n_1 ^2 k^2d^2}}.\]
Moreover, for $l\in\{j,j+1\}$
\[\begin{split}
\Big\lvert \frac{1}{1-(\theta  v)^2}-\frac{1}{ 1-\frac{\sigma _{l} ^2}{n_1 ^2 k^2d^2}}\Big\rvert &\leq \Big\lvert \frac{j+1}{N}\theta -\frac{j-1}{n_1 k d}\pi\Big\rvert \frac{2\theta}{(1-\theta^2)^2}\\
&\leq \frac{K}{N}.
\end{split}\]
Consequently,
\[  \sup_{j\in\{1 ,\dots,N-[N^\alpha]-1\}} \sup_{v \in \left[\frac{j}{N},\frac{j+1}{N}\right]} \Big\lvert \frac{1}{1-(\theta  v)^2}-\frac{1}{\sqrt{ \left( 1-\frac{\sigma_{j}^2}{n_1 ^2 k^2d^2}\right) \left( 1-\frac{\sigma_{j+1}^2}{n_1 ^2 k^2d^2}\right)}}\Big\rvert=\mathcal{O}\left(\frac{1}{N}\right).\]
Next,
\[\begin{split}
\Big\lvert n_1 k \Big( \sqrt{1-\frac{\sigma_{j}^2}{n_1 ^2 k^2d^2} }&-\sqrt{1 -\frac{\sigma_{j+1}^2}{n_1 ^2 k^2d^2}}  \Big)-\frac{\pi}{d}\frac{\frac{\sigma_{j}}{n_1  k d}}{\sqrt{1-\frac{\sigma_{j}^2}{n_1 ^2 k^2d^2}}}\Big\rvert \\
&\leq  K\Big\lvert n_1 k \Big(\sqrt{1-\frac{\sigma_{j}^2}{n_1 ^2 k^2d^2} }-\sqrt{1 -\frac{\sigma_{j+1}^2}{n_1 ^2 k^2d^2}} \Big) -  \frac{1}{d}(\sigma_{j+1} -\sigma_{j}) \frac{\frac{\sigma_{j}}{n_1  k d}}{\sqrt{1-\frac{\sigma_{j}^2}{n_1 ^2 k^2 d^2}}} \Big\rvert\\
& \quad+ K\Big\lvert (\sigma_{j+1} -\sigma_{j}-\pi) \frac{\frac{\sigma_{j}}{n_1  k d}}{\sqrt{1-\frac{\sigma_{j}^2}{n_1 ^2 k^2 d^2}}} \Big\rvert\\
&\leq \frac{K}{N}\frac{1}{(1-\theta^2)^{3/2}}(\sigma_{j+1}-\sigma_j)^2+K\big\lvert \sigma_{j+1} -\sigma_{j}-\pi\big\rvert\frac{\theta}{\sqrt{1-\theta^2}}\leq K N^{\frac{1}{2}-\frac{3}{2}\alpha}
\end{split}\]
and then
\[\sup_{j\in \{1,\dots,N-[N^\alpha ]-1\}}\Big\lvert \frac{1}{a^2+\left(\beta_j -\beta_{j+1} \right)^2 }-\frac{1}{a^2+\frac{\pi^2 }{ d^2}\frac{\frac{\sigma_{j}^2}{n_1 ^2 k^2 d^2}}{1-\frac{\sigma_{j}^2}{n_1 ^2 k^2 d^2}}} \Big\rvert =\mathcal{O}  \Big(N^{\frac{1}{2}-\frac{3}{2}\alpha } \Big).\]
Moreover $\forall j\in \big\{1,\dots,N-[N^\alpha]-1\big\}$ and $\forall v\in \Big[\frac{j}{N}, \frac{j+1}{N}\Big]$, we have
\[\Big\lvert \frac{\frac{\sigma_{j}}{n_1 k d}}{\sqrt{1-\frac{\sigma_{j}^2}{n_1 ^2 k^2 d^2}}}- \frac{\theta v}{\sqrt{1-( \theta v)^2}} \Big\rvert  \leq K \Big\lvert \frac{\sigma_{j}}{n_1 k d} - \theta v \Big\rvert\leq \frac{K}{N},\]
and finally  
\[\sup_{j\in \{1,\dots,N-[N^\alpha ]-1\}}\sup_{v\in \big[\frac{j}{N},\frac{j+1}{N} \big]} \Big\lvert \frac{1}{a^2+(\beta_j -\beta_{j+1} )^2}-\frac{1}{a^2+\frac{\pi^2 }{ d^2}\frac{(\theta v)^2}{1-( \theta x)^2}} \Big\rvert=\mathcal{O}\Big(N^{\frac{1}{2}-\frac{3}{2}\alpha}\Big).\]
This concludes the proof of Lemma \ref{point1P2}.$\square$
\end{preuve}

\begin{lem}\label{point2P2}
\[\lim_{N} \sup_{j\in\{1 ,\dots,N-[N^\alpha]-1\}} \sup_{v \in \left[\frac{j}{N},\frac{j+1}{N}\right]}\Big\lvert N\Big( B^N_j  -B^N_{j+1}\Big) - \frac{n_1 ^2 d^2 \theta ^2}{\pi^2}\frac{d}{dv}a_\infty(v)\Big\rvert =0.\]
\end{lem}

\begin{preuve}[of Lemma \ref{point2P2}]
We separate the proof of this lemma into two step. First, for each $j \in \big\{1,\dots,N-[N^\alpha]-2\big\}$ let
\[C^N_j=N\left( \frac{1}{\sqrt{1-\frac{\sigma_{j+1}^2}{n_1 ^2 k^2 d^2}}\sqrt{1-\frac{\sigma_{j+2}^2}{n_1 ^2 k^2 d^2}}}-\frac{1}{\sqrt{1-\frac{\sigma_{j}^2}{n_1 ^2 k^2 d^2}}\sqrt{1-\frac{\sigma_{j+1}^2}{n_1 ^2 k^2 d^2}}}\right).\]
We can write $\forall v\in \big[\frac{j}{N}, \frac{j+1}{N}\big]$
\[\begin{split}
C^N_j-&\frac{2\theta^2 v}{(1-(\theta v)^2)^2} = \frac{1}{\sqrt{1-\frac{\sigma_{j+1}^2}{n_1 ^2 k^2d^2}}}N\int_{\frac{\sigma_j}{n_1 k d}}^{\frac{\sigma_{j+2}}{n_1 k d}} \frac{w}{(1-w^{2})^{\frac{3}{2}}} dw - \frac{2\theta^2 v}{(1-(\theta v)^2)^2}\\
                         & =  \frac{1}{\sqrt{1-\frac{\sigma_{j+1}^2}{n_1 ^2 k^2 d^2}}}\Big( N\int_{\frac{\sigma_j}{n_1 k d}}^{\frac{\sigma_{j+2}}{n_1 k d}}\frac{w}{(1-w^{2})^{\frac{3}{2}}}dw  -N\Big(\frac{\sigma_{j+2}}{n_1 k d}-\frac{\sigma_j}{n_1 k d}\Big)\frac{\theta v}{(1-(\theta v)^{2})^{\frac{3}{2}}} \Big) \\
                         & +  N\Big(\frac{\sigma_{j+2}}{n_1 k d}-\frac{\sigma_j}{n_1 k d}\Big)\frac{\theta v}{(1-(\theta v)^{2})^{\frac{3}{2}}} \Big( \frac{1}{\sqrt{1-\frac{\sigma_{j+1}^2}{n_1 ^2 k^2 d^2}}}-\frac{1}{\sqrt{1-(\theta v)^2}} \Big) \\
                         & +  \Big(N\Big(\frac{\sigma_{j+2}}{n_1 k d}-\frac{\sigma_j}{n_1 k d}\Big)-2\theta  \Big)\frac{\theta  v}{(1-(\theta v)^2)^2}. 
 \end{split}\]                        
We can check that the function $v\mapsto \frac{\theta v}{(1-(\theta v)^2)^2}$ is bounded on $[0,1]$ and 
\[\Big\lvert N\Big(\frac{\sigma_{j+2}}{n_1 k d}-\frac{\sigma_j}{n_1 k d}\Big)-2 \theta  \Big\rvert  \leq  \frac{N}{n_1 k d} \Big\lvert \sigma _{j+2}- \sigma _{j}-2\pi \Big\rvert + 2\theta \Big\lvert \frac{N \pi}{n_1 d k\theta }-1  \Big\rvert \leq K\,N^{\frac{1}{2}-\frac{3}{2}\alpha}.\]
Moreover, $v\mapsto \frac{\theta v}{(1-(\theta v)^2)^2}$ is bounded on $[0,1]$ and 
\[\Big\lvert  \frac{1}{\sqrt{1-(\theta v)^2}}-\frac{1}{\sqrt{1-\frac{\sigma_{j+1}^2}{n_1 ^2 k^2}}} \Big\rvert  \leq   \frac{K}{N}\frac{\theta}{(1-\theta^2)^{3/2}}.\]
Finally, $0\leq\frac{1}{\sqrt{1-\frac{\sigma_{j+1}^2}{n_1 ^2 k^2 d^2}}}\leq \frac{1}{\sqrt{1-\theta^2}}$ and
\[\begin{split}
\Big\lvert N\int_{\frac{\sigma_j}{n_1 k d}}^{\frac{\sigma_{j+2}}{n_1 k d}}\frac{w}{(1-w^{2})^{\frac{3}{2}}}dw  -N\Big(\frac{\sigma_{j+2}}{n_1 k d}-\frac{\sigma_j}{n_1 k d}\Big)&\frac{\theta v}{(1-(\theta v)^{2})^{\frac{3}{2}}}\Big\rvert\\
&\leq N\int_{\frac{\sigma_j}{n_1k d}}^{\frac{\sigma_{j+2}}{n_1 k d}}\lvert w- \theta v\rvert dw \frac{2\theta^2 +1}{(1-\theta^2)^{\frac{5}{2}}}\\
&\leq K N \Big[ \Big(\theta v -\frac{\sigma_j}{n_1 k d}\Big)^2 + \Big(\theta v -\frac{\sigma_{j+2}}{n_1 k d}\Big)^2  \Big]\\
&\leq \frac{K}{N}.
\end{split}\]
Consequently,
\[  \sup_{j\in\{1,\dots,N-[N^\alpha ]-1\}}\sup_{v\in \big[\frac{j}{N},\frac{j+1}{N} \big]} \Big\lvert C^N_j-\frac{2\theta^2 v}{(1-(\theta v)^2)^2} \Big\rvert =\mathcal{O}(N^{\frac{1}{2}-\frac{3}{2}\alpha}).\]
Second, for $j\in \big\{1,\dots,N-[N^\alpha]-1\big\}$ and $v\in\big[\frac{j}{N},\frac{j+1}{N}\big]$, we have
\[\begin{split}
\Big\lvert N&\Big( \frac{1}{a^2+(\beta_{j+1}-\beta_{j+2})^2 }-\frac{1}{a^2+(\beta_{j}-\beta_{j+1})^2}\Big)+\frac{\frac{\pi^2 }{d^2}\frac{2\theta^2 v}{(1-(\theta v)^2)^2}}{\left( a^2+\frac{\pi^2 }{ d^2}\frac{(\theta v)^2}{1-(\theta v)^2} \right)^2}\Big\rvert \\
&\leq\Big\lvert N \Big( \frac{1}{a^2+(\beta_{j+1}-\beta_{j+2})^2}-\frac{1}{a^2+(\beta_{j}-\beta_{j+1})^2}\Big)\\
&\hspace{4cm}- N \big( (\beta_{j+1}-\beta_{j+2})-(\beta_{j}-\beta_{j+1})\big)\frac{-2   (\beta_{j}-\beta_{j+1})}{\big(a^2+(\beta_{j}-\beta_{j+1})^2 \big)^2}\Big\rvert\\
&+ \Big\lvert  N \big( (\beta_{j+1}-\beta_{j+2})-(\beta_{j}-\beta_{j+1})\big)\frac{-2  (\beta_{j}-\beta_{j+1})}{\big(a^2+(\beta_{j}-\beta_{j+1})^2 \big)^2} + \frac{  \frac{\pi^2 }{ d^2}  \frac{2\theta^2 v}{(1-(\theta v)^2)^2}  }{  \Big( a^2+\frac{\pi^2 }{d^2}  \frac{(\theta v)^2}{1-(\theta v)^2}   \Big)^2  }\Big\rvert.
\end{split}\]
For the first term on the right of the previous inequality, we have
\[\begin{split} 
\Big\lvert N & \Big( \frac{1}{a^2+(\beta_{j+1}-\beta_{j+2})^2}-\frac{1}{a^2+(\beta_{j}-\beta_{j+1})^2 }\Big)\\
&\hspace{3.5cm}- N \big( (\beta_{j+1}-\beta_{j+2})-(\beta_{j}-\beta_{j+1})\big)\frac{-2  (\beta_{j}-\beta_{j+1})}{\big(a^2+(\beta_{j}-\beta_{j+1})^2 \big)^2}\Big\rvert\\
&\leq K\,N \big((\beta_{j+1}-\beta_{j+2})-(\beta_{j}-\beta_{j+1})\big)^2,
\end{split}\]
and we shall see just below that 
\[ \sup_{j\in\{1,\dots,N-[N^\alpha ]-2\}}\lvert\beta_{j+2}-2\beta_{j+1}-\beta_{j}\rvert=\mathcal{O}\Big(\frac{1}{N}\Big).\]
Now, for the second term we have previously get
\[ \sup_{j\in\{1,\dots,N-[N^\alpha ]-1\}}\sup_{v\in \big[\frac{j}{N},\frac{j+1}{N} \big]} \Big\lvert \beta_{j}-\beta_{j+1}-\frac{\pi}{d}\frac{\theta v}{\sqrt{1-(\theta v)^2}}\Big\rvert=\mathcal{O}(N^{\frac{1}{2}-\frac{3}{2}\alpha}). \]
Then, to finish the proof of this lemma it suffices to show that 
\[\sup_{j\in\{1,\dots,N-[N^\alpha ]-2\}}\sup_{v\in \big[\frac{j}{N},\frac{j+1}{N} \big]} \Big\lvert N \big( \beta _j -2\beta _{j+1} +\beta _{j+2} \big)+  \frac{\frac{\pi}{d}\theta }{(1-(\theta v)^2)^{\frac{3}{2}}} \Big\rvert=\mathcal{O}(N^{2-3\alpha}).\]
To show this relation we shall use the following decompositions. For $l\in\{j,j+1\}$
\[ \begin{split}\sqrt{1-\frac{\sigma_l ^2}{n_1 ^2 k^2 d^2}}-\sqrt{1-\frac{\sigma_{l+1}^2}{n_1 ^2 k^2d^2}}&-\frac{1}{n_1 k d}(\sigma_{l+1}-\sigma_l)\frac{\frac{\sigma_l}{n_1 k d}}{\sqrt{1-\frac{\sigma_l ^2}{n_1 ^2k^2 d^2}}} \\
&= \int_{\frac{\sigma_l}{n_1 k d}}^{\frac{\sigma_{l+1}}{n_1 k d}}\left( \frac{\sigma_{l+1}}{n_1 k d}-w\right)\frac{1}{(1-w^2)^{\frac{3}{2}}}dw ,\end{split}\]
and
\[\begin{split}   
N n_1 k \Big( &\sqrt{1-\frac{\sigma_j ^2}{n_1 ^2 k^2 d^2}}-2\sqrt{1-\frac{\sigma_{j+1}^2}{n_1 ^2 k^2d^2}}+ \sqrt{1-\frac{\sigma_{j+2} ^2}{n_1 ^2 k^2 d^2}} \Big)+\frac{\frac{\pi}{d}\theta }{\big( 1-(\theta v)^2\big)^{\frac{3}{2}}}\\
&=\frac{N}{d}\Big( (\sigma _{j+1}-\sigma_j)\frac{\frac{\sigma_j }{n_1 k d}}{\sqrt{1-\frac{\sigma_j^2 }{n_1 ^2 k^2 d^2}}}-(\sigma _{j+2}-\sigma_{j+1})\frac{\frac{\sigma_{j+1} }{n_1 k d}}{\sqrt{1-\frac{\sigma_{j+1}^2 }{n_1 ^2 k^2 d^2}}} \Big)+\frac{\frac{\pi}{d}\theta }{\left( 1-(\theta v)^2\right)^{\frac{3}{2}}}\\
&\quad+N n_1 k \Big( \int_{\frac{\sigma_j }{n_1 k d}}^{\frac{\sigma_{j+1} }{n_1 k d}} \Big(\frac{\sigma_{j+1} }{n_1 k d}-w\Big)\frac{1}{(1-w^2)^{\frac{3}{2}}}dw - \int_{\frac{\sigma_{j+1} }{n_1 k d}}^{\frac{\sigma_{j+2} }{n_1 k d}} \Big(\frac{\sigma_{j+2} }{n_1 k d}-w\Big)\frac{1}{(1-w^2)^{\frac{3}{2}}}dw \Big).
\end{split}\]
First, using Lemma \ref{coefgP2} we have
\[
\int_{\frac{\sigma_j }{n_1 k d}}^{\frac{\sigma_{j+1} }{n_1 k d}} \Big(\frac{\sigma_{j+1} }{n_1 k d}-w\Big)\frac{1}{(1-w^2)^{\frac{3}{2}}}dw - \int_{\frac{\sigma_{j+1} }{n_1 k d}}^{\frac{\sigma_{j+2} }{n_1 k d}} \Big(\frac{\sigma_{j+2} }{n_1 k d}-w\Big)\frac{1}{(1-w^2)^{\frac{3}{2}}}dw=\mathcal{O}(N^{\frac{1}{2}-\frac{3}{2}\alpha-2}).\]
Second, we have
\[\begin{split}
\frac{N}{d}\Big( (\sigma _{j+1}-\sigma_j)\frac{\frac{\sigma_j }{n_1 k}}{\sqrt{1-\frac{\sigma_j^2 }{n_1 ^2 k^2}}}-&(\sigma _{j+2}-\sigma_{j+1})\frac{\frac{\sigma_{j+1} }{n_1 k}}{\sqrt{1-\frac{\sigma_{j+1}^2 }{n_1 ^2 k^2}}} \Big)+\frac{\frac{\pi}{d}\theta }{\left( 1-(\theta v)^2\right)^{\frac{3}{2}}}\\
&=\frac{N}{d}\big(\sigma_{j+1}-\sigma_j-\pi\big) \Big[\frac{\frac{\sigma_j }{n_1 k d}}{\sqrt{1-\frac{\sigma_j^2 }{n_1 ^2 k^2 d^2}}}-\frac{\frac{\sigma_{j+1} }{n_1 k d}}{\sqrt{1-\frac{\sigma_{j+1}^2 }{n_1 ^2 k^2 d^2}}}\Big]\\
&\quad-\frac{N}{d}(\sigma_{j+2}-2\sigma_{j+1}+\sigma_{j})\frac{\frac{\sigma_{j+1} }{n_1 k d}}{\sqrt{1-\frac{\sigma_{j+1}^2 }{n_1 ^2 k^2 d^2}}}\\
&\quad+ \frac{\pi}{d}\Big( N \Big[\frac{\frac{\sigma_j }{n_1 k d}}{\sqrt{1-\frac{\sigma_j^2 }{n_1 ^2 k^2 d^2}}}-\frac{\frac{\sigma_{j+1} }{n_1 k d}}{\sqrt{1-\frac{\sigma_{j+1}^2 }{n_1 ^2 k^2 d^2}}}\Big] +\frac{\theta}{(1-(\theta v)^2)^{\frac{3}{2}}}\Big),
\end{split}\]
where, according to Lemma \ref{coefgP2}, the first and the third term are $\mathcal{O}(N^{\frac{1}{2}-\frac{3}{2}\alpha})$, and the second term is $\mathcal{O}(N^{2-3\alpha})$. That concludes the proof of Lemma \ref{point2P2} for $\alpha\in(2/3,1)$.$\square$
\end{preuve}

Consequently, thanks to Lemma \ref{point1P2} and Lemma \ref{point2P2}, we get
\[ \sup_{v\in \left[\frac{1}{N},\frac{N-[N^\alpha ]}{N}\right]}\lvert {\mathcal{L}}^{N} \varphi(v)-{\mathcal{L}_{a_\infty}}\varphi(v) \rvert=\mathcal{O}\big(N^{(2-3\alpha)\vee(\alpha -1)}\big),\]
this concludes the proof of Lemma \ref{generatorP2}.$\square$
\end{preuve}

\begin{lem}\label{tightP2}
$\mathbb{P}^{N,\tau} _{d(N)}$ is tight on ${\cal{D}}([0,+\infty ),\mathbb{R})$.
\end{lem}

\begin{preuve}[of Lemma \ref{tightP2}]
Let ${\mathcal{M}}_t =\sigma( x(u) ,\, 0\leq u\leq t )$. According to Theorem $3$ in \cite[Chapter 3]{kushner}, we have to show the two following points. First,   
\[\lim_{K\to +\infty } \overline{\lim_{N}} \,\,\mathbb{P}^{N,\tau} _{d(N)} \left( \sup_{t\geq 0}\lvert x(t) \rvert \geq K \right)=0.\]
The first point is satisfied since we have $\forall N$, $\mathbb{P}^{N,\tau} _{d(N)} \left( \sup_{t\geq 0}\lvert x(t) \rvert\leq 1 \right)=1$.
Second, for each $N$, $h \in (0,1)$, $s\in[0,h]$ and $t\geq 0$,
\[ \E^{\mathbb{P}^{N,\tau} _{d(N)} }\big( (x(t+s) -x(t) )^2 \vert \mathcal{M} _t \big)\leq K\, h.\]
Concerning the second point, letting $\varphi\in\mathcal{C}^\infty_b(\mathbb{R})$ such that $\varphi(s)=s$ if $\lvert s\rvert \leq 1$, we have
\[\begin{split}
\E^{\mathbb{P}^{N,\tau} _{d(N)} } \big( (x(t+s) -x(t) )^2 \vert \mathcal{M}_t\big) & \leq  2\,\E^{\mathbb{P}^{N,\tau} _{d(N)} }\big( (M^N_{\varphi}(t+s)  - M^N _{\varphi}(t) )^2 \vert \mathcal{M}_t\big)\\
                                                                                          & + 2\,\E^{\mathbb{P}^{N,\tau} _{d(N)}}\left( \left(\int_{t}^{t+s} {\mathcal{L}}^{N}\varphi(x(w))dw\right) ^2\Big\vert \mathcal{M}_t\right),
\end{split} \]                                                                                         
with 
\[M^N _{\varphi} (t) =\varphi(x(t)) -\varphi(x(0)) -\int_{0}^{t} {\mathcal{L}}^{N} \varphi (x(s))ds,\] 
which is a $( \mathcal{M}_t)_{t\geq 0}$-martingale under $\mathbb{P}^{N} _{d(N)}$ and we know that
\[{\mathbb{P}^{N,\tau} _{d(N)}}\left(  \sup_{t\geq 0} \lvert x(t)\rvert\leq \frac{N-[N^\alpha]}{N} \right)=1.\]
Moreover,  by Lemma \ref{generatorP2}
\[\sup_N \sup_{v\in \left[-\frac{N-[N^\alpha ]}{N},-\frac{1}{N}  \right]\cup\left[\frac{1}{N} ,\frac{N-[N^\alpha ]}{N} \right]} \lvert {\mathcal{L}}^{N}\varphi(v)\lvert <+\infty \]
and the fact that $\mathcal{L}^{N}\varphi(0)=0$, we get
\[\E^{\mathbb{P}^{N,\tau} _{d(N)}}\left( \left(\int_{t}^{t+s} {\mathcal{L}}^{N}\varphi(x(w) )dw\right) ^2\Big \vert \mathcal{M}_t \right)  \leq C h^2 .\]
We recall that
\[< M^N _{\varphi} > _t =\int_{0}^{t}\left( {\mathcal{L}}^{N}{\varphi}^2-2\varphi{\mathcal{L}}^{N}\varphi \right)(x(s))ds.\]                                                                                         
Then, using the martingale property of $(M^N _{\varphi} (t))_{t\geq 0}$, we have
\[\begin{split}
\E^{\mathbb{P}^{N,\tau} _{d(N)}} \big( (M^N _{\varphi} (t+s)  - M^N _{\varphi} (t) )^{2}\vert& \mathcal{M}_t \big) = \E^{\mathbb{P}^{N} _{d(N)}}\left(( {M^N _{\varphi}}((t+s)\wedge\tau^\alpha _N)  - {M^N_{\varphi}}(t\wedge\tau^\alpha _N))^{2}\vert \mathcal{M}_t\right)\\
&= \E^{\mathbb{P}^{N} _{d(N)}}\left( {M^N_{\varphi}}((t+s)\wedge\tau^\alpha _N) ^2  - {M^N_{\varphi}}(t\wedge\tau^\alpha _N)^{2}\vert \mathcal{M}_t\right) \\
&=\E^{\mathbb{P}^{N} _{d(N)}}\left( < M^N_{\varphi} > _{(t+s)\wedge\tau^\alpha _N}  - < M^N _{\varphi} > _{t\wedge\tau^\alpha _N} \vert \mathcal{M}_t\right)\\
&= \E^{\mathbb{P}^{N} _{d(N)}}\left( \int_{t\wedge\tau^\alpha _N}^{(t+s)\wedge\tau^\alpha _N}\left( {\mathcal{L}}^{N}\varphi^2-2\varphi{\mathcal{L}}^{N}\varphi \right)(x(w))dw \Big\vert \mathcal{M}_t\right) \\
&\leq C\, h.
\end{split}\] 
In fact, by Lemma \ref{generatorP2} we have
\[\begin{split} 
\sup_N \sup_{v\in \left[-\frac{N-[N^\alpha ]}{N},-\frac{1}{N}  \right]\cup\left[\frac{1}{N} ,\frac{N-[N^\alpha ]}{N} \right]} \lvert {\mathcal{L}}^{N}\varphi (v) \rvert & <  +\infty ,\\
\sup_N \sup_{v\in \left[-\frac{N-[N^\alpha ]}{N},-\frac{1}{N}  \right]\cup\left[\frac{1}{N} ,\frac{N-[N^\alpha ]}{N} \right]} \lvert {\mathcal{L}}^{N}\varphi^2 (v) \rvert & <  +\infty,
\end{split}\] 
in addition to ${\mathcal{L}}^{N}\varphi (0) = 0$ and $\sup_{N}{\mathcal{L}}^{N}\varphi^2(0) =\frac{\Gamma^c _{1\,2}}{N^2}<+\infty$.
$\square$
\end{preuve}

\begin{lem}\label{martingaleP22}
Let $\mathbb{Q}_{u}$ be a limit point of the relatively compact sequence $\left( \mathbb{P}^{N,\tau }_{d(N)} \right) _N$. Then, $\forall \varphi \in \mathcal{C}^{\infty}_{0} (\mathbb{R})$ and $\forall r\in(0,1)$, $(M_\varphi (t\wedge\tau_r))_{t\geq 0}$ is a $(\mathcal{M})_t$-martingale under $\mathbb{Q}_{u}$.
 \end{lem} 

\begin{preuve}[of Lemma \ref{martingaleP22}]
Let $\left({\mathbb{P}^{\Ns ,\tau} _{d(\Ns )}}\right)_{\Ns}$ be a converging subsequence. Throughout this proof we will take $N$ for $\Ns$ to simplify the notations. Let $0\leq t_1 < t_2 $ and $\Phi$ be a bounded continuous $\mathcal{M}_{t_1}$-measurable function. We have
\[\E^{\mathbb{P}^{N,\tau} _{d(N)}} \left( M^N _\varphi (t_2\wedge\tau _r) \Phi \right)=\E^{\mathbb{P}^{N,\tau} _{d(N)}} \left( M^N _\varphi (t_1\wedge\tau _r) \Phi \right).\]
Furthermore, $\forall t \geq 0$ 
\[\begin{split}
\E^{\mathbb{P}^{N,\tau} _{d(N)}} \left( \int_{0}^{t\wedge\tau _r} {\mathcal{L}}^{N} \varphi(x(s))ds \Phi  \right) 
& =  \E^{\mathbb{P}^{N,\tau} _{d(N)}} \left( \int_{0}^{t\wedge\tau_r } {\mathcal{L}}^{N} \varphi(x(s))\textbf{1}_{\left( x(s)\in I^\alpha _{N} \right) }ds \Phi  \right)\\
                                                                                                                            & +  \E^{\mathbb{P}^{N,\tau} _{d(N)}}\left( \int_{0}^{t\wedge\tau_r} {\mathcal{L}}^{N} \varphi(x(s))\textbf{1}_{( x(s)=0 ) }ds \Phi  \right) ,
\end{split}\]
with $I^\alpha _N= \left[-(N-[N^\alpha ])/N,-1/N \right]\cup\left[1/N ,(N-[N^\alpha ])/N \right]$. Using Lemma \ref{generatorP2}
\[   \lim_N \left \lvert\E^{\mathbb{P}^{N,\tau} _{d(N)}} \left( \int_{0}^{t\wedge\tau _r } \left({\mathcal{L}}^{N} \varphi(x(s))-\mathcal{L}_{\overline{a}_\infty} \varphi(x(s))\right)\textbf{1}_{\left( x(s)\in I^\alpha _{N} \right) }ds \Phi  \right)\right \rvert =0.\]
Consequently, we have to prove the two following points:
\begin{itemize}
\item  $ \lim_N \E^{\mathbb{P}^{N,\tau} _{d(N)}} \left(M _\varphi (t \wedge\tau _r) \Phi  \right)=\E^{\mathbb{Q} _{u}} \left(M _\varphi (t \wedge\tau _r) \Phi  \right)$. \\
\item $ \lim_N \E^{\mathbb{P}^{N,\tau} _{d(N)}}\left( \int_{0}^{t\wedge\tau _r} \textbf{1}_{( x(s)=0 ) }ds  \right)=0$.
\end{itemize}

We prove the first point as follows. The problem is to apply the mapping theorem to the functional $M_\varphi (t \wedge \tau_r)$ and to do this we must have $\mathbb{Q}_u (D_{M_\varphi (t \wedge \tau_r)} )=0$, where $D_{M_\varphi (t \wedge \tau_r)}$ is the set of discontinuities of $M_\varphi (t \wedge \tau_r)$ for the Skorokhod topology. While $M_\varphi (t)$ is continuous for this topology, it is not necessarily true for $\tau_r$. However, we can follow the proof of Lemma 11.1.3 in \cite{stroock} and then use a family of stopping times for which we can apply the mapping theorem.

We know that the size of the jumps of $(X^N _t)_t$ is constant equal to $1/N$, therefore we have $\mathbb{Q}_u (\mathcal{C}([0,+\infty),\mathbb{R}))=1$ (see Theorem 13.4 in \cite{billingsley} for instance). Then
\[ \mathbb{Q}_u \big(D_{M_\varphi (t \wedge \tau_r)} \big)=\mathbb{Q}_u \big(D_{M_\varphi (t \wedge \tau_r)} \cap \mathcal{C}([0,+\infty),\mathbb{R})\big).\]
We recall that  the Skorokhod topology on $\mathcal{C}([0,+\infty),\mathbb{R})$ coincides with the usual topology defined on this space. Therefore, $D_{M_\varphi (t \wedge \tau_r)} \cap \mathcal{C}([0,+\infty),\mathbb{R})$ is the set of discontinuities of $M_\varphi (t \wedge \tau_r)$ under the topology of  $\mathcal{C}([0,+\infty),\mathbb{R})$, and $\tau_r$ restrict to $\mathcal{C}([0,+\infty),\mathbb{R})$ is lower semi-continuous. Consequently, according to the proof of lemmas 11.1.2 in \cite{stroock}, there exists a sequence $(r_n)_n$ such that $r_n \nearrow r$ and
\[\mathbb{Q}_u\left((\tau_{r_n}<+\infty)\cap D_{ \tau_{r_n}}\cap \mathcal{C}([0,+\infty),\mathbb{R})  \right)= \mathbb{Q}_u\left((\tau_{r_n}<+\infty)\cap D_{ \tau_r}\big)  \right)=0.\]
Then, $\mathbb{Q}_u (D_{M_\varphi (t \wedge \tau_{r_n})} )=0$ and we can apply the mapping theorem to $M_\varphi (t \wedge \tau_{r_n})$, i.e
 \[  \lim_N \E^{\mathbb{P}^{N,\tau} _{d(N)}} \left(M _\varphi (t \wedge\tau _{r_n}) \Phi  \right)=\E^{\mathbb{Q} _{u}} \left(M _\varphi (t \wedge\tau _{r_n}) \Phi  \right) .\]
Finally, we obtain
\[\E^{\mathbb{Q} _{u}} \left(M _\varphi (t_2 \wedge\tau _{r_n}) \Phi  \right)=\E^{\mathbb{Q} _{u}} \left(M _\varphi (t_1 \wedge\tau _{r_n}) \Phi  \right),\] 
and 
\[\lim_n \E^{\mathbb{Q} _{u}} \left(M _\varphi (t \wedge\tau _{r_n}) \Phi  \right) =  \E^{\mathbb{Q} _{u}} \left(M _\varphi (t \wedge\tau _{r}) \Phi  \right)\]
because $\tau_{r_n} \nearrow \tau_{r}$. Consequently, 
\[ \E^{\mathbb{Q} _{u}} \left(M _\varphi (t_2 \wedge\tau _{r}) \Phi  \right)=\E^{\mathbb{Q} _{u}} \left(M _\varphi (t_1 \wedge\tau _{r}) \Phi  \right) .\]

For the second point, we have 
\[\E^{\mathbb{P}^{N,\tau} _{d(N)}}\left( \int_{0}^{t\wedge\tau _r} \textbf{1}_{( x(s)=0 ) }ds  \right)= \E_{d(N)}\left[ \int_{0}^{t\wedge\tau _r} \textbf{1}_{( X^N _s =0 ) }  ds\right] \leq  \E_{0}\left[ \int_{0}^t \textbf{1}_{( X^N _s =0 ) }  ds\right], \] 
since the stopped process spends less time in $0$ than the original process and the last inequality is given by the Markov property. We denote by $N^0 _t$ the number of returns in $0$ during the time interval $[0,t]$ and by $(Y_j)_{j\geq0}$ the renewal process associated with the return times in $0$, $(\sigma^{(i)}_0)_{i\geq1}$, of the process $(X^N _t)_t$, with $Y_0 =\sigma^{(0)}_0=0$. Moreover, for $\alpha' \in (0,1)$
\[ \begin{split}
\E_{0}\left[ \int_{0}^t \textbf{1}_{( X^N _s =0 ) }  ds\right]& \leq t \mathbb{P}_{0}\left(N^0 _t\geq[ N^{1+\alpha'}]\right)+ \E_{0}\left[\sum_{j=0}^{[N^{1+\alpha'}]} \int_{\sigma^{(j)}_0}^{\sigma^{(j+1)}_0} \textbf{1}_{( X^N _s =0 ) }ds \right] \\
& \leq \frac{t}{[N^{1+\alpha'}]}\mathbb{E}_0 [ N^0 _t] + \frac{[N^{1+\alpha'}]+1}{\Gamma^c _{2\,1}},
\end{split}\]
since $\Big(  \int_{\sigma^{(j)}_0}^{\sigma^{(j+1)}_0} \textbf{1}_{( X^N _s =0 ) }ds\Big)_j $ is an i.i.d sequence with mean $1/\Gamma^c_{2\,1}$. We recall that $N^0 _t +1$ is a stopping time for $(Y_j)_{j\geq1}$. Then,
\[ \E_0 \left[ \sigma^{(N^0 _t +1)}_0 \right]= \E_0 \left[ \sum_{j=1}^{N^0 _t +1} Y_{j} \right] =\left(\E_0\left[ N^0 _t\right]+1 \right)\E_0\left[ \sigma_0 ^{(1)}\right].\]
Furthermore, 
\[ \begin{split}
\E_0 \left[ \sigma^{(N^0 _t +1)}_0 \right]&=  \E_0 \left[ \sigma^{(N^0 _t +1)}_0\Big( \textbf{1}_{\big(X^N_t=0\big)}+\textbf{1}_{\big(X^N_t\not=0\big)}\Big)\right]\\
&=\E_0 \left[ \inf\big(s>T^N_t,\,X^N_{t+s}=0\big) \textbf{1}_{\big(X^N_t=0\big)}\right]+ \E_0 \left[ \inf\big(s>0,\,X^N_{t+s}=0\big)\textbf{1}_{\big(X^N_t\not=0\big)}\right]
\end{split}\]
where $T_t^N=\inf\big(s>0,\, X^N_{t+s}\not=0\big)$.
Then, using the Markov property we get
\[\begin{split} \E_0 \left[ \inf\big(s>T^N_t,\,X^N_{t+s}=0\big) \textbf{1}_{\big(X^N_t=0\big)}\right]&=\Big(t+\mathbb{E}_0\left[ \sigma_0 ^{(1)}\right]\Big)\mathbb{P}_0\big(X^N_t=0\big)\\
&\leq t\mathbb{P}_0\big(X^N_t=0\big)+ \frac{2N-1}{\Gamma^c _{1\,2}}\end{split}\]
and
\[\begin{split}
\E_0 \left[ \inf\big(s>0,\,X^N_{t+s}=0\big)\textbf{1}_{\big(X^N_t\not=0\big)}\right]&=\sum_{\substack{j=-(N-1)\\j\not=0}}^{N-1} \E_0 \left[ \inf\big(s>0,\,X^N_{t+s}=0\big)\textbf{1}_{\big(X^N_t=j\big)}\right]\\
&=\sum_{\substack{j=-(N-1)\\j\not=0}}^{N-1}\Big(t+ \E_j \left[ \tau^{(1)}_0\right]\Big)\mathbb{P}_0\big(X^N_t=j\big)\\
&=\sum_{j=1}^N \sum_{l=1}^j \frac{N-l}{\Gamma^c_{l\,l+1}}\mathbb{P}_0\left(\lvert X^N_t\rvert =j\right)+t\mathbb{P}_0(X^N_t \not=0)\\
&\leq K\frac{N-1}{N^2}\E_0 \left[\lvert X^N_t\rvert \textbf{1}_{\big(X^N_t\not=0\big)} \right]+t\mathbb{P}_0(X^N_t \not=0),
\end{split}\]
where $K$ is a constant independent of $N$. 
Consequently, 
\[ \E_0 [ N^0 _t] \leq \tilde{K}\frac{\Gamma^c_{1\,2}}{2N-1} -\frac{2N-1}{\Gamma^c_{1\,2}}=\mathcal{O}(N),\]
and
\[ \E_{0}\left[ \int_{0}^t \textbf{1}_{( X^N _s =0 ) }  ds\right] =\mathcal{O}\left( \frac{1}{N^{\alpha' \wedge (1-\alpha')} }\right).\]
$\square$
\end{preuve}

From Lemma \ref{martingaleP22}, we have $\forall r\in (0,1)$, $\mathbb{Q}_u =\overline{\mathbb{P}}_u$ on $\mathcal{M}_{\tau_r}$. From this relation and the fact that $ \mathbb{Q}_u(\mathcal{C}([0,+\infty),\mathbb{R})) =\overline{\mathbb{P}}_u (\mathcal{C}([0,+\infty),\mathbb{R})) =1$, $\mathbb{Q}_u =\overline{\mathbb{P}}_u$ on $\mathcal{M}_{\tau_1}$ since $\tau_r \nearrow \tau_1$ as $r\nearrow 1$.

Let $f \in \mathcal{C}^\infty ([0,1])$ with compact support included in $[0,1)$ and let $\left({\mathbb{P}^{\Ns,\tau} _{d(\Ns)}}\right)_{N'}$ be a converging subsequence as in the previous proof. We have 
\begin{equation}\label{egaliteP2}
\mathcal{T}^{l(\Ns )} _f ( \omega , t )= \E_{d(\Ns)}\left[ f\left( \lvert X^{\Ns} _t\rvert +\frac{1}{\Ns} \right) \textbf{1}_{\left(t<  \tau^{\alpha}_{\Ns}\right)} \right] + r(\Ns),\end{equation}
with
\[ \begin{split}
r(N) = \E_{d(N)}& \left[  e^{-\Lambda^c_{N} \int_0 ^t \textbf{1}_{\left(\lvert X^N_v \rvert = \frac{N-1}{N}\right)}dv -\Lambda^c_{N-1} \int_0 ^t \textbf{1}_{\left(\lvert X^N_v \rvert = \frac{N-2}{N}\right)}dv }\right.  \\
& \times \left.f\left( \lvert X^N _t\rvert +\frac{1}{N} \right)\left( \textbf{1}_{(  \tau^{\alpha}_{N} \leq t < \tau^0 _{N} +\lambda )} + \textbf{1}_{(t\geq \tau^0 _{N} +\lambda )}\right)  \right],\end{split}\]
where $\tau^0 _{N}=\tau^ {(1)}_{ (N-1)/N}\wedge \tau^ {(1)}_{- (N-1)/N}$ and $\lambda\in (0,t)$.
Using Lemma \ref{tightP2} and Lemma \ref{martingaleP22}, we can study the first term on the right in \eqref{egaliteP2}.
\[\begin{split} 
 \E_{d(\Ns)}\left[ f\left( \lvert X^{\Ns} _t\rvert +\frac{1}{\Ns} \right) \textbf{1}_{\left(t<  \tau^{\alpha}_{\Ns}\right)} \right] &= \E^{\mathbb{P}^{\Ns}_{d(\Ns)}}\left[ f\left( \lvert x(t)\rvert +\frac{1}{\Ns} \right) \textbf{1}_{\left(t<  \tau^{\alpha}_{\Ns}\right)} \right]\\
 & =\E^{\mathbb{P}^{\Ns,\tau}_{d(\Ns)}} \left[ f\left( \lvert x(t) \rvert +\frac{1}{\Ns} \right) \textbf{1}_{\left(t<  \tau^{\alpha}_{\Ns}\right)} \right] ,
 \end{split}\] 
since $\left(t<  \tau^{\alpha}_{N}\right) \in \mathcal{M}_{\tau ^\alpha _N}$ and $ \mathbb{P}^{N} _{d(N)}= \mathbb{P}^{N,\tau} _{d(N)} $ on $ \mathcal{M}_{\tau ^\alpha _N}$. Moreover, 
\[\begin{split} 
\E^{\mathbb{P}^{\Ns,\tau}_{d(\Ns)}}& \left[ f\left( \lvert x(t) \rvert +\frac{1}{\Ns} \right) \textbf{1}_{\left(t<  \tau^{\alpha}_{\Ns}\right)} \right] \\
&= \E^{\mathbb{P}^{\Ns,\tau}_{d(\Ns)}} \left[ f\left( \lvert x(t) \rvert +\frac{1}{\Ns} \right)  \right] - f\left( \frac{\Ns-[\Ns ^{\alpha}] +1}{\Ns} \right) \mathbb{P}^{\Ns,\tau}_{d(\Ns)}\left(t\geq  \tau^{\alpha}_{\Ns}\right)\\
&=\E^{\mathbb{P}^{\Ns,\tau}_{d(\Ns)}} \left[ f( \lvert x(t) \rvert)  \right] +o(1).
\end{split} \]   
Consequently,
\[ \lim_{\Ns}\E_{d(\Ns)}\left[ f\left( \lvert X^{\Ns} _t\rvert +\frac{1}{\Ns} \right) \textbf{1}_{\left(t<  \tau^{\alpha}_{\Ns}\right)} \right] =\E^{\mathbb{Q}_u}[f(\lvert x(t)\rvert)] .\]
However,
\[ \E^{\mathbb{Q}_u}\left[f(\lvert x(t)\rvert)\textbf{1}_{( \tau_1 \leq t )}\right] =0.\]
In fact, let
$\tau_s = \inf (t\geq 0, \forall v> t, x(v)=x(t) )$ be the first time for which the process becomes constant. From the Portmanteau theorem 
\[ 1=\overline{\lim _{\Ns}} \,\, \mathbb{P} ^{\Ns ,\tau}_{d(\Ns)}\left( \overline{\left( \tau_s \leq \tau_1 \right)}\right) \leq \mathbb{Q}_u\left( \overline{\left( \tau_s \leq \tau_1 \right)} \right),\]
where $\overline{A}$ denote the closure under the Skorokhod topology of  a subset $A$ of $\D ([0,+\infty),\mathbb{R})$.
Moreover, we have 
\[\begin{split}  \overline{\left( \tau_s \leq \tau_1 \right) \cap \left(\tau_1 \leq t \right)}&\cap(x(0)=u) \cap\mathcal{C}([0,+\infty),\mathbb{R})\\
&= \left( \tau_s \leq \tau_1 \right) \cap \left(\tau_1 \leq t \right)\cap(x(0)=u)\cap \mathcal{C}([0,+\infty),\mathbb{R}).\end{split} \]
Then,
\[ \begin{split}  \mathbb{Q}_u \big( \lvert x(t)\rvert \in supp(f),  \tau_1 \leq t  \big)& \leq \mathbb{Q}_u \big( \lvert x(t)\rvert \in supp(f) ,\tau_s \leq \tau_1 \leq t \big)\\
&\leq  \mathbb{Q}_u \big( \lvert x(t)\rvert \in supp(f) , \lvert x(t) \rvert =1\big)=0, \end{split}\]
and 
\[ \lim_{\Ns}\E_{d(\Ns)}\left[ f\left( \lvert X^{\Ns} _t\rvert +\frac{1}{\Ns} \right) \textbf{1}_{\left(t<  \tau^{\alpha}_{\Ns}\right)} \right]=  \E^{\mathbb{Q}_u}\left[f(\lvert x(t)\rvert)\textbf{1}_{(t< \tau_1)}\right]=\E^{\overline{\mathbb{P}}_u}\left[f(\lvert x(t)\rvert)\textbf{1}_{(t< \tau_1)}\right].\]
Finally, by the following lemma we get
\[ \lim_{\Ns} \mathcal{T}^{l(\Ns )} _f ( \omega , t )= \E^{\overline{\mathbb{P}}_u}\left[f(\lvert x(t)\rvert)\textbf{1}_{(t< \tau_1)}\right].\] 
We can remark that this limit does not depend on the subsequence $(\Ns)$. The following lemma represents the loss of energy from the propagating modes produced by the coupling between the propagating and the radiating modes.  Moreover, this lemma implies the absorbing condition at the boundary $1$ in Theorem \ref{hfapproxP2}, which implies the dissipation behavior in Theorem \ref{hfexpdecP2}.   
\begin{lem}\label{resteP2}
$\lim_{\Ns} r(\Ns)=0.$
\end{lem}
\begin{preuve}
\[ \begin{split} \lvert r(\Ns)\rvert \leq \|f\|_{\infty}& \left( \E_{d(\Ns)}\left[  e^{-\Lambda^c_{\Ns} \int_0 ^t \textbf{1}_{\left(\lvert X^{\Ns}_s \rvert = \frac{\Ns-1}{\Ns}\right)}ds } \textbf{1}_{(t\geq\tau^0 _{\Ns} +\lambda )} \right] \right. \\
&+\left. \mathbb{P}_{d(\Ns)}\left( \lvert X^{\Ns} _t\rvert +\frac{1}{\Ns} \in supp(f),  \tau^{\alpha}_{\Ns} \leq t < \tau^0 _{\Ns} +\lambda   \right) \right). \end{split} \]
First, let $\alpha'\in (3/4,1)$ and $N^{N} _t$ the number of passages in $(N-1)/ N$ during the time interval $[0,t]$.
\[ \begin{split} \E_{d(\Ns)}\left[  e^{-\Lambda^c_{\Ns} \int_0 ^t \textbf{1}_{\left(\lvert X^{\Ns}_s \rvert = \frac{\Ns-1}{\Ns}\right)}ds } \textbf{1}_{(t\geq \tau^0 _{\Ns} +\lambda )} \right] &\leq \E_{d(\Ns)}\left[  e^{-\Lambda^c_{\Ns} \int_0 ^t \textbf{1}_{\left(\lvert X^{\Ns}_s \rvert = \frac{\Ns-1}{\Ns}\right)}ds }\right. \\
&\times\left.\left( \textbf{1}_{(t\geq \tau^{(1)} _{(\Ns-1)/ \Ns} +\lambda )}+\textbf{1}_{(t\geq \tau^{(1)}_{-(\Ns-1)/ \Ns} +\lambda )} \right)\right].  \end{split}\]
We shall work only with  $\E_{d(\Ns)}\left[  e^{-\Lambda^c_{\Ns} \int_0 ^t \textbf{1}_{\left(\lvert X^{\Ns}_s \rvert = \frac{\Ns-1}{\Ns}\right)}ds } \textbf{1}_{(t\geq \tau^{(1)}_{(\Ns-1)/\Ns} +\lambda )} \right]$ but the same proof works for the other term. 
\[\begin{split}
\E_{d(\Ns)}&\left[  e^{-\Lambda^c_{\Ns} \int_0 ^t \textbf{1}_{\left(\lvert X^{\Ns}_s \rvert = \frac{\Ns-1}{\Ns}\right)}ds } \textbf{1}_{(t\geq \tau^{(1)}_{(\Ns-1)/\Ns} +\lambda )} \right] \\
& \leq   \E_{d(\Ns)}\left[  e^{-\Lambda^c_{\Ns} \int_0 ^t \textbf{1}_{\left(X^{\Ns}_s = \frac{\Ns-1}{\Ns}\right)}ds } \textbf{1}_{(N^{\Ns} _t\geq [{N}^{\alpha'}]+1)} \right]\\
&+\mathbb{P} _{d(\Ns)}\left(  N^{\Ns} _t\leq [{N}^{\alpha'}], t -\tau^{(1)}_{(\Ns-1)/\Ns} \geq \lambda \right).
\end{split} \]
On  $(N^{\Ns} _t\geq [{\Ns}^{\alpha'}]+1)$, we have
\[ \begin{split}
e^{-\Lambda^c_{\Ns} \int_0 ^t \textbf{1}_{\left(\lvert X^{\Ns}_s  \rvert =  \frac{\Ns-1}{\Ns}\right)}ds }& \leq e^{-\Lambda^c_{\Ns} \sum_{j=1}^{N^{\Ns} _t -1}  \int_{\tau^{(j)}_{(\Ns-1)/ \Ns}}^{\tau^{(j+1)}_{(\Ns-1)/\Ns}} \textbf{1}_{\left( X^{\Ns} _s = \frac{\Ns-1}{\Ns} \right) }ds}\\
&\leq \prod_{j=1}^{[{\Ns}^{\alpha'}]} e^{ -\Lambda^c_{\Ns}\int_{\tau^{(j)}_{(\Ns-1)/\Ns}}^{\tau^{(j+1)}_{(\Ns-1)/\Ns}} \textbf{1}_{\left( X^{\Ns} _s = \frac{\Ns-1}{\Ns}\right ) }ds}.\end{split}\]
We denote by $\sigma^{(1)}_{(N-1)/N}$ the time of the first return in ${(N-1)/N}$, then
\[  \begin{split}
\E_{d(\Ns)}&\left[  e^{-\Lambda^c_{\Ns} \int_0 ^t \textbf{1}_{\left(X^{\Ns} _s =  \frac{\Ns-1}{\Ns}\right)}ds } \textbf{1}_{(N^{\Ns}  _t\geq [{\Ns} ^{\alpha'}]+1)} \right] \\
&\leq \prod_{j=1}^{[{\Ns} ^{\alpha'}]} \E_{d(\Ns)}\left[ e^{-\Lambda^c_{\Ns}  \int_{\tau^{(j)}_{(\Ns-1)/\Ns}}^{\tau^{(j+1)}_{(\Ns-1)/\Ns}} \textbf{1}_{\left( X^{\Ns} _s = \frac{\Ns-1}{\Ns} \right) }ds}\right]\\
 &\leq \left( \E_{\frac{\Ns-1}{\Ns}}\left[ e^{-\Lambda^c_{\Ns}  \int_{0}^{\sigma^{(1)}_{(\Ns-1)/\Ns}} \textbf{1}_{\left( X^{\Ns} _s = \frac{\Ns-1}{\Ns} \right) }ds}\right] \right)^{[{\Ns} ^{\alpha'}]}
\end{split}\]
since $\left(   \int_{\tau^{(j)}_{(N-1)/N}}^{\tau^{(j+1)}_{(N-1)/N}} \textbf{1}_{\left( X^N _s = \frac{N-1}{N} \right) }ds \right)_j $ is an i.i.d sequence.
Moreover, we can check that $\Lambda^c_N\geq C N^{3/2}$ and then 
\[\E_{\frac{\Ns-1}{\Ns}}\left[ e^{-\Lambda^c_{\Ns}  \int_{0}^{\sigma^{(1)}_{(\Ns-1)/\Ns}} \textbf{1}_{\left( X^{\Ns} _s = \frac{\Ns-1}{\Ns} \right) }ds}\right]\leq \E_{\frac{\Ns-1}{\Ns}}\left[ e^{-C\,\Ns^{3/2}  \int_{0}^{\sigma^{(1)}_{(\Ns-1)/\Ns}} \textbf{1}_{\left( X^{\Ns} _s = \frac{\Ns-1}{\Ns} \right) }ds}\right].\]
In fact, a computation gives
\[\begin{split}
\Lambda^c_N = \frac{ a k^4   A^2_N}{16 \pi \beta_N}\int_{n_1 k d \theta }^{n_1 kd } & \frac{\eta \sqrt{\eta^2-(n_1 k d\theta)^2}  }{\sqrt{(n_1 kd)^2 -\eta^2}\left( 1+(\beta_N-\frac{1}{d}\sqrt{(n_1 k d)^2 -\eta^2})^2   \right)} \\&\times \frac{S(\eta -\sigma_N,\eta-\sigma_N)}{ (\eta^2 -(n_1 kd \theta)^2)\sin^2 (\eta )+\eta^2 \cos^2 (\eta ) }d\eta.
\end{split}\]
However, we recall that the support of $S$ lies in the square $\left[- \frac{3\pi}{2},\frac{3\pi}{2} \right]\times \left[- \frac{3\pi}{2},\frac{3\pi}{2} \right]$, then we can restrict the integration over $\left[n_1 k d \theta,n_1 k d \theta+\frac{3\pi}{2}\right]$. Moreover,
\[\begin{split}
\left( \beta_N -\sqrt{(n_1 k)^2 -\eta^2/d^2} \right)^2 &=(\eta -\sigma)^2 \frac{y^2}{1-y^2}, \quad \text{ for some }y \in \left[ \frac{\sigma_N}{n_1 k d} ,\frac{\eta }{n_1 k d} \right] \\
&\leq \left( \frac{3\pi}{2}\right)^2 \frac{\eta ^2/(n_1 k d)^2}{1-\eta^2/(n_1 k d)^2 }\leq K,
\end{split}\]
where $K$ stands for a constant independent of $N$, because $\theta <1$ and $k\gg1$.
Therefore,
\[\Lambda^c _N \geq K' \int _{n_1 k d \theta}^{n_1 k  d \theta+\frac{3\pi}{2}}\eta  \sqrt{\eta^2-(n_1 k d \theta)^2}  d\eta \geq K'' N^{3/2}, \]
where we assume that the function $S$ has a positive minimum, and then $K''>0$. 

Now, let us remark that $\forall v\in [0,\ln(N^{1/4})]$
\[e^{-v}\leq 1-\frac{1}{N^{1/4}}v\]
and
\[ \E_{\frac{\Ns-1}{\Ns}}\left[   \int_{0}^{\sigma^{(1)}_{(\Ns-1)/\Ns}} \textbf{1}_{\left( X^{\Ns} _s = \frac{\Ns-1}{\Ns} \right) }ds\right]=\frac{1}{\Gamma^c _{\Ns-1\, \Ns}}.\]
Then, we get
\[ \E_{\frac{\Ns-1}{\Ns}}\left[ e^{-K''\,\Ns^{3/2}  \int_{0}^{\sigma^{(1)}_{(\Ns-1)/\Ns}} \textbf{1}_{\left( X^{\Ns} _s = \frac{\Ns-1}{\Ns} \right) }ds}\right]\leq 1-\frac{K''}{\Ns^{3/4}}\left(1-\frac{1}{\ln(\Ns^{1/4})}\right)\]
and
\[\E_{d(\Ns)}\left[  e^{-\Lambda_{\Ns} \int_0 ^t \textbf{1}_{\left(X^{\Ns} _s =  \frac{\Ns-1}{\Ns}\right)}ds } \textbf{1}_{(N^{\Ns}  _t\geq [{\Ns} ^{\alpha'}]+1)} \right] \leq e^{[\Ns^{\alpha'}] \ln\left[1-\frac{K''}{\Ns^{3/4}}\left(1-\frac{1}{\ln(\Ns^{1/4})}\right)\right]}.\]
Moreover,
\[\begin{split}
\mathbb{P}_{d(N)}\left(  N^N _t\leq [N^{\alpha'}], t -\tau^{(1)}_{(N-1)/N} \geq \lambda \right)&\leq \mathbb{P}_{d(N)}\left(\tau^{([N^{\alpha'}] )}_{(N-1)/N} -\tau^{(1)}_{(N-1)/N} \geq \lambda \right)\\
&\leq \frac{1}{\lambda}\E_{d(N)}\left(\tau^{([N^{\alpha'}] )}_{(N-1)/N} -\tau^{(1)}_{(N-1)/N}\right)\\
&\leq \frac{N^{\alpha'}}{\lambda}\E_{(N-1)/N} \left[ \sigma_{(N-1)/N}^{(1)} \right] \\
&\leq \frac{K}{N^{1-\alpha'}}.
\end{split} \]
Consequently,
\[\lim_{\Ns} \E_{d(\Ns)}\left[  e^{-\Lambda^c _{\Ns} \int_0 ^t \textbf{1}_{\left(\lvert X^{\Ns}_s \rvert = \frac{\Ns-1}{\Ns}\right)}ds } \textbf{1}_{(t\geq \tau^{0}_{\Ns} +\lambda )}  \right] =0.\]
Second, let $c_f \in (0,1)$ such that $supp(f)\subset [0, c_f -1/\Ns]$ and $x\in [0,c_f)$, then
\[ \begin{split}
\mathbb{P}_{d(N)}&\left( \lvert X^N _t\rvert +\frac{1}{N} \in supp(f),  \tau^{\alpha}_{N} \leq t < \tau^{0}_{N} +\lambda   \right) \leq  \mathbb{P}_{d(N)}\left( X^N _t \in [-c_f , c_f],  \tau^{\alpha}_{N} \leq t < \tau^{0}_{N} +\lambda   \right) \\
&\leq  \mathbb{P}_{d(N)}\left( X^N _t \in [-c_f , c_f], X^N _{\tau^{\alpha}_{N}}= \frac{N-[N^\alpha]}{N}, \tau^{\alpha}_{N} \leq t < \tau^{0}_{N} +\lambda   \right) \\
&+ \mathbb{P}_{d(N)}\left( X^N _t \in [-c_f , c_f],  X^N _{\tau^{\alpha}_{N}}=-\frac{N-[N^\alpha]}{N},  \tau^{\alpha}_{N} \leq t < \tau^{0}_{N} +\lambda   \right). 
\end{split} \] 
We shall treat only the case where $X^N _{\tau^{\alpha}_{N}}= (N-[N^\alpha])/N$, but the following proof works also in the other case. Let $\tilde{c}_f \in (c_f,1)$, $\rho\in(0,1)$ such that $[\tilde{c}_f -\rho,\tilde{c}_f+\rho]\subset (c_f ,1)$ and $\lambda' \in (0,1)$. Using the strong Markov property we have
\[ \begin{split}
\mathbb{P}_{d(N)}& \left( X^N _t \in [-c_f , c_f], X^N _{\tau^{\alpha}_{N}}= \frac{N-[N^\alpha]}{N}, \tau^{\alpha}_{N} \leq t < \tau^{0}_{N} +\lambda   \right)\\
& \leq \mathbb{P}_{\frac{N-[N^{\alpha}]}{N}}\left( \tau^{(1)}_{(N-1)/N} > \lambda' \right) + \mathbb{P}^N _{\frac{[N\tilde{c}_f]}{N}} \left( \tau_{ \tilde{c}_f \pm \rho} \leq \lambda +\lambda' \right),
\end{split}\]
where $\tau_{ \tilde{c}_f \pm \rho}=\inf (t \geq 0, \lvert x(t)- \tilde{c}_f  \rvert \geq \rho )$.
First, a computation gives
\[  \mathbb{P}_{\frac{N-[N^{\alpha}]}{N}}\left( \tau^{(1)}_{(N-1)/N} > \lambda' \right)  \leq\frac{1}{\lambda'} \E_{\frac{N-[N^{\alpha}]}{N}}\left[ \tau^{(1)}_{(N-1)/N}  \right]=\frac{1}{\lambda'}\sum_{l=N-[N^\alpha]}^{N-2} \frac{N+1 +l}{\Gamma^c _{l+1\,l+2}}\leq \frac{K}{N^{1-\alpha}}.\]
Second, the sequence $(r(\Ns))_{\Ns}$ is bounded. Let $(r(\Nss))_{\Nss}$ be a converging subsequence. We recall that  $ \mathbb{P}^{N}_{c(N)} = \mathbb{P}^{N,\tau}_{c(N)} $ on $ \mathcal{M}_{\tau ^\alpha _N}$, where $c(N)=[N\tilde{c}_f]/N$, and by Lemma \ref{tightP2} the sequence $\left(\mathbb{P}^{\Nss,\tau}_{c(\Nss) }\right)_{\Nss}$ is tight. Let $\left(\mathbb{P}^{\Nsss,\tau}_{c(\Nsss) }\right)_{\Nsss}$ be a converging subsequence to $\mathbb{Q}_{\tilde{c}_f}$. Moreover, $\tau_{ \tilde{c}_f \pm \rho} \leq \tau^\alpha _N$ and therefore $\left( \tau_{ \tilde{c}_f \pm \rho}  \leq \lambda + \lambda' \right) \in \mathcal{M}_{\tau ^\alpha _N}$. Consequently, by the Portmanteau theorem 
\[\begin{split}
\overline{\lim_{\Nsss}}\,\,\mathbb{P}^{\Nsss}_{c(\Nsss)} \left(\tau_{ \tilde{c}_f \pm \rho} \leq \lambda + \lambda' \right)&=\overline{\lim_{\Nsss}}\,\,\mathbb{P}^{\Nsss,\tau}_{c(\Nsss)} \left( \tau_{ \tilde{c}_f \pm \rho} \leq \lambda + \lambda' \right) \\
&\leq \overline{\lim_{\Nsss}}\,\,\mathbb{P}^{\Nsss,\tau}_{c(\Nsss)} \left( {\overline{\left(\tau_{ \tilde{c}_f \pm \rho}\leq \lambda + \lambda' \right) }}\right)\\
&\leq \mathbb{Q}_{\tilde{c}_f} \left( {\overline{\left(\tau_{ \tilde{c}_f \pm \rho} \leq \lambda+\lambda' \right) }}\right).
\end{split}\] 
We recall that $\mathbb{Q}_{\tilde{c}_f} \left( \mathcal{C}([0,+\infty),\mathbb{R} )\right) =1$ and we can show that
\[ {\overline{\left(\tau_{ \tilde{c}_f \pm \rho} \leq\lambda+ \lambda' \right) }}\cap  \mathcal{C}([0,+\infty),\mathbb{R} )={\left(\tau_{ \tilde{c}_f \pm \rho} \leq \lambda +\lambda' \right) }\cap  \mathcal{C}([0,+\infty),\mathbb{R} ).\]
Then,
\[ \overline{\lim_{\Nsss}}\,\mathbb{P}^{\Nsss,\tau}_{c(\Nsss)} \left(\tau_{ \tilde{c}_f \pm \rho}\leq \lambda' \right)\leq \mathbb{Q}_{\tilde{c}_f} \left(\tau_{ \tilde{c}_f \pm \rho} \leq \lambda + \lambda' \right),  \] 
and
\[\overline{\lim_{\Nsss}} \, r(\Nsss)\leq  \mathbb{Q}_{\tilde{c}_f} \left(\tau_{ \tilde{c}_f \pm \rho} \leq \lambda + \lambda' \right).\]
Finally, $\lim_{\Nsss} r(\Nsss)=0$ and the limit of all subsequences $(r(\Nss))_{\Nss}$ of $(r(\Ns))_{\Ns}$ is $0$. $\square$
\end{preuve}

To finish, $\left(\mathcal{T}^{l(N )} _f ( \omega , t )\right)_N$ is a bounded sequence. Let $\left(\mathcal{T}^{l(\Ns )} _f ( \omega , t )\right)_{\Ns}$ be a converging subsequence.  By the previous work, there exists an another subsequence such that 
 \[   \lim_{\Nss}\mathcal{T}^{l(\Nss )} _f ( \omega , t ) =  \E^{\overline{\mathbb{P}}_u}\left[f(\lvert x(t)\rvert)\textbf{1}_{(t< \tau_1)}\right],\] 
where the limit does not depend on the particular subsequence, then all subsequence limits of $\left(\mathcal{T}^{l(N )} _f ( \omega , t )\right)_N$ are equal to $ \E^{\overline{\mathbb{P}}_u}\left[f(\lvert x(t)\rvert)\textbf{1}_{(t< \tau_1)}\right]$.  Consequently,
 \[   \lim_{N}\mathcal{T}^{l(N )} _f ( \omega , t ) =  \E^{\overline{\mathbb{P}}_u}\left[f(\lvert x(t)\rvert)\textbf{1}_{(t< \tau_1)}\right].\] 

Now, we have to show that this equality holds even for a sequence $(l(N))_N$ such that $l(N)/N \to u=1$, i.e $\lim_{N}\mathcal{T}^{l(N )} _f ( \omega , t ) =0$. To do this, we write for $\lambda\in(0,t)$,
\[\begin{split}
 \mathcal{T}^{l(N )} _f ( \omega , t ) &\leq \|f\|_{\infty}\Big(\mathbb{P}_{d(N)}\left( t < \tau^{(1)}_{(N-1)/N)} +\lambda \right) \\
 &\quad+ \E_{d(N)}\left[ e^{-\Lambda_N \int_0 ^t \textbf{1}_{\left(X^N _u=\frac{N-1}{N} \right)}  du    }\textbf{1}_{\left(t \geq \tau^{(1)}_{(N-1)/N} +\lambda \right)}   \right]\Big).
\end{split}\]
We have already shown in Lemma \ref{resteP2} that the second term on the right in the previous inequality goes to $0$. The proof did not depend of the sequence $(d(N))_N$. Moreover, we have
\[\mathbb{P}_{d(N)}\left( t < \tau^{(1)}_{(N-1)/N)}+\lambda\right) \leq \frac{1}{t-\lambda} \E_{d(N)}\left[ \tau^{(1)}_{(N-1)/N)} \right] ,\]
and 
\[ \E_{d(N)}\left[ \tau^{(1)}_{(N-1)/N)} \right] =\sum_{j=l(N)-1}^{N-2}\frac{N+1+j}{\Gamma_{j+1,j+2}}\leq K\left(1-\frac{l(N)}{N} \right). \]
Consequently, we have $\forall u \in [0,1]$ and $\forall (l(N))_N$ such that $l(N)/N \to u$,
\begin{equation}\label{prorephfaP2}   \lim_{N}\mathcal{T}^{l(N )} _f ( \omega , t ) =  \E^{\overline{\mathbb{P}}_u}\left[f(\lvert x(t)\rvert)\textbf{1}_{(t< \tau_1)}\right],\end{equation} 
where the limit satisfies the required conditions.
Finally, from the decomposition used in the proof of Theorem \ref{hfexpdecP2}, we have $\forall \varphi \in L^2(0,1)$ and $\tilde{\varphi}$ a smooth function with compact support
\[
\| \mathcal{T}^N_{\varphi}(L,.)- \mathcal{T}_{\varphi}(L,.)\|_{L^2(0,1)}\leq 2\| \varphi- \tilde{\varphi}\|_{L^2(0,1)}+ \|\mathcal{T}^N_{\tilde{\varphi}}(L,.)- \mathcal{T}_{\tilde{\varphi}}(L,.)\|_{L^2(0,1)}.\]
Using the density of the smooth functions with compact support in $L^2(0,1)$ for $\|.\|_{L^2(0,1)}$ and the dominated convergence theorem we get the first point of Theorem \ref{hfapproxP2}.
The second point is a direct consequence of the probabilistic representation \eqref{prorephfaP2}
and the density for the sup norm over $[0,1]$ in $\{\varphi\in\mathcal{C}^0([0,1]),\,\varphi(1)=0\}$ of the smooth functions with compact support included in $[0,1)$.

\subsection{Proof of Theorem \ref{hfapprox0P2}}

As in the proof of Theorem \ref{hfapproxP2}, we use a probabilistic representation of $\mathcal{T}_{j} ^{0,l} (z)$ by using the Feynman-Kac formula. However, we introduce the jump Markov process which is a symmetric version with respect to reflecting barrier $(N-1)/N$ of that used in the proof of Theorem \ref{hfapproxP2}.

Let $\big( X^{N}_t \big)_{t\geq 0}$ be a jump Markov process  with state space $\big\{-(N-1)/N,\dots,(N-1)/N ,\dots, 3(N-1)/N\big\}$ and generator given by
\[{\mathcal{L}}^{N}\phi \left(\frac{l}{N}\right)=\Gamma^c _{\lvert l\rvert +2\,\lvert l\rvert +1} \left( \phi \left(\frac{l-1}{N}\right) -\phi \left(\frac{l}{N}\right) \right)+\Gamma^c _{\lvert l\rvert \,\lvert l\rvert +1}\left( \phi \left(\frac{l+1}{N}\right) -\phi \left(\frac{l}{N}\right) \right)\]
for $l\in \{-(N-2),\dots,-1 \} $,
\[{\mathcal{L}}^{N}\phi \left(\frac{l}{N}\right)= \Gamma^c _{l\,l+1}\left( \phi \left(\frac{l-1}{N}\right) -\phi \left(\frac{ l}{N}\right) \right)+\Gamma ^c_{ l +2\,l +1}\left( \phi \left(\frac{ l +1}{N}\right) -\phi \left(\frac{ l }{N}\right) \right)\]
for $l\in \{1,\dots,N-2 \} $,
\[\begin{split}{\mathcal{L}}^{N}\phi \left(\frac{l}{N}\right)&= \Gamma^c _{\lvert l-2(N-1)\rvert+2\,\,\lvert l-2(N-1)\rvert+1}\left( \phi \left(\frac{l-1}{N}\right) -\phi \left(\frac{ l}{N}\right) \right)\\
&\quad +\Gamma ^c_{ \lvert l-2(N-1)\rvert \,\,\lvert l-2(N-1)\rvert+1}\left( \phi \left(\frac{ l +1}{N}\right) -\phi \left(\frac{ l }{N}\right) \right)\end{split}\]
for $l\in \{N,\dots,2N-3 \} $,
\[\begin{split}{\mathcal{L}}^{N}\phi \left(\frac{l}{N}\right)&=\Gamma^c _{l+2-2(N-1)\,\,l+1-2(N-1)} \left( \phi \left(\frac{l+1}{N}\right) -\phi \left(\frac{l}{N}\right) \right)\\
&\quad+\Gamma^c _{l-2(N-1) \,\, l+1-2(N-1)}\left( \phi \left(\frac{l-1}{N}\right) -\phi \left(\frac{l}{N}\right) \right)\end{split}\]
for $l\in \{ 2N-1,\dots,3N-2\} $,
\[{\mathcal{L}}^{N}\phi \left(-\frac{ N-1}{N}\right)=\Gamma^{c} _{N-1\,N}\left( \phi \left(-\frac{ N-2}{N}\right) -\phi \left(-\frac{ N-1}{N}\right) \right),\]
\[{\mathcal{L}}^{N}\phi \left(\frac{ 3N-3}{N}\right)=\Gamma^{c} _{N-1\,N}\left( \phi \left(\frac{3N-4}{N}\right) -\phi \left(\frac{3N-3}{N}\right) \right),\]
\[{\mathcal{L}}^{N}\phi (0)=\frac{\Gamma^c _{2\,1}}{2}\left( \phi \left(\frac{1}{N}\right) -\phi (0) \right)+\frac{\Gamma^{c} _{2\,1}}{2}\left( \phi \left(\frac{-1}{N}\right) -\phi (0) \right),\]
\[{\mathcal{L}}^{N}\phi \left(\frac{N-1}{N}\right)\!=\!\frac{\Gamma^c _{N-1\,N}}{2}\left( \phi \left(\frac{N-2}{N}\right) -\phi \left(\frac{N-1}{N}\right) \right)+\frac{\Gamma^{c} _{N-1\,N}}{2}\left( \phi \left(\frac{N}{N}\right) -\phi \left(\frac{N-1}{N}\right) \right),\]
\[{\mathcal{L}}^{N}\phi \left(\frac{2N-2}{N}\right)\!=\!\frac{\Gamma^c _{2\,1}}{2}\left( \phi \left(\frac{2N-3}{N}\right) -\phi \left(\frac{2N-2}{N}\right) \right)+\frac{\Gamma^{c} _{2\,1}}{2}\left( \phi \left(\frac{2N-1}{N}\right) -\phi \left(\frac{2N-2}{N}\right)\right).\]
We recall that $\mathcal{T}^{0,l} ( z )$ can be viewed as a probability measure on $[0,1]$ by setting
\[\mathcal{T}^{0,l}_f ( z )=\sum_{j=1}^N f\left(\frac{j}{N}\right)\mathcal{T}^{0,l}_j ( z )\] for all bounded continuous function $f$ on $[0,1]$. Let $0<r\ll1$ and $f$ be a smooth function with support included in $[0,1-r)$. 
In order to make the link between $\mathcal{T}^{0,l} (z)$ and the process $X^N$, let us introduce an extension of $f$ by setting
\[ f^{N,s}(v)=\left\{\begin{array}{ccc} 
f\left(-v+1/N\right) &\text{ if }& v\in\big[-(N-1)/N,0\big]\\  
f\left(v+1/N\right) &\text{ if }& v\in\big[0,(N-1)/N\big]\\
f\left(-v+(2N-1)/N\right) &\text{ if }& v\in\big[(N-1)/N,2(N-1)/N\big]\\ 
f\left(v-(2N-3)/N\right) &\text{ if }& v\in\big[2(N-1)/N,(3N-3)/N\big].
\end{array}\right.\]
With these two functions we get the following representation. $\forall l\in\{1,\dots,N\}$,
\[\mathcal{T}^{0,l} _f (z) =  \mathbb{E}_{\frac{l-1}{N}}\left[ f^{N,s}(X^N _z) \right]. \]
Moreover, we have
\[\mathcal{T}^{0,l} _f (z) =  \mathbb{E}_{\frac{l-1}{N}}\left[ f^{s}(X^N _z) \right]+\mathcal{O}\left(\frac{1}{N}\right)=  \mathbb{E}_{\frac{l-1}{N}}\left[ f^{s}(g_r(X^N _z)) \right]+\mathcal{O}\left(\frac{1}{N}\right), \]
where
\[g_r(v)=\left\{\begin{array}{ccc}
v & \text{ if } &v\in (-(1-r),1-r) \cup (1+r,3-r) \\
v_s & \text{ elsewhere,}
\end{array}\right.\]
with $v_s \in (1-r,1-r/2)$, and where
\[ f^{s}(v)=\left\{\begin{array}{ccc} 
f(-v) &\text{ if }& v\in[-1,0]\\  
f(v) &\text{ if }& v\in[0,1]\\
f(-v+2) &\text{ if }& v\in[1,2]\\ 
f(v-2) &\text{ if }& v\in[2,3].
\end{array}\right.\]
Let $u\in[0,1)$ such that $l(N)/N \to u$. One can assume $u\in [0,1-r)$ by changing $r$ if necessary. 
As in the proof of Theorem \ref{hfapproxP2}, we have the following lemma. 
\begin{lem}\label{generatorP2p} 
$\forall\varphi\in {\mathcal{C}}^{\infty  }_{0}(\mathbb{R})$.
\[\lim_{N \to +\infty } \sup_{v\in I_N}\Big\lvert {\mathcal{L}}^{N} \varphi\left(\frac{[Nv]}{N}\right)-{\mathcal{L}_{a_{r,\infty}}}\varphi(v) \Big\rvert=0,\]
where 
\[\begin{split}
I_N&= \left[-\frac{N-1-[Nr]}{N},-\frac{1}{N}\right]\cup \left[\frac{1}{N},\frac{N-1-[Nr]}{N}\right]\\
&\quad \cup\left[\frac{N-1+[Nr]}{N},\frac{2N-3}{N}\right]\cup\left[\frac{2N-1}{N},\frac{3N-3-[Nr]}{N}\right],\end{split}\]
and $a_{r,\infty}$ is a $\mathcal{C}^1$-extended version of $a_\infty$ such that 
\[a_{r,\infty}(v)=\left\{\begin{array}{ccc} 
a_\infty(-v) &\text{ if }& v\in(-(1-r),0]\\  
a_\infty(v) &\text{ if }& v\in[0,1-r)\\
a_\infty(-v+2) &\text{ if }& v\in(1+r,2]\\ 
a_\infty(v-2) &\text{ if }& v\in[2,3-r),
\end{array}\right.\]  
and the martingale problem associated to $\mathcal{L}_{a_{r,\infty}}$ and starting from $u$ is well-posed.
\end{lem}

\begin{lem}\label{tightP2p}
The law of the process $(g_r(X^N))_N$ starting from $d(N)=(l(N)-1)/N$ is tight on $\mathcal{D}([0,+\infty),\mathbb{R})$.
\end{lem}
\begin{preuve}[of Lemma \ref{tightP2p}]
Let $\mathcal{F}^N_t=\sigma(X^N_s,\,\,s\leq t)$. According to Theorem $3$ in \cite[Chapter 3]{kushner}. We have to show only the two following points. First, we have  
\[\lim_{K\to +\infty } \overline{\lim_{N}} \,\,\mathbb{P}_{d(N)} \left( \sup_{t\geq 0}\lvert g_r(X^N_t) \rvert \geq K \right)=0,\]
since $\forall N$, $\sup_{t\geq 0}\lvert g_r(X^N_t) \rvert \leq 3$.
Second, we have for each $N$, $h \in (0,1)$, $s\in[0,h]$ and $t\geq 0$,
\[ \E _{d(N)} \big( (g_r(X^N_{t+s}) -g_r(X^N_{t}) )^2 \vert \mathcal{F}^N _t \big)\leq K\, h.\]
In fact, we have
\[\begin{split}
\E_{d(N)}  \big( (g_r(X^N_{t+s}) -g_r(X^N_{t}) )^2 \vert \mathcal{F}^N_t\big) & \leq  2\,\E _{d(N)} \big( (M^N_{g_r}(t+s)  - M^N _{g_r}(t) )^2 \vert \mathcal{F}^N_t\big)\\
                                                                                          & + 2\,\E _{d(N)}\left( \left(\int_{t}^{t+s} {\mathcal{L}}^{N}g_r(X^N_{w})dw\right) ^2\Big\vert \mathcal{F}^N_t\right),
\end{split} \]                                                                                         
with 
\[M^N _{g_r} (t) =g_r(X^N_t) -g_r(X^N_0) -\int_{0}^{t} {\mathcal{L}}^{N} g_r (X^N_s)ds,\] 
which is a $( \mathcal{F}^N_t)_{t\geq 0}$-martingale. We also have   
\[\sup_N \sup_{v\in \left[-\frac{N-1}{N},3\frac{N-1}{N}\right]\setminus\big\{0,2\frac{N-1}{N}\big\}} \lvert {\mathcal{L}}^{N}g_r(v)\lvert <+\infty\]
since by Lemma \ref{generatorP2p}
\[\sup_N\sup_{v\in I_N\cup\{v_s\}}\big\lvert \mathcal{L}^N g_r(v)\big\rvert<+\infty.\]
Moreover, $\mathcal{L}^{N}g_r(0)=\mathcal{L}^{N}g_r(2(N-1)/N)=0$. Then, we get
\[\E _{d(N)}\left( \left(\int_{t}^{t+s} {\mathcal{L}}^{N}g_r( X^N_w )dw\right) ^2\Big \vert \mathcal{F}^N_t \right)  \leq C h^2 .\]
We recall that
\[< M^N _{g_r} > _t =\int_{0}^{t}\left( {\mathcal{L}}^{N}{g_r}^2-2g_r{\mathcal{L}}^{N}g_r \right)(X^N_s)ds.\]                                                                                         
Consequently, by the martingale property of $(M^N _{g_r} (t))_{t\geq 0}$, 
\[\begin{split}
\E_{d(N)} \big( (M^N _{g_r} (t+s)  - M^N _{g_r} (t) )^{2}\vert \mathcal{F}^N_t \big) &= \E_{d(N)}\left(( {M^N _{g_r}}(t+s)  - {M^N_{g_r}}(t))^{2}\vert \mathcal{F}^N_t\right)\\
&= \E_{d(N)}\left( {M^N_{g_r}}(t+s) ^2  - {M^N_{g_r}}(t)^{2}\vert \mathcal{F}^N_t\right) \\
&=\E_{d(N)}\left( < M^N_{g_r} > _{t+s}  - < M^N _{Id} > _{t} \vert \mathcal{F}^N_t\right)\\
&= \E_{d(N)}\left( \int_{t}^{t+s}\left( \mathcal{L}^{N}g_r ^2-2g_r\mathcal{L}^{N}g_r \right)(X^N _w)dw \Big\vert \mathcal{F}^N_t\right) \\
&\leq C\, h.
\end{split}\] 
In fact, in addition to the previous arguments, we also have 
\[\sup_N \sup_{v\in I_N\cup\{v_s\}} \lvert {\mathcal{L}}^{N}g_r^2 (v) \rvert  <  +\infty,\] 
$\sup_{N}{\mathcal{L}}^{N}g_r^2(0) =\frac{\Gamma^c _{1\,2}}{N^2}<+\infty$, and $\sup_{N}{\mathcal{L}}^{N}g_r^2(2(N-1)/2) =2\frac{\Gamma^c _{1\,2}}{N^2}<+\infty$. That concludes the proof Lemma \ref{tightP2p}.
$\square$
\end{preuve} 

Now, let us introduce some notations. $\forall j \in \mathbb{N}^{\ast}$, let
\[\begin{split}
\tau_{r}^{(j)}&=\inf\big(t>\tau_{r,c}^{(j-1)},\quad x(t)\in[-1,-(1-r))\cup(1-r,1+r)\cup(3-r,3] \big)\\
\tau_{r,c}^{(j)}&=\inf\big(t>\tau_r^{(j)},\quad x(t)\in (-(1-r),1-r)\cup(1+r,3-r)\big),
\end{split}\]
with $\tau_{r,c}^{(0)}=0$. Using the previous lemma, there exists $(N')$ such that 
\[\lim_{N'\to +\infty}\mathbb{E}_{d(N')}\big[ f^s (g_r(X^{N'}_z))\big]=\mathbb{E}^{\mathbb{Q}_u}\big[ f^s (x(z))\big].\]
Moreover,
\[\begin{split}
\mathbb{E}^{\mathbb{Q}_u}\big[ f^s (x(z))\big]&=\sum_{j\geq 1}\mathbb{E}^{\mathbb{Q}_u}\Big[ f^s (x(z))\textbf{1}_{(\tau_{r,c}^{(j-1)}\leq z<\tau_r^{(j)})}\Big]\\
&=\sum_{j\geq 1}\mathbb{E}^{\mathbb{Q}_u}\Big[\mathbb{E}^{\mathbb{Q}_u}\Big[ f^s (x(z))\textbf{1}_{(\tau_{r,c}^{(j-1)}\leq z<\tau_r^{(j)})}\Big\vert \mathcal{M}_{\tau_{r,c}^{(j-1)}} \Big] \Big],
\end{split}\]
where $\mathcal{M}_t=\sigma(x(s),\,\, 0\leq s\leq t)$. With the following lemma we can identify each excursion between $\tau_{r,c}^{(j-1)}$ and $\tau_{r}^{(j)}$.
\begin{lem}\label{caractP2p}
$\forall j\in \mathbb{N}^{\ast}$, the conditional law ${\mathbb{Q}_{u}}\big(\cdot\big\vert \mathcal{M}_{\tau_{r,c}^{(j-1)}}\big)$ coincide up to the stopping time $\tau_r^{(j)}$ with the conditional law ${\overline{\mathbb{P}}^r_u}\big(\cdot \big\vert \mathcal{M}_{\tau_{r,c}^{(j-1)}}\big)$, where $\overline{\mathbb{P}}^r_u$ is the unique solution of the martingale problem associated to $\mathcal{L}_{a_{r,\infty}}$ and starting from $u$.
\end{lem}
\begin{preuve}[of Lemma \ref{caractP2p}]
This proof is a conditional version of Lemma \ref{martingaleP22}. Moreover, this lemma follows from Lemma \ref{generatorP2p} and the fact that we are studying excursions between $\tau_{r,c}^{(j-1)}$ and $\tau_{r}^{(j)}$. By Lemma \ref{generatorP2p}, in addition to $g_r(X^N_z)=X^N_z$ for $\tau_{r,c}^{(j-1)}\leq z<\tau_r^{(j)}$,
\[\lim_N \E _{d(N)} \left[ \int_{\tau_{r,c}^{(j-1)}}^{t\wedge\tau _r^{(j)} } \left \lvert{\mathcal{L}}^{N} \varphi(X^N_s)-\mathcal{L}_{a_{r,\infty}} \varphi(X^N_s)\right \rvert ds \Big\vert \mathcal{M}_{\tau_{r,c}^{(j-1)}}\right]=0,\]
and we also have 
\[\begin{split}
 \E_{0}\left[ \int_{0}^t \textbf{1}_{( X^N _s =0 ) }  ds\right] &=\mathcal{O}\left( \frac{1}{N^{\alpha' \wedge (1-\alpha')} }\right),\\
 \E_{2(N-1)/N}\left[ \int_{0}^t \textbf{1}_{( X^N _s =2(N-1)/N ) }  ds\right] &=\mathcal{O}\left( \frac{1}{N^{\alpha' \wedge (1-\alpha')} }\right)
 \end{split}\]
by symmetry of the process $X^N$. As in the proof of Lemma \ref{martingaleP22}, we get that
$\forall \varphi \in \mathcal{C}^{\infty}_{0} (\mathbb{R})$
\[ \varphi(x(t\wedge\tau_{r}^{(j)}))-\varphi(x(\tau_{r,c}^{(j-1)}))-\int_{\tau_{r,c}^{(j-1)}}^{t\wedge\tau_{r}^{(j)}}\mathcal{L}_{a_{r,\infty}}\varphi(x(s))ds \]
is a martingale under the conditional law ${\mathbb{Q}_{u}}\big(\cdot\big\vert \mathcal{M}_{\tau_{r,c}^{(j-1)}}\big)$.
Finally, from the uniqueness of the martingale problem associated to $\mathcal{L}_{a_{r,\infty}}$, ${\mathbb{Q}_{u}}\big(\cdot\big\vert \mathcal{M}_{\tau_{r,c}^{(j-1)}}\big)$ coincide up to the stopping time $\tau_r^{(j)}$ with ${\overline{\mathbb{P}}^r_u}\big(\cdot\big\vert \mathcal{M}_{\tau_{r,c}^{(j-1)}}\big)$ (see Theorem 6.2.2 in \cite{stroock}). That concludes the proof of Lemma \ref{caractP2p}.$\square$
\end{preuve}
From the previous lemma, $\forall j\in\mathbb{N}^{\ast}$, we have
\[ \mathbb{E}^{\mathbb{Q}_u}\Big[ f^s (x(z))\textbf{1}_{(\tau_{r,c}^{(j-1)}\leq z<\tau_r^{(j)})}\Big]=\mathbb{E}^{\overline{\mathbb{P}}^r_u}\Big[ f^s (x(z))\textbf{1}_{(\tau_{r,c}^{(j-1)}\leq z<\tau_r^{(j)})}\Big] \]
and then
\[\lim_{N'\to +\infty}\mathbb{E}_{d(N')}\big[ f^s (g_r(X^{N'}_z))\big]=\mathbb{E}^{\mathbb{Q}_u}\big[ f^s (x(z))\big]=\mathbb{E}^{\overline{\mathbb{P}}^r_v}\big[ f^s (x(z))\big],  \]
where the limit does not depend to $(N')$. Consequently,
\[\lim_{N\to +\infty}\mathcal{T}^{0,l(N)} _f (z)=\mathbb{E}^{\overline{\mathbb{P}}^r_u}\big[ f^s (x(z))\big]=\mathcal{T}_f (z,u),\]
with
\[\frac{\partial}{\partial z}\mathcal{T}_f (z,u)=\mathcal{L}_{a_{r,\infty}}\mathcal{T}_f (z,u)=\mathcal{L}_{a_{\infty}}\mathcal{T}_f (z,u).\]
For the boundary conditions, first let $h\in(0,1)$ such that $0<h\ll 1$, we have
\[ \frac{1}{h}\big(\mathcal{T}_f (z,h)-\mathcal{T}_f (z,-h)\big)=\frac{1}{h}\lim_{N\to +\infty} \big(\mathbb{E}_{\frac{[Nh]}{N}}\big[ f^s (X^{N}_z)\big] -\mathbb{E}_{-\frac{[Nh]}{N}}\big[ f^s (X^{N}_z)\big] \big)=0,\] 
because of the symmetry of the process $X^N$ and $f^s$, and therefore, 
\[2\frac{\partial}{\partial u}\mathcal{T}_f (z,0)=0.\]
Second, in the same way, let $h\in(0,1)$ such that $h\ll 1$. Moreover, one can assume $r<h$ by changing $r$ if necessary. Then, we have
\[ \frac{1}{h}\big(\mathcal{T}_f (z,1-h)-\mathcal{T}_f (z,1+h)\big)=\frac{1}{h}\lim_{N\to +\infty} \big(\mathbb{E}_{\frac{[N(1-h)]}{N}}\big[ f^s (X^{N}_z)\big] -\mathbb{E}_{\frac{[N(1+h)]}{N}}\big[ f^s (X^{N}_z)\big] \big)=0,\]
and therefore,
\[2\frac{\partial}{\partial u}\mathcal{T}_f (z,1)=0.\]
As a result, using the density of the smooth functions with compact support in $L^2(0,1)$ for $\|.\|_{L^2(0,1)}$ and the dominated convergence theorem we get the first point of Theorem \ref{hfapprox0P2}. The second point is a consequence of the maximum principle and the density for the sup norm over $[0,1]$ in $\{\varphi\in\mathcal{C}^0([0,1]),\,\varphi(1)=0\}$ of the smooth functions with compact support included in $[0,1)$.
$\blacksquare$

\section*{Acknowledgments}  I wish to thank my Ph.D. supervisor Josselin Garnier for his suggestions, his help and his support.


\begin{thebibliography}{20} 

\bibitem{adler}
{\sc R.J.~Adler}, 
{\em The geometry of random fields}, Wiley, London, 1981.

\bibitem{adlertaylor}
{\sc R.J~ Adler and J.~Taylor}, 
{\em Random fields and geometry}, Springer, New York, 2007.

\bibitem{billingsley}
{\sc P.~Billingsley}, {\em
Convergence of probability measure}, $2^{nd}$ ed.,
Wiley InterScience, 1999. 

\bibitem{carmona}
{\sc R.~Carmona and J.-P.~Fouque}, {\em Diffusion-approximation for the advection-diffusion of a passive scalar by a space-time gaussian velocity field}, Seminar on Stochastic Analysis, Random Fields and Applications. (Edited by E. Bolthausen, M. Dozzi and F. Russo). Birkhauser, Basel (1994) 37-50. 

\bibitem{donvar} 
{\sc M.D.~Donsker and S.R.S.~Varadhan}, {\em Asymptotic evaluation of certain Markov process expectations for large time, I-IV}, Comm. Pure Appl. Math, 28 (1975), pp.~1--47, 279--301, 29 (1979), pp.~368--461, 36 (1983) pp.~183--212.

\bibitem{fouque}
{\sc J.-P.~Fouque}, {\em La convergence en loi pour les processus \`a valeur dans un espace nucl{\'e}aire}, Ann. Inst. Henri Poincar\'e, 20 (1984), pp.~225--245.

\bibitem{book}
{\sc J.-P.~Fouque, J.~Garnier, G.~Papanicolaou, and K.~S{\o}lna}, {\em Wave propagation and time reversal in randomly layered  media}, Springer, New York, 2007.

\bibitem{garniereva}
{\sc J.~Garnier}, {\em The role of evanescent modes in randomly perturbed single-mode waveguides }, Discrete and Continuous Dynamical Systems-Series B, 8 (2007), pp.~455--472.


\bibitem{papa}
{\sc J.~Garnier and G.~Papanicolaou}, {\em Pulse propagation and time reversal in random waveguides}, SIAM J. Appl. Math, 67 (2007), pp.~1718--1739.

\bibitem{garnier2}
{\sc J.~Garnier and K.~S{\o}lna}, {\em Effective transport equations and enhanced backscattering in random waveguides }, SIAM J. Appl. Math, 68 (2008), pp.~1574--1599.

\bibitem{gomez}
{\sc C.~Gomez}, {\em Time-reversal superresolution in random waveguides}, SIAM Multiscale Model. Simul. 7 (2009), pp.~1348--1386.

\bibitem{gomez2}
{\sc C.~Gomez}, {\em Wave propagation in shallow-water acoustic random waveguides}, preprint (2009).

\bibitem{inf2}
{\sc G.~Kallianpur and J.~Xiong}, {\em Stochastic differential equations in infinite dimensional spaces}, IMS Lecture notes-monograph series, 1995.

\bibitem{papanicolaou}
{\sc W.~Kohler and G.~Papanicolaou}, {\em Wave propagation in randomly inhomogeneous ocean}, Lecture Notes in Physics, Vol. 70, J. B. Keller and J. S. Papadakis, eds., Wave Propagation and Underwater acoustics, Springer-Verlag, Berlin, 1977.

\bibitem{khasminskii}
{\sc R.Z.~Khasminskii}, {\em A limit theorem for solutions of differential equations with random right hand side}, Theory Probab. Appl. 11 (1966), pp.~390--406.


\bibitem{kushner}
{\sc H.-J.~Kushner}, {\em Approximation and weak convergence methods for random processes}, MIT press, Cambridge, 1984.

\bibitem{magnanini}
{\sc R.~Magnanini and F.~Santosa}, {\em Wave propagation in a 2-D optical waveguide}, SIAM J. Appl. Math, 61 (2000),
pp.~1237--1252.


\bibitem{marcuse}
{\sc D.~Marcuse}, {\em Theory of dielectric optical waveguides}, $2^{nd}$ ed., Academic press, New York, 1991.

\bibitem{mellen}
{\sc R.H.~Mellen, D.G.~Browning, and J.M.~Ross}, {\em Attenuation in randomly inhomogeneous sound chanels},  J. Acoust. Soc. Am., 56 (1974), no. 3, pp.~80--82.


\bibitem{metivier}
{\sc M.~Metivier}, {\em Stochastic partial differential equations in infinite dimensional spaces}, Scuola normale superiore, Pisa, 1988.




\bibitem{mitoma}
{\sc I.~Mitoma}, {\em On the sample continuity of $\mathcal{S}'$-processes}, J. Math. Soc. Japan, 35 (1983), pp.~629--636.

\bibitem{papakohler}
{\sc G.~Papanicolaou and W.~Kohler}, {\em Asymptotic theory of mixing stochastic ordinary differential equations}, Comm. Pure Appl. Math, 27 (1974), pp.~641–668.  

\bibitem{pekeris}
{\sc C.I.~Pekeris}, {\em Theory of propagation of explosive sound in shallow water}, propagation of sound in the ocean, geological society of America, memoir 27, pp.~1--117. 

\bibitem{perrey}
{\sc E.~Perrey-Debain and I.D.~Abrahams}, {\em A diffusion analysis approach to TE mode propagation in randomly perturbed optical waveguides}, SIAM J. Appl. Math, 68 (2007), pp.~523--543.

\bibitem{rowe}
{\sc H.E.~Rowe}, {\em Electromagnetic propagation in multi-mode random media}, Wiley, New York, 1999.

\bibitem{schaefer}
{\sc H.~Schaefer}, {\em Topological vector spaces}, Springer-Verlag, New York, 1971.


\bibitem{stroock}
{\sc D.W.~Stroock and S.R.S.~Varadhan}, {\em Multidimensional diffusion processes}, Springer-Verlag, Berlin, 1979.

\bibitem{wilcox}
{\sc C.~Wilcox}, {\em Spectral analysis of the Pekeris operator in the theory of acoustic wave propagation in shallow water},  Arch. Rational Mech. Anal, 60 (1975/76), no. 3, pp.~259--300.

\bibitem{wilcox2}
{\sc C.~Wilcox}, {\em Transient electromagnetic wave propagation in a dielectric waveguide}, in Symposia Mathematica, vol. XVIII (Convegno sulla Teoria Matematica dell'Elettromagnetismo, INDAM, Rome, 1974), Academic Press, London, 1976, pp.~239--277. 

\bibitem{inf}
{\sc M.~Yor}, {\em Existence et unicit{\'e} de diffusion {\`a} valeurs dans un espace de Hilbert}, Ann. Inst. Henri Poincar{\'e}, 10 (1974), pp.~55--88.


\end{thebibliography}
\end{document}